\documentclass[11pt,leqno]{amsart}
\usepackage{amsmath, amssymb}
\usepackage{graphicx,color,hyperref}
\usepackage{verbatim,enumitem}

\def\smallskip{\vskip3pt plus1pt minus1pt}
\def\medskip{\vskip5pt plus2pt minus1pt}

\def\plainsubsection#1|{%
  \par\vskip0.25cm\penalty -100
  \centerline{{\sc #1}}
  \vskip3pt plus 1pt minus 0pt
  \penalty 500}

\long\def\claim#1|#2\endclaim{\par\vskip 4pt\noindent 
{\bf #1.}\ {\em #2}\par\vskip 5pt}

\def\plainproof{\noindent{\em Proof}}

\def\today{\ifcase\month\or
January\or February\or March\or April\or May\or June\or July\or August\or
September\or October\or November\or December\fi \space\number\day,
\number\year}

\catcode`\@=11
\newcount\@tempcnta \newcount\@tempcntb 
\def\timeofday{{%
\@tempcnta=\time \divide\@tempcnta by 60 \@tempcntb=\@tempcnta
\multiply\@tempcntb by -60 \advance\@tempcntb by \time
\ifnum\@tempcntb > 9 \number\@tempcnta:\number\@tempcntb
  \else\number\@tempcnta:0\number\@tempcntb\fi}}
\catcode`\@=12

\def\Bibitem#1&#2&#3&#4&%
{\hangindent=2cm\hangafter=1
\noindent\rlap{\hbox{\bf #1}}\kern2cm{\rm #2}{\it #3}{\rm #4.}} 

\def\bB{{\mathbb B}}
\def\bC{{\mathbb C}}
\def\bD{{\mathbb D}}
\def\bG{{\mathbb G}}
\def\bH{{\mathbb H}}
\def\bK{{\mathbb K}}
\def\bN{{\mathbb N}}
\def\bP{{\mathbb P}}
\def\bQ{{\mathbb Q}}

\def\bZ{{\mathbb Z}}
\def\bOne{{\mathchoice {\rm 1\mskip-4mu l} {\rm 1\mskip-4mu l}
{\rm 1\mskip-4.5mu l} {\rm 1\mskip-5mu l}}}

\def\cA{{\mathcal A}}

\def\cE{{\mathcal E}}
\def\cF{{\mathcal F}}
\def\cG{{\mathcal G}}
\def\cI{{\mathcal I}}
\def\cJ{{\mathcal J}}
\def\cL{{\mathcal L}}

\def\cO{{\mathcal O}}

\def\cR{{\mathcal R}}

\def\cV{{\mathcal V}}
\def\cW{{\mathcal W}}
\def\cX{{\mathcal X}}
\def\cY{{\mathcal Y}}

\def\gm{{\frak m}}

\def\bfk{{\bf k}}
\def\bfe{{\bf e}}

\def\abu{a_{\scriptscriptstyle\bullet}}
\def\bbu{b_{\scriptscriptstyle\bullet}}
\def\cbu{c_{\scriptscriptstyle\bullet}}
\def\onebu{1_{\scriptscriptstyle\bullet}}

\def\ii{i}
\def\ld{,\ldots,}
\def\bu{{\scriptstyle\bullet}}
\def\ort{\mathop{\hbox{\kern1pt\vrule width4pt height0.4pt depth0pt
    \vrule width0.4pt height7pt depth0pt\kern3pt}}}
\def\bigglp#1({\raise-#1\hbox{$\bigg($}}
\def\biggrp#1){\raise-#1\hbox{$\bigg)$}}

\def\qedsquare{\hbox{
\vrule height 1.5ex  width 0.1ex  depth 0ex\kern-0.1ex
\vrule height 1.5ex  width 1.5ex  depth -1.4ex\kern-1.5ex
\vrule height 0.1ex  width 1.5ex  depth 0ex\kern-0.1ex
\vrule height 1.5ex  width 0.1ex  depth 0ex}\kern0.5pt}
\def\qed{~\hfill\qedsquare\vskip6pt plus2pt minus1pt}
\def\lambdawedge{\mathop{\raise1.5pt\hbox{$\scriptstyle\bigwedge$}}\nolimits}

\let\ssm\smallsetminus

\let\le\leqslant
\let\leq\leqslant

\let\ge\geqslant
\let\geq\geqslant

\let\ol=\overline

\let\wt=\widetilde

\def\swt#1{\smash{\widetilde#1}}
\def\swh#1{\smash{\widehat#1}}
\def\build#1^#2_#3{\mathrel{\mathop{\null#1}\limits^{#2}_{#3}}}

\def\mertorelbar{\vrule width0.6ex height0.65ex depth-0.55ex}
\def\merto{\mathrel{\mertorelbar\kern1.3pt\mertorelbar\kern1.3pt\mertorelbar
    \kern1.3pt\mertorelbar\kern-1ex\raise0.28ex\hbox{${\scriptscriptstyle>}$}}}
\let\lra=\longrightarrow
\let\lra=\longrightarrow
\def\lhra{\lhook\joinrel\longrightarrow}

\def\lraww{\mathrel{\rlap{$\longrightarrow$}\kern-1pt\longrightarrow}}
\def\vdasharrow{\rotatebox{-90}{$\dashrightarrow$}}

\catcode`\@=11
\newdimen\@rrowlength \@rrowlength=6ex
\def\ssrelbar{\vrule width\@rrowlength height0.64ex depth-0.56ex\kern-4pt}
\def\llra#1{\@rrowlength=#1\ssrelbar\rightarrow}
\def\vlra#1{\hbox to#1mm{\rightarrowfill}}
\catcode`\@=12

\def\lcm{\mathop{\rm lcm}\nolimits}
\def\card{\mathop{\rm card}\nolimits}

\def\Id{\mathop{\rm Id}\nolimits}
\def\Ker{\mathop{\rm Ker}\nolimits}

\def\Hom{\mathop{\rm Hom}\nolimits}

\def\Aut{\mathop{\rm Aut}\nolimits}
\def\Tr{\mathop{\rm Tr}\nolimits}
\def\GL{\mathop{\rm GL}\nolimits}
\def\PGL{\mathop{\rm PGL}\nolimits}

\def\Pic{\mathop{\rm Pic}\nolimits}

\def\Proj{\mathop{\rm Proj}\nolimits}
\def\Supp{\mathop{\rm Supp}\nolimits}
\def\Vol{\mathop{\rm Vol}\nolimits}
\def\Ricci{\mathop{\rm Ricci}\nolimits}
\def\Vect{\mathop{\rm Vect}\nolimits}
\def\codim{\mathop{\rm codim}\nolimits}
\def\rank{\mathop{\rm rank}\nolimits}
\def\div{\mathop{\rm div}\nolimits}

\def\mod{\mathop{\rm mod}\nolimits}
\def\pr{\mathop{\rm pr}\nolimits}
\def\Gr{\mathop{\rm Gr}\nolimits}

\def\Bs{\mathop{\rm Bs}\nolimits}
\def\dbar{{\overline\partial}}
\def\ddbar{{\partial\overline\partial}}
\def\bddK{{{}^b\kern-1pt K}}

\def\IEL{\mathop{\rm IEL}\nolimits}
\def\ECL{\mathop{\rm ECL}\nolimits}
\def\Mordell{\mathop{\rm Mordell}\nolimits}
\def\GL{\mathop{\rm GL}\nolimits}
\def\Sing{\mathop{\rm Sing}}
\def\Reg{\mathop{\rm Reg}}

\def\reg{{\rm reg}}
\def\sing{{\rm sing}}
\def\Sing{{\rm Sing}}

\def\GG{{\rm GG}}
\def\Kob{{\rm Kob}}
\def\loc{{\rm loc}}

\def\Const{{\rm Const}}

\newdimen\plainitemindent \plainitemindent=18pt
\def\plainitem#1{\vskip3pt\noindent
\hangindent\plainitemindent\hbox to\plainitemindent{#1\hss}\ignorespaces}

\catcode`\@=11
\def\openup{\afterassignment\@penup\dimen@=}
\def\@penup{\advance\lineskip\dimen@
  \advance\baselineskip\dimen@
  \advance\lineskiplimit\dimen@}
\newdimen\jot \jot=3pt
\newskip\plaincentering \plaincentering=0pt plus 1000pt minus 1000pt
\def\ialign{\everycr{}\tabskip\z@skip\halign}
\def\eqalign#1{\null\,\vcenter{\openup\jot\m@th
  \ialign{\strut\hfil$\displaystyle{##}$&$\displaystyle{{}##}$\hfil
      \crcr#1\crcr}}\,}
\newif\ifdt@p
\def\displ@y{\global\dt@ptrue\openup\jot\m@th
  \everycr{\noalign{\ifdt@p \global\dt@pfalse \ifdim\prevdepth>-1000\p@
      \vskip-\lineskiplimit \vskip\normallineskiplimit \fi
      \else \penalty\interdisplaylinepenalty \fi}}}
\def\@lign{\tabskip\z@skip\everycr{}} % restore inside \displ@y
\def\displaylines#1{\displ@y \tabskip\z@skip
  \halign{\hbox to\displaywidth{$\@lign\hfil\displaystyle##\hfil$}\crcr
    #1\crcr}}
\def\eqalignno#1{\displ@y \tabskip\plaincentering
  \halign to\displaywidth{\hfil$\@lign\displaystyle{##}$\tabskip\z@skip
    &$\@lign\displaystyle{{}##}$\hfil\tabskip\plaincentering
    &\llap{$\@lign##$}\tabskip\z@skip\crcr
    #1\crcr}}
\def\leqalignno#1{\displ@y \tabskip\plaincentering
  \halign to\displaywidth{\hfil$\@lign\displaystyle{##}$\tabskip\z@skip
    &$\@lign\displaystyle{{}##}$\hfil\tabskip\plaincentering
    &\kern-\displaywidth\rlap{$\@lign##$}\tabskip\displaywidth\crcr
    #1\crcr}}
\def\plaincases#1{\left\{\,\vcenter{\normalbaselines\m@th
    \ialign{$##\hfil$&\quad##\hfil\crcr#1\crcr}}\right.}
\def\plainmatrix#1{\null\,\vcenter{\normalbaselines\m@th
    \ialign{\hfil$##$\hfil&&\quad\hfil$##$\hfil\crcr
      \mathstrut\crcr\noalign{\kern-\baselineskip}
      #1\crcr\mathstrut\crcr\noalign{\kern-\baselineskip}}}\,}

\newcommand\myleaders{\leavevmode\leaders\hbox{\kern1pt.}\hfill}
\renewcommand{\tocsection}[3]{%
  \indentlabel{\@ifnotempty{#2}{\ignorespaces#1 
      \rlap{\kern0.1em #2.}\hphantom{00.\ \ }}}#3\myleaders}
\catcode`\@=12

\title[Kobayashi and Green-Griffiths-Lang conjectures]{Recent 
results on the Kobayashi and\\ Green-Griffiths-Lang 
conjectures}

\author{Jean-Pierre Demailly}
\date{}
\thanks{Expanded version of the Takagi lectures delivered 
at the University of Tokyo on November 6, 2015}
\thanks{Work supported by the advanced ERC grant ALKAGE No.\ 670846 
started in September 2015}

\begin{document}

\begin{abstract}
The study of entire holomorphic curves contained in
projective algebraic varieties is intimately related to fascinating
questions of geometry and number theory -- especially through the
concepts of curvature and positivity which are central themes in
Kodaira's contributions to mathematics. The aim of these lectures is to
present recent results concerning the geometric side of the problem. The
Green-Griffiths-Lang conjecture stipulates that for every projective
variety $X$ of general type over~$\bC$, there exists a proper
algebraic subvariety $Y$ of $X$ containing all entire curves
$f:\bC\to X$.  Using the formalism of directed varieties and jet bundles, 
we show that this assertion holds true in case $X$ satisfies a strong 
general type
condition that is related to a certain jet-semi-stability property of
the tangent bundle $T_X$. It is possible to exploit similar techniques
to investigate a famous conjecture of Shoshichi Kobayashi (1970),
according to which a generic algebraic hypersurface of dimension $n$
and of sufficiently large degree $d\ge d_n$ in the complex projective
space $\bP^{n+1}$ is hyperbolic: in the early 2000's, Yum-Tong Siu 
proposed a strategy that led in 2015 to a proof based on a clever 
use of slanted vector fields on jet spaces, combined with Nevanlinna
theory arguments. In 2016, the conjecture has been settled in a different
way by Damian Brotbek, making a more direct use of Wronskian differential
operators and associated multiplier ideals; shortly afterwards, Ya Deng 
showed how the proof could be modified to yield an explicit value of $d_n$.
We give here a short proof based on a substantial simplification of their
ideas, producing a bound very similar to Deng's original estimate, namely
$d_n=\lfloor\frac{1}{3}(en)^{2n+2}\rfloor$.
\vskip10pt
\noindent
{\bf MSC~Classification (2010):} 32H20, 32L10, 53C55, 14J40
\end{abstract}

\maketitle
\hbox to \textwidth{\hfill\it Contribution to the $16^{th}$ Takagi lectures
\hfill}
\hbox to \textwidth{\hfill\it in celebration of the $100^{th}$ anniversary of
K.~Kodaira's birth\hfill}
\vskip25pt

\tableofcontents

\setcounter{section}{-1}
\section{Introduction}

The goal of these lectures is to study the conjecture of Kobayashi
[Kob70] on the hyperbolicity of generic hypersurfaces of high
degree in projective space, and the related conjecture by
Green-Griffiths [GrGr80] and Lang [Lang86] on the structure of
entire curve loci. 
\medskip

Let us recall that a complex space $X$ is said to be hyperbolic in the
sense of Kobayashi if analytic disks $f:\bD\to X$ through a given
point form a normal family. By a well known result of Brody [Bro78], a compact
complex space is Kobayashi hyperbolic iff it does not contain any entire
holomorphic curve $f:\bC\to X$ (``Brody hyperbolicity'').

In this paper entire holomorphic curves are assumed to be
non-constant and simply called entire curves. If $X$ is not hyperbolic,
a basic question
is thus to analyze the geometry of entire holomorphic curves $f:\bC\to X$,
and especially to understand the {\em entire curve locus of~$X$}, defined as
the Zariski closure
$$
\ECL(X)=\overline{\bigcup_f f(\bC)}^{\rm Zar}.\leqno(0.1)
$$
The Green-Griffiths-Lang conjecture, in its strong form, can be stated 
as follows.

\claim 0.2. GGL conjecture|Let $X$ be a projective variety of general type. 
Then $Y=\ECL(X)$ is a proper algebraic subvariety $Y\subsetneq X$.
\endclaim

\noindent Equivalently, there exists $Y\subsetneq X$ such that every 
entire curve $f:\bC\to X$ satisfies $f(\bC)\subset Y$.
A~weaker form of the GGL conjecture states that entire curves are
algebraically degenerate, i.e.\ that
$f(\bC)\subset Y_f\subsetneq X$ where $Y_f$ may depend on $f$. 

If $X\subset\bP^N_\bC$ is defined over a number field $\bK_0$ (i.e.\ by
polynomial equations with coefficients in $\bK_0$), one defines
the Mordell locus, denoted $\Mordell(X)$, to be the smallest complex subvariety
$Y$ in $X$ such that the set of $\bK$-points $X(\bK)\smallsetminus Y$ is
finite for every number field~$\bK\supset\bK_0$. Lang [Lang86]
conjectured that one should always have $\Mordell(X)=\ECL(X)$
in this situation. This conjectural arithmetical statement would
be a vast generalization of the Mordell-Faltings theorem, and is 
one of the strong motivations to study the geometric GGL conjecture
as a first~step. S.~Kobayashi [Kob70] had earlier made
the following tantalizing conjecture.

\claim 0.3.~Conjecture {\rm(Kobayashi)}|\vskip0pt
\plainitem{\rm(a)}A $($very$)$ generic hypersurface
$X\subset\bP^{n+1}$ of degree $d\ge d_n$ large enough is hyperbolic,
especially it does not possess any entire holomorphic curve $f:\bC\to X$.
\vskip0pt
\plainitem{\rm(b)}The complement $\bP^n\ssm H$ of a $($very$)$ generic hypersurface
$H\subset\bP^n$ of degree $d\ge d'_n$ large enough is hyperbolic.
\vskip0pt
\endclaim

It should be noticed that the existence of a smooth hyperbolic hypersurface
$X\subset\bP^{n+1}$ in 0.3~(a), or a hyperbolic complement
$\bP^n\ssm H$ with $H$ smooth irreducible in 0.3~(b), is already a
hard problem; many efforts were initially concentrated on this
problem. As Zaidenberg observed, a smooth deformation of a union of
$2n + 1$ hyperplanes in $\bP^n$ is not necessarily Kobayashi
hyperbolic, and the issue is non-trivial at all. The existence problem
was initially solved for sufficiently high degree hypersurfaces through
a number of examples:\vskip2pt

\noindent
$\scriptstyle\bullet$ case (a) for $n = 2$ and degree $d\geq 50$
by Brody–Green [BrGr77]$\,$;\vskip2pt
\noindent
$\scriptstyle\bullet$ case (b) for $n = 2$ by [AzSu80]
(as a consequence of [BrGr77])$\,$;\vskip2pt
\noindent
$\scriptstyle\bullet$ cases (a) and (b) for $n \geq 3$
by Masuda–Noguchi [MaNo96].\vskip2pt
\noindent
Improvements in the degree estimates were later obtained in 
[Shi98], [Fuj01], [ShZa02], in addition to many
other papers dealing with low dimensional varieties ($n=2,3$).\medskip

We now describe a number of known results concerning the question of
generic hyperbolicity, according to the Kobayashi conjectures 0.3~(a,b).
M.~Zaidenberg observed in [Zai87] that the complement of a general 
hypersurface of degree $2n$ in $\bP^n$ is not hyperbolic; as a consequence,
one must take $d'_n\ge 2n+1$ in 0.3~(b). This observation, along with
Fujimoto's classical result that the complement
of $2n + 1$ hyperplanes of $\bP^n$ in general position is hyperbolic and
hyperbolically embedded in $\bP^n$ ([Fuj72]) led Zaidenberg to propose
the bound $d'_n=2n+1$ for $n\ge 1$.
Another famous result due to Clemens 
[Cle86], Ein [Ein88, Ein91] and Voisin [Voi96], states that every
subvariety $Y$ of a generic algebraic hypersurface $X\subset \bP^{n+1}$ 
of degree $d\ge 2n+1$ is of general type for $n\geq 2$ (for surfaces
$X\subset\bP^3$, Geng Xu [Xu94] also obtained some refined information
for the genera of algebraic curves drawn in~$X$). The bound was
subsequently improved to $d\ge d_n=2n$ for $n\geq 5$ by Pacienza [Pac04].
That the same bound $d_n$ holds for Kobayashi hyperbolicity would then
be a consequence of the Green-Griffiths-Lang conjecture.
By these results, one can hope in the compact case that the optimal
bound $d_n$ is $d_1=4$, $d_n=2n+1$
for $n=2,3,4$ and $d_n=2n$ for $n\geq 5$. The case of complements
$\bP^n\ssm H$ (the so-called ``logarithmic case'') is a priori
somewhat easier to deal with: in fact, on can then exploit the fact that
the hyperbolicity of
the hypersurface $X=\{w^d=P(z)\}\subset\bP^{n+1}$ implies the hyperbolicity
of the complement $\bP^n\ssm H$, when $H=\{P(z)=0\}$.
Pacienza and Rousseau [PaRo07] proved that for $H$ very general
of degree $d\geq 2n+2-k$, any $k$-dimensional log-subvariety $(Y,D)$
of $(\bP^n,H)$ is of log-general type, i.e.\ any log-resolution
$\mu:\widetilde Y\to Y$ of $(Y,D)$ has a big log-canonical bundle
$K_{\widetilde Y}(\mu^*D)$.
\medskip

One of the early important result in the direction of Conjecture 0.2 is the
proof of the Bloch conjecture, as proposed by Bloch [Blo26a] and Ochiai
[Och77]: this is the special case of the conjecture when the
irregularity of $X$ satisfies $q=h^0(X,\Omega^1_X)>\dim X$. Various
solutions have then been obtained in fundamental papers of Noguchi
[Nog77a,~Nog81a,~Nog81b], Kawamata [Kaw80], Green-Griffiths [GrGr80],
McQuillan [McQ96],
and the book of Noguchi-Winkelmann [NoWi13], by means of different
techniques. Especially, assuming X to be of (log-) general type, it is
now known by [NWY07, NWY13] and [LuWi12] 
that if the (log-) irregularity is $q \geq\dim X$, then
no entire curve $f :\bC \to X$ has a Zariski dense image, and the GGL
conjecture holds in the compact (i.e.\ non logarithmic) case.
In the case of complex surfaces, major
progress was achieved by Lu, Miyaoka and Yau [LuYa90], [LuMi95, 96], [Lu96];
McQuillan [McQ96] extended these results to the case of all surfaces
satisfying $c_1^2>c_2$, in a situation where there are many symmetric 
differentials,
e.g.\ sections of $H^0(X,S^mT^*_X\otimes\cO(-1))$, $m\gg 1$ (cf.\ also 
[McQ99], [DeEG00] for applications to hyperbolicity). A more recent
result is the deep statement due to Diverio, Merker and Rousseau [DMR10],
confirming Conjecture~0.2 when  $X\subset\bP^{n+1}$ is a generic
non-singular hypersurface of 
sufficiently large degree~$d\ge 2^{n^5}$ (cf.~\S10); in the
case $n=2$ of surfaces in $\bP^3$, we are here in the more difficult 
situation where symmetric differentials do not exist (we have $c_1^2<c_2$
in this case). Conjecture~0.2 was also considered by S.~Lang
[Lang86, Lang87] in view of arithmetic counterparts of the
above geometric statements.
\medskip

Although these optimal conjectures are still unsolved at present, 
substantial progress was achieved in the meantime, for a large part via the
technique of producing jet differentials. This is done either by direct
calculations or by various indirect methods: Riemann-Roch
calculations, vanishing theorems~... Vojta [Voj87] and McQuillan [McQ98] 
introduced the ``diophantine approximation'' method, which was soon 
recognized to be an important tool in the study of holomorphic foliations, 
in parallel with Nevanlinna theory and the construction of Ahlfors currents.
Around 2000, Siu [Siu02, 04] showed that generic hyperbolicity results
in the direction of the Kobayashi conjecture could be investigated by
combining the algebraic techniques of Clemens, Ein and Voisin with the
existence of certain ``vertical'' meromorphic vector fields on the jet
space of the universal hypersurface of high degree; these vector
fields are actually used to differentiate the global sections of the
jet bundles involved, so as to produce new sections with a better
control on the base locus.  Also, during the years 2007--2010, it 
was realized [Dem07a, 07b, Dem11] that
holomorphic Morse inequalities could be used to prove the existence
of jet differentials; in 2010, Diverio, Merker and Rousseau [DMR10] were 
able in that way to prove the
Green-Griffiths conjecture for generic hypersurfaces of high degree in
projective space, e.g.\ for $d\ge 2^{n^5}$ -- their proof 
makes an essential use of Siu's
differentiation technique via meromorphic vector fields, as improved
by P\u{a}un [Pau08] and Merker [Mer09] in 2008. The present study 
will be focused on the holomorphic Morse inequality technique; as
an application, a partial answer to the Kobayashi and Green-Griffiths-Lang
conjecture can be obtained in a very wide context~: the basic general
result achieved in [Dem11] consists of showing that for every
projective variety of general type~$X$, there exists a global
algebraic differential operator $P$ on $X$ (in fact many such
operators $P_j$) such that every entire curve $f:\bC\to X$ must
satisfy the differential equations $P_j(f;f',\ldots,f^{(k)})=0$. One
also recovers from there the result of Diverio-Merker-Rousseau on the
generic Green-Griffiths conjecture (with an even better bound
asymptotically as the dimension tends to infinity), as well as a
result of Diverio-Trapani [DT10] on the hyperbolicity of
generic 3-dimensional hypersurfaces in $\bP^4$. Siu [Siu04, Siu15] has 
introduced a more explicit but more computationally involved approach
based on the use of ``slanted vector fields'' on jet spaces, extending
ideas of Clemens [Cle86] and Voisin [Voi96] (cf.\ section 10 for
details); [Siu15] explains how this strategy can be used
to assert the Kobayashi conjecture for $d\ge d_n$, with a very 
large bound and non-effective bound $d_n$ instead of~$2n+1$.
\medskip

As we will see here, it is useful to work in a more general context 
and to consider the category of directed varieties. When the problems
under consideration are birationally invariant, as is the case of the
Green-Griffiths-Lang conjecture, varieties can be replaced by non
singular models; for this reason, we will mostly restrict ourselves
to the case of non-singular varieties in the rest of the
introduction. A~{\em directed projective
manifold} is a pair $(X,V)$ where $X$ is a projective manifold equipped
with an analytic linear subspace $V\subset T_X$, i.e.\ a closed
irreducible complex analytic subset $V$ of the total space of~$T_X$,
such that each fiber $V_x=V\cap T_{X,x}$ is a complex vector
space. If $X$ is not connected, $V$ should rather be assumed to be
irreducible merely over each connected component of $X$, but we will
hereafter assume that our manifolds are connected. A morphism
$\Phi:(X,V)\to(Y,W)$ in the category of directed manifolds is an
analytic map $\Phi:X\to Y$ such that $\Phi_*V\subset W$. We refer to
the case $V=T_X$ as being the {\it absolute case}, and to the case
$V=T_{X/S}=\Ker d\pi$ for a fibration $\pi:X\to S$, as being the 
{\it relative case}; $V$ may also be taken to be the tangent space 
to the leaves of a singular analytic foliation on~$X$, or maybe even 
a non-integrable linear subspace of $T_X$.
We are especially interested in {\it entire curves} that are tangent to $V$,
namely non-constant holomorphic morphisms $f:(\bC,T_\bC)\to (X,V)$ of directed manifolds. In the absolute case, these are just arbitrary entire curves $f:\bC\to X$.

\claim 0.4. Generalized GGL conjecture|Let $(X,V)$ be a 
projective directed manifold. Define the entire curve locus of $(X,V)$ 
to be the Zariski closure of the locus of entire curves tangent to $V$, i.e.
$$
\ECL(X,V)=\overline{\mathop{\kern-40pt\bigcup}_{f:(\bC,T_\bC)\to(X,V)}f(\bC)}^{\rm Zar}.
$$
Then, if $(X,V)$ is of general type in the sense that the canonical sheaf 
sequence $K_V^\bullet$ is big $($cf.\ Prop~$2.11$ below$)$, 
$Y=\ECL(X,V)$ is a proper algebraic subvariety $Y\subsetneq X$.
\endclaim

\noindent [We will say that $(X,V)$ is {\it Brody hyperbolic} if 
$\ECL(X,V)=\emptyset\,$; by Brody's reparametrization technique, 
this is equivalent to Kobayashi hyperbolicity whenever
$X$ is compact.]
\medskip

In case $V$ has no singularities, the {\it canonical sheaf} $K_V$ is 
defined to be $(\det(\cO(V)))^*$ where $\cO(V)$ is the sheaf of holomorphic
sections of $V$, but in general this naive definition would not work.
Take for instance a generic pencil of
elliptic curves $\lambda P(z)+\mu Q(z)=0$ of degree $3$ in $\bP_\bC^2$,
and the linear space $V$ consisting of the tangents to the fibers
of the rational map $\bP_\bC^2\merto\bP_\bC^1$
defined by $z\mapsto Q(z)/P(z)$. Then $V$ is given by
$$
0\longrightarrow \cO(V)\longrightarrow\cO(T_{\bP_\bC^2})
\build\llra{8ex}^{PdQ-QdP}_{}\cO_{\bP^2_\bC}(6)\otimes\cJ_S
\longrightarrow 0
$$
where $S=\Sing(V)$ consists of the 9 points
$\{P(z)=0\}\cap\{Q(z)=0\}$, and $\cJ_S$ is the corresponding ideal
sheaf of~$S$. Since $\det(\cO(T_{\bP^2}))=\cO(3)$, we see that
$(\det(\cO(V))^*=\cO(3)$ is ample, thus the generalized GGL conjecture 0.4
would not have a
positive answer (all leaves are elliptic or singular rational curves
and thus covered by entire curves). An even more ``degenerate''
example is obtained with a generic pencil of conics, in which case
$(\det(\cO(V))^*=\cO(1)$ and $\# S=4$.
\medskip

If we want to get a positive answer to the generalized GGL conjecture 0.4, it is
therefore indispensable to give a definition of $K_V$ that incorporates
in a suitable way the singularities of $V\,;$ this will be done in
Def.~2.12 (see also Prop.~2.11). The goal is then to give a positive answer 
to the conjecture under some possibly more restrictive conditions for the 
pair $(X,V)$. These conditions will be expressed in terms of the tower
of Semple jet bundles
$$
(X_k,V_k)\to(X_{k-1},V_{k-1})\to\ldots\to(X_1,V_1)\to(X_0,V_0):=(X,V)
\leqno(0.5)
$$
which we define more precisely in Section~1, following [Dem95]. 
It is constructed inductively by setting $X_k=P(V_{k-1})$
(projective bundle of {\it lines} of $V_{k-1}$), and all $V_k$ have the 
same rank $r=\rank V$, so that $\dim X_k=n+k(r-1)$ where $n=\dim X$.
Entire curve loci have their counterparts for all stages of the Semple tower,
namely, one can define
$$
\ECL_k(X,V)=\overline{\mathop{\kern-40pt\bigcup}_{f:(\bC,T_\bC)\to(X,V)}f_{[k]}(\bC)}^{\rm Zar},\leqno(0.6)
$$
where $f_{[k]}:(\bC,T_\bC)\to(X_k,V_k)$ is the $k$-jet of $f$.
These are by definition algebraic subvarieties of $X_k$,
and if we denote by $\pi_{k,\ell}:X_k\to X_\ell$ the natural
projection from $X_k$ to $X_\ell$, $0\le\ell\le k$, we get immediately
$$
\pi_{k,\ell}(\ECL_k(X,V))=\ECL_\ell(X,V),\qquad \ECL_0(X,V)=\ECL(X,V).
\leqno(0.7)
$$
Let $\cO_{X_k}(1)$ be the tautological line bundle over $X_k$ associated
with the projective structure. We define the $k$-stage Green-Griffiths 
locus of $(X,V)$ to be
$$
\GG_k(X,V)=\overline{(X_k\smallsetminus\Delta_k)\cap
\bigcap_{m\in\bN}\left(\hbox{base locus of }\cO_{X_k}(m)\otimes 
\pi_{k,0}^*A^{-1}\right)}
\leqno(0.8)
$$
where $A$ is any ample line bundle on $X$ and $\Delta_k=\bigcup_{2\le \ell\le k}
\pi_{k,\ell}^{-1}(D_\ell$) is the union of ``vertical divisors'' 
(see (6.9) and (7.17); 
the vertical divisors play no role and have to be removed in this context;
for this, one uses the fact that $f_{[k]}(\bC)$ is not contained in any
component of $\Delta_k$, cf.~[Dem95]).
Clearly, $\GG_k(X,V)$ does not depend on the choice of~$A$.

\claim 0.9. Basic vanishing theorem for entire curves|
Let $(X,V)$ be an arbitrary directed variety with $X$ non-singular, and
let $A$ be an ample line bundle on~$X$. Then any entire 
curve $f:(\bC,T_\bC)\to(X,V)$ satisfies the differential equations
$P(f;f',\ldots,f^{(k)})=0$ arising from sections 
$\sigma\in H^0(X_k,\cO_{X_k}(m)\otimes\pi_{k,0}^*A^{-1})$. As a 
consequence, one has
$$
\ECL_k(X,V)\subset \GG_k(X,V).
$$
\endclaim

The main argument goes back to [GrGr80]. We will give here a complete proof 
of Theorem~0.9, based only on the arguments [Dem95], namely on 
the Ahlfors-Schwarz lemma (the alternative proof given in [SiYe96b] uses 
Nevanlinna theory and is analytically more involved). By (0.7) and Theorem~0.9
we infer that
$$
\ECL(X,V)\subset\GG(X,V),
\leqno(0.10)
$$
where $\GG(X,V)$ is the global Green-Griffiths locus of $(X,V)$ defined by
$$
\GG(X,V)=\bigcap_{k\in\bN}\pi_{k,0}\left(\GG_k(X,V)\right).
\leqno(0.11)
$$
The main result of [Dem11] (Theorem~2.37 and Cor.~3.4) implies the 
following useful information:

\claim 0.12. Theorem|Assume that $(X,V)$ is of ``general type'', i.e.\ that 
the pluricanonical sheaf sequence $K_V^\bullet$ is big on $X$.
Then there exists an integer $k_0$ such that $\GG_k(X,V)$ is a proper algebraic
subset of $X_k$ for $k\ge k_0$ $[\,$though $\pi_{k,0}(\GG_k(X,V))$ might still 
be equal to $X$ for all $k\,]$.
\endclaim

In fact, if $F$ is an invertible sheaf on $X$ such that
$K_V^\bullet\otimes F$ is big (cf.\ Prop.~2.11), the probabilistic estimates 
of [Dem11, Cor.~2.38 and Cor.~3.4] produce global sections of
$$
\cO_{X_k}(m)\otimes\pi_{k,0}^*\cO\Big({m\over kr}\Big(
1+{1\over 2}+\ldots+{1\over k}\Big)F\Big)\leqno(0.13)
$$
for $m\gg k\gg 1$. The (long and elaborate) proof uses a curvature computation
and singular holomorphic Morse inequalities to show that the line bundles 
involved in (0.11) are big on $X_k$ for $k\gg 1$. One applies this to 
$F=A^{-1}$ with $A$ ample on $X$ to produce sections and conclude that 
$\GG_k(X,V)\subsetneq X_k$. 
\medskip

Thanks to (0.10), the GGL conjecture is satisfied whenever 
$\GG(X,V)\subsetneq X$. By [DMR10], this happens for instance in the
absolute case when $X$ is a generic hypersurface of 
degree~\hbox{$d\ge 2^{n^5}$ in $\bP^{n+1}$} (see also [Pau08]
for better bounds in low dimensions, and [Siu02, Siu04]).
However, as already mentioned in [Lang86], very simple examples show that 
one can have $\GG(X,V)=X$ even when $(X,V)$ is of general type, and 
this already occurs in the absolute case as soon as $\dim X\ge 2$. 
A typical example is a product of directed manifolds
$$(X,V)=(X',V')\times(X'',V''),\qquad
V=\pr^{\prime\,*}V'\oplus\pr^{\prime\prime\,*}V''.\leqno(0.14)$$
The absolute case $V=T_X$, $V'=T_{X'}$, $V''=T_{X''}$ on a product of curves
is the simplest instance. It is then easy to check that $\GG(X,V)=X$, cf.\
(3.2). Diverio and Rousseau [DR15] have given many more such examples, 
including the case of indecomposable varieties $(X,T_X)$, e.g.\ Hilbert modular 
surfaces, or more generally compact quotients of bounded symmetric domains
of rank${}\ge 2$. 
\medskip

The problem here is the failure of some sort of stability 
condition that is introduced in Remark~11.10. This leads us to make
the assumption that the directed pair $(X,V)$ is
{\it strongly of general type}: by this, we mean that the
induced directed structure $(Z,W)$ on each non-vertical subvariety
$Z\subset X_k$ that projects onto $X$ either has $\rank W=0$ or is
of general type modulo $X_\bu\to X$, in the sense that
$K_{W_\ell}^\bullet\otimes\cO_{Z_\ell}(p)_{\restriction Z_\ell}$ is big for some
stage of the Semple tower of $(Z,W)$ and some $p\geq 0$
(see~Section~11 for details -- one may
have to replace $Z_\ell$ by a suitable modification).
Our main result can be stated as follows:

\claim 0.15. Theorem (partial solution to the generalized GGL conjecture)|Let 
$(X,V)$ be a directed pair that is strongly of general type.
Then the Green-Griffiths-Lang conjecture holds true for~$(X,V)$, namely
$\ECL(X,V)$ is a proper algebraic subvariety of $X$.
\endclaim

The proof proceeds through a complicated induction on $n=\dim X$ and 
$k=\rank V$, which is the main reason why we have to introduce
directed varieties, even in the absolute case. An interesting feature of 
this result is that the conclusion
on $\ECL(X,V)$ is reached without having to know anything about 
the Green-Griffiths locus $\GG(X,V)$, even a posteriori. Nevertheless,
this is not yet enough to confirm the GGL conjecture. Our hope is
that pairs $(X,V)$ that are of general type without being strongly
of general type -- and thus exhibit some sort of ``jet-instability'' -- can
be investigated by different methods, e.g.\ by the diophantine
approximation techniques of McQuillan [McQ98]. However, Theorem~0.15
provides a sufficient criterion for Kobayashi hyperbolicity [Kob70],
thanks to the following concept of algebraic 
jet-hyperbolicity.

\claim 0.16. Definition|A directed variety $(X,V)$ will be said to be 
algebraically jet-hyperbolic if the induced directed variety structure
$(Z,W)$ on every non-vertical irreducible algebraic variety $Z$
of~$X_k$ with $\rank W\ge 1$ is such that
$K_{W_\ell}^\bullet\otimes\cO_{Z_\ell}(p)_{\restriction Z_\ell}$ is big for some
stage of the Semple tower of $(Z,W)$ and some $p\geq 0$
$[$possibly after taking a suitable modification of $Z_\ell\,$;
see Sections~$11, 12$ for the definition of induced directed 
structures and further details$]$. We also say that a projective
manifold $X$ is algebraically jet-hyperbolic if $(X,T_X)$ is.
\endclaim

In this context, Theorem 0.15 yields the following connection between
algebraic jet-hyperbolicity and the analytic concept of Kobayashi 
hyperbolicity.

\claim 0.17. Theorem|Let $(X,V)$ be a directed variety structure on a
projective manifold $X$. Assume that $(X,V)$ is algebraically jet-hyperbolic.
Then $(X,V)$ is Kobayashi hyperbolic.
\endclaim

The following conjecture would then make a bridge between these theorems
and the GGL and Kobayashi conjectures.

\claim 0.18. Conjecture|Let $X\subset\bP^{n+c}$ be a complete intersection
of hypersurfaces of respective degrees $d_1,\ldots, d_c$, $\codim X=c$.
\plainitem{\rm(a)} If $X$ is non-singular and of general type, i.e.\ if
$\sum d_j\ge n+c+2$, then $X$ is in fact strongly of general type.
\vskip0pt
\plainitem{\rm(b)} If $X$ is $($very$)$ generic and
$\sum d_j\ge 2n+c$, then $X$ is algebraically jet-hyperbolic.
\vskip0pt
\endclaim

Since Conjecture 0.18 only deals with algebraic statements, our hope is that
a proof can be obtained through a suitable deepening of the techniques
introduced by Clemens, Ein, Voisin and Siu. Under the slightly 
stronger condition $\sum d_j\ge 2n+c+1$, Voisin showed indeed 
that every subvariety $Y\subset X$ is of general type, if $X$ is generic.
To prove the Kobayashi conjecture in its optimal incarnation, we would 
need to show that such $Y$'s are strongly of general type.
\medskip

In the direction of getting examples of low degrees,
Dinh Tuan Huynh [DTH16a] showed that there are families of
hyperbolic hypersurfaces of degree $2n+2$ in $\bP^{n+1}$ for
$2\leq n\leq 5$, and in [DTH16b] he showed that certain small
deformations (in Euclidean topology) of a union of
\hbox{$\lceil(n+3)^2/4\rceil$} hyperplanes in general position
in $\bP^{n+1}$ are hyperbolic. In [Ber18], G.~B\'erczi stated a positivity
conjecture for Thom polynomials of Morin singularities (see also [BeSz12]),
and announ\-ced that it would imply a polynomial bound $d_n=2\,n^9+1$ for the
generic hyperbolicity of hypersurfaces.
By using the ``technology'' of Semple towers and following new ideas
introduced by D.~Brotbek [Brot17]
and Ya Deng [Deng16], we prove here the following effective 
(although non-optimal) version of the Kobayashi conjecture on
generic hyperbolicity.

\claim 0.19. Theorem|Let $Z$ be a projective $(n+1)$-dimensional manifold
and $A$ a very ample line bundle on $Z$. Then, for a general section
$\sigma\in H^0(Z,A^d)$ and $d\geq d_n$, the hypersurface
$X_\sigma=\sigma^{-1}(0)$ is Kobayashi hyperbolic and, in fact, satisfies
the stronger property of being algebraically jet hyperbolic.
The bound $d_n$ for the degree can be taken to be
$d_n:=\lfloor\frac{1}{3}(en)^{2n+2}\rfloor$.
\endclaim

I would like to thank Damian Brotbek, Ya Deng, Simone Diverio, Gianluca
Pacienza, Erwan Rousseau, Mihai P\u{a}un and Mikhail Zaidenberg for
very stimulating discussions on these questions. These notes also owe a lot
to their work. I~also with to thank the unknown referees for a large number of
corrections and very useful suggestions.

\section{Basic hyperbolicity concepts}

\plainsubsection 1.A. Kobayashi hyperbolicity|
We first recall a few basic facts concerning the concept of
hyperbolicity, according to S.~Kobayashi [Kob70, Kob76, Kob98]. 
Let $X$ be a complex space. Given two points $p,q\in X$, let us
consider a {\em chain  of analytic disks} from $p$ to $q$, that is a 
sequence of holomorphic maps $f_0,f_1,\dots,f_k:\Delta\to X$ from
the unit disk $\Delta=D(0,1)\subset\bC$ to $X$, together with
pairs of points $a_0,b_0,\dots,a_k,b_k$ of $\Delta$ such that
$$
p=f_0(a_0),\quad q=f_k(b_k),\quad f_i(b_i)=f_{i+1}(a_{i+1}),\qquad 
i=0,\dots,k-1.
$$
Denoting this chain by $\alpha$, we define its length $\ell(\alpha)$ to be
$$
\ell(\alpha)=d_P(a_1,b_1)+\cdots+d_P(a_k,b_k)\leqno(1.1')
$$ 
where $d_P$ is the Poincar\'e distance on $\Delta$, and the {\em Kobayashi 
pseudodistance} $d^K_X$ on $X$ to be
$$
d^K_X(p,q)=\inf_{\alpha}\ell(\alpha).\leqno(1.1'')
$$
A {\em Finsler metric} (resp.\ {\em pseudometric}\/) on a vector bundle $E$
is a homogeneous positive (resp.\ nonnegative) function $N$ on the
total space $E$, that is,
$$
N(\lambda\xi)=|\lambda|\,N(\xi)\qquad
\hbox{for all $\lambda\in\bC$ and $\xi\in E$,}
$$
but in general $N$ is not assumed to be subadditive 
(i.e.\ convex) on the fibers of~$E$. A~Finsler (pseudo-)metric on $E$ 
is thus nothing but a Hermitian (semi-)norm
on the tautological line bundle $\cO_{P(E)}(-1)$ of lines of~$E$ over the
projectivized bundle $Y=P(E)$. 
The {\em Kobayashi-Royden infinitesimal pseudometric} on $X$ is the
Finsler pseudometric on the tangent bundle $T_X$ defined by
$$
\bfk_X(\xi)=\inf\big\{\lambda>0\,;\,\exists f:\Delta\to X,\,f(0)=x,\,
\lambda f'(0)=\xi\big\},\qquad x\in X,~\xi\in T_{X,x}.\leqno(1.2)
$$
Here, if $X$ is not smooth at $x$, we take $T_{X,x}=(\gm_{X,x}/\gm^2_{X,x})^*$ 
to be the Zariski tangent space, i.e.\ the tangent space of a minimal smooth
ambient vector space containing the germ $(X,x)$; all tangent vectors
may not be reached by analytic disks and in those cases we put 
$\bfk_X(\xi)=+\infty$. When $X$ is a smooth manifold, it follows from the
work of H.L.~Royden ([Roy71], [Roy74]) that $\bfk_X$ is upper-continuous
on $T_X$ and that $d^K_X$ is the integrated
pseudodistance associated with the pseudometric, i.e.
$$
d^K_X(p,q)=\inf_\gamma\int_\gamma\bfk_X(\gamma'(t))\,dt,
$$
where the infimum is taken over all piecewise smooth curves joining $p$ to 
$q\,$; in the case of complex spaces, a similar formula holds, involving
jets of analytic curves of arbitrary order, cf.\ S.~Venturini [Ven96].
When $X$ is a non-singular projective variety, it has been shown in
[DeLS94] that the Kobayashi pseudodistance and the Kobayashi-Royden
infinitesimal pseudometric can be computed by looking only at analytic
disks that are contained in algebraic curves.

\claim 1.3.~Definition|A complex space $X$ is said to be
{\em hyperbolic} $($in the sense of Kobayashi$)$ if $d^K_X$ is actually a
distance, namely if $d^K_X(p,q)>0$ for all pairs of distinct points
$(p,q)$ in~$X$.
\endclaim

When $X$ is hyperbolic, it is interesting to investigate when the Kobayashi
metric is complete: one then says that $X$ is a {\em complete hyperbolic} space.
However, we will be mostly concerned with compact spaces here, so completeness
is irrelevant in that case.

Another important property is the {\em monotonicity} of the Kobayashi
pseudometric with respect to holomorphic mappings. In fact, if $\Phi:X\to Y$
is a holomorphic map, it is easy to see from the definition that
$$
d^K_Y(\Phi(p),\Phi(q))\le d^K_X(p,q),\qquad\hbox{for all $p,q\in X$.}\leqno(1.4)
$$
The proof merely consists of taking the composition $\Phi\circ f_i$ for all
chains of analytic disks connecting $p$ and $q$ in~$X$. Clearly the Kobayashi 
pseudodistance $d^K_\bC$ on $X=\bC$ is identically zero, as one can see by looking at
arbitrarily large analytic disks $\Delta\to\bC$, $t\mapsto\lambda t$. Therefore,
if there is an {\em entire curve} $\Phi:\bC\to X$, namely a non-constant
holomorphic map defined on the whole complex plane~$\bC$, then by monotonicity
$d^K_X$ is identically zero on the image $\Phi(\bC)$ of the curve, and therefore 
$X$ cannot be hyperbolic. When $X$ is hyperbolic, it follows that $X$ cannot contain
rational curves $C\simeq\bP^1$, or elliptic curves $\bC/\Lambda$, or more generally
any non-trivial image $\Phi:W=\bC^p/\Lambda\to X$ of a $p$-dimensional complex torus
(quotient of $\bC^p$ by a lattice). The only case where hyperbolicity is easy to assess
is the case of curves $(\dim_\bC X=1)$. 

\claim 1.5.~Case of complex curves|Up to bihomorphism, any smooth complex curve
$X$ belongs to one $($and only one$)$ of the following three types$\,:$
\vskip-2pt
\plainitem{\rm(a)} $($rational curve$)$ $X\simeq\bP^1\,;$
\vskip-2pt
\plainitem{\rm(b)} $($parabolic type$)$ $\widehat X\simeq\bC$, $X\simeq \bC,\,\bC^*$ or $X\simeq
\bC/\Lambda$ $($elliptic curve$)\,;$
\vskip-2pt
\plainitem{\rm(c)} $($hyperbolic type$)$ $\widehat X\simeq\Delta$. All compact curves
$X$ of genus $g\ge 2$ enter in this category, as well as $X=\bP^1\ssm\{a,b,c\}\simeq
\bC\ssm\{0,1\}$, or $X=\bC/\Lambda\ssm\{a\}$ $($elliptic curve minus one
point$)$.\vskip0pt
\endclaim

In fact, as the disk is simply connected,
every holomorphic map $f:\Delta\to X$ lifts to the universal cover 
$\widehat f:\Delta\to\widehat X$, so that $f=\rho\circ\widehat f$ where
$\rho:\widehat X\to X$ is the projection map, and the conclusions (a,b,c)
follow easily from the Poincar\'e-Koebe uniformization theorem: every simply
connected Riemann surface is biholomorphic to $\bC$, the unit disk
$\Delta$ or the complex projective line $\bP^1$. 

In some rare cases, the one-dimensional case can be used to study the 
case of higher 
dimensions. For instance, it is easy to see by looking at projections
that the Kobayashi pseudodistance on a product 
$X\times Y$ of complex spaces is given by
$$
\leqalignno{
&d^K_{X\times Y}((x,y),(x',y'))=\max\big(d^K_X(x,x'),d^K_Y(y,y')\big),&(1.6)\cr
&\bfk_{X\times Y}(\xi,\xi')=\max\big(\bfk_X(\xi),\bfk_Y(\xi')\big),&(1.6')\cr}
$$
and from there it follows that a product of hyperbolic spaces is hyperbolic. As a
consequence $(\bC\ssm\{0,1\})^2$, which is also a complement of five lines in $\bP^2$,
is hyperbolic. 

\plainsubsection 1.B. Brody criterion for hyperbolicity|
Throughout this subsection, we assume that $X$ is a complex
manifold. In this context, we have the following well-known result
of Brody [Bro78]. Its main interest is to relate hyperbolicity to
the non-existence of entire curves.

\claim 1.7.~Brody reparametrization lemma|Let $\omega$ be a Hermitian
metric on~$X$ and let $f:\Delta\to X$ be a holomorphic map. For every
$\varepsilon>0$, there exists a radius $R\ge(1-\varepsilon)\|f'(0)\|_\omega$
and a homographic transformation $\psi$ of the disk $D(0,R)$ onto
$(1-\varepsilon)\Delta$ such that
$$
\|(f\circ\psi)'(0)\|_\omega=1,\qquad \|(f\circ\psi)'(t)\|_\omega
\le{1\over 1-|t|^2/R^2}\quad\hbox{for every $t\in D(0,R)$.}
$$
\endclaim

\plainproof. Select $t_0\in\Delta$ such that $(1-|t|^2)\|f'((1-\varepsilon)t)\|_\omega$
reaches its maximum for $t=t_0$. The reason for this choice is that
$(1-|t|^2)\|f'((1-\varepsilon)t)\|_\omega$
is the norm of the differential $f'((1-\varepsilon)t):T_\Delta\to T_X$
with respect to the Poincar\'e metric $|dt|^2/(1-|t|^2)^2$ on $T_\Delta$,
which is conformally invariant under $\Aut(\Delta)$. One then adjusts
$R$ and $\psi$ so that $\psi(0)=(1-\varepsilon)t_0$ and
$|\psi'(0)|\,\|f'(\psi(0))\|_\omega=1$. As
$|\psi'(0)|={1-\varepsilon\over R}(1-|t_0|^2)$, the only possible
choice for $R$ is
$$
R=(1-\varepsilon)(1-|t_0|^2)\|f'(\psi(0))\|_\omega\ge (1-\varepsilon)
\|f'(0)\|_\omega.
$$
The inequality for $(f\circ\psi)'$ follows from the fact that the Poincar\'e
norm is maximum at the origin, where it is equal to $1$ by the choice
of~$R$.\qed

Using the Ascoli-Arzel\`a theorem we obtain immediately:

\claim 1.8.~Corollary {\rm(Brody)}|Let $(X,\omega)$ be a compact complex Hermitian manifold.
Given a sequence of holomorphic
mappings \hbox{$f_\nu:\Delta\to X$} such that
$\lim\|f_\nu'(0)\|_\omega=+\infty$, one can find a sequence of homographic
transformations $\psi_\nu:D(0,R_\nu)\to(1-1/\nu)\Delta$ with
$\lim R_\nu=+\infty$, such that, after passing possibly to a subsequence,
$(f_\nu\circ\psi_\nu)$ converges uniformly on every compact subset of $\bC$
towards a non-constant holomorphic map $g:\bC\to X$ with $\|g'(0)\|_\omega=1$
and $\sup_{t\in\bC}\|g'(t)\|_\omega\le 1$.
\endclaim

An entire curve $g:\bC\to X$ such that $\sup_\bC\|g'\|_\omega=M<+\infty$ is called a
{\em Brody curve}; this concept does not depend on the choice of $\omega$ when $X$ is 
compact, and one can always assume $M=1$ by rescaling the parameter~$t$. 

\claim 1.9.~Brody criterion|Let $X$ be a compact complex manifold. The following properties
are equivalent$\,:$
\plainitem{\rm(a)} $X$ is hyperbolic$\,;$
\plainitem{\rm(b)} $X$ does not possess any entire curve $f:\bC\to X\,;$
\plainitem{\rm(c)} $X$ does not possess any Brody curve $g:\bC\to X\,;$
\plainitem{\rm(d)} The Kobayashi infinitesimal metric $\bfk_X$ is uniformly bounded below, namely
$$
\bfk_X(\xi)\ge c\|\xi\|_\omega,\qquad c>0,
$$
for any Hermitian metric $\omega$ on $X$.\vskip0pt
\endclaim

\plainproof. (a)${}\Rightarrow{}$(b). If $X$ possesses an entire curve
$f:\bC\to X$, then by looking at arbitrary large disks
$D(0,R)\subset\bC$, it is easy to see that the Kobayashi distance of
any two points in $f(\bC)$ is zero, so $X$ is not
hyperbolic.\smallskip

\noindent
(b)${}\Rightarrow{}$(c). This is trivial.
\smallskip 

\noindent
(c)${}\Rightarrow{}$(d). If (d) does not hold, there exists a sequence
of tangent vectors $\xi_\nu\in T_{X,x_\nu}$ with
$\Vert\xi_\nu\|_\omega=1$ and $\bfk_X(\xi_\nu)\to 0$. By definition,
this means that there exists an analytic curve $f_\nu:\Delta\to X$
with $f(0)=x_\nu$ and
$\Vert f'_\nu(0)\Vert_\omega\ge (1-{1\over \nu})/\bfk_X(\xi_\nu)\to+\infty$.
One can then produce a Brody curve $g=\bC\to X$ by Corollary~1.8,
contradicting~(c).\smallskip

\noindent
(d)${}\Rightarrow{}$(a). In fact (d) implies after integrating that
$d^K_X(p,q)\ge c\,d_\omega(p,q)$ where $d_\omega$ is the geodesic
distance associated with~$\omega$, so $d^K_X$ must be non
degenerate.\qed

Notice also that if $f:\bC\to X$ is an entire curve such that
$\|f'\|_\omega$ is unbounded, one can apply the Corollary~1.8 to
$f_\nu(t):=f(t+a_\nu)$ where the sequence $(a_\nu)$ is chosen such
that $\|f'_\nu(0)\|_\omega= \|f(a_\nu)\|_\omega\to+\infty$. Brody's
result then produces reparametrizations
$\psi_\nu:D(0,R_\nu)\to D(a_\nu,1-1/\nu)$ and a Brody curve
$g=\lim f\circ\psi_\nu:\bC\to X$ such that $\sup\|g'\|_\omega=1$ and
$g(\bC)\subset\overline{f(\bC)}$. It may happen that the image
$g(\bC)$ of such a limiting curve is disjoint from $f(\bC)$. In fact
Winkelmann [Win07] has given a striking example, actually a projective
$3$-fold $X$ obtained by blowing-up a 3-dimensional abelian
variety~$Y$, such that every Brody curve $g:\bC\to X$ lies in the
exceptional divisor $E\subset X\,$; however, entire curves
$f:\bC\to X$ can be dense, as one can see by taking $f$ to be the
lifting of a generic complex line embedded in the abelian
variety~$Y$. For further precise information on the localization of
Brody curves, we refer the reader to the remarkable results of
[Duv08].

The absence of entire holomorphic curves in a given complex manifold
is often referred to as {\em Brody hyperbolicity}. Thus, in the
compact case, Brody hyperbolicity and Kobayashi hyperbolicity
coincide (but Brody hyperbolicity is in general a strictly weaker
property when $X$ is non-compact).

\plainsubsection 1.C. Geometric applications|

We give here two immediate consequences of the Brody criterion: the
openness property of hyperbolicity and a hyperbolicity criterion for
subvarieties of complex tori.  By definition, a {\em holomorphic
  family} of compact complex manifolds is a holomorphic proper
submersion $\cX\to S$ between two complex manifolds.

\claim 1.10. Proposition|Let $\pi:\cX\to S$ be a holomorphic family of
compact complex manifolds. Then the set of $s\in S$ such that the
fiber $X_s=\pi^{-1}(s)$ is hyperbolic is open in the Euclidean
topology.
\endclaim

\plainproof. Let $\omega$ be an arbitrary Hermitian metric on $\cX$,
$(X_{s_\nu})_{s_\nu\in S}$ a sequence of non-hyperbolic fibers, and
$s=\lim s_\nu$. By the Brody criterion, one obtains a sequence of
entire maps $f_\nu:\bC\to X_{s_\nu}$ such that
$\|f'_\nu(0)\|_\omega=1$ and $\|f'_\nu\|_\omega\le 1$. Ascoli's theorem
shows that there is a subsequence of $f_\nu$ converging uniformly to a
limit $f:\bC\to X_s$, with $\|f'(0)\|_\omega=1$. Hence
$X_s$ is not hyperbolic and the collection of non-hyperbolic fibers is
closed in $S$.\qed
\medskip

Consider now an $n$-dimensional complex torus~$W$, i.e.\ an additive quotient 
$W=\bC^n/\Lambda$,  where $\Lambda\subset\bC^n$ is a (cocompact) lattice.
By taking a composition of entire curves $\bC\to\bC^n$ with
the projection $\bC^n\to W$ we obtain an infinite dimensional space 
of entire curves in~$W$.

\claim 1.11.~Theorem|Let $X\subset W$ be a compact complex submanifold of a complex torus. Then $X$ is hyperbolic if and only if it does not contain any translate of a subtorus.
\endclaim

\plainproof. If $X$ contains some translate of a subtorus, then it contains lots of 
entire curves and so $X$ is not hyperbolic.

Conversely, suppose that $X$ is not hyperbolic. Then by the Brody
criterion there exists an entire curve $f:\bC\to X$ such that
$\|f'\|_\omega\le \|f'(0)\|_\omega=1$, where $\omega$ is the flat
metric on $W$ inherited from $\bC^n$. This means that any lifting
$\wt f=(\wt f,\dots,\wt f_\nu):\bC\to\bC^n$ is
such that
$$
\sum_{j=1}^n|f_j'|^2\le 1.
$$
Then, by Liouville's theorem, $\wt f'$ is constant and therefore $\wt f$
is affine linear. But then the closure of the image of $f$ is a
translate $a+H$ of a connected 
(possibly real) subgroup $H$ of $W$. We conclude that $X$ contains the analytic 
Zariski closure of $a+H$, namely $a+H^\bC$ where $H^\bC\subset W$ is the smallest 
closed complex subgroup of $W$ containing $H$.\qed
\medskip

\section{Directed manifolds}

\plainsubsection 2.A. Basic definitions concerning directed manifolds|
Let us consider a pair $(X,V)$ consisting of an $n$-dimensional complex
manifold $X$ equipped with a {\em linear subspace $V\subset T_X$}: if
we assume $X$ to be connected, this is by definition an irreducible
closed analytic subspace of the total space of $T_X$ such that each fiber
$V_x=V\cap T_{X,x}$ is a vector subspace of~$T_{X,x}$.
If $\cW\subset\Omega^1_X$ is the sheaf of $1$-forms vanishing
on~$V$, then $\cW$ is coherent (this follows from the direct image theorem
by looking at the proper morphism $P(V)\subset P(T_X)\to X$), and $V$
is locally defined by
$$
V_x=\{\xi\in T_{X,x}\,;\;\alpha_j(x)\cdot\xi=0,\,
1\leq j\leq N\},\quad \alpha_j\in H^0(U,\Omega^1_X),\quad x\in U,
$$
where $(\alpha_1,\ldots,\alpha_N)$, is a local family of generators of $\cW$
on a small open set $U$.
We can also associate to $V$ a coherent sheaf
$\cV:=\cW^\perp=\Hom(\Omega^1_X/\cW,\cO_X)\subset \cO(T_X)$, which
is a {\it saturated subsheaf}
of $\cO(T_X)$, i.e.\ such that $\cO(T_X)/\cV$ has no torsion; then 
$\cV$ is also reflexive, i.e.\ $\cV^{**}=\cV$. We will refer to 
such a pair as being a (complex) {\em directed manifold}, and we will in
general think of $V$ as a linear space (rather than considering the
associated saturated subsheaf $\cV\subset \cO(T_X)$).
A~morphism
$\Phi:(X,V)\to (Y,W)$ in the category of complex directed manifolds is
a holomorphic map such that $\Phi_*(V)\subset W$.

Here, the rank $x\mapsto \dim_\bC V_x$
is Zariski lower semi-continuous, and it may a priori jump.
The rank $r:=\rank(V)\in\{0,1,\ldots,n\}$ of $V$ is by definition the dimension of 
$V_x$ at a generic point. The dimension may be larger at non-generic points;
this happens e.g.\ on $X=\bC^n$ for the rank~$1$ linear space $V$
generated by the Euler vector field: $V_z=\bC\sum_{1\le j\le n} 
z_j{\partial\over\partial z_j}$ for $z\ne 0$, and $V_0=\bC^n$. 
Our philosophy is that directed manifolds are also useful to study the
``absolute case'', i.e.\ the case $V=T_X$, because there are certain functorial
constructions which are quite natural in the category of directed manifolds
(see e.g.\ \S$\,$5,~6,~7). We think of directed manifolds as a kind
of ``relative situation'', covering e.g.\ the case when $V$ is the
relative tangent space to a holomorphic map $X\to S$. It is important
to notice that the local sections of $\cV$ need not generate the fibers
of $V$ at singular points, as one 
sees already in the case of the Euler vector field when $n\ge 2$.
We also want to stress that no assumption need be made 
on the Lie bracket tensor $[\bu,\bu]:\cV\times \cV\to \cO(T_X)/\cV$, i.e.\ we do not 
assume any kind of integrability for~$\cV$ or $\cW$.

The singular set $\Sing(V)$ is by definition the set of points where
$\cV$ is not locally free, it can also be defined as the indeterminacy set of 
the (meromorphic) classifying map $\alpha:X\merto G_r(T_X)$, $z\mapsto
V_z$ to the Grassmannian of $r$ dimensional subspaces of~$T_X$. We thus have 
$V_{\restriction X\ssm \Sing(V)}=\alpha^*S$ where $S\to G_r(T_X)$ is the tautological 
subbundle of $G_r(T_X)$. The singular set $\Sing(V)$ is an
analytic subset of $X$ of codim${}\ge 2$, and hence $V$ is always a
holomorphic subbundle outside of codimension~$2$. Thanks to this remark,
one can most often treat linear spaces as vector bundles (possibly modulo
passing to the Zariski closure along $\Sing(V)$).

\plainsubsection 2.B. Hyperbolicity properties of directed manifolds|
\noindent
Most of what we have done in \S1 can be extended to the category of directed manifolds.

\claim 2.1.~Definition|Let $(X,V)$ be a complex directed manifold.
\plainitem{\rm(i)} The Kobayashi-Royden infinitesimal metric
of $(X,V)$ is the Finsler metric on $V$ defined for any $x\in X$ and
$\xi\in V_x$ by
$$
\bfk_{(X,V)}(\xi)=\inf\big\{\lambda>0\,;\,\exists f:\Delta\to X,\,f(0)=x,\,
\lambda f'(0)=\xi,\,f'(\Delta)\subset V\big\}.
$$
Here $\Delta\subset\bC$ is the unit disk and the map $f$ is an arbitrary
holomorphic map which is tangent to~$V$, i.e., such that $f'(t)\in V_{f(t)}$
for all $t\in\Delta$. We say that $(X,V)$ is infinitesimally hyperbolic if
$\bfk_{(X,V)}$ is positive definite on every fiber~$V_x$ and satisfies
a uniform lower bound $\bfk_{(X,V)}(\xi)\ge\varepsilon\|\xi\|_\omega$
in terms of any smooth Hermitian metric $\omega$ on $X$, when $x$
describes a compact subset of~$X$.
\plainitem{\rm(ii)} More generally, the Kobayashi-Eisenman infinitesimal
pseudometric of $(X,V)$ is the pseudometric defined on all decomposable
$p$-vectors $\xi=\xi_1\wedge\cdots\wedge\xi_p\in\Lambda^pV_x$, $1\le p\le
r=\rank(V)$, by
$$
\bfe^p_{(X,V)}(\xi)=\inf\big\{\lambda>0\,;\,\exists f:\bB_p\to X,\,f(0)=x,\,
\lambda f_*(\tau_0)=\xi,\,f_*(T_{\bB_p})\subset V\big\},
$$
where $\bB_p$ is the unit ball in $\bC^p$ and $\tau_0=\partial/\partial t_1
\wedge\cdots\wedge\partial/\partial t_p$ is the unit $p$-vector of $\bC^p$
at the origin. We say that $(X,V)$ is infinitesimally $p$-measure hyperbolic
if $\bfe^p_{(X,V)}$ is positive definite on every fiber~$\Lambda^pV_x$
and satisfies a locally uniform lower bound in terms of any smooth metric.
\vskip0pt
\endclaim

If $\Phi:(X,V)\to (Y,W)$ is a morphism of directed manifolds, it is
immediate to check that we have the monotonicity property
$$
\leqalignno{
&\bfk_{(Y,W)}(\Phi_*\xi)\le \bfk_{(X,V)}(\xi),
\qquad\forall\xi\in V,&(2.2)\cr
&\bfe^p_{(Y,W)}(\Phi_*\xi)\le \bfe^p_{(X,V)}(\xi),
\qquad\forall\xi=\xi_1\wedge\cdots\wedge\xi_p\in\Lambda^pV.&(2.2^p)\cr}
$$
The following proposition shows that virtually all reasonable definitions
of the hyperbolicity property are equivalent if $X$ is compact (in
particular, the additional assumption that there is locally uniform
lower bound for $\bfk_{(X,V)}$ is not needed). We merely say in that case
that $(X,V)$ is {\em hyperbolic}.

\claim 2.3.~Proposition|For an arbitrary directed manifold $(X,V)$,
the Kobayashi-Royden infinitesimal metric $\bfk_{(X,V)}$ is upper
semi-continuous on the total space of~$V$. If $X$ is compact, $(X,V)$ is
infinitesimally hyperbolic if and only if there are no entire
curves $g:\bC\to X$ tangent to~$V$. In that case, $\bfk_{(X,V)}$ is a
continuous $($and positive definite$)$ Finsler metric on~$V$.
\endclaim

\plainproof. The proof is almost identical to the standard proof for $\bfk_X$,
for which we refer to Royden [Roy71, Roy74]. One of the main ingredients is that
one can find a Stein neighborhood of the graph of any analytic disk (thanks to
a result of [Siu76], cf.\ also [Dem90a] for more general results). This allows
to obtain ``free'' small deformations of any given analytic disk, as there are many 
holomorphic vector fields on a Stein manifold.\qed

Another easy observation is that the concept of $p$-measure hyperbolicity gets 
weaker and weaker as $p$ increases (we leave it as an exercise to the reader, this
is mostly just linear algebra).

\claim 2.4. Proposition|If $(X,V)$ is $p$-measure hyperbolic, then it is
$(p+1)$-measure hyperbolic for all $p\in\{1,\ldots,\rank(V)-1\}$.
\endclaim

Again, an argument extremely similar to the proof of 1.10 shows that relative hyperbo\-licity 
is an open property.

\claim 2.5. Proposition|Let $(\cX,\cV)\to S$ be a holomorphic family
of compact directed manifolds $($by this, we mean a proper holomorphic
map $\cX\to S$ together with an analytic linear subspace $\cV\subset T_{\cX/S}\subset
T_{\cX}$ of the relative tangent bundle, defining a deformation $(X_s,V_s)_{s\in S}$
of the fibers$)$. Then the set of $s\in S$ such that the fiber $(X_s,V_s)$
is hyperbolic is open in $S$ with respect to the Euclidean topology.
\endclaim

Let us mention here an impressive result proved by Marco Brunella
[Bru03, Bru05, Bru06] concerning the behavior of the Kobayashi metric
on foliated varieties.

\claim 2.6. Theorem {\rm (Brunella)}|Let $X$ be a compact K\"ahler manifold
equipped with a $($possibly singular$)$ rank $1$ holomorphic foliation 
which is not a foliation by rational curves. Then the canonical bundle 
$K_{\cF}=\cF^*$ of the foliation is pseudoeffective $($i.e.\ the curvature
of $K_{\cF}$ is${}\ge 0$ in the sense of currents$)$.
\endclaim

The proof is obtained by putting on $K_{\cF}$ precisely the metric induced
by the Kobayashi metric on the leaves whenever they are generically 
hyperbolic (i.e.\ covered by the unit disk). The case of parabolic leaves
(covered by $\bC$) has to be treated separately.

\plainsubsection 2.C. Pluricanonical sheaves of a directed variety|

Let $(X,V)$ be a directed projective manifold where $V$ is possibly singular,
and let $r=\rank V$. If $\mu:\widehat X\to X$ is a proper modification
(a composition of blow-ups with smooth centers, say), we get a directed
manifold $(\widehat X,\widehat V)$ by taking $\widehat V$ 
to be the closure of $\mu_*^{-1}(V')$, where $V'=V_{\restriction X'}$ is the
restriction of $V$ over a Zariski open set $X'\subset X\smallsetminus \Sing(V)$ 
such that $\mu:\mu^{-1}(X')\to X'$ is a biholomorphism. We say that
$(\widehat X,\widehat V)$ is a {\em modification} of
$(X,V)$ and write $\swh{V}=\mu^*V$.

We will be interested in taking modifications realized by iterated blow-ups of 
certain non-singular subvarieties of the singular set $\Sing(V)$, so 
as to eventually ``improve'' the singularities of $V\,$; outside of
$\Sing(V)$ the effect of blowing-up will be irrelevant. The 
canonical sheaf $K_V$, resp.\
the pluricanonical sheaf sequence  $K^{[m]}_V$, will be defined here
in several steps, using the concept of bounded pluricanonical forms 
that was already introduced in [Dem11].

\claim 2.7. Definition|For a directed pair $(X,V)$ with $X$ non-singular,
we define $\bddK_V$, resp.\ $\bddK^{[m]}_V$, for any integer $m\ge 0$, 
to be the rank~$1$ analytic sheaves such that
$$
\eqalign{
\bddK_V(U)&=\hbox{sheaf of locally bounded sections 
of}~~\cO_X\big(\Lambda^r V^{\prime *}\big)(U\cap X'),\cr
\bddK^{[m]}_V(U)&=\hbox{sheaf of locally bounded sections 
of}~~\cO_X\big((\Lambda^r V^{\prime *})^{\otimes m}\big)(U\cap X'),\cr}
$$
where $r=\rank(V)$, $X'=X\smallsetminus \Sing(V)$, $V'=V_{\restriction X'}$,
and ``locally bounded'' means bounded with respect to a smooth Hermitian 
metric $h$ on $T_X$, on every set $W\cap X'$ such that $W$ is relatively
compact in $U$.
\endclaim

In the trivial case $r=0$, we simply set $\bddK_V^{[m]}=\cO_X$ for all~$m$;
clearly $\ECL(X,V)=\emptyset$ in that case, so there is not much to
say.  The above definition of $\bddK_V^{[m]}$ may look like an
analytic one, but it can easily be turned into an equivalent algebraic
definition:

\claim 2.8. Proposition|Consider the natural morphism 
$\cO(\Lambda^rT_X^*)\to \cO(\Lambda^r V^*)$ where $r=\rank V$ and
$\cO(\Lambda^r V^*)$ is defined as the quotient of 
$\cO(\Lambda^rT_X^*)$ by $r$-forms that have zero restrictions 
to $\cO(\Lambda^rV^*)$ on $X\smallsetminus \Sing(V)$. The 
bidual $\cL_V=\cO_X(\Lambda^r V^*)^{**}$ is an invertible sheaf,
and our natural morphism can be written
$$
\cO(\Lambda^rT_X^*)\to \cO(\Lambda^rV^*)=\cL_V\otimes\cJ_V\subset \cL_V
\leqno(2.8_1)
$$
where $\cJ_V$ is a certain ideal sheaf of $\cO_X$ whose zero set is
contained in $\Sing(V)$, and the arrow on the left is surjective by 
definition. Then 
$$
\bddK_V^{[m]}=\cL_V^{\otimes m}\otimes\overline{\cJ^m_V}\leqno(2.8_2)
$$
where $\overline{\cJ^m_V}$ is the integral closure of $\cJ_V^m$ in $\cO_X$.
In particular, $\bddK^{[m]}_V$ is always a coherent sheaf.
\endclaim

\plainproof. Let $(u_k)$ be a set of generators of $\cO(\Lambda^rV^*)$ obtained
(say) as the images of a basis $(dz_I)_{|I|=r}$ of $\Lambda^rT^*_X$
in some local coordinates near a point $x\in X$. Write $u_k=g_k\ell$ 
where $\ell$ is a local generator of $\cL_V$ at~$x$. Then $\cJ_V=(g_k)$ by
definition. The boundedness condition expressed in Def.~2.7 means that we 
take sections of the form $f\ell^{\otimes m}$ where $f$ is a holomorphic function 
on~$U\cap X'$ (and $U$ a neighborhood of $x$), such that 
$$
|f|\le C\Big(\sum|g_k|\Big)^m\leqno(2.8_3)
$$
for some constant $C>0$. But then $f$ extends holomorphically to $U$ into 
a function that lies in the integral closure ${\overline\cJ_V^m}$ (it is
well known that the latter is characterized analytically 
by condition~$(2.8_3)$). This proves Prop.~2.8.\qed

\claim 2.9. Lemma|Let $(X,V)$ be a directed variety.
\plainitem{\rm(a)} For any modification 
$\mu:(\widehat X,\widehat V)\to (X,V)$, there are always well defined
injective natural morphisms of rank $1$ sheaves
$$
\bddK_V^{[m]}\hookrightarrow \mu_*\big(\bddK_{\widehat V}^{[m]}\big)
\hookrightarrow \cL_V^{\otimes m}.
$$
\plainitem{\rm(b)} The direct image $\mu_*\big(\bddK_{\widehat V}^{[m]}\big)$ 
may only increase when we replace $\mu$ by a ``higher'' modification
$\widetilde\mu=\mu'\circ\mu:\widetilde X\to\widehat X\to X$
and $\swh{V}=\mu^*V$ by $\swt{V}=\swt{\mu}^*V$, i.e.\ there are injections
$$
\mu_*\big(\bddK_{\widehat V}^{[m]}\big)\hookrightarrow
\widetilde \mu_*\big(\bddK_{\widetilde V}^{[m]}\big)\hookrightarrow
\cL_V^{\otimes m}.
$$
We~refer to this property as the {\rm monotonicity principle}.\vskip0pt
\endclaim

\plainproof. (a) The existence of the first arrow is seen as follows:
the differential $\mu_*=d\mu:\swh{V}\to \mu^*V$ is smooth, so it is bounded
with respect to ambient Hermitian metrics on $X$ and $\swh{X}$. Going
to the duals reverses the arrows while preserving boundedness with
respect to the metrics. We thus get an arrow
$$
\mu^*({}^bV^\star)\hookrightarrow{}^b\widehat V^\star.
$$
By taking the top exterior power, followed by the $m$-th tensor product and 
the integral closure of the ideals involved, we get an injective arrow
$\mu^*\big(\bddK_V^{[m]}\big)\hookrightarrow\bddK_{\widehat V}^{[m]}$.
Finally we apply the direct image functor $\mu_*$ and the canonical
morphism $\cF\to\mu_*\mu^*\cF$ to get the first inclusion morphism.
The second arrow comes from the fact that $\mu^*\big(\bddK_V^{[m]}\big)$
coincides with $\cL_V^{\otimes m}$ (and with $\det(V^*)^{\otimes m}$) on 
the complement of
the codimension~$2$ set $S=\Sing(V)\cup\mu({\rm Exc}(\mu))$, and the fact
that for every open set $U\subset X$, sections of $\cL_V$ defined on 
$U\ssm S$ automatically extend to $U$ by Riemann's extension theorem
(or Hartog's extension theorem~$\ldots$), even without any boundedness
assumption.
\smallskip

\noindent (b) Given $\mu':\swt{X}\to\swh{X}$, we argue as in (a) that there is a
bounded morphism $d\mu':\swt{V}\to\swh{V}$.\qed
\medskip

By the monotonicity principle and the strong Noetherian 
property of coherent ideals, we infer that there exists a maximal 
direct image when $\mu:\swh{X}\to X$ runs over all non-singular 
modifications of~$X$. The following definition is thus legitimate.

\claim 2.10. Definition|We define the pluricanonical sheaves
$K_V^m$ of $(X,V)$ to be the inductive limits
$$
K^{[m]}_V:=\lim_{{\scriptstyle\lra\atop\scriptstyle\mu}}
\mu_*\big(\bddK_{\widehat V}^{[m]}\big)=
\max_\mu\mu_*\big(\bddK_{\widehat V}^{[m]}\big)
$$
taken over the family of all modifications $\mu:(\swh{X},\swh{V})\to (X,V)$,
with the trivial $($filtering$)$ partial order. The canonical sheaf $K_V$ 
itself is defined to be the same as $K_V^{[1]}$. By construction, 
we have for every $m\ge 0$ inclusions
$$
\bddK_V^{[m]}\hookrightarrow K_V^{[m]}\hookrightarrow \cL_V^{\otimes m},
$$
and $K_V^{[m]}=\cJ_V^{[m]}\cdot\cL_V^{\otimes m}$ for a certain sequence 
of integrally closed ideals $\cJ_V^{[m]}\subset\cO_X$.
\endclaim

It is clear from this construction that $K^{[m]}_V$ is birationally invariant,
i.e.\ that $K^{[m]}_V=\mu_*(K^{[m]}_{V'})$ for every modification
$\mu:(X',V')\to (X,V)$. Moreover the sequence is submultiplicative, i.e.\
there are injections
$$
K^{[m_1]}_V\otimes K^{[m_2]}_V\hookrightarrow K^{[m_1+m_2]}_V
$$
for all non-negative integers $m_1,\,m_2\,$; the corresponding sequence of
ideals $\cJ_V^{[m]}$ is thus also submultiplicative.
By blowing up $\cJ_V^{[m]}$ and taking a desingularization $\widehat X$
of the blow-up, one  can always find a {\it log-resolution} of 
$\cJ_V^{[m]}$, i.e.\ a 
modification $\mu_m:\widehat X_m\to X$ such 
that $\mu_m^*\cJ_V^{[m]}\subset\cO_{\widehat X_m}$ is an invertible ideal
sheaf; it~follows that 
$$
\mu_m^* K_V^{[m]}=\mu_m^*\cJ_V^{[m]}\cdot(\mu_m^*\cL_V)^{\otimes m}
$$
is an invertible sheaf on $\widehat X_m$.
We do not know whether $\mu_m$ can be taken independent of $m$, nor
whether the inductive limit introduced in Definition~2.10 is reached
for a $\mu$ that is independent of~$m$. If such a ``uniform'' 
$\mu$ exists, it could be thought of as some sort of replacement for
the resolution of singularities of directed structures (which do
not exist in the naive sense that $V$ could be made non-singular).
By means of a standard Serre-Siegel argument, one can easily show

\claim 2.11. Proposition|Let $(X,V)$ be a directed variety $(X,V)$ and $F$
be an invertible sheaf on~$X$. The following properties are equivalent~$:$
\plainitem{\rm(a)} there exists a constant $c>0$ and $m_0>0$ such that
$h^0(X,K_V^{[m]}\otimes F^{\otimes m})\ge c\,m^n$ for $m\ge m_0$, where
$n=\dim X\,;$
\plainitem{\rm(b)} the space of sections $H^0(X,K_V^{[m]}\otimes F^{\otimes m})$ 
provides a 
generic embedding of $X$ in projective space for sufficiently large $m\,;$
\plainitem{\rm(c)} there exists $m>0$ and a log-resolution 
$\mu_m:\widehat X_m\to X$ 
of $K_V^{[m]}$ such that $\mu_m^*(K_V^{[m]}\otimes F^{\otimes m})$ is a big 
invertible sheaf on $\widehat X_m\,;$
\plainitem{\rm(d)} there exists $m>0$, a modification
$\swt{\mu}_m:(\widetilde X_m,\widetilde V_m)\to (X,V)$ and a log-resolution 
$\mu'_m:\widehat X_m\to \widetilde X$ of $\bddK_{\widetilde V_m}^{[m]}$ such 
that $\mu_m^{\prime\,*}(\bddK_{\widetilde V_m}^{[m]}\otimes
\swt{\mu}_m^*F^{\otimes m})$ is a big invertible sheaf on
$\widehat X_m$.\vskip0pt
\noindent
We will express any of these equivalent properties by saying that the
twisted pluricanonical sheaf sequence $K_V^\bullet \otimes F^\bullet$ is big.
\endclaim

In the special case $F=\cO_X$, we introduce

\claim 2.12. Definition|We say that $(X,V)$ is of
general type if $K_V^\bullet$ is big.
\endclaim

\claim 2.13. Remarks|{\rm \smallskip\noindent
{\rm(a)} At this point, it is important to stress the difference between
``our'' canonical sheaf $K_V$, and the sheaf $\cL_V$, which
is considered by some experts as ``the canonical sheaf of the foliation'' 
defined by $V$, in the integrable case. Notice that $\cL_V$ can also be
obtained as the direct image $\cL_V=i_*\cO(\det V^*)$ 
associated with the injection $i:X\ssm \Sing(V)\hookrightarrow X$. 
The discrepancy already occurs
with the rank~$1$ linear space $V\subset T_{\bP^n_\bC}$ consisting at each
point $z\ne 0$ of the tangent to the line $(0z)$ (so that
necessarily $V_0=T_{\bP^n_\bC,0}$). As a sheaf (and not as a linear space), 
$i_*\cO(V)$ is the invertible sheaf generated by the vector field 
$\xi=\sum z_j\partial/\partial z_j$ on the affine open set
$\bC^n\subset\bP^n_\bC$, and therefore $\cL_V:=i_*\cO(V^*)$ is generated
over $\bC^n$ by the unique $1$-form $u$ such that $u(\xi)=1$. 
Since $\xi$ vanishes at $0$, the generator $u$ is {\em unbounded} 
with respect to a smooth metric $h_0$
on  $T_{\bP^n_\bC}$, and it is easily seen that $K_V$ is the non-invertible
sheaf $K_V=\cL_V\otimes \gm_{\bP^n_\bC,0}$. We can make it invertible
by considering the blow-up $\mu:\wt X\to X$ of $X=\bP^n_\bC$
at $0$, so that $\mu^*K_V$ is isomorphic to $\mu^*\cL_V\otimes\cO_{\wt X}(-E)$ 
where $E$ is the exceptional divisor. The integral curves $C$ of $V$ are of
course lines through~$0$, and when a standard parametrization is used,
their derivatives do not vanish at~$0$, while the sections of 
$i_*\cO(V)$ do -- a first sign
that $i_*\cO(V)$ and $i_*\cO(V^*)$ are the {\em wrong objects} to consider. 

\noindent
{\rm (b)} When $V$ is of rank 1 , we get a foliation by curves on $X$. 
If $(X,V)$ is of general type (i.e.\ $K_V^\bu$ is big), we will see in 
Prop.~4.9 that almost all leaves of $V$ are hyperbolic, i.e.\ covered 
by the unit 
disk. This would not be true if $K_V^\bu$ was replaced by $\cL_V$, 
In fact, the examples of pencils of conics or cubic curves in $\bP^2$
already produce this phenomenon, as we have seen in the introduction,
right after the generalized GGL conjecture~0.4. For this second reason,
we believe that $K_V^\bu$  is a more appropriate concept
of ``canonical sheaf'' than $\cL_V$~is.

\noindent
{\rm(c)} When $\dim X=2$, a singularity of a (rank $1$) foliation $V$ is 
said to be {\em simple} if the linear part of the local vector field generating
$\cO(V)$  has two distinct eigenvalues $\lambda\ne 0$, $\mu\ne 0$ such that
the quotient $\lambda/\mu$ is not a positive rational number. Seidenberg's
theorem [Sei68] says there always exists a composition of blow-ups 
$\mu:\swh{X}\to X$ such that $\swh{V}=\mu^*V$ only has simple 
singularities. It is easy to check that the inductive limit canonical sheaf 
$K_V^{[m]}=\mu_*(\bddK_{\widehat V}^{[m]})$ is reached whenever
$\swh{V}=\mu^*V$ has simple singularities.}
\endclaim

\section{Algebraic hyperbolicity}

In the case of projective algebraic varieties, hyperbolicity is
expected to be related to other properties of a more algebraic
nature. Theorem~3.1 below is a first step in this direction.

\claim 3.1.~Theorem|Let $(X,V)$ be a compact complex directed manifold
and let $\sum\omega_{jk}dz_j\otimes d\ol z_k$ be a Hermitian metric on~$T_X$,
with associated positive $(1,1)$-form
$\omega={i\over 2}\sum\omega_{jk}dz_j\wedge d\ol z_k$.
Consider the following
three properties, which may or not be satisfied by~$(X,V)\,:$
\plainitem{\rm(i)} $(X,V)$ is hyperbolic.
\plainitem{\rm(ii)} There exists $\varepsilon>0$ such that every compact
irreducible curve $C\subset X$ tangent to $V$ satisfies
$$
-\chi(\ol C)=2g(\ol C)-2\ge\varepsilon\,\deg_\omega(C),
$$
where $\deg_\omega(C)=\int_C\omega$,  and where $g(\ol C)$ is the genus of
the normalization $\ol C$ of~$C$ and $\chi(\ol C)$ its Euler charac\-teristic
$($the degree coincides with the usual concept of degree if $X$ is projective,
embedded in $\bP^N$ via a very ample line bundle $A$, and
$\omega=\Theta_{A,h_A}>0\;;$ such an estimate is of course independent 
of the choice of~$\omega$, provided that $\varepsilon$ is changed
accordingly$.)$
\plainitem{\rm(iii)} There does not exist any non-constant holomorphic map
$\Phi:Z\to X$ from an abelian variety $Z$ to $X$ such that
$\Phi_*(T_Z)\subset V$.
\smallskip\noindent
Then {\rm(i)$\,\Rightarrow\,$(ii)$\,\Rightarrow\,$(iii)}.
\endclaim

\plainproof. (i)$\,\Rightarrow\,$(ii). If $(X,V)$ is hyperbolic, there is a
constant $\varepsilon_0>0$ such that $\bfk_{(X,V)}(\xi)\ge\varepsilon_0
\|\xi\|_\omega$ for all $\xi\in V$. Now, let $C\subset X$ be a compact
irreducible curve tangent to $V$ and let $\nu:\ol C\to C$ be its
normalization. As $(X,V)$ is hyperbolic, $\ol C$ cannot be a rational
or elliptic curve. Hence $\ol C$ admits the disk as its universal
covering $\rho:\Delta\to\ol C$.

The Kobayashi-Royden metric $\bfk_\Delta$ is the Finsler metric
$|dz|/(1-|z|^2)$ associated with the Poincar\'e metric
$|dz|^2/(1-|z|^2)^2$ on~$\Delta$, and $\bfk_{\ol C}$ is such that
$\rho^* \bfk_{\ol C}=\bfk_\Delta$. In other words, the metric
$\bfk_{\ol C}$ is induced by the unique Hermitian metric on $\ol C$
of constant Gaussian curvature~$-4$. If $\sigma_\Delta=
{i\over 2}dz\wedge d\ol z/(1-|z|^2)^2$ and $\sigma_{\ol C}$ are
the corresponding area measures, the Gauss-Bonnet formula
(integral of the curvature${}=2\pi\,\chi(\ol C)$) yields
$$
\int_{\ol C}d\sigma_{\ol C}=-{1\over 4}\int_{\ol C}{\rm curv}
(\bfk_{\ol C})=-{\pi\over 2}\chi(\ol C)
$$
On the other hand, if $j:C\to X$ is the inclusion, the monotonicity
property (2.2) applied to the holomorphic map $j\circ\nu:\ol C\to X$
shows that
$$
\bfk_{\ol C}(t)\ge \bfk_{(X,V)}\big((j\circ\nu)_* t\big)
\ge \varepsilon_0\big\|(j\circ\nu)_* t\big\|_\omega,
\qquad\forall t\in T_{\ol C}.
$$
From this, we infer $d\sigma_{\ol C}\ge\varepsilon_0^2(j\circ\nu)^*
\omega$, thus
$$
-{\pi\over 2}\chi(\ol C)=\int_{\ol C}d\sigma_{\ol C}\ge
\varepsilon_0^2\int_{\ol C}(j\circ\nu)^*\omega=
\varepsilon_0^2\int_C\omega.
$$
Property (ii) follows with $\varepsilon=2\varepsilon_0^2/\pi$.
\smallskip

\noindent (ii)$\,\Rightarrow\,$(iii). First observe that (ii) excludes
the existence of elliptic and rational curves tangent to~$V$. Assume that
there is a non-constant holomorphic map $\Phi:Z\to X$ from an abelian
variety $Z$ to $X$ such that $\Phi_*(T_Z)\subset V$. We must have
$\dim\Phi(Z)\ge 2$, otherwise $\Phi(Z)$ would be a
curve covered by images of holomorphic maps $\bC\to\Phi(Z)$, and so
$\Phi(Z)$ would be elliptic or rational, contradiction. Select a
sufficiently general curve $\Gamma$ in $Z$ (e.g., a curve obtained as an
intersection of very generic divisors in a given very ample linear system 
$|L|$ in~$Z$). Then all isogenies $u_m:Z\to Z$, $s\mapsto ms$ map $\Gamma$ 
in a $1:1$ way to curves $u_m(\Gamma)\subset Z$, except maybe for
finitely many double points of $u_m(\Gamma)$ when $\dim Z=2\;$: we leave this
as an exercise to the reader, using Bertini type arguments). It follows
that the normalization of $u_m(\Gamma)$ is isomorphic to $\Gamma$. If
$\Gamma$ is general enough and $\tau_a:Z\to Z$, $w\mapsto w+a$ denote
translations of~$Z$, similar arguments show that for general $a\in Z$ the images
$$
C_{m,a}:=\Phi(\tau_a(u_m(\Gamma)))\subset X
$$
are also generically $1:1$ images of $\Gamma$, thus $\ol C_{m,a}\simeq\Gamma$
and $g(\ol C_{m,a})=g(\Gamma)$. We claim that on average $C_{m,a}$ has 
degree${}\geq\Const\,m^2$. In fact, if $\mu$ is the translation invariant
probability measure on $Z$
$$
\int_{C_{m,a}}\omega=
\int_\Gamma u_m^*(\tau_a^*\Phi^*\omega),\quad\hbox{and hence}\quad
\int_{a\in Z}\bigg(\int_{C_{m,a}}\omega\bigg)d\mu(a)
=\int_\Gamma u_m^*\beta
$$
where $\beta=\int_{a\in Z}(\tau_a^*\Phi^*\omega)\,d\mu(a)$ is a translation
invariant $(1,1)$-form on $Z$. Therefore $\beta$ is a constant
coefficient $(1,1)$-form, so $u_m^*\beta=m^2\beta$ and the right hand side is
$cm^2$ with $c=\int_\Gamma\beta>0$. For a suitable choice of $a_m\in Z$,
we have $\deg_\omega C_{m,a_m}\geq cm^2$ and
$(2g(\overline C_{m,a_m})-2)/\deg_\omega(C_{m,a_m})\to 0$, contradiction.
\qed

\claim 3.2.~Definition|We say that a projective directed manifold
$(X,V)$ is ``algebraically hyperbolic'' if it satisfies Property
$3.1$~{\rm(ii)}, namely, if there exists $\varepsilon>0$ such that
every algebraic curve $C\subset X$ tangent to $V$ satisfies
$$
2g(\ol C)-2\ge\varepsilon\,\deg_\omega(C).
$$
\endclaim

A nice feature of algebraic hyperbolicity is that it satisfies an algebraic
analogue of the openness property.

\claim 3.3.~Proposition|Let $(\cX,\cV)\to S$ be an algebraic family
of projective algebraic directed manifolds $($given by a projective
morphism $\cX\to S)$. Then the set of $t\in S$ such that the fiber
$(X_t,V_t)$ is algebraically hyperbolic is open with respect to
the ``countable Zariski topology'' of $S$ $($by definition, this is
the topology for which closed sets are countable unions of algebraic
sets$)$.
\endclaim

\plainproof. After replacing $S$ by a Zariski open subset, we may assume that
the total space $\cX$ itself is quasi-projective. Let $\omega$ be the
K\"ahler metric on~$\cX$ obtained by pulling back the Fubini-Study
metric via an embedding in a projective space. If the integers $d>0$,
$g\ge 0$ are fixed, the set $A_{d,g}$ of $t\in S$ such that
$X_t$ contains an algebraic $1$-cycle $C=\sum m_jC_j$ tangent to $V_t$ with
$\deg_\omega(C)=d$ and $g(\ol C)=\sum m_j\,g(\ol C_j)\le g$ is a closed
algebraic subset of $S$ (this follows from the existence of a relative
cycle space of curves of given degree, and from the fact that the
geometric genus is Zariski lower semi-continuous). Now, the set of non
algebraically hyperbolic fibers is by definition
$$
\bigcap_{k>0}~~\bigcup_{2g-2<d/k}~A_{d,g}.
$$
This concludes the proof (of course, one has to know that the countable
Zariski topology is actually a topology, namely that the class
of countable unions of algebraic sets is stable under arbitrary
intersections; this can be easily checked by an induction on
dimension).\qed

\claim 3.4.~Remark|{\rm More explicit versions of the openness
property have been dealt with in the literature. H.~Clemens
([Cle86] and [CKM88]) has shown that on a very generic surface of degree
$d\ge 5$ in $\bP^3$, the curves of type $(d,k)$ are of genus $g>kd(d-5)/2$
(recall that a very generic surface $X\subset\bP^3$ of degree${}\ge 4$ has
Picard group generated by $\cO_X(1)$ thanks to the Noether-Lefschetz theorem;
thus any curve on the surface is a complete intersection with another
hypersurface of degree $k\,$; such a curve is said to be of type
$(d,k)\,$; genericity is taken here in the sense of the countable
Zariski topology).  Improving on this result of Clemens, Geng Xu
[Xu94] proved that every curve contained in a very generic surface of
degree $d\ge 5$ satisfies the sharp bound $g\ge d(d-3)/2-2$. In
April~2018, I.~Coskun and E.~Riedl improved the above bounds, and got
the more precise bound $g\ge 1+(dk(d-5)+k)/2\,$; this result actually
shows that a very generic surface of degree $d\ge 5$ is algebraically
hyperbolic in the sense of Definition~3.2.
In higher dimension, L.~Ein ([Ein88], [Ein91]) proved that every subvariety
of a very generic hypersurface $X\subset\bP^{n+1}$ of degree $d\ge 2n+1$ $(n\ge 2)$,
is of general type. This was reproved by a simple efficient technique by
C.~Voisin in [Voi96], along with other improvements.}
\endclaim

\claim 3.5.~Remark|{\rm In view of Proposition 1.10, it would be
interesting to know whether algebraic hyperbolicity is open with
respect to the Euclidean topology$\,$; still more interesting would
be to know whether Kobayashi hyperbolicity is open for the countable
Zariski topology (of course, both properties would follow
immediately if one knew that algebraic hyperbolicity and Kobayashi
hyperbolicity  coincide, but they seem otherwise highly non-trivial to
establish). The latter openness property has raised an important amount
of work around the following more particular question: is a (very) generic
hypersurface $X\subset\bP^{n+1}$ of degree $d$ large enough
(say $d\ge 2n+1$) Kobayashi hyperbolic$\,$? Again, ``very generic'' is to
be taken here in the sense of the countable Zariski topology. Brody-Green
[BrGr77] and Nadel [Nad89] produced examples of hyperbolic
surfaces in $\bP^3$ for all degrees $d\ge 50$, and Masuda-Noguchi
[MaNo96] gave examples of such hypersurfaces in $\bP^n$ for
arbitrary $n\ge 2$, of degree $d\ge d_0(n)$ large enough. The 
hyperbolicity of complements $\bP^n\ssm H$ of generic
divisors may be inferred from the compact case; in fact if
$H=\{P(z_0\ld z_n)=0\}$ is a smooth generic divisor of degree $d$,
one may look at the hypersurface
$$
X=\big\{z_{n+1}^d=P(z_0\ld z_n)\big\}\subset\bP^{n+1}
$$
which is a cyclic $d\,{:}\,1$ covering of $\bP^n$. Since any holomorphic
map $f:\bC\to\bP^n\ssm H$ can be lifted to $X$, it is clear that the
hyperbolicity of $X$ would imply the hyperbolicity of $\bP^n\ssm H$.
The hyperbolicity of complements of divisors in $\bP^n$ has been
investigated by many authors. In the case $n=2$, Huynh, Vu and Xie
[HVX17, Theorem~1.2] have announced that $\bP^2\ssm C$ is hyperbolic
for a very general curve~$C$ of degree $d\geq 11$ (and that a very
general surface $X\subset\bP^3$ of degree $d\geq 15$ is hyperbolic,
[HVX17, Theorem~1.5]).
The reader can also consult [CFZ17, \S4] for more details and references
in these directions.
\qed}
\endclaim 

In the ``absolute case''~$V=T_X$, it seems reasonable to expect that
Properties 3.1~(i), (ii) are equivalent, i.e.\ that
Kobayashi and algebraic hyperbolicity coincide. However, it was observed
by Serge Cantat [Can00] that Property 3.1~(iii) is not sufficient to imply
the hyperbolicity of $X$, at least when $X$ is a general complex
surface: a general (non-algebraic) K3 surface is known to have no elliptic
 curves and does not admit either any surjective map from an abelian
variety; however such a surface is not Kobayashi hyperbolic. We are uncertain
about the sufficiency of Property 3.1~(iii) when $X$ is assumed to be
projective.

\section{The Ahlfors-Schwarz lemma for metrics of negative curvature}

One of the most basic ideas is that hyperbolicity should somehow be
related with suitable negativity properties of the curvature. For
instance, it is a standard fact already observed in Kobayashi [Kob70]
that the negativity of $T_X$ (or the ampleness of $T^*_X$) implies
the hyperbolicity of~$X$. There are many ways of improving or
generalizing this result. We present here a few simple examples of such
generalizations. 

\plainsubsection 4.A. Exploiting curvature via potential theory|

If $(V,h)$ is a holomorphic vector bundle equipped with
a smooth Hermitian metric, we denote by $\nabla_h=\nabla'_h+\nabla''_h$
the associated Chern connection and by $\Theta_{V,h}={i\over 2\pi}
\nabla_h^2$ its Chern curvature tensor.

\claim 4.1.~Proposition|Let $(X,V)$ be a compact directed
manifold. Assume that $V$ is non-singular and that
$V^*$ is ample. Then $(X,V)$ is hyperbolic.
\endclaim

\plainproof\ (from an original idea of [Kob75]). Recall that a vector bundle
$E$ is said to be ample if $S^mE$ has enough global sections
$\sigma_1\ld\sigma_N$ so as to generate $1$-jets of sections at any
point, when $m$ is large. One obtains a Finsler metric $N$ on $E^*$
by putting
$$
N(\xi)=\Big(\sum_{1\le j\le N}|\sigma_j(x)\cdot\xi^m|^2\Big)^{1/2m},\qquad
\xi\in E^*_x,
$$
and $N$ is then a strictly plurisubharmonic function on the total space
of $E^*$ minus the zero section (in other words, the line
bundle $\cO_{P(E^*)}(1)$ has a metric of positive curvature). By
the ampleness assumption on $V^*$, we thus have a Finsler
metric $N$ on $V$ which is strictly plurisubharmonic outside the zero
section. By the Brody lemma, if $(X,V)$ is not hyperbolic, there is an
entire curve $g:\bC\to X$ tangent to $V$ such that
$\sup_\bC\|g'\|_\omega\le 1$ for some given Hermitian metric $\omega$
on~$X$. Then $N(g')$ is a bounded subharmonic function on $\bC$ which
is strictly subharmonic on $\{g'\ne 0\}$. This is a contradiction, for
any bounded subharmonic function on $\bC$ must be constant.\qed
\medskip

\plainsubsection 4.B. Ahlfors-Schwarz lemma|

Proposition 4.1 can be generalized a little bit further by means of
the Ahlfors-Schwarz lemma (see e.g.\ [Lang87]; we refer to [Dem95]
for the generalized version presented here; the proof is merely an 
application of the maximum principle plus a regularization argument).

\claim 4.2.~Ahlfors-Schwarz lemma|Let $\gamma(t)=\gamma_0(t)\,i\,dt
\wedge d\ol t$ be a Hermitian metric on $\Delta_R$ where $\log\gamma_0$ is
a subharmonic function such that $i\,\ddbar\log\gamma_0(t)\ge A\,\gamma(t)$
in the sense of currents, for some positive constant $A$. Then $\gamma$ can
be compared with the Poincar\'e metric of $\Delta_R$ as follows:
$$
\gamma(t)\le{2\over A}{R^{-2}|dt|^2\over(1-|t|^2/R^2)^2}.
$$
More generally, let $\gamma=i\sum\gamma_{jk}dt_j\wedge d\ol t_k$ be an
almost everywhere positive Hermitian form on the ball $B(0,R)\subset\bC^p$,
such that $-\Ricci(\gamma):=i\,\ddbar\log\det(\gamma)\ge A\gamma$ in the
sense of currents, for some constant $A>0$ $($this means in particular
that $\det(\gamma)=\det(\gamma_{jk})$ is such that $\log\det(\gamma)$ is
plurisubharmonic$)$. Then the $\gamma$-volume form is controlled by
the Poincar\'e volume form~:
$$
\det(\gamma)\le\Big({p+1\over AR^2}\Big)^p{1\over(1-|t|^2/R^2)^{p+1}}.
$$
\endclaim

\plainsubsection 4.C. Applications of the Ahlfors-Schwarz lemma to hyperbolicity|

Let $(X,V)$ be a {\em projective} directed variety. We assume throughout this subsection that $X$ is {\em non-singular}.

\claim 4.3.~Proposition|Assume that $V$ itself is non-singular and that the dual bundle $V^*$ is ``very big'' in the following sense: there
exists an ample line bundle $L$ and a sufficiently large integer $m$
such that the global sections in $H^0(X,S^mV^*\otimes L^{-1})$
generate all fibers over $X\ssm Y$, for some analytic subset $Y\subsetneq X$.
Then all entire curves $f:\bC\to X$ tangent to $V$ satisfy
$f(\bC)\subset Y$.
\endclaim

\plainproof. Let $\sigma_1\ld\sigma_N\in H^0(X,S^mV^*\otimes L^{-1})$ be
a basis of sections generating $S^mV^*\otimes L^{-1}$ over $X\ssm Y$.
If $f:\bC\to X$ is tangent to~$V$, we define a semi-positive Hermitian form
$\gamma(t)=\gamma_0(t)\,|dt|^2$ on $\bC$ by putting
$$
\gamma_0(t)=\sum\|\sigma_j(f(t))\cdot f'(t)^m\|_{L^{-1}}^{2/m}
$$
where $\|~~\|_L$ denotes a Hermitian metric with positive
curvature on~$L$. If $f(\bC)\not\subset Y$, the form $\gamma$ is not
identically $0$ and we then find 
$$
i\,\ddbar\log\gamma_0\ge {2\pi\over m}f^*\Theta_L
$$
where $\Theta_L$ is the curvature form. The positivity assumption combined
with an obvious homogeneity argument yield
$$
{2\pi\over m}f^*\Theta_L\ge\varepsilon\|f'(t)\|_\omega^2\,|dt|^2\ge
\varepsilon'\,\gamma(t)
$$
for any given Hermitian metric $\omega$ on~$X$. Now, for any $t_0$ with
$\gamma_0(t_0)>0$, the Ahlfors-Schwarz lemma shows that $f$ can only exist 
on a disk $D(t_0,R)$ such that $\gamma_0(t_0)\le {2\over\varepsilon'}R^{-2}$,
contradiction.\qed

There are similar results for $p$-measure hyperbolicity, see e.g.\ [Carl72]
and [Nog77b]:

\claim 4.4.~Proposition|Assume that $V$ is non-singular and that
$\Lambda^pV^*$ is ample. Then $(X,V)$ is infinitesimally
$p$-measure hyperbolic. More generally, assume that $\Lambda^pV^*$
is very big with base locus contained in $Y\subsetneq X$ $($see
Proposition~$3.3)$.
Then $\bfe^p$ is non-degenerate over $X\ssm Y$.
\endclaim

\plainproof. By the ampleness assumption, there is a smooth Finsler metric $N$ on
$\Lambda^pV$ which is strictly plurisubharmonic outside the zero section.
We select also a Hermitian metric $\omega$ on~$X$. For any holomorphic
map $f:\bB_p\to X$ we define a semi-positive Hermitian metric
$\swt{\gamma}$ on $\bB_p$ by putting $\swt{\gamma}=f^*\omega$. Since
$\omega$ need not have any good curvature estimate, we introduce
the function $\delta(t)=N_{f(t)}(\Lambda^pf'(t)\cdot\tau_0)$, where
$\tau_0=\partial/\partial t_1\wedge\cdots\wedge\partial/\partial t_p$,
and select a metric $\gamma=\lambda\swt{\gamma}$ conformal to
$\swt{\gamma}$ such that $\det(\gamma)=\delta$. Then $\lambda^p$ is equal
to the ratio $N/\Lambda^p\omega$ on the element $\Lambda^pf'(t)\cdot\tau_0
\in\Lambda^pV_{f(t)}$. Since $X$ is compact, it is clear that the
conformal factor $\lambda$ is bounded by an absolute constant
independent of~$f$. From the curvature assumption we then get
$$
i\,\ddbar\log\det(\gamma)=i\,\ddbar\log\delta\ge (f,\Lambda^pf')^*
(i\,\ddbar\log N)\ge\varepsilon f^*\omega\ge\varepsilon'\,\gamma.
$$
By the Ahlfors-Schwarz lemma we infer that $\det(\gamma(0))\le C$ for some
constant $C$, i.e., $N_{f(0)}(\Lambda^pf'(0)\cdot\tau_0)\le C'$. This means
that the Kobayashi-Eisenman pseudometric $\bfe^p_{(X,V)}$ is positive
definite everywhere and uniformly bounded from below. In the
case $\Lambda^pV^*$ is very big with base locus $Y$, we use essentially
the same arguments, but we then only have $N$ being positive definite
on $X\ssm Y$.\qed

\claim 4.5.~Corollary {\rm([Gri71], KobO71])}|If $X$ is a projective
variety of general type, the Kobayashi-Eisenmann volume form $\bfe^n$,
$n=\dim X$, can degenerate only along a proper algebraic set
$Y\subsetneq X$.
\endclaim

The converse of Corollary 4.5 is expected to be true, namely, the generic
non-degeneracy of $\bfe^n$ should imply that $X$ is of general type; 
this is only known for surfaces (see [GrGr80] and [MoMu82]):

\claim 4.6.~General Type Conjecture {\rm(Green-Griffiths [GrGr80])}|
A projective algebraic variety $X$ is measure hyperbolic $($i.e.\ $\bfe^n$
degenerates only along a proper algebraic subvariety$)$ if and only if
$X$ is of general type.
\endclaim

An essential step in the proof of the necessity of having general type
subvarieties would be to show that manifolds of Kodaira
dimension $0$ (say, Calabi-Yau manifolds and holomorphic symplectic
manifolds, all of which have $c_1(X)=0$) are not measure hyperbolic, 
e.g.\ by exhibiting enough families of curves $C_{s,\ell}$ covering $X$
such that $(2g(\ol C_{s,\ell})-2)/\deg(C_{s,\ell})\to 0$. 

\claim 4.7.~Conjectural corollary {\rm(Lang)}|A projective algebraic
variety $X$ is hyperbolic if and only if all its algebraic subvarieties
$($including $X$ itself\/$)$ are of general type.
\endclaim

\claim 4.8.~Remark|{\rm The GGL conjecture implies the ``if'' part 
of~4.7, and the General Type Conjecture~4.6 implies the ``only if'' part
of 4.7. In fact if the GGL conjecture holds and every subvariety $Y$ of $X$ is
of general type, then it is easy to infer that every entire curve
$f:\bC\to X$ has to be constant by induction on $\dim X$, because
in fact $f$ maps $\bC$ to a certain subvariety~$Y\subsetneq X$. Therefore
$X$ is hyperbolic. Conversely, if Conjecture 4.6 holds and $X$ has 
a certain subvariety $Y$ which is not of general type, then $Y$ is not measure
hyperbolic. However Proposition 2.4 shows that hyperbolicity implies
measure hyperbolicity. Therefore $Y$ is not hyperbolic
and so $X$ itself is not hyperbolic either.}
\endclaim

We end this section by another easy application of the Ahlfors-Schwarz 
lemma for the case of rank 1 (possibly singular) foliations.

\claim 4.9.~Proposition|Let $(X,V)$ be a projective directed manifold. 
Assume that $V$ is of rank $1$ and that $K_V^\bullet$ is big. Then 
$S$ be the union of the singular set $\Sing(V)$ and of the base locus
of $K_V^\bullet$ $($namely the intersection of the images $\mu_m(B_m)$
of the base loci $B_m$ of the invertible sheaves 
$\mu_m^*K_V^{[m]}$, $m>0$, obtained by
taking log-resolutions$)$. Then $\ECL(X,V)\subset S$, in other words,
all non-hyperbolic leaves of $V$ are contained in~$S$.
\endclaim

\plainproof. By Prop.~2.11~(d), we can take a blow-up $\swt{\mu}_m:\swt{X}_m\to X$ and a 
log-resolution $\mu'_m:\swh{X}_m\to \swt{X}_m$ such that
$F_m=\mu_m^{\prime\,*}(\bddK_{\widetilde V_m}^{[m]})$ is a big invertible sheaf.
This means that (after possibly increasing~$m$) we can find sections 
$\sigma_1,\ldots\sigma_N\in H^0(\swh{X}_m,F_m)$ that define a (singular)
Hermitian metric with strictly positive curvature on $F_m$, cf.\ Def.~8.1
below. Now, for every entire curve $f:(\bC,T_\bC)\to(X,V)$ not contained 
in $S$, we can choose $m$ and a lifting 
$\swt{f}:(\bC,T_\bC)\to(\swt{X},\swt{V})$
such that $\swt{f}(\bC)$ is not contained in the base locus of our sections. 
Again, we can define a semi-positive Hermitian form
$\gamma(t)=\gamma_0(t)\,|dt|^2$ on $\bC$ by putting
$$
\gamma_0(t)=\sum\|\sigma_j(f(t))\cdot f'(t)^m\|_{L^{-1}}^{2/m}.
$$
Then $\gamma$ is not identically zero and we have $i\ddbar\log\gamma_0\ge
\varepsilon\gamma$ by the strict positivity of the curvature. One should also
notice that $\gamma_0$ is locally bounded from above by the assumption that
the $\sigma_j$'s come from {\em locally bounded} sections on $\swt{X}_m$.
This contradicts the Ahlfors-Schwarz lemma, and thus it cannot happen that
$f(\bC)\not\subset S$.\qed

\section{Projectivization of a directed manifold}

\plainsubsection 5.A. The $1$-jet functor|
The basic idea is to introduce a functorial process which produces a new
complex directed manifold $(\swt{X},\swt{V})$ from a given one~$(X,V)$.
The new structure $(\swt{X},\swt{V})$ plays the role of a space of $1$-jets
over~$X$.  Fisrt assume that $V$ is {\em non-singular}. We let
$$
\swt{X}=P(V),\qquad \swt{V}\subset T_{\swt{X}}
$$
be the projectivized bundle of lines of $V$, together with a subbundle
$\swt{V}$ of $T_{\swt{X}}$ defined as follows: for every point $(x,[v])\in
\swt{X}$ associated with a vector $v\in V_x\ssm\{0\}$,
$$
\swt{V}_{(x,[v])}=\big\{\xi\in T_{\swt{X},\,(x,[v])}\,;\,\pi_*\xi\in
\bC v\big\},\qquad\bC v\subset V_x\subset T_{X,x},\leqno(5.1)
$$
where $\pi:\swt{X}=P(V)\to X$ is the natural projection and $\pi_*:
T_{\swt{X}}\to\pi^* T_X$ is its differential. On $\swt{X}=P(V)$
we have a tautological line bundle $\cO_{\swt{X}}(-1)\subset\pi^* V$
such that $\cO_{\swt{X}}(-1)_{(x,[v])}=\bC v$. The bundle $\swt{V}$ is
characterized by the two exact sequences
$$
\leqalignno{
&0\lra T_{\swt{X}/X}\lra\swt{V}\build\lra^{\pi_*}_{}\cO_{\swt{X}}(-1)
\lra 0,&(5.2)\cr
&0\lra\cO_{\swt{X}}\lra \pi^* V\otimes\cO_{\swt{X}}(1)
\lra T_{\swt{X}/X}\lra 0,&(5.2')\cr}
$$
where $T_{\swt{X}/X}$ denotes the relative tangent bundle of the fibration
$\pi:\swt{X}\to X$. The first sequence is a direct consequence of the
definition of $\swt{V}$, whereas the second is a relative version of the
Euler exact sequence describing the tangent bundle of the fibers
$P(V_x)$. From these exact sequences we infer
$$
\dim\swt{X}=n+r-1,\qquad \rank\swt{V}=\rank V=r,\leqno(5.3)
$$
and by taking determinants we find $\det(T_{\swt{X}/X})=
\pi^*\det(V)\otimes\cO_{\swt{X}}(r)$. Thus
$$
\det(\swt{V})=\pi^*\det(V)\otimes\cO_{\swt{X}}(r-1).\leqno(5.4)
$$
By definition, $\pi:(\swt{X},\swt{V})\to(X,V)$ is a morphism of
complex directed manifolds. Clearly, our construction is functorial, i.e.,
for every morphism of directed manifolds $\Phi:(X,V)\to(Y,W)$, there
is a commutative diagram
$$
\plainmatrix{(\swt{X},\swt{V})&\build\lra^{\textstyle\pi}_{}&(X,V)\cr
\wt\Phi~\smash{\raise 1.2em\hbox{$\vdasharrow$}}&&\big\downarrow\Phi\cr
(\swt{Y},\swt{W})&\build\lra^{\textstyle\pi}_{}&(Y,W)\cr}\leqno(5.5)
$$
where the left vertical arrow is the meromorphic map $P(V)\merto P(W)$
induced by the differential $\Phi_*:V\to\Phi^* W$ ($\swt{\Phi}$ is
actually holomorphic if $\Phi_*:V\to\Phi^* W$ is injective).

\plainsubsection 5.B. Lifting of curves to the $1$-jet bundle|

Suppose that we are given a holomorphic curve $f:\Delta_R\to X$
parametrized by the disk $\Delta_R$ of centre $0$ and radius $R$
in the complex plane, and that $f$ is a tangent curve of the
directed manifold, i.e., $f'(t)\in V_{f(t)}$ for every $t\in \Delta_R$.
If $f$ is non-constant, there is a well defined and unique tangent line
$[f'(t)]\in P(V_{f(t)})$ for every~$t$, even at stationary points, and the map
$$
\swt{f}:\Delta_R\to\swt{X},\qquad
t\mapsto\swt{f}(t):=(f(t),[f'(t)])\leqno(5.6)
$$
is holomorphic; in fact, at a stationary point $t_0$, we can write
$f'(t)=(t-t_0)^su(t)$ with $s\in\bN^*$ and $u(t_0)\ne 0$,
and we define the tangent line at $t_0$ to be $[u(t_0)]$, so that
$\swt{f}(t)=(f(t),[u(t)])$ near $t_0\,$; even for $t=t_0$, we still denote
$[f'(t_0)]=[u(t_0)]$ for simplicity of notation. By definition
$f'(t)\in\cO_{\swt{X}}(-1)_{\swt{f}(t)}=\bC\,u(t)$, so the derivative
$f'$ defines a section
$$
f':T_{\Delta_R}\to\swt{f}^*\cO_{\swt{X}}(-1).\leqno(5.7)
$$
Moreover $\pi\circ\swt{f}=f$, and thus
$$
\pi_*\swt{f}'(t)=f'(t)\in\bC u(t)\Longrightarrow
\swt{f}'(t)\in\swt{V}_{(f(t),u(t))}=\swt{V}_{\swt{f}(t)}
$$
and we see that $\swt{f}$ is a tangent trajectory of $(\swt{X},\swt{V})$.
We say that $\swt{f}$ is the {\em canonical lifting} of $f$ to~$\swt{X}$.
Conversely, if $g:\Delta_R\to\swt{X}$ is a tangent trajectory
of $(\swt{X},\swt{V})$, then by definition of $\swt{V}$ we see that
$f=\pi\circ g$ is a tangent trajectory of $(X,V)$ and that $g=\swt{f}$
(unless $g$ is contained in a vertical fiber $P(V_x)$, in which case
$f$ is constant).

For any point $x_0\in X$, there are local coordinates $(z_1\ld z_n)$ on a
neighborhood $\Omega$ of $x_0$ such that the fibers $(V_z)_{z\in\Omega}$
can be defined by linear equations
$$
V_z=\Big\{\xi=\sum_{1\le j\le n}\xi_j{\partial\over\partial z_j}\,;\,
\xi_j= \sum_{1\le k\le r}a_{jk}(z)\xi_k~\hbox{\rm for $j=r+1\ld n$}\Big\},
\leqno(5.8)
$$
where $(a_{jk})$ is a holomorphic $(n-r)\times r$ matrix. It follows that
a vector $\xi\in V_z$ is completely determined by its first $r$ components
$(\xi_1\ld\xi_r)$, and the affine chart $\xi_j\ne 0$ of $P(V)_{\restriction
\Omega}$ can be described by the coordinate system
$$
\Big(z_1\ld z_n;{\xi_1\over\xi_j}\ld{\xi_{j-1}\over\xi_j},
{\xi_{j+1}\over\xi_j}\ld{\xi_r\over\xi_j}\Big).\leqno(5.9)
$$
Let $f\simeq(f_1\ld f_n)$ be the components of $f$ in the coordinates
$(z_1\ld z_n)$ (we suppose here $R$ so small that $f(\Delta_R)\subset\Omega$).
It should be observed that $f$ is uniquely determined by its initial value
$x$ and by the first $r$ components $(f_1\ld f_r)$. Indeed, as $f'(t)\in
V_{f(t)}\,$, we can recover the other components by integrating the system
of ordinary differential equations
$$
f_j'(t)=\sum_{1\le k\le r}a_{jk}(f(t))f_k'(t),\qquad j>r,\leqno(5.10)
$$
on a neighborhood of~$0$, with initial data $f(0)=x$.
We denote by $m=m(f,t_0)$ the {\em multiplicity} of $f$ at any point
$t_0\in\Delta_R$, that is, $m(f,t_0)$ is the smallest integer $m\in\bN^*$
such that $f_j^{(m)}(t_0)\ne 0$ for some~$j$. By (5.10), we can always
suppose $j\in\{1\ld r\}$, for example $f_r^{(m)}(t_0)\ne 0$. Then
$f'(t)=(t-t_0)^{m-1}u(t)$ with $u_r(t_0)\ne 0$, and the lifting $\swt{f}$
is described in the coordinates of the affine chart $\xi_r\ne 0$ of
$P(V)_{\restriction\Omega}$ by
$$
\swt{f}\simeq\Big(f_1\ld f_n;{f_1'\over f_r'}\ld{f_{r-1}'\over f_r'}\Big).
\leqno(5.11)
$$

\plainsubsection 5.C. Curvature properties of the $1$-jet bundle|

We end this section with a few curvature computations.
Assume that $V$ is non singular and equipped with a smooth Hermitian
metric~$h$. Denote by
$\nabla_h=\nabla'_h+\nabla''_h$ the associated Chern connection and
by $\Theta_{V,h}={i\over 2\pi}\nabla_h^2$ its Chern curvature tensor.
For every point $x_0\in X$, there exists a ``normalized'' holomorphic frame
$(e_\lambda)_{1\le\lambda\le r}$ on a neighborhood of~$x_0$, such that
$$
\langle e_\lambda,e_\mu\rangle_h=\delta_{\lambda\mu}-
\sum_{1\le j,k\le n}c_{jk\lambda\mu}z_j\ol z_k+O(|z|^3),\leqno(5.12)
$$
with respect to any holomorphic coordinate system $(z_1\ld z_n)$ centered
at~$x_0$. A~computation of $d'\langle e_\lambda,e_\mu\rangle_h=
\langle\nabla'_h e_\lambda,e_\mu\rangle_h$ and
$\nabla^2_h e_\lambda=d''\nabla'_h e_\lambda$ then gives
$$
\leqalignno{\nabla'_h e_\lambda
&=-\sum_{j,k,\mu}c_{jk\lambda\mu}\ol z_k\,dz_j\otimes e_\mu+O(|z|^2),\cr
\Theta_{V,h}(x_0)
&={i\over 2\pi}\sum_{j,k,\lambda,\mu}c_{jk\lambda\mu}dz_j\wedge
d\ol z_k\otimes e_\lambda^*\otimes e_\mu.&(5.13)\cr}
$$
The above curvature tensor can also be viewed as a Hermitian form
on $T_X\otimes V$. In fact, one associates with $\Theta_{V,h}$ the Hermitian
form $\langle\Theta_{V,h}\rangle$ on $T_X\otimes V$ defined for all
$(\zeta,v)\in T_X\times_X V$ by
$$
\langle\Theta_{V,h}\rangle(\zeta\otimes v)=
\sum_{1\le j,k\le n,\,1\le\lambda,\mu\le r}
c_{jk\lambda\mu}\zeta_j\ol\zeta_kv_\lambda\ol v_\mu.\leqno(5.14)
$$
Let $h_1$ be the Hermitian metric on the tautological line bundle
$\cO_{P(V)}(-1)\subset\pi^* V$ induced by the metric $h$ of~$V$. We
compute the curvature $(1,1)$-form $\Theta_{h_1}(\cO_{P(V)}(-1))$
at an arbitrary point $(x_0,[v_0])\in P(V)$, in terms of $\Theta_{V,h}$. For
simplicity, we suppose that the frame $(e_\lambda)_{1\le\lambda\le r}$ has
been chosen in such a way that $[e_r(x_0)]=[v_0]\in P(V)$ and $|v_0|_h=1$.
We get holomorphic local coordinates $(z_1\ld z_n\,;\,\xi_1\ld\xi_{r-1})$
on a neighborhood of $(x_0,[v_0])$ in $P(V)$ by assigning
$$
(z_1\ld z_n\,;\,\xi_1\ld\xi_{r-1})\longmapsto(z,[\xi_1e_1(z)+\cdots
+\xi_{r-1}e_{r-1}(z)+e_r(z)])\in P(V).
$$
Then the function
$$
\eta(z,\xi)=\xi_1e_1(z)+\cdots+\xi_{r-1}e_{r-1}(z)+e_r(z)
$$
defines a holomorphic section of $\cO_{P(V)}(-1)$ in a neighborhood
of~$(x_0,[v_0])$. By using the expansion (5.12) for $h$, we find
$$
\leqalignno{
|\eta|_{h_1}^2=|\eta|_h^2=1+|\xi|^2&{}-\sum_{1\le j,k\le n}c_{jkrr}
z_j\ol z_k+O((|z|+|\xi|)^3),\cr
\Theta_{h_1}(\cO_{P(V)}(-1))_{(x_0,[v_0])}
&=-{i\over 2\pi}\ddbar\log|\eta|_{h_1}^2\cr
&={i\over 2\pi}\Big(\sum_{1\le j,k\le n}c_{jkrr}dz_j\wedge d\ol z_k
-\sum_{1\le\lambda\le r-1}d\xi_\lambda\wedge d\ol\xi_\lambda\Big).
&(5.15)}
$$

\section{Jets of curves and Semple jet bundles}

\plainsubsection 6.A. Semple tower of non-singular directed varieties|

Let $X$ be a complex $n$-dimensional manifold. Following ideas of 
\hbox{Green-Griffiths} [GrGr80], we let $J_kX\to X$ be the bundle of $k$-jets
of germs of parametrized curves in~$X$, that is, the set of equivalence
classes of holomorphic maps $f:(\bC,0)\to(X,x)$, with the equivalence
relation $f\sim g$ if and only if all derivatives $f^{(j)}(0)=g^{(j)}(0)$
coincide for $0\le j\le k$, when computed in some local coordinate system
of $X$ near~$x$. The projection map $J_kX\to X$ is simply $f\mapsto f(0)$.
If $(z_1\ld z_n)$ are local holomorphic coordinates on an open set
$\Omega\subset X$, the elements $f$ of any fiber $J_kX_x$,
$x\in\Omega$, can be seen as $\bC^n$-valued maps $$f=(f_1\ld
f_n):(\bC,0)\to\Omega\subset\bC^n,$$ and they are completely
determined by their Taylor expansion of order $k$ at~$t=0$
$$
f(t)=x+t\,f'(0)+{t^2\over 2!}f''(0)+\cdots+{t^k\over k!}f^{(k)}(0)+
O(t^{k+1}).
$$
In these coordinates, the fiber $J_kX_x$ can thus be identified with the
set of $k$-tuples of vectors 
$(\xi_1,\ldots,\xi_k)=(f'(0)\ld f^{(k)}(0))\in(\bC^n)^k$.
It follows that $J_kX$ is a holomorphic fiber bundle with typical fiber
$(\bC^n)^k$ over $X$ (however, $J_kX$ is not a vector bundle for $k\ge 2$,
because of the nonlinearity of coordinate changes; see formula (7.2)
in~\S$\,$7).

According to the philosophy developed throughout this paper, we describe
the concept of jet bundle in the general situation of complex directed
manifolds. If $X$ is equipped with a holomorphic subbundle $V\subset T_X$,
we associate to $V$ a $k$-jet bundle $J_kV$ as follows, assuming 
$V$ {\em non-singular} throughout subsection 6.A.

\claim 6.1.~Definition|Let $(X,V)$ be a complex directed manifold.
We define $J_kV\to X$ to be the bundle of $k$-jets of curves
\hbox{$f:(\bC,0)\to X$} which are tangent to $V$, i.e., such that
$f'(t)\in V_{f(t)}$ for all $t$ in a neighborhood of~$0$, together with
the projection map $f\mapsto f(0)$ onto~$X$.
\endclaim

It is easy to check that $J_kV$ is actually a subbundle of~$J_kX$. In
fact, by using (5.8) and (5.10), we see that the fibers $J_kV_x$ are
parametrized by 
$$
\big((f_1'(0)\ld f_r'(0));(f_1''(0)\ld f_r''(0));\ldots;
(f_1^{(k)}(0)\ld f_r^{(k)}(0))\big)\in(\bC^r)^k
$$
for all $x\in\Omega$, and hence $J_kV$ is a locally trivial
$(\bC^r)^k$-subbundle of~$J_kX$. Alternatively, we can pick a
local holomorphic connection $\nabla$ on $V$ such that for
any germs $w=\sum_{1\leq j\leq n}w_j{\partial\over \partial z_j}
\in\cO(T_{X,x})$ and $v=\sum_{1\leq\lambda\leq r}v_\lambda e_\lambda\in\cO(V)_x$
in a local trivializing frame $(e_1,\ldots,e_r)$ of~$V_{\restriction \Omega}$ we have
$$
\nabla_wv(x)=\sum_{1\leq j\leq n,\,1\leq\lambda\leq r}
w_j{\partial v_\lambda\over\partial z_j}e_\lambda(x)+
\sum_{1\leq j\leq n,\,1\leq\lambda,\mu\leq r}
\Gamma_{j\lambda}^\mu(x)w_jv_\lambda\,e_\mu(x).
$$
We can of course take the frame obtained from (5.8) by lifting
the vector fields $\partial/\partial z_1,\ldots,\partial/\partial z_r$,
and the ``trivial connection'' given by the zero Christoffel symbols
$\Gamma=0$. One then obtains a
trivialization $J^kV_{\restriction \Omega}\simeq V^{\oplus k}_{\restriction \Omega}$ by
considering
$$
J_kV_x\ni f\mapsto(\xi_1,\xi_2,\ldots,\xi_k)=(\nabla f(0),\nabla^2f(0),\ldots,\nabla^k f(0))\in V_x^{\oplus k}
$$
and computing inductively the successive derivatives $\nabla f(t)=f'(t)$ and
$\nabla^sf(t)$ via
$$
\nabla^sf=(f^*\nabla)_{d/dt}(\nabla^{s-1}f)
=\sum_{1\leq\lambda\leq r}
{d\over dt}\Big(\nabla^{s-1}f\Big)_\lambda e_\lambda(f)+
\sum_{1\leq j\leq n,\,1\leq\lambda,\mu\leq r}
\Gamma_{j\lambda}^\mu(f)f'_j\Big(\nabla^{s-1}f\Big)_\lambda e_\mu(f).
$$
This identification depends of course on the choice of $\nabla$ and cannot 
be defined globally in general (unless we are in the rare situation
where $V$ has a global holomorphic connection.\qed

We now describe a convenient process for constructing ``projectivized
jet bundles'', which will later appear as natural quotients of our jet
bundles~$J_kV$ (or~rather, as suitable desingularized compactifications
of the quotients). Such spaces have already been considered since
a long time, at least in the special case $X=\bP^2$, $V=T_{\bP^2}$ (see
Gherardelli [Ghe41], Semple [Sem54]), and they have been mostly used as a 
tool for establishing enumerative formulas dealing with the order of contact
of plane curves (see [Coll88], [CoKe94]); the article [ASS97] is also
concerned with such generalizations of jet bundles, as well as  [LaTh96] 
by Laksov and Thorup.

We define inductively the {\em projectivized $k$-jet bundle $X_k$}
(or {\em Semple $k$-jet bundle}) and the associated subbundle 
$V_k\subset T_{X_k}$ by
$$
(X_0,V_0)=(X,V),\qquad (X_k,V_k)=(\swt{X}_{k-1},\swt{V}_{k-1}).
\leqno(6.2)
$$
In other words, $(X_k,V_k)$ is obtained from $(X,V)$ by iterating
$k$-times the lifting construction $(X,V)\mapsto(\swt{X},\swt{V})$ described
in~\S$\,$5. By (5.2--5.7), we find
$$
\dim X_k=n+k(r-1),\qquad\rank V_k=r,\leqno(6.3)
$$
together with exact sequences
$$
\leqalignno{
&0\lra T_{X_k/X_{k-1}}\lra V_k\build{\vlra{11}}^{(\pi_k)_*}_{}
\cO_{X_k}(-1)\lra 0,&(6.4)\cr
&0\lra\cO_{X_k}\lra\pi_k^* V_{k-1}\otimes\cO_{X_k}(1)
\lra T_{X_k/X_{k-1}}\lra 0,
&(6.4')\cr}
$$
where $\pi_k$ is the natural projection $\pi_k:X_k\to X_{k-1}$ and
$(\pi_k)_*$ its differential. Formula (5.4) yields
$$
\det V_k=\pi_k^*\det V_{k-1}\otimes\cO_{X_k}(r-1).\leqno(6.5)
$$
Every non-constant tangent trajectory
\hbox{$f:\Delta_R\to X$} of $(X,V)$ lifts to a well defined and
unique tangent trajectory $f_{[k]}:\Delta_R\to X_k$ of $(X_k,V_k)$.
Moreover, the derivative $f_{[k-1]}'$ gives rise to a section
$$
f_{[k-1]}':T_{\Delta_R}\to f_{[k]}^*\cO_{X_k}(-1).\leqno(6.6)
$$
In coordinates, one can compute $f_{[k]}$ in terms of its components in
the various affine charts (5.9) occurring at each step: we get inductively
$$
f_{[k]}=(F_1\ld F_N),\qquad f_{[k+1]}=\Big(F_1\ld F_N,
{F_{s_1}'\over F_{s_r}'}\ld{F_{s_{r-1}}'\over F_{s_r}'}\Big)\leqno(6.7)
$$
where $N=n+k(r-1)$ and $\{s_1\ld s_r\}\subset\{1\ld N\}$. If $k\ge 1$,
$\{s_1\ld s_r\}$ contains the last $r-1$ indices of $\{1\ld N\}$
corresponding to the ``vertical'' components of the projection
$X_k\to X_{k-1}$, and in general, $s_r$ is an index such that
$m(F_{s_r},0)=m(f_{[k]},0)$, that is, $F_{s_r}$ has the smallest
vanishing order among all components $F_s$ ($s_r$ may be vertical
or not, and the choice of $\{s_1\ld s_r\}$ need not be unique).

By definition, there is a canonical injection $\cO_{X_k}(-1)
\hookrightarrow\pi_k^* V_{k-1}$, and a composition with the
projection $(\pi_{k-1})_*$ (analogue for order $k-1$ of the
arrow~$(\pi_k)_*$ in sequence (6.4)) yields for any $k\ge 2$
a natural line bundle morphism
$$
\cO_{X_k}(-1)\lhra\pi_k^* V_{k-1}\build{\vlra{16}}^
{(\pi_k)^*(\pi_{k-1})_*}_{}\pi_k^*\cO_{X_{k-1}}(-1),\leqno(6.8)
$$
which admits precisely \hbox{$D_k=P(T_{X_{k-1}/X_{k-2}})
\subset P(V_{k-1})=X_k$} as its zero divisor (clearly, $D_k$ is a
hyperplane subbundle of~$X_k$). Hence we find
$$
\cO_{X_k}(1)=\pi_k^*\cO_{X_{k-1}}(1)\otimes\cO(D_k).\leqno(6.9)
$$
Now, we consider the composition of projections
$$
\pi_{j,k}=\pi_{j+1}\circ\cdots\circ\pi_{k-1}\circ\pi_k:X_k\lra X_j.
\leqno(6.10)
$$
Then $\pi_{0,k}:X_k\to X_0=X$ is a locally trivial holomorphic fiber bundle
over~$X$, and the fibers $X_{k,x}=\pi_{0,k}^{-1}(x)$ are $k$-stage
towers of $\bP^{r-1}$-bundles. Since we have (in both directions) morphisms 
$(\bC^r,T_{\bC^r})\leftrightarrow(X,V)$ of directed manifolds which
are bijective on the level of bundle morphisms, the fibers are
all isomorphic to a ``universal'' non-singular projective algebraic 
variety of dimension $k(r-1)$ which we will denote by~$\cR_{r,k}\,$; it is
not hard to see that $\cR_{r,k}$ is rational (as will indeed follow from
the proof of Theorem~7.11 below). 

\plainsubsection 6.B. Semple tower of singular directed varieties|

Let $(X,V)$ be a directed variety. We assume $X$ non-singular, but here $V$ is
allowed to have singularities. We are going to give a natural definition
of the Semple tower$(X_k,V_k)$ in that case.

Let us take $X'=X\smallsetminus \Sing(V)$ and $V'=V_{\restriction X'}$. By subsection 6.A, we have a well defined Semple tower $(X'_k,V'_k)$ over the
Zariski open set $X'$. We also have an ``absolute''
Semple tower $(X^a_k,V^a_k)$ obtained from $(X^a_0,V^a_0)=(X,T_X)$,
which is non-singular. The injection $V'\subset T_X$ induces by functoriality
(cf.\ (5.5)) an injection 
$$(X'_k,V'_k)\subset (X^a_k,V^a_k).\leqno(6.11)$$
\claim 6.12. Definition|Let $(X,V)$ be a directed variety, with $X$
non-singular. When $\Sing(V)\neq\emptyset$, we
define $X_k$ and $V_k$ to be the respective closures of $X'_k$, $V'_k$ 
associated with $X'=X\smallsetminus \Sing(V)$ and $V'=V_{\restriction X'}$, where
the closure is taken in the non-singular absolute
Semple tower $(X^a_k,V^a_k)$ obtained from $(X^a_0,V^a_0)=(X,T_X)$.
\endclaim

We leave the reader check that the following functoriality property still 
holds.

\claim 6.13. Fonctoriality|If $\Phi:(X,V)\to(Y,W)$ is a morphism of
directed varieties such that $\Phi_*:T_X\to\Phi^*T_Y$ is injective
$($i.e.\ $\Phi$ is an immersion$\,)$, then there is a corresponding
natural morphism $\Phi_{[k]}:(X_k,V_k)\to(Y_k,W_k)$
at the level of Semple bundles. If one merely assumes that the differential
$\Phi_*:V\to\Phi^*W$ is non-zero, there is still a natural
meromorphic map $\Phi_{[k]}:(X_k,V_k)\merto(Y_k,W_k)$ for all~$k\ge 0$.
\endclaim

In case $V$ is singular, the $k$-th stage $X_k$ of the Semple tower will also
be singular, but we can replace $(X_k,V_k)$ by a suitable modification 
$(\widehat X_k,\widehat V_k)$ if we want to work with 
a non-singular model $\widehat X_k$ of~$X_k$. The exceptional
set of $\widehat X_k$ over $X_k$ can be chosen to lie above
$\Sing(V)\subset X$, and proceeding inductively with respect to $k$, we can 
also arrange the modifications in such a way that we get a tower 
structure $(\widehat X_{k+1},\widehat V_{k+1})\to
(\widehat X_k,\widehat V_k)\,$; however,
in general, it will not be possible to achieve that 
$\widehat V_k$ is a subbundle of $T_{\widehat X_k}$.

\section{Jet differentials}

\plainsubsection 7.A. Green-Griffiths jet differentials|

We first introduce the concept of jet differentials in the sense of Green-Griffiths 
[GrGr80]. The goal is to provide an intrinsic geometric description of
holomorphic differential equations that a germ of curve $f:(\bC,0)\to X$
may satisfy. In the sequel, we fix a directed manifold $(X,V)$ and suppose
implicitly that all germs of curves $f$ are tangent to~$V$. 

Let $\bG_k$ be the group of germs of $k$-jets of biholomorphisms of $(\bC,0)$,
that is, the group of germs of biholomorphic maps
$$
t\mapsto\varphi(t)=a_1t+a_2t^2+\cdots+a_kt^k,\qquad
a_1\in\bC^*,~a_j\in\bC,~j\ge 2,
$$
in which the composition law is taken modulo terms $t^j$ of degree $j>k$.
Then $\bG_k$ is a $k$-dimensional nilpotent complex Lie group,
which admits a natural fiberwise right action on $J_kV$. The action
consists of reparametrizing $k$-jets of maps $f:(\bC,0)\to X$
by a biholomorphic change of parameter $\varphi:(\bC,0)\to(\bC,0)$, that is,
$(f,\varphi)\mapsto f\circ\varphi$. There is an exact sequence of groups
$$
1\to \bG'_k\to \bG_k\to\bC^*\to 1,
$$
where $\bG_k\to\bC^*$ is the obvious morphism $\varphi\mapsto\varphi'(0)$,
and $\bG'_k=[\bG_k,\bG_k]$ is the group of $k$-jets of biholomorphisms tangent
to the identity. Moreover, the subgroup $\bH\simeq\bC^*$ of homotheties
$\varphi(t)=\lambda t$ is a (non-normal) subgroup of $\bG_k$, and we have a
semi-direct decomposition $\bG_k=\bG'_k\ltimes\bH$. The corresponding
action on $k$-jets is described in coordinates by
$$
\lambda\cdot(f',f'',\ldots,f^{(k)})=
(\lambda f',\lambda^2f'',\ldots,\lambda^kf^{(k)}).
$$

Following [GrGr80], we introduce the vector bundle $E^\GG_{k,m}V^*\to X$
whose fibers are complex valued polynomials $Q(f',f'',\ldots,f^{(k)})$ on
the fibers of $J_kV$, of weighted degree $m$ with respect to the
$\bC^*$ action defined by $\bH$, that is, such that
$$
Q(\lambda f',\lambda^2 f'',\ldots,\lambda^k f^{(k)})=\lambda^m
Q(f',f'',\ldots,f^{(k)})\leqno(7.1)
$$
for all $\lambda\in\bC^*$ and $(f',f'',\ldots,f^{(k)})\in J_kV$.
Here we view $(f',f'',\ldots,f^{(k)})$ as indeterminates with components
$$
\big((f_1'\ld f_r');(f_1''\ld f_r'');\ldots;(f_1^{(k)}\ld f_r^{(k)})\big)
\in(\bC^r)^k.
$$
Notice that the concept of polynomial on the fibers of $J_kV$ makes sense,
for all coordinate changes $z\mapsto w=\Psi(z)$ on $X$ induce polynomial
transition automorphisms on the fibers of $J_kV$, given by a formula
$$
(\Psi\circ f)^{(j)}=\Psi'(f)\cdot f^{(j)}+\sum_{s=2}^{s=j}{~}
\sum_{j_1+j_2+\cdots+j_s=j}c_{j_1\ldots j_s}\Psi^{(s)}(f)\cdot
(f^{(j_1)}\ld f^{(j_s)})\leqno(7.2)
$$
with suitable integer constants $c_{j_1\ldots j_s}$ (this is easily
checked by induction on~$s$). In the ``absolute case'' $V=T_X$, we simply
write $E^\GG_{k,m}T^*_X=E^\GG_{k,m}$. If \hbox{$V\subset V'\subset V^a:=T_X$} 
are holomorphic subbundles, there are natural inclusions
$$
J_kV\subset J_kV'\subset J_kV^a,\qquad X_k\subset X_k'\subset X_k^a.
$$
The restriction morphisms induce surjective arrows
$$
E^\GG_{k,m}T_X^*\to E^\GG_{k,m}V^{\prime*}\to E^\GG_{k,m}V^*,
$$
and in particular $E^\GG_{k,m}V^*$ can be seen as a quotient of
$E^\GG_{k,m}T_X^*$.
(The notation $V^*$ is used here to make the contravariance property
implicit from the notation). Another useful consequence of these inclusions
is that one can extend the
definition of $J_kV$ and $X_k$ to the case where $V$ is an arbitrary
linear space, simply by taking the closure of $J_kV_{X\ssm\Sing(V)}$
and $X_{k|X\ssm\Sing(V)}$ in the smooth bundles 
$J_kX$ and $X_k^a$, respectively.

If $Q\in E^\GG_{k,m}V^*$ is decomposed into multihomogeneous
components of multidegree $(\ell_1,\ell_2\ld\ell_k)$ in $f',f''\ld
f^{(k)}$ (the decomposition is of course coordinate dependent), these
multidegrees must satisfy the relation
$$
\ell_1+2\ell_2+\cdots+k\ell_k=m.
$$
The bundle $E^\GG_{k,m}V^*$ will be called the {\em bundle of jet
differentials of order $k$ and weighted degree~$m$}. It is clear from (7.2)
that a coordinate change $f\mapsto\Psi\circ f$ transforms every monomial
$(f^{(\bullet)})^\ell=(f')^{\ell_1}(f'')^{\ell_2}\cdots(f^{(k)})^{\ell_k}$
of partial weighted degree $|\ell|_s:=\ell_1+2\ell_2+\cdots+s\ell_s$,
$1\le s\le k$, into a polynomial $((\Psi\circ f)^{(\bullet)})^\ell$ in
$(f',f''\ld f^{(k)})$ whose non-zero monomials have the same partial
weighted degree of order
$s$ if $\ell_{s+1}=\cdots=\ell_k=0$, and a larger or equal partial degree
of order $s$ otherwise. Hence, for each $s=1\ld k$, we get a well defined
(i.e., coordinate invariant) decreasing filtration $F_s^\bullet$ on
$E^\GG_{k,m}V^*$ as follows:
$$
F^p_s(E^\GG_{k,m}V^*)=\left\{
{\displaystyle
\hbox{\rm $Q(f',f''\ld f^{(k)})\in E^\GG_{k,m}V^*$ involving}
\atop\hbox{\rm only monomials $(f^{(\bullet)})^\ell$ with $|\ell|_s\ge p$}
\hfill}\right\},\qquad
\forall p\in\bN.\leqno(7.3)
$$
The graded terms $\Gr^p_{k-1}(E^\GG_{k,m}V^*)$ associated with the
filtration $F^p_{k-1}(E^\GG_{k,m}V^*)$ are precisely the
homogeneous polynomials $Q(f'\ld f^{(k)})$ whose monomials
$(f^{\bullet})^\ell$ all have partial weighted degree $|\ell|_{k-1}=p$
(hence their degree $\ell_k$ in~$f^{(k)}$ is such that $m-p=k\ell_k$,
and $\Gr^p_{k-1}(E^\GG_{k,m}V^*)=0$ unless $k|m-p$).  The transition
automorphisms of the graded bundle are induced by coordinate changes
$f\mapsto\Psi\circ f$, and they are described by substituting the
arguments of~$Q(f'\ld f^{(k)})$ according to formula (7.2), namely
$f^{(j)}\mapsto(\Psi\circ f)^{(j)}$ for $j<k$, and
$f^{(k)}\mapsto\Psi'(f)\circ f^{(k)}$ for $j=k$ (when $j=k$, the other
terms fall in the next stage $F^{p+1}_{k-1}$ of the filtration).
Therefore $f^{(k)}$ behaves as an element of $V\subset T_X$ under
coordinate changes.  We thus find
$$
G_{k-1}^{m-k\ell_k}(E^\GG_{k,m}V^*)=E^\GG_{k-1,m-k\ell_k}V^*\otimes
S^{\ell_k}V^*.
\leqno(7.4)
$$
Combining all filtrations $F_s^\bullet$ together, we find inductively a
filtration $F^\bullet$ on $E^\GG_{k,m}V^*$ such that the graded terms are
$$
\Gr^\ell(E^\GG_{k,m}V^*)=S^{\ell_1}V^*\otimes S^{\ell_2}V^*\otimes
\cdots\otimes S^{\ell_k}V^*,\qquad\ell\in\bN^k,\quad
|\ell|_k=m.\leqno(7.5)
$$

The bundles $E^\GG_{k,m}V^*$ have other interesting properties. In fact,
$$
E^\GG_{k,\bu}V^*:=\bigoplus_{m\ge 0}E^\GG_{k,m}V^*
$$
is in a natural way a bundle of graded algebras (the product is
obtained simply by taking the product of polynomials). There are
natural inclusions \hbox{$E^\GG_{k,\bu}V^*\subset E^\GG_{k+1,\bu}
V^*$} of algebras, and hence $E^\GG_{\infty,\bu}V^*=\bigcup_{k\ge 0}
E^\GG_{k,\bu}V^*$ is also an algebra. Moreover, the sheaf of
holomorphic sections $\cO(E^\GG_{\infty,\bu} V^*)$ admits a
canonical derivation $D^\GG$ given by a collection of $\bC$-linear maps
$$
D^\GG:\cO(E^\GG_{k,m}V^*)\to\cO(E^\GG_{k+1,m+1}V^*),
$$
constructed in the following way. A holomorphic section of
$E^\GG_{k,m}V^*$ on a coordinate open set $\Omega\subset X$ can be
seen as a differential operator on the space of germs
$f:(\bC,0)\to\Omega$ of the form
$$
Q(f)=\sum_{|\alpha_1|+2|\alpha_2|+\cdots+k|\alpha_k|=m}
a_{\alpha_1\ldots\alpha_k}(f)\,(f')^{\alpha_1}(f'')^{\alpha_2}\cdots
(f^{(k)})^{\alpha_k}\leqno(7.6)
$$
in which the coefficients $a_{\alpha_1\ldots\alpha_k}$ are holomorphic
functions on $\Omega$. Then $D^\GG Q$ is given by the formal derivative
$(D^\GG Q)(f)(t)=d(Q(f))/dt$ with respect to the
$1$-dimensional parameter $t$ in~$f(t)$. For example, in dimension 2,
if $Q\in H^0(\Omega,\cO(E^\GG_{2,4}))$ is the section of weighted
degree $4$
$$
Q(f)=a(f_1,f_2)\,f_1^{\prime 3}f_2'+b(f_1,f_2)\,f_1^{\prime\prime 2},
$$
we find that $D^\GG Q\in H^0(\Omega,\cO(E^\GG_{3,5}))$ is given by
$$
\eqalign{(D^\GG Q)&(f)=
{\partial a\over\partial z_1}(f_1,f_2)\,f_1^{\prime 4}f_2'+
{\partial a\over\partial z_2}(f_1,f_2)\,f_1^{\prime 3}f_2^{\prime 2}+
{\partial b\over\partial z_1}(f_1,f_2)\,f_1'f_1^{\prime\prime 2}\cr
&{}+{\partial b\over\partial z_2}(f_1,f_2)\,f_2'f_1^{\prime\prime 2}
+a(f_1,f_2)\,\big(3f_1^{\prime 2}f_1''f_2'+f_1^{\prime 3}f_2'')+
b(f_1,f_2)\,\,2f_1''f_1'''.\cr}
$$
Associated with the graded algebra bundle $E^\GG_{k,\bu}V^*$, we 
have an analytic fiber bundle
$$X_k^\GG:=\Proj(E^\GG_{k,\bu}V^*)=(J_kV\ssm\{0\})/\bC^*\leqno(7.7)$$
over $X$, which has weighted projective spaces
$\bP(1^{[r]},2^{[r]}\ld k^{[r]})$ as fibers (these
weighted projective spaces are singular for $k>1$, but they only have
quotient singularities, see [Dol81]$\,$; here $J_kV\ssm\{0\}$ is the set
of non-constant jets of order~$k\,$; we refer e.g.\ to Hartshorne's book
[Har77] for a definition of the $\Proj$ functor). 
As such, it possesses a canonical sheaf $\cO_{X^\GG_k}(1)$ such that
$\cO_{X^\GG_k}(m)$ is invertible when $m$ is a multiple of
$\lcm(1,2,\ldots,k)$. Under the natural projection 
$\pi_k:X^\GG_k\to X$, the direct image 
$(\pi_k)_*\cO_{X^\GG_k}(m)$ coincides with polynomials
$$
P(z\,;\,\xi_1,\ldots,\xi_k)=\sum_{\alpha_\ell\in\bN^r,\,1\le\ell\le k} 
a_{\alpha_1\ldots\alpha_k}(z)\,\xi_1^{\alpha_1}\ldots\xi_k^{\alpha_k}
\leqno(7.8)
$$
of weighted degree $|\alpha_1|+2|\alpha_2|+\ldots+k|\alpha_k|=m$ on 
$J^kV$ with holomorphic coefficients; in other words, we obtain precisely
the sheaf of sections of the bundle $E_{k,m}^\GG V^*$ of jet differentials 
of order $k$ and degree $m$.

\claim 7.9. Proposition|By construction, if $\pi_k:X_k^\GG\to X$ is the
natural projection, we have the direct image formula
$$
(\pi_k)_*\cO_{X^\GG_k}(m)=\cO(E_{k,m}^\GG V^*)
$$
for all $k$ and $m$.
\endclaim

\plainsubsection 7.B. Invariant jet differentials|

In the geometric context, we are not really interested in the 
bundles $(J_kV\ssm\{0\})/\bC^*$ themselves, but
rather in their quotients $(J_kV\ssm\{0\})/\bG_k$ (would such nice complex
space quotients exist!). We will see that the Semple bundle $X_k$ 
constructed in \S$\,$6 plays the role of such a quotient. First we 
introduce a canonical subalgebra of the bundle algebra~$E^\GG_{k,\bu}V^*$.

\claim 7.10.~Definition|We introduce a subbundle
$E_{k,m}V^*\subset E^\GG_{k,m}V^*$, called the bundle of
invariant jet differentials of order~$k$ and degree $m$, defined as
follows: $E_{k,m}V^*$ is the set of polynomial differential
operators $Q(f',f'',\ldots,f^{(k)})$ which are invariant under
arbitrary changes of parametrization, i.e., for every $\varphi\in \bG_k$
$$
Q\big((f\circ\varphi)',(f\circ\varphi)'',\ldots,(f\circ\varphi)^{(k)})=
\varphi'(0)^m Q(f',f'',\ldots,f^{(k)}).
$$
\endclaim

Alternatively, $E_{k,m}V^*=(E^\GG_{k,m}V^*)^{\bG'_k}$ is the set of
invariants of $E^\GG_{k,m}V^*$ under the action of~$\bG'_k$. Clearly,
$E_{\infty,\bu}V^*=\bigcup_{k\ge 0}\bigoplus_{m\ge0}E_{k,m}V^*$
is a subalgebra of $E^\GG_{k,m}V^*$ (observe however that this algebra is
not invariant under the derivation~$D^\GG$, since e.g.\ 
$f_j''=D^\GG f_j$ is not an invariant polynomial). 
In addition to this, there are natural induced filtrations
$F^p_s(E_{k,m}V^*)=E_{k,m}V^*\cap F^p_s(E^\GG_{k,m}V^*)$
(all locally trivial over~$X$). These induced filtrations will play an
important role later on.

\claim 7.11.~Theorem|Suppose that $V$ has rank $r\ge 2$. Let $\pi_{0,k}:
X_k\lra X$ be the Semple jet bundles constructed in section~$6$,
and let $J_kV^\reg$ be the bundle of regular $k$-jets of maps
$f:(\bC,0)\to X$, that is, jets $f$ such that $f'(0)\ne 0$. 
\plainitem{\rm(i)} The quotient $J_kV^\reg/\bG_k$ has the structure of a
locally trivial bundle over~$X$, and there is a holomorphic embedding
$J_kV^\reg/\bG_k\hookrightarrow X_k$ over $X$, which identifies
$J_kV^\reg/\bG_k$ with $X_k^\reg$ $($thus $X_k$ is a relative
compactification of $J_kV^\reg/\bG_k$ over~$X)$.
\plainitem{\rm(ii)} The direct image sheaf
$$
(\pi_{0,k})_*\cO_{X_k}(m)\simeq\cO(E_{k,m}V^*)
$$
can be identified with the sheaf of holomorphic sections of
$E_{k,m}V^*$.
\plainitem{\rm(iii)} For every $m>0$, the relative base locus of the linear
system $|\cO_{X_k}(m)|$ is equal to the set $X_k^\sing$
of singular $k$-jets. Moreover, $\cO_{X_k}(1)$ is relatively big
over~$X$.
\vskip0pt
\endclaim

\plainproof. (i) For $f\in J_kV^\reg$, the lifting $\swt{f}$ is obtained by taking
the derivative $(f,[f'])$ without any cancellation of zeroes in~$f'$,
and hence we get a uniquely defined $(k-1)$-jet $\swt{f}:(\bC,0)\to\swt{X}$.
Inductively, we get a well defined $(k-j)$-jet $f_{[j]}$ in~$X_j$, and
the value $f_{[k]}(0)$ is independent of the choice of the
representative $f$ for the $k$-jet. As the lifting process commutes
with reparametrization, i.e., $(f\circ\varphi)^\sim=\swt{f}\circ\varphi$
and more generally $(f\circ\varphi)_{[k]}=f_{[k]}\circ\varphi$, we
conclude that there is a well defined set-theoretic map
$$
J_kV^\reg/\bG_k\to X_k^\reg,\qquad f~\mod~\bG_k\mapsto f_{[k]}(0).
$$
This map is better understood in coordinates as follows. Fix
coordinates $(z_1\ld z_n)$ near a point $x_0\in X$, such that
$V_{x_0}=\Vect(\partial/\partial z_1\ld\partial/\partial z_r)$. Let
$f=(f_1\ld f_n)$ be a regular $k$-jet tangent to~$V$. Then there exists
$i\in\{1,2\ld r\}$ such that $f_i'(0)\ne 0$, and there is a unique
reparametrization $t=\varphi(\tau)$ such that $f\circ\varphi=g
=(g_1,g_2\ld g_n)$ with $g_i(\tau)=\tau$ (we just express the curve as
a graph over the $z_i$-axis, by means of a change of parameter $\tau=f_i(t)$,
i.e.\ $t=\varphi(\tau)=f_i^{-1}(\tau)$). Suppose $i=r$ for the simplicity of
notation. The space $X_k$ is a $k$-stage tower of $\bP^{r-1}$-bundles.
In~the corresponding inhomogeneous coordinates on these $\bP^{r-1}$'s,
the point $f_{[k]}(0)$ is given by the collection of derivatives
$$
\big((g_1'(0)\ld g_{r-1}'(0));(g_1''(0)\ld g_{r-1}''(0));\ldots;
(g_1^{(k)}(0)\ld g_{r-1}^{(k)}(0))\big).
$$
[Recall that the other components $(g_{r+1}\ld g_n)$ can be recovered from
$(g_1\ld g_r)$ by integrating the differential system (5.10)].
Thus the map $J_kV^\reg/\bG_k\to X_k$ is a bijection onto $X_k^\reg$, and the
fibers of these isomorphic bundles can be seen as unions of $r$ affine
charts ${}\simeq(\bC^{r-1})^k$, associated with each choice of the axis
$z_i$ used to describe the curve as a graph. The change of parameter formula
${d\over d\tau}={1\over f_r'(t)}{d\over dt}$ expresses all derivatives
$g_i^{(j)}(\tau)=d^jg_i/d\tau^j$ in terms of the derivatives
$f_i^{(j)}(t)=d^jf_i/dt^j$
$$
\leqalignno{
(g_1'\ld g_{r-1}')&=\Big({f'_1\over f'_r}\ld{f'_{r-1}\over f'_r}\Big);\cr
(g_1''\ld g_{r-1}'')&=\Big({f''_1f'_r-f''_rf'_1\over f^{\prime 3}_r}\ld
{f''_{r-1}f'_r-f''_rf'_{r-1}\over f^{\prime 3}_r}\Big);~\ldots~;&(7.12)\cr
\qquad(g_1^{(k)}\ld g_{r-1}^{(k)})&=\Big({f^{(k)}_1f'_r-f^{(k)}_rf'_1\over
f^{\prime k+1}_r}\ld{f^{(k)}_{r-1}f'_r-f^{(k)}_rf'_{r-1}\over
f^{\prime k+1}_r}\Big)+(\hbox{order}<k).\cr}
$$
Also, it is easy to check that $f_r^{\prime 2k-1}g_i^{(k)}$ is an
invariant polynomial in $f'$, $f''\ld f^{(k)}$ of total degree $2k-1$, i.e.,
a section of $E_{k,2k-1}$. 

\noindent (ii) Since the bundles $X_k$ and $E_{k,m}V^*$ are both locally
trivial over $X$, it is sufficient to identify sections $\sigma$ of
$\cO_{X_k}(m)$ over a fiber $X_{k,x}=\pi_{0,k}^{-1}(x)$ with the fiber
$E_{k,m}V^*_x$, at any point  $x\in X$. Let $f\in J_kV_x^\reg$
be a regular $k$-jet at~$x$. By (6.6), the derivative $f_{[k-1]}'(0)$
defines an element of the fiber of $\cO_{X_k}(-1)$ at $f_{[k]}(0)\in X_k$.
Hence we get a well defined complex valued operator
$$
Q(f',f''\ld f^{(k)})=\sigma(f_{[k]}(0))\cdot(f_{[k-1]}'(0))^m.
\leqno(7.13)
$$
Clearly, $Q$ is holomorphic on $J_kV_x^\reg$ (by the holomorphicity
of~$\sigma$), and the $\bG_k$-invariance condition of Definition~7.10 
is satisfied
since $f_{[k]}(0)$ does not depend on reparametrization and
$$(f\circ\varphi)_{[k-1]}'(0)=f_{[k-1]}'(0)\varphi'(0).$$
Now, $J_kV_x^\reg$ is the complement of a
linear subspace of codimension $r$ in $J_kV_x$, and hence $Q$ extends
holomorphically to all of $J_kV_x\simeq(\bC^r)^k$ by Riemann's
extension theorem (here we use the hypothesis $r\ge 2\,$; if $r=1$, the
situation is anyway not interesting since $X_k=X$ for all~$k$). Thus
$Q$ admits an everywhere convergent power series
$$
Q(f',f''\ld f^{(k)})=\sum_{\alpha_1,\alpha_2\ld\alpha_k\in\bN^r}
a_{\alpha_1\ldots\alpha_k}\,(f')^{\alpha_1}(f'')^{\alpha_2}\cdots
(f^{(k)})^{\alpha_k}.
$$
The $\bG_k$-invariance asserted in Definition~7.10 implies in
particular that $Q$ must be
multihomogeneous in the sense of (7.1), and thus $Q$ must be a
polynomial. We conclude that $Q\in E_{k,m}V^*_x$, as desired.

Conversely, for all $w$ in a neighborhood of any given
point $w_0\in X_{k,x}$, we can find a holomorphic family
of germs $f_w:(\bC,0)\to X$ such that $(f_w)_{[k]}(0)=w$ and
$(f_w)_{[k-1]}'(0)\ne 0$ (just take the projections to $X$ of
integral curves of $(X_k,V_k)$
integrating a nonvanishing local holomorphic section of $V_k$ near~$w_0$).
Then every $Q\in E_{k,m}V^*_x$ yields a
holomorphic section $\sigma$ of $\cO_{X_k}(m)$ over the fiber $X_{k,x}$
by putting
$$
\sigma(w)=Q(f_w',f_w''\ld f_w^{(k)})(0)\,\big((f_w)_{[k-1]}'(0)\big)^{-m}.
\leqno(7.14)
$$

\noindent (iii) By what we saw in (i)--(ii), every section $\sigma$ of
$\cO_{X_k}(m)$ over the fiber $X_{k,x}$ is given by a polynomial
$Q\in E_{k,m}V^*_x$, and this polynomial can be expressed
on the Zariski open chart $f'_r\ne 0$ of $X_{k,x}^\reg$ as
$$
Q(f',f''\ld f^{(k)})=f_r^{\prime m}\swh{Q}(g',g''\ld g^{(k)}),
\leqno(7.15)
$$
where $\swh{Q}$ is a polynomial and $g$ is the reparametrization of $f$ such
that $g_r(\tau)=\tau$. In fact $\swh{Q}$ is obtained from $Q$ by substituting
$f'_r=1$ and $f^{(j)}_r=0$ for~$j\ge 2$, and conversely $Q$ can be recovered
easily from $\swh{Q}$ by using the substitutions (7.12).

In this context, the jet differentials $f\mapsto f'_1\ld f\mapsto f'_r$ can
be viewed as sections of $\cO_{X_k}(1)$ on a neighborhood
of the fiber $X_{k,x}$. Since these sections vanish exactly
on $X_k^\sing$, the relative base locus of $\cO_{X_k}(m)$ is
contained in $X_k^\sing$ for every~$m>0$. We see that $\cO_{X_k}(1)$ is 
big by considering the sections of $\cO_{X_k}(2k-1)$ associated with the
polynomials $Q(f'\ld f^{(k)})=f_r^{\prime 2k-1}g_i^{(j)}$, $1\le i\le r-1$,
$1\le j\le k$; indeed, these sections separate all points in the open
chart $f_r'\ne 0$ of~$X_{k,x}^\reg$. 

Now, we check that every section $\sigma$ of $\cO_{X_k}(m)$ over $X_{k,x}$
must vanish on $X_{k,x}^\sing$. Pick an arbitrary element $w\in X_k^\sing$
and a germ of curve \hbox{$f:(\bC,0)\to X$} such that $f_{[k]}(0)=w$,
$f_{[k-1]}'(0)\ne 0$ and $s=m(f,0)\gg 0$ (such an $f$ exists by
[Dem95, Corollary~5.14]). There are local coordinates $(z_1\ld z_n)$ on $X$ such
that $f(t)=(f_1(t)\ld f_n(t))$ where $f_r(t)=t^s$. Let $Q$, $\swh{Q}$ be 
the polynomials associated with $\sigma$ in these coordinates and let
$(f')^{\alpha_1}(f'')^{\alpha_2}\cdots (f^{(k)})^{\alpha_k}$ be a
monomial occurring in $Q$, with $\alpha_j\in\bN^r$, $|\alpha_j|=\ell_j$,
\hbox{$\ell_1+2\ell_2+\cdots+k\ell_k=m$}. Putting $\tau=t^s$, the curve
$t\mapsto f(t)$ becomes a Puiseux expansion $\tau\mapsto
g(\tau)=(g_1(\tau)\ld g_{r-1}(\tau),\tau)$ in which $g_i$ is a power
series in~$\tau^{1/s}$, starting with exponents of $\tau$ at least
equal to~$1$. The derivative $g^{(j)}(\tau)$ may involve negative
powers of $\tau$, but the exponent is always${}\ge 1+{1\over s}-j$.
Hence the Puiseux expansion of $\swh{Q}(g',g''\ld g^{(k)})$ can only
involve powers of $\tau$ of exponent
\hbox{$\ge~{}-\max_\ell((1-{1\over s})\ell_2+\cdots+(k-1-{1\over s})\ell_k)$}.
Finally $f_r'(t)=st^{s-1}=s\tau^{1-1/s}$, and so the lowest exponent of
$\tau$ in $Q(f'\ld f^{(k)})$ is at least equal to
$$
\eqalign{
\Big(1-{1\over s}\Big)m-\max_\ell\Big(&\Big(1-{1\over s}\Big)\ell_2+\cdots+
\Big(k-1-{1\over s}\Big)\ell_k\Big)\cr
&\ge\min_\ell\Big(1-{1\over s}\Big)\ell_1+\Big(1-{1\over s}\Big)\ell_2
+\cdots+\Big(1-{k-1\over s}\Big)\ell_k,\cr}
$$
where the minimum is taken over all monomials
$(f')^{\alpha_1}(f'')^{\alpha_2}\cdots (f^{(k)})^{\alpha_k}$,
$|\alpha_j|=\ell_j$, occurring in $Q$. Choosing $s\ge k$, we already find
that the minimal exponent is positive, and hence\break
\hbox{$Q(f'\ld f^{(k)})(0)=0$}, so that $\sigma(w)=0$ by (7.14).\qed

Theorem 7.11~(iii) shows that $\cO_{X_k}(1)$ is never relatively ample
over $X$ for $k\ge 2$. In order to overcome this difficulty, we define for
every \hbox{$\abu=(a_1\ld a_k)\in\bZ^k$} a line bundle
$\cO_{X_k}(\abu)$ on $X_k$ such that
$$
\cO_{X_k}(\abu)=\pi_{1,k}^*\cO_{X_1}(a_1)\otimes
\pi_{2,k}^*\cO_{X_2}(a_2)\otimes\cdots\otimes\cO_{X_k}(a_k).
\leqno(7.16)
$$
By (6.9), we have $\pi_{j,k}^*\cO_{X_j}(1)=
\cO_{X_k}(1)\otimes\cO_{X_k}
(-\pi_{j+1,k}^* D_{j+1}-\cdots-D_k)$. Therefore by putting $D^*_j=
\pi_{j+1,k}^* D_{j+1}$ for $1\le j\le k-1$ and $D^*_k=0$, we find
an identity
$$
\leqalignno{
&\cO_{X_k}(\abu)=\cO_{X_k}(b_k)\otimes\cO_{X_k}(-\bbu\cdot D^*),
\qquad\hbox{\rm where}&(7.17)\cr
&\bbu=(b_1\ld b_k)\in\bZ^k,\quad b_j=a_1+\cdots+a_j,&\cr
&\bbu\cdot D^*=\sum_{1\le j\le k-1}b_j\,\pi_{j+1,k}^* D_{j+1}.&\cr}
$$
In particular, if $\bbu\in\bN^k$, i.e., $a_1+\cdots+a_j\ge 0$, we get a
morphism
$$
\cO_{X_k}(\abu)=\cO_{X_k}(b_k)\otimes\cO_{X_k}(-\bbu\cdot D^*)
\to\cO_{X_k}(b_k).\leqno(7.18)
$$
The following result gives a sufficient condition for the relative nefness
or ampleness of weighted jet bundles. Let us recall that a line
bundle $L\to X$ on a projective variety $X$ is
said to be nef if \hbox{$L\cdot C\geq 0$} for~all~irreducible algebraic curves
$C\subset X$, and that a vector bundle $E\to X$ is said to be nef if
$\cO_{\bP(E)}(1)$ is nef on $\bP(E):=P(E^*)\,$; any vector bundle
generated by global sections is nef (cf.~[DePS94] for more details).

\claim 7.19.~Proposition|Take a very ample line bundle $A$ on $X$, and
consider on $X_k$ the line bundle
$$
L_k=\cO_{X_k}(3^{k-1},3^{k-2},\ldots,3,1)\otimes
\pi_{k,0}^*A^{\otimes 3^k}
$$
defined inductively by $L_0=A$ and
$L_k=\cO_{X_k}(1)\otimes\pi_{k,k-1}^*L_{k-1}^{\otimes 3}$.
Then $V_k^*\otimes L_k^{\otimes 2}$ is a nef vector bundle on $X_k$, which is
in fact generated by its global sections, for all $k\geq 0$. Equivalently,
 for all $k\geq 1$,
$$
L'_k=\cO_{X_k}(1)\otimes\pi_{k,k-1}^*L_{k-1}^{\otimes 2}=
\cO_{X_k}(2\cdot 3^{k-2},2\cdot 3^{k-3},\ldots,6,2,1)\otimes
\pi_{k,0}^*A^{\otimes 2\cdot 3^{k-1}}
$$
is nef over~$X_k$ and generated by sections.
\endclaim

The statement concerning $L'_k$ is obtained by projectivizing the vector
bundle $E=V_{k-1}^*\otimes L_{k-1}^{\otimes 2}$ on $X_{k-1}$, whose associated
tautological line bundle is $\cO_{\bP(E)}(1)=L'_k$ on
\hbox{$\bP(E)=P(V_{k-1})=X_k$}. Also one gets inductively that
$$
L_k=\cO_{\bP(V_{k-1}\otimes L_{k-1}^{\otimes 2})}(1)\otimes\pi_{k,k-1}^*L_{k-1}
\quad\hbox{\it is very ample on $X_k$}.\leqno(7.20)
$$
\smallskip

\plainproof. Let $X\subset \bP^N$ be the embedding provided by $A$, so
that $A=\cO_{\bP^N}(1)_{\restriction X}$. As is well known,
if $Q$ is the tautological quotient vector bundle on $\bP^N$,
the twisted cotangent bundle 
$$
T^*_{\bP^N}\otimes\cO_{\bP^N}(2)=\Lambda^{N-1}Q
$$
is nef; hence its quotients $T_X^*\otimes A^{\otimes 2}$ and
$V_0^*\otimes L_0^{\otimes 2}=V^*\otimes A^{\otimes 2}$ are nef
(any tensor power of nef vector bundles is nef, and so is any quotient).
We now proceed by induction, assuming $V_{k-1}^*\otimes L_{k-1}^{\otimes 2}$
to be nef, $k\geq 1$. By taking the second wedge power of the central
term in $(6.4')$, we get an injection
$$
0\lra T_{X_k/X_{k-1}}\lra \Lambda^2\big(\pi_k^\star V_{k-1}\otimes
\cO_{X_k}(1)\big).
$$
By dualizing and twisting with $\cO_{X_{k-1}}(2)\otimes\pi_k^\star
L_{k-1}^{\otimes 2}$, we find a surjection
$$
\pi_k^\star\Lambda^2(V_{k-1}^\star\otimes L_{k-1})\lra T_{X_k/X_{k-1}
}^\star\otimes\cO_{X_k}(2)\otimes\pi_k^\star L_{k-1}^{\otimes 2}\lra 0.
$$
By the induction hypothesis, we see that $T_{X_k/X_{k-1}}^\star\otimes
\cO_{X_k}(2)\otimes\pi_k^\star L_{k-1}^{\otimes 2}$ is nef. Next, the
dual of (6.4) yields an exact sequence
$$
0\lra\cO_{X_k}(1)\lra V_k^\star\lra T_{X_k/X_{k-1}}^\star\lra 0.
$$
As an extension of nef vector bundles is nef, the nefness of
$V_k^*\otimes L_k^{\otimes 2}$ will follow if we check that
$\cO_{X_k}(1)\otimes L_k^{\otimes 2}$ and
$T^\star_{X_k/X_{k-1}}\otimes L_k^{\otimes 2}$ are both nef. However, this follows
again from the induction hypothesis if we observe that the latter implies
$$
L_k\geq \pi_{k,k-1}^* L_{k-1}\quad\hbox{and}\quad
L_k\geq \cO_{X_k}(1)\otimes\pi_{k,k-1}^* L_{k-1},
$$
in the sense that $L''\geq L'$ if the ``difference'' $L''\otimes (L')^{-1}$
is nef. All statements remain valid if we replace ``nef'' with
``generated by sections'' in the above arguments.\qed

\claim 7.21. Corollary|A $\bQ$-line bundle $\cO_{X_k}(\abu)\otimes
\pi_{k,0}^*A^{\otimes p}$ with $\abu\in\bQ^k$, $p\in\bQ$, is nef,
resp.\ ample, on $X_k$ as soon as
$$
\hbox{$a_j\geq 3a_{j+1}$ for $j=1,2,\ldots,k-2$ and $a_{k-1}\geq 2a_k\geq 0$,
$p\geq 2\sum a_j$},
$$
resp.\
$$
\hbox{$a_j\geq 3a_{j+1}$ for $j=1,2,\ldots,k-2$ and $a_{k-1}>2a_k>0$,
$p> 2\sum a_j$}.
$$
\endclaim

\plainproof. This follows easily by taking convex combinations of the
$L_j$ and $L'_j$ and applying Proposition~7.19 and our 
observation (7.20).\qed

\claim 7.22.~Remark|{\rm As in Green-Griffiths [GrGr80], Riemann's
extension theorem shows that for every meromorphic map $\Phi:X\merto Y$
there are well-defined pull-back morphisms
$$
\Phi^*:H^0(Y,E^\GG_{k,m}T^*_Y)\to H^0(X,E^\GG_{k,m}T^*_X),\qquad
\Phi^*:H^0(Y,E_{k,m}T^*_Y)\to H^0(X,E_{k,m}T^*_X).
$$
In particular the dimensions $h^0(X,E^\GG_{k,m}T^*_X)$ and
$h^0(X,E_{k,m}T^*_X)$ are bimeromorphic invariants of~$X$.}
\endclaim

\claim 7.23.~Remark|{\rm As $\bG_k$ is a non-reductive group, it is not
a priori clear that the graded ring $\cA_{n,k,r}=\bigoplus_{m\in\bZ}
E_{k,m}V^\star$ (even pointwise over $X$) is finitely generated.
This can be checked by hand
([Dem07a], [Dem07b]) for $n=2$ and $k\le 4$. Rousseau [Rou06] also checked
the case $n=3$, $k=3$, and then Merker [Mer08, Mer10] proved the finiteness
for $n=2,3,4$, $k\leq 4$ and $n=2$, $k=5$. Recently, B\'erczi and 
Kirwan [BeKi12] made an attempt to prove the finiteness in full generality,
but it appears that the general case is still unsettled.}
\endclaim

\plainsubsection 7.C. Semple tower of a directed variety of general type|

Even if $(X,V)$ is of general type, it is not true that $(X_k,V_k)$ is 
of general type: the fibers of $X_k\to X$ are towers of $\bP^{r-1}$ bundles,
and the canonical bundles of projective spaces are always negative~! However,
a twisted version holds true.

\claim 7.24. Lemma|If $(X,V)$ is of general type, there is a
modification $(\widehat X,\widehat V)$ such that all pairs
$\smash{(\widehat X_k,\widehat V_k)}$ of the associated Semple tower
have a twisted canonical bundle $K_{\widehat V_k}\otimes
\cO_{\widehat X_k}(p)$ that is big when one multiplies $K_{\widehat V_k}$ by 
a suitable $\bQ$-line bundle $\cO_{\widehat X_k}(p)$, $p\in\bQ_+$.
\endclaim

\plainproof. First assume that $V$ has no singularities.
The exact sequences (6.4) and $(6.4')$ provide
$$
K_{V_k}:=\det V_k^*=\det(T^*_{X_k/X_{k-1}})\otimes\cO_{X_k}(1)=
\pi^*_{k,k-1}K_{V_{k-1}}\otimes\cO_{X_k}(-(r-1)),
$$
where $r=\hbox{rank}(V)$. Inductively we get
$$
K_{V_k}=\pi_{k,0}^*K_V\otimes\cO_{X_k}(-(r-1)\onebu),\qquad
\onebu=(1,\ldots,1)\in\bN^k.\leqno(7.25)
$$
We know by [Dem95] that $\cO_{X_k}(\cbu)$ is relatively ample over
$X$ when we take the special weight
$\cbu=(2\,3^{k-2},...,2\,3^{k-j-1},...,6,2,1)\,$; hence
$$
K_{V_k}\otimes\cO_{X_k}((r-1)\onebu+\varepsilon\cbu)
=\pi_{k,0}^* K_V\otimes\cO_{X_k}(\varepsilon\cbu)
$$
is big over $X_k$ for any sufficiently small positive rational number 
$\varepsilon\in\bQ^*_+$. Thanks to Formula~(1.9), we can in fact
replace the weight $(r-1)\onebu+\varepsilon\cbu$ by its total degree
$p=(r-1)k+\varepsilon|\cbu|\in\bQ_+$. The general case of a singular 
linear space follows by considering suitable ``sufficiently high'' 
modifications $\widehat X$ of $X$, the related directed structure
$\widehat V$ on $\widehat X$, and embedding 
$(\widehat X_k,\widehat V_k)$ in the absolute Semple tower
$(\widehat X^a_k,\widehat V^a_k)$ of $\widehat X$.
We still have a well defined morphism of rank $1$ sheaves
$$
\pi_{k,0}^*K_{\widehat V}\otimes\cO_{\widehat X_k}(-(r-1)\onebu)\to
K_{\widehat V_k}
\leqno(7.26)
$$
because the multiplier ideal sheaves involved at each stage behave according
to the monoto\-nicity principle applied to the projections
$\pi^a_{k,k-1}:\widehat X^a_k\to \widehat X^a_{k-1}$ and their
differentials $(\pi^a_{k,k-1})_*$, which yield well-defined transposed
morphisms from the $(k-1)$-st stage to the $k$-th stage at the level of
exterior differential forms. Our contention follows.\qed

\plainsubsection 7.D. Induced directed structure on a subvariety of a jet bundle|

We discuss here the concept of induced directed structure for subvarieties
of the Semple tower of a directed variety $(X,V)$. This will be very 
important to proceed inductively with the base loci of jet differentials.
Let $Z$ be an irreducible algebraic subset of some $k$-jet bundle
$X_k$ over~$X$, $k\ge 0$. We define the linear subspace 
\hbox{$W\subset T_Z\subset T_{X_k|Z}$} to be the closure
$$
W:=\overline{T_{Z'}\cap V_k}\leqno(7.27)
$$
taken on a suitable Zariski open set $Z'\subset Z_{\rm reg}$ where 
the intersection $T_{Z'}\cap V_k$ has a constant rank and 
is a subbundle of $T_{Z'}$. Alternatively, we could also take 
$W$ to be the closure of $T_{Z'}\cap V_k$ in the $k$-th stage
$(X^a_k,V^a_k)$ of the absolute Semple tower, which has the advantage 
of being non-singular. We say that $(Z,W)$ is the 
{\it induced} directed variety structure; this concept of induced
structure already applies of course in the case $k=0$. 
If $f:(\bC,T_{\bC})\to(X,V)$ satisfies $f_{[k]}(\bC)\subset Z$,
then 
$$\hbox{either}~~f_{[k]}(\bC)\subset Z_\alpha\quad\hbox{or}\quad
f'_{[k]}(\bC)\subset W,\leqno(7.28)$$
where $Z_\alpha$ is one of the connected components of $Z\smallsetminus Z'$
and $Z'$ is chosen as in (7.27); especially, if $W=0$, we conclude that
$f_{[k]}(\bC)$ must be contained in one of the $Z_\alpha$'s.
In the sequel, we always consider such a subvariety $Z$ of $X_k$ as a 
directed pair $(Z,W)$ by
taking the induced structure described above. By (7.28), if we proceed
by induction on $\dim Z$, the study of curves tangent to $V$ that have
an order $k$ lifting $f_{[k]}(\bC)\subset Z$ is reduced to the study of
curves tangent to $(Z,W)$. Let us first quote the 
following easy observation.

\claim 7.29. Observation|For $k\ge 1$, let 
$Z\subsetneq X_k$ be an irreducible algebraic subset that projects
onto~$X_{k-1}$, i.e.\ $\pi_{k,k-1}(Z)=X_{k-1}$.
Then the induced directed variety \hbox{$(Z,W)\subset(X_k,V_k)$},
satisfies
$$1\le \rank W<r:=\rank(V_k).$$
\endclaim

\plainproof. Take a Zariski open subset $Z'\subset Z_{\rm reg}$
such that $W'=T_{Z'}\cap V_k$ is a vector bundle over~$Z'$. Since 
$X_k\to X_{k-1}$ is a $\bP^{r-1}$-bundle, $Z$ has codimension at most 
$(r-1)$ in $X_k$. Therefore $\rank W\ge \rank V_k-(r-1)\ge 1$.
On the other hand, if we had $\rank W=\rank V_k$ generically, then
$T_{Z'}$ would contain $V_{k|Z'}$, in particular it would contain all
vertical directions $T_{X_k/X_{k-1}}\subset V_k$ that are tangent
to the fibers of $X_k\to X_{k-1}$. By taking the flow along vertical 
vector fields, we would
conclude that $Z'$ is a union of fibers of $X_k\to X_{k-1}$ up to an
algebraic set of smaller dimension, but this is excluded since $Z$
projects onto $X_{k-1}$ and $Z\subsetneq X_k$.\qed

We introduce the following definition that slightly differs from the one 
given in [Dem14] -- it is actually more flexible and more general.

\claim 7.30. Definition|For $k\ge 1$, let $Z\subset X_k$ be an
irreducible algebraic subset of $X_k$ and $(Z,W)$ the induced
directed structure. We assume moreover that 
$Z\not\subset D_k=P(T_{X_{k-1}/X_{k-2}})$ 
$($and~put $D_1=\emptyset$ in what follows to avoid to have to 
single out the case $k=1)$. In this situation
we say that $(Z,W)$ is of general type modulo the Semple tower
$X_\bu\to X$ if either $W=0$, or $\rank W\ge 1$ and there exists $\ell\geq 0$
and $p\in\bQ_{\geq 0}$ such that 
$$
K_{\widehat W_\ell}^\bullet\otimes\cO_{\widehat Z_\ell}(p)
=K_{\widehat W_\ell}^\bullet\otimes
\cO_{\widehat X_{k+\ell}}(p)_{\restriction \widehat Z_\ell}
\quad\hbox{is big over~$\widehat Z_\ell$},
\leqno(7.31)
$$
possibly after replacing $(Z_\ell,W_\ell)$ by a suitable $($non-singular$)$
modification $(\widehat Z_\ell,\widehat W_\ell)$ obtained via
an embedded resolution of singularities
$$
\mu_\ell:(\widehat Z_\ell\subset \widehat X_{k+\ell})\to
(Z_\ell\subset X_{k+\ell}).
$$
\endclaim

Notice that by (7.26), Condition (7.31) is satisfied if we assume the
existence of $p\geq 0$ such that
$$
\pi_{k+\ell}^*K_{\widehat W}^\bullet\otimes
\cO_{\widehat X_{k+\ell}}(p)_{\restriction \widehat Z_\ell}
\quad\hbox{is big over~$\widehat Z_\ell\subset \widehat X_{k+\ell}$}.
\leqno(7.32)
$$
In fact we infer (7.31) with $\cO_{\widehat Z_\ell}(p)$ replaced by
$\cO_{\widehat Z_\ell}((0,...,0,p)+(r_W-1)\onebu)\subset
\cO_{\widehat Z_\ell}(p+(r_W-1)\ell)$. As a consequence, (7.31) is
satisfied if $K_{\widehat W}^\bu$ is big (i.e.\ $(Z,W)$ is of 
general type), or if $\cO_{\widehat Z_\ell}(1)$ is big on
some $\widehat Z_\ell$, $\ell\geq 1$, but (7.32) is weaker than these
two bigness conditions, since 
we only require that some combination is big. Also, we have the 
following easy observation.

\claim 7.33. Proposition|Let $(X,V)$ be a projective directed variety.
Assume that there exist $\ell\geq 1$ and a weight
$\abu\in\bQ_{>0}^\ell$ such that $\cO_{X_\ell}(\abu)$ is ample over $X_\ell$.
Then every induced directed variety $(Z,W)\subset (X_k,V_k)$
is if general type modulo $X\bu\to X$ for every~$k\geq 1$.
\endclaim

\plainproof. Corollary 7.21 shows that for $\ell'>\ell$ and a suitable weight
$\bbu\in\bQ_{>0}^{\ell'}$, the line bundle $\cO_{X_{\ell'}}(\bbu)$ is relatively
ample with respect to the projection $X_{\ell'}\to X_\ell$. From this, one
deduces that the assumption also holds for arbitrary $\ell'>\ell$ and
a suitable weight $\abu'\in\bQ_{>0}^{\ell'}$. Now,
we use (7.32), in combination with Lemma~2.9~(b); in fact,
$\cO_{\widehat X_{k+\ell}}(1)_{\restriction\widehat Z_\ell}$ is big over 
$\widehat Z_\ell\subset\widehat X_{k+\ell}$ for $\ell\gg 1$,
since we get many sections by pulling back the sections of
$\cO_{\widehat X_{\ell'}}(m\abu')$, $\ell'=k+\ell$, and by restricting them
to $\widehat Z_\ell$.\qed

\plainsubsection 7.E. Relation between invariant and non-invariant jet 
differentials|

We show here that the existence of $\bG_k$-invariant global jet 
differentials is essentially equivalent to the existence of
non-invariant ones. We have seen that the direct image sheaf
$$
\pi_{k,0}\cO_{X_k}(m):=E_{k,m}V^*\subset E_{k,m}^\GG V^*
$$
has a stalk at point $x\in X$ that consists of algebraic differential 
operators $P(f_{[k]})$ acting on germs
of $k$-jets $f:(\bC,0)\to(X,x)$ tangent to $V$, satisfying the 
invariance property
$$
P((f\circ\varphi)_{[k]})=
(\varphi')^m P(f_{[k]})\circ \varphi\leqno(7.34)
$$
whenever $\varphi\in \bG_k$ is in the group of $k$-jets of biholomorphisms
$\varphi:(\bC,0)\to(\bC,0)$. The right action $J_kV\times\bG_k\to J_kV$, 
$(f,\varphi)\mapsto f\circ\varphi$ induces a dual left action of $\bG_k$ on 
$\bigoplus_{m'\le m}E^\GG_{k,m'}V^*$ by 
$$
\bG_k\times\bigoplus_{m'\le m}E^\GG_{k,m'}V^*_x\to 
\bigoplus_{m'\le m}E^\GG_{k,m'}V^*_x,~~
(\varphi,P)\mapsto \varphi^*P,\qquad
(\varphi^*P)(f_{[k]})=P((f\circ\varphi)_{[k]}),
\leqno(7.35)
$$
so that $\psi^*(\varphi^*P)=(\psi\circ\varphi)^*P$. 
Notice that for a global curve $f:(\bC,T_\bC)\to(X,V)$ and a global operator
$P\in H^0(X,E^\GG_{k,m}V^*\otimes F)$ we have to modify a little bit the
definition to consider germs of curves at points $t_0\in \bC$ other than $0$.
This leads to putting
$$
\varphi^*P(f_{[k]})(t_0)=P((f\circ\varphi_{t_0})_{[k]})(0)
\qquad\hbox{where $\varphi_{t_0}(t)=t_0+\varphi(t)$,~~~$t\in D(0,\varepsilon)$}.
$$
The $\bC^*$-action on a homogeneous polynomial of degree $m$ is simply 
$h_\lambda^*P=\lambda^mP$ for a dilation $h_\lambda(t)=\lambda t$, 
$\lambda\in\bC^*$, but since $\varphi\circ h_\lambda\ne h_\lambda\circ\varphi$
in general, $\varphi^*P$ is no longer homogeneous when $P$~is. However,
by expanding the derivatives of $t\mapsto f(\varphi(t))$ at $t=0$, we 
find an expression
$$
(\varphi^*P)(f_{[k]})=\sum_{\alpha\in\bN^k,\,|\alpha|_w=m}\varphi^{(\alpha)}(0)\,
P_\alpha(f_{[k]}),\leqno(7.36)
$$
where $\alpha=(\alpha_1,\ldots,\alpha_k)\in\bN^k$,
$\varphi^{(\alpha)}=(\varphi')^{\alpha_1}(\varphi'')^{\alpha_2}\ldots
(\varphi^{(k)})^{\alpha_k}$,
$|\alpha|_w=\alpha_1+2\alpha_2+\ldots+k\alpha_k$ is the weighted degree
and $P_\alpha$ is a homogeneous polynomial. Since any additional derivative 
taken on $\varphi'$ means one less derivative left for $f$, it is easy to see
that for $P$ homogeneous of degree $m$ we have
$$
m_\alpha:=\deg(P_\alpha)=m-(\alpha_2+2\alpha_3+\ldots+(k-1)\alpha_k)=
\alpha_1+\alpha_2+\ldots+\alpha_k,
$$
in particular $\deg(P_\alpha)<m$ unless $\alpha=(m,0,\ldots,0)$, in which case
$P_\alpha=P$. Let us fix a non-zero global section 
$P\in H^0(X,E^\GG_{k,m}V^*\otimes F)$ for some line bundle $F$ over~$X$,
and pick a non-zero component $P_{\alpha_0}$ of minimum degree $m_{\alpha_0}$ in
the decomposition of~$P$ (of course $m_{\alpha_0}=m$ if and only if $P$ is 
already invariant). We have by construction
$$
P_{\alpha_0}\in H^0(X,E^\GG_{k,m_{\alpha_0}}V^*\otimes F).
$$
We claim that $P_{\alpha_0}$ is $\bG_k$-invariant. Otherwise, there is for
each $\alpha$ a decomposition
$$
(\psi^*P_\alpha)(f_{[k]})=\sum_{\beta\in\bN^k,\,|\beta|_w=m_\alpha}
\psi^{(\beta)}(0)\,P_{\alpha,\beta}(f_{[k]}),\leqno(7.37)
$$
and the non-invariance of $P_{\alpha_0}$ would yield some non-zero 
term $P_{\alpha_0,\beta_0}$ of degree 
$$
\deg(P_{\alpha_0,\beta_0})<\deg(P_{\alpha_0})\le\deg(P)=m.
$$
By replacing $f$ with
$f\circ\psi$ in (7.36) and plugging (7.37) into it, we would get
an identity of the form
$$
(\psi\circ\varphi)^*P(f_{[k]})
=\sum_{\alpha\in\bN^k}(\psi\circ\varphi)^{(\alpha)}(0)\,P_{\alpha}(f_{[k]})=
\sum_{\alpha,\beta\in\bN^k}
\varphi^{(\alpha)}(0)\psi^{(\beta)}(0)\,P_{\alpha,\beta}(f_{[k]}),
$$
but the term in the middle would have all components of 
degree${}\ge m_{\alpha_0}$,
while the term on the right possesses a component 
of degree${}<m_{\alpha_0}$ for a sufficiently generic
choice of $\varphi$ and $\psi$, contradiction. Therefore, we have
shown the existence of a non-zero invariant section
$$
P_{\alpha_0}\in H^0(X,E_{k,m_{\alpha_0}}V^*\otimes F),\quad m_{\alpha_0}\le m.
\eqno\qedsquare
$$

\section{$k$-jet metrics with negative curvature}

The goal of this section is to show that hyperbolicity is closely
related to the existence of $k$-jet metrics with suitable negativity
properties of the curvature. The connection between these properties is
in fact a simple consequence of the Ahlfors-Schwarz lemma. Such ideas have
been already developed long ago by Grauert-Reckziegel [GRe65],
Kobayashi [Kob75] for $1$-jet metrics (i.e., Finsler metrics on~$T_X$)
and by Cowen-Griffiths [CoGr76], Green-Griffiths [GrGr80] and Grauert
[Gra89] for higher order jet metrics. 

\plainsubsection 8.A. Definition of $k$-jet metrics|

Even in the standard case $V=T_X$, the definition given below differs from 
that of [GrGr80], in which the $k$-jet metrics are not supposed to be 
$\bG'_k$-invariant. We prefer to deal here with $\bG'_k$-invariant objects, because
they reflect better the intrinsic geometry. Grauert [Gra89] actually deals
with $\bG'_k$-invariant metrics, but he apparently does not take care of
the way the quotient space $J_k^\reg V/\bG_k$ can be compactified;
also, his metrics are always induced by the Poincar\'e metric, and it is
not at all clear whether these metrics have the expected curvature
properties (see Problem~8.14 below). In the present situation, it is
important to allow also Hermitian metrics possessing some singularities
(``singular Hermitian metrics'' in the sense of [Dem90b]).

\claim 8.1.~Definition|Let $L\to X$ be a holomorphic line bundle over
a complex manifold~$X$. We say that $h$ is a singular metric on~$L$ if
for any trivialization $L_{\restriction U}\simeq U\times\bC$ of~$L$,
the metric is given by $|\xi|_h^2=|\xi|^2e^{-\varphi}$ for some real
valued weight function $\varphi\in L^1_\loc(U)$. The curvature current
of $L$ is then defined to be the closed $(1,1)$-current
$\Theta_{L,h}={i\over 2\pi}\ddbar\varphi$, computed in the sense of
distributions. We say that $h$ admits a closed subset $\Sigma\subset
X$ as its degeneration set if $\varphi$ is locally bounded on
$X\ssm\Sigma$ and is unbounded on a neighborhood of any point
of~$\Sigma$.
\endclaim

An especially useful situation is the case when the curvature of $h$
is positive definite. By this, we mean that there exists a smooth
positive definite Hermitian metric $\omega$ and a continuous positive
function $\varepsilon$ on $X$ such that $\Theta_{L,h}\ge\varepsilon
\omega$ in the sense of currents, and we write in this case
$\Theta_{L,h}\gg0$. We need the following basic fact (quite standard
when $X$ is projective algebraic); however we want to avoid
any algebraicity assumption here, so as to cover potential
applications to non algebraic complex tori.

\claim 8.2.~Proposition|Let $L$ be a holomorphic line bundle on a compact
complex manifold~$X$.
\plainitem{\rm(i)} $L$ admits a singular Hermitian metric $h$ with positive
  definite curvature current $\Theta_{L,h}\gg0$ if and only if $L$ is~big.
Now, define $B_m$ to be the base locus of the linear system 
$|H^0(X,L^{\otimes m})|$ and let 
$$
\Phi_m:X\ssm B_m\to\bP^N
$$ 
be the corresponding meromorphic map. Let $\Sigma_m$ be the closed 
analytic set equal to the union of $B_m$ and of the set of points 
$x\in X\ssm B_m$ such that the fiber $\Phi_m^{-1}(\Phi_m(x))$ is 
positive dimensional.
\plainitem{\rm(ii)} If $\Sigma_m\ne X$ and $G$ is any line bundle,
  the base locus of $L^{\otimes k}\otimes G^{-1}$ is contained in~$\Sigma_m$
  for~$k$ large. As a consequence, $L$ admits a singular Hermitian
  metric $h$ with degeneration set~$\Sigma_m$ and with $\Theta_{L,h}$
  positive definite on~$X$. 
\plainitem{\rm(iii)} Conversely, if $L$ admits a Hermitian metric $h$ with
  degeneration set~$\Sigma$ and positive definite curvature
  current~$\Theta_{L,h}$, there exists an integer~$m>0$ such that the
  base locus $B_m$ is contained in~$\Sigma$ and
  $\Phi_m:X\ssm\Sigma\to\bP_m$ is an embedding.
\vskip0pt 
\endclaim

\plainproof. (i) is proved e.g.\ in [Dem90b,~92], and (ii) and (iii) are
well-known results in the basic theory of linear systems.\qed
\medskip

We now come to the main definitions. By (6.6), every regular $k$-jet
$f\in J_kV$ gives rise to an element $f_{[k-1]}'(0)\in \cO_{X_k}(-1)$.
Thus, measuring the ``norm of $k$-jets'' is the same as taking a
Hermitian metric on $\cO_{X_k}(-1)$. 

\claim 8.3.~Definition|A smooth, $($resp.\ continuous, resp.\ singular$)$
$k$-jet metric on a complex directed manifold $(X,V)$ is a
Hermitian metric $h_k$ on the line bundle $\cO_{X_k}(-1)$ over~$X_k$
$($i.e.\ a Finsler metric on the vector bundle $V_{k-1}$ over $X_{k-1})$,
such that the weight functions $\varphi$ representing the metric are smooth
$($resp.\ conti\-nuous, $L^1_\loc)$. We let $\Sigma_{h_k}\subset X_k$ be
the singularity set of the metric, i.e., the closed subset of points in a 
neighborhood of which the weight $\varphi$ is not locally bounded.
\endclaim

We will always assume here that the weight function $\varphi$ is quasi psh.
Recall that a function $\varphi$ is said to be quasi psh if $\varphi$
is locally the sum of a plurisubharmonic function and of a smooth
function (so that in particular $\varphi\in L^1_\loc$). Then the curvature 
current
$$
\Theta_{h_k^{-1}}(\cO_{X_k}(1))={i\over 2\pi}\ddbar\varphi.
$$
is well defined as a current and is locally bounded from below by a 
negative $(1,1)$-form with constant coefficients.

\claim 8.4.~Definition|Let $h_k$ be a $k$-jet metric on~$(X,V)$. We say that
$h_k$ has negative jet curvature $($resp.\ negative total jet curvature$)$
if $\Theta_{h_k}(\cO_{X_k}(-1))$ is negative definite along the
subbundle $V_k\subset T_{X_k}$ $($resp.\ on all of $T_{X_k})$, i.e.,
if there is $\varepsilon>0$ and a smooth Hermitian metric $\omega_k$ on
$T_{X_k}$ such that
$$
\langle\Theta_{h_k^{-1}}(\cO_{X_k}(1))\rangle(\xi)
\ge\varepsilon|\xi|^2_{\omega_k},\qquad\forall\xi\in V_k\subset T_{X_k}
\quad(\hbox{resp.}\quad\forall\xi\in T_{X_k}).
$$
$($If the metric $h_k$ is not smooth, we suppose that its weights $\varphi$
are quasi psh, and the curvature inequality is taken in the sense of
distributions.$)$
\endclaim

It is important to observe that for $k\ge 2$ there cannot exist any smooth
Hermitian metric $h_k$ on $\cO_{X_k}(1)$ with positive definite curvature
along $T_{X_k/X}$, since $\cO_{X_k}(1)$ is not relatively ample over~$X$.
However, it is relatively big, and Prop.~8.2~(i) shows that $\cO_{X_k}(-1)$
admits a singular Hermitian metric with negative total jet curvature 
(whatever the singularities of the metric are) if and only if $\cO_{X_k}(1)$
is big over $X_k$. It is therefore crucial to allow singularities in the
metrics in Def.~8.4.

\plainsubsection 8.B. Special case of $1$-jet metrics|

A $1$-jet metric $h_1$ 
on $\cO_{X_1}(-1)$ is the same as a Finsler metric $N=\sqrt{h_1}$
on~$V\subset T_X$. Assume until the end of this paragraph that $h_1$ is
smooth. By the well known Kodaira embedding theorem, the
existence of a smooth metric $h_1$ such that
$\Theta_{h_1^{-1}}(\cO_{X_1}(1))$ is positive on all of $T_{X_1}$ is
equivalent to $\cO_{X_1}(1)$ being ample, that is, $V^*$~ample.

\claim 8.5 Remark|{\rm In the absolute case $V=T_X$, there are only few examples of varieties
$X$ such that $T^*_X$ is ample, mainly quotients of the ball
$\bB_n\subset\bC^n$ by a discrete cocompact group of automorphisms.}
\endclaim

The $1$-jet negativity condition considered in Definition~8.4 is much
weaker. For example, if the Hermitian metric $h_1$ comes from a 
(smooth) Hermitian metric $h$ on $V$, then formula (5.15) implies that
$h_1$ has negative total jet curvature (i.e.\ $\Theta_{h_1^{-1}}
(\cO_{X_1}(1))$ is positive) if and only if $\langle\Theta_{V,h}\rangle
(\zeta\otimes v)<0$ for all $\zeta\in T_X\ssm\{0\}$, $v\in V\ssm\{0\}$, 
that is, if $(V,h)$ is {\em negative in the sense of Griffiths}. 
On the other hand, $V_1\subset T_{X_1}$ consists by definition of 
tangent vectors $\tau\in T_{X_1,(x,[v])}$ whose horizontal projection
${}^H\!\tau$ is proportional to~$v$. Thus $\Theta_{h_1}(\cO_{X_1}(-1))$
is negative definite on $V_1$ if and only if $\Theta_{V,h}$ satisfies
the much weaker condition that the {\em holomorphic sectional curvature}
$\langle\Theta_{V,h}\rangle(v\otimes v)$ is negative on every complex
line.\qed

\plainsubsection 8.C. Vanishing theorem for invariant jet differentials|

We now come back to the general situation of jets of arbitrary order~$k$.
Our first observation is the fact that the $k$-jet negativity property
of the curvature becomes actually weaker and weaker as $k$ increases.

\claim 8.6.~Lemma|Let $(X,V)$ be a compact complex directed
manifold. If $(X,V)$ has a $(k-1)$-jet metric $h_{k-1}$ with negative 
jet curvature, then there is a $k$-jet metric $h_k$ with negative jet
curvature such that $\Sigma_{h_k}\subset\pi_k^{-1}(\Sigma_{h_{k-1}})\cup D_k$.
$($The same holds true for negative total jet curvature$)$.
\endclaim

\plainproof. Let $\omega_{k-1}$, $\omega_k$ be given smooth Hermitian metrics on
$T_{X_{k-1}}$ and $T_{X_k}$. The hypothesis implies
$$
\langle\Theta_{h_{k-1}^{-1}}(\cO_{X_{k-1}}(1))\rangle(\xi)\ge
\varepsilon|\xi|_{\omega_{k-1}}^2,\qquad\forall\xi\in V_{k-1}
$$
for some constant~$\varepsilon>0$. On the other hand, as
$\cO_{X_k}(D_k)$ is relatively ample over $X_{k-1}$ ($D_k$ is a hyperplane
section bundle), there exists a smooth metric $\wt h$ on $\cO_{X_k}(D_k)$
such that
$$
\langle\Theta_{\wt h}(\cO_{X_k}(D_k))\rangle(\xi)\ge
\delta|\xi|_{\omega_k}^2-C|(\pi_k)_*\xi|_{\omega_{k-1}}^2,\qquad
\forall\xi\in T_{X_k}
$$
for some constants $\delta,C>0$. Combining both inequalities (the second one
being applied to $\xi\in V_k$ and the first one to $(\pi_k)_*\xi\in
V_{k-1}$), we get
$$
\eqalign{
\langle\Theta_{(\pi_k^* h_{k-1})^{-p}\wt h\,}(\pi_k^*\cO_{X_{k-1}}(p)
\otimes{}&\cO_{X_k}(D_k))\rangle(\xi)\cr
&\ge\delta|\xi|_{\omega_k}^2+(p\varepsilon-C)
|(\pi_k)_*\xi|_{\omega_{k-1}}^2,\qquad\forall\xi\in V_k.\cr}
$$
Hence, for $p$ large enough, $(\pi_k^* h_{k-1})^{-p}\wt h$ has positive
definite curvature along~$V_k$. Now, by (6.9), there is a sheaf injection
$$
\cO_{X_k}(-p)=\pi_k^*\cO_{X_{k-1}}(-p)\otimes\cO_{X_k}(-pD_k)
\hookrightarrow\big(\pi_k^*\cO_{X_{k-1}}(p)\otimes\cO_{X_k}(D_k)
\big)^{-1}
$$
obtained by twisting with $\cO_{X_k}((p-1)D_k)$.  Therefore
$h_k:=((\pi_k^* h_{k-1})^{-p}\wt h)^{-1/p}=(\pi_k^*
h_{k-1})\wt h^{-1/p}$ induces a singular metric on $\cO_{X_k}(1)$ in
which an additional degeneration divisor $p^{-1}(p-1)D_k$ appears.
Hence we get $\Sigma_{h_k}=\pi_k^{-1}\Sigma_{h_{k-1}}\cup D_k$ and
$$
\Theta_{h_k^{-1}}(\cO_{X_k}(1))={1\over p}\Theta_{(\pi_k^* 
h_{k-1})^{-p}\wt h}+{p-1\over p}[D_k]
$$
is positive definite along~$V_k$. The same proof works in the case of
negative total jet curvature.\qed

One of the main motivations for the introduction of $k$-jets metrics is the
following list of algebraic sufficient conditions.

\claim 8.7.~Algebraic sufficient conditions|{\rm We suppose here that $X$
is projective algebraic, and we make one of the additional assumptions (i),
(ii) or (iii) below.

\noindent (i) Assume that there exist integers $k,m>0$ and $\bbu\in\bN^k$
such that the line bundle $L:=\cO_{X_k}(m)\otimes\cO_{X_k}(-\bbu\cdot D^*)$
is ample over~$X_k$. Then there is a smooth Hermitian metric $h_L$ on $L$ 
with positive definite curvature on~$X_k$. By means of the morphism
$\mu:\cO_{X_k}(-m)\to L^{-1}$, we get an induced metric $h_k=(\mu^*
h_L^{-1})^{1/m}$ on $\cO_{X_k}(-1)$ which is degenerate on the support of
the zero divisor $\div(\mu)=\bbu\cdot D^*$. Hence
$\Sigma_{h_k}=\Supp(\bbu\cdot D^*)\subset X_k^\sing$ and
$$
\Theta_{h_k^{-1}}(\cO_{X_k}(1))={1\over m}\Theta_{h_L}(L)+{1\over m}
[\bbu\cdot D^*]\ge{1\over m}\Theta_{h_L}(L)>0.
$$
In particular $h_k$ has negative total jet curvature.

\noindent (ii) Assume more generally that there exist integers $k,m>0$ and
an ample line bundle $A$ on $X$ such that
$H^0(X_k,\cO_{X_k}(m)\otimes\pi_{0,k}^* A^{-1})$ has non-zero
sections $\sigma_1\ld\sigma_N$. Let $Z\subset X_k$ be the base locus of
these sections; necessarily $Z\supset X_k^\sing$ by 7.11~(iii). By taking a
smooth metric $h_A$ with positive curvature on $A$, we get a singular
metric $h_k'$ on $\cO_{X_k}(-1)$ such that
$$
h_k'(\xi)=\Big(\sum_{1\le j\le N}|\sigma_j(w)\cdot\xi^m|_{h_A^{-1}}^2
\Big)^{1/m},\qquad w\in X_k,\quad\xi\in\cO_{X_k}(-1)_w.
$$
Then $\Sigma_{h_k'}=Z$, and by computing ${i\over 2\pi}\ddbar\log h_k'(\xi)$
we obtain
$$
\Theta_{h_k^{\prime\,-1}}(\cO_{X_k}(1))\ge{1\over m}\pi_{0,k}^*\Theta_A.
$$
By (7.17) and an induction on $k$, there exists $\bbu\in\bQ^k_+$ such
that $\cO_{X_k}(1)\otimes\cO_{X_k}(-\bbu\cdot D^*)$ is relatively ample
over~$X$. Hence $L=\cO_{X_k}(1)\otimes\cO_{X_k}(-\bbu\cdot D^*)\otimes
\pi_{0,k}^* A^{\otimes p}$ is ample on $X$ for $p\gg 0$. The
arguments used in (i) show that there is a $k$-jet metric $h_k''$ on
$\cO_{X_k}(-1)$ with $\Sigma_{h_k''}=\Supp(\bbu\cdot D^*)=X_k^\sing$ and
$$
\Theta_{h_k^{\prime\prime\,-1}}(\cO_{X_k}(1))=\Theta_L+[\bbu\cdot D^*]-
p\,\pi_{0,k}^*\Theta_A,
$$
where $\Theta_L$ is positive definite on $X_k$. The metric
$h_k=(h_k^{\prime\,mp}h_k'')^{1/(mp+1)}$ then satisfies 
$\Sigma_{h_k}=\Sigma_{h_k'}=Z$ and
$$
\Theta_{h_k^{-1}}(\cO_{X_k}(1))\ge{1\over mp+1}\Theta_L>0.
$$

\noindent (iii) If $E_{k,m}V^*$ is ample, there is an ample
line bundle $A$ and a sufficiently high symmetric power such that
$S^p(E_{k,m}V^*)\otimes A^{-1}$ is generated by sections. These
sections can be viewed as sections of $\cO_{X_k}(mp)\otimes\pi_{0,k}^*
A^{-1}$ over $X_k$, and their base locus is exactly $Z=X_k^\sing$ by
7.11~(iii). Hence the $k$-jet metric $h_k$ constructed in (ii) has negative
total jet curvature and satisfies $\Sigma_{h_k}=X_k^\sing$.\qed}
\endclaim

An important fact, first observed by [GRe65] for $1$-jet metrics and by
[GrGr80] in the higher order case, is that $k$-jet negativity implies
hyperbolicity. In particular, the existence of enough global jet
differentials implies hyperbolicity.

\claim 8.8.~Theorem|Let $(X,V)$ be a compact complex directed
manifold. If $(X,V)$ has a $k$-jet metric $h_k$ with negative jet
curvature, then every entire curve $f:\bC\to X$ tangent to $V$ is such
that $f_{[k]}(\bC)\subset\Sigma_{h_k}$. In particular, if $\Sigma_{h_k}
\subset X_k^\sing$, then $(X,V)$ is hyperbolic.
\endclaim

\plainproof. The main idea is to use the Ahlfors-Schwarz lemma, following the
approach of [GrGr80]. However we will give here all necessary details
because our setting is slightly different. Assume that there is a $k$-jet
metric $h_k$ as in the hypotheses of Theorem~8.8. Let 
$\omega_k$ be a smooth Hermitian metric on $T_{X_k}$. By 
hypothesis, there exists $\varepsilon>0$ such that
$$
\langle\Theta_{h_k^{-1}}(\cO_{X_k}(1))\rangle(\xi)\ge\varepsilon
|\xi|_{\omega_k}^2\qquad\forall\xi\in V_k.
$$
Moreover, by (6.4), $(\pi_k)_*$ maps $V_k$ continuously to
$\cO_{X_k}(-1)$ and the weight $e^\varphi$ of $h_k$ is locally bounded from
above. Hence there is a constant $C>0$ such that
$$
|(\pi_k)_*\xi|_{h_k}^2\le C|\xi|_{\omega_k}^2,\qquad\forall\xi\in V_k.
$$
Combining these inequalities, we find
$$
\langle\Theta_{h_k^{-1}}(\cO_{X_k}(1))\rangle(\xi)\ge
{\varepsilon\over C}|(\pi_k)_*\xi|_{h_k}^2,\qquad\forall\xi\in V_k.
$$
Now, let $f:\Delta_R\to X$ be a non-constant holomorphic map tangent to~$V$
on the disk~$\Delta_R$. We use the line bundle morphism (6.6)
$$
F=f_{[k-1]}':T_{\Delta_R}\to f_{[k]}^*\cO_{X_k}(-1)
$$
to obtain a pull-back metric
$$
\gamma=\gamma_0(t)\,dt\otimes d\ol t=F^* h_k\qquad
\hbox{on $T_{\Delta_R}$}.
$$
If $f_{[k]}(\Delta_R)\subset\Sigma_{h_k}$ then $\gamma\equiv 0$. Otherwise,
$F(t)$ has isolated zeroes at all singular points of $f_{[k-1]}$ and so
$\gamma(t)$ vanishes only at these points and at points of the degeneration
set $(f_{[k]})^{-1}(\Sigma_{h_k})$ which is a polar set in~$\Delta_R$.
At other points, the Gaussian curvature of $\gamma$ satisfies
$$
{i\,\ddbar\log\gamma_0(t)\over\gamma(t)}
={-2\pi\,(f_{[k]})^*\Theta_{h_k}(\cO_{X_k}(-1))\over F^* h_k}
={\langle\Theta_{h_k^{-1}}(\cO_{X_k}(1))\rangle(f_{[k]}'(t))\over
|f_{[k-1]}'(t)|_{h_k}^2}\ge {\varepsilon\over C},
$$
since $f_{[k-1]}'(t)=(\pi_k)_* f_{[k]}'(t)$. The Ahlfors-Schwarz lemma
4.2 implies that $\gamma$ can be compared with the Poincar\'e
metric as follows:
$$
\gamma(t)\le{2C\over\varepsilon}{R^{-2}|dt|^2\over(1-|t|^2/R^2)^2}\quad
\Longrightarrow\quad|f_{[k-1]}'(t)|_{h_k}^2\le
{2C\over\varepsilon}{R^{-2}\over(1-|t|^2/R^2)^2}.
$$
If $f:\bC\to X$ is an entire curve tangent to~$V$ such that
$f_{[k]}(\bC)\not\subset\Sigma_{h_k}$, the above estimate implies
as $R\to+\infty$ that $f_{[k-1]}$ must be a constant, and hence also~$f$. Now,
if $\Sigma_{h_k}\subset X_k^\sing$, the inclusion $f_{[k]}(\bC)\subset
\Sigma_{h_k}$ implies $f'(t)=0$ at every point. Therefore $f$ is a constant
and $(X,V)$ is hyperbolic.\qed

Combining Theorem~8.8 with 8.7~(ii) and (iii), we get the following
consequences.

\claim 8.9.~Vanishing theorem|Assume that there exist integers $k,m>0$ and an
ample line bundle $L$ on $X$ such that $H^0(X_k,\cO_{X_k}(m)\otimes
\pi_{0,k}^* L^{-1})\simeq H^0(X,E_{k,m}V^*\otimes L^{-1})$
has non-zero sections $\sigma_1\ld\sigma_N$. Let $Z\subset X_k$ be the
base locus of these sections. Then every entire curve $f:\bC\to X$
tangent to $V$ is such that $f_{[k]}(\bC)\subset Z$. In other words,
for every global $\bG_k$-invariant polynomial differential operator
$P$ with values in $L^{-1}$, every entire curve $f$ must satisfy
the algebraic differential equation $P(f_{[k]})=0$.
\endclaim

\claim 8.10.~Corollary|Let $(X,V)$ be a compact complex directed
manifold. If $E_{k,m}V^*$ is ample for some positive integers
$k,m$, then $(X,V)$ is hyperbolic.
\endclaim

\claim 8.11.~Remark|{\rm Green and Griffiths [GrGr80] stated that 
Theorem~8.9 is even true for sections $\sigma_j\in H^0(X,E^\GG_{k,m}(V^*)
\otimes L^{-1})$, in the special case $V=T_X$ they consider. This is
proved below in \S8.D; the reader is also referred to Siu and Yeung 
[SiYe97] for a proof based on a use of Nevanlinna theory and 
the logarithmic derivative lemma (the original proof given in [GrGr80] 
does not seem to be complete, as it relies on an unsettled pointwise version 
of  the Ahlfors-Schwarz lemma for general jet differentials); other proofs
seem to have been circulating in the literature in the last years.
Let us first give a very short proof in the case where $f$ is supposed 
to have a bounded derivative (thanks to the Brody criterion, 
this is enough if one is merely interested in proving
hyperbolicity; thus Corollary 8.10 will be valid with $E^\GG_{k,m}V^*$
in place of $E_{k,m}V^*$). In fact, if $f'$ is bounded, one can apply
the Cauchy inequalities to all components $f_j$ of $f$ with respect to
a finite collection of coordinate patches covering~$X$. As $f'$ is bounded,
we can do this on sufficiently small discs $D(t,\delta)\subset\bC$ 
of constant radius $\delta>0$. Therefore
all derivatives $f'$, $f''$, $\ldots\,f^{(k)}$ are bounded. From this
we conclude that $\sigma_j(f)$ is a bounded section of $f^*
L^{-1}$. Its norm $|\sigma_j(f)|_{L^{-1}}$ (with respect to any
positively curved metric $|~~|_L$ on $L$) is a bounded subharmonic
function, which is moreover strictly subharmonic at all points where
$f'\ne 0$ and $\sigma_j(f)\ne 0$. This is a contradiction unless $f$ is
constant or $\sigma_j(f)\equiv 0$.}
\endclaim

The above results justify the following definition and problems.

\claim 8.12.~Definition|We say that $X$, resp.\ $(X,V)$, has non
degenerate negative $k$-jet curvature if there exists a $k$-jet metric
$h_k$ on $\cO_{X_k}(-1)$ with negative jet curvature such that
$\Sigma_{h_k}\subset X_k^\sing$.
\endclaim

\claim 8.13.~Conjecture|Let $(X,V)$ be a compact directed manifold.
Then $(X,V)$ is hyperbolic if and only if $(X,V)$ has nondegenerate 
negative $k$-jet curvature for $k$ large enough.
\endclaim

This is probably a hard problem. In fact, it is shown in [Dem97, 
Section~8] that the smallest admissible integer $k$ must depend on the
geometry of $X$ and need not be uniformly bounded as soon as $\dim X\ge 2$
(even in the absolute case~$V=T_X$). On the other hand, if $(X,V)$ is
hyperbolic, we get for each integer $k\ge 1$ a generalized
Kobayashi-Royden metric $\bfk_{(X_{k-1},V_{k-1})}$ on $V_{k-1}$ (see
Definition~2.1), which can be also viewed as a $k$-jet metric $h_k$ on
$\cO_{X_k}(-1)\,$; we will call it the {\em Grauert $k$-jet metric}
of $(X,V)$, although it formally differs from the jet metric considered 
in [Gra89] (see also [DGr91]). By looking at the projection
$\pi_k:(X_k,V_k)\to(X_{k-1},V_{k-1})$, we see that the sequence $h_k$
is monotonic, namely $\pi_k^*h_k\le h_{k+1}$ for every~$k$. If
$(X,V)$ is hyperbolic, then $h_1$ is nondegenerate and therefore by
monotonicity $\Sigma_{h_k}\subset X_k^\sing$ for $k\ge 1$. Conversely,
if the Grauert metric satisfies $\Sigma_{h_k} \subset X_k^\sing$, it 
is easy to see that $(X,V)$ is hyperbolic. The following problem is
thus especially meaningful.

\claim 8.14.~Problem|Estimate the $k$-jet curvature
$\Theta_{h_k^{-1}}(\cO_{X_k}(1))$ of the Grauert metric $h_k$
on $(X_k,V_k)$ as $k$ tends to $+\infty$.
\endclaim

\plainsubsection 8.D. Vanishing theorem for non-invariant jet differentials|

As an application of the arguments developed in \S7.E, we 
indicate here a proof of the basic vanishing theorem for non-invariant
jet differentials. This version has been first proved in full generality 
by Siu [Siu97]  (cf.\ also [Dem97]), with a different and more involved
technique based on Nevanlinna theory and the logarithmic derivative lemma.

\claim 8.15.~Theorem|Let $(X,V)$ be a projective directed and $A$ an
ample divisor on~$X$.  Then one has $P(f\,;\,f',f'',\ldots,f^{(k)})=0$
for every entire curve $f:(\bC,T_\bC)\to (X,V)$ and every global
section $P\in H^0(X,E^\GG_{k,m}V^*\otimes \cO(-A))$.
\endclaim

\noindent{\em Sketch of proof}. In general, we know by Theorem~8.9 that the result 
is true when $P$ is invariant,
i.e.\ for $P\in H^0(X,E_{k,m}V^*\otimes \cO(-A))$. Now, we prove Theorem
8.15 by induction on $k$ and $m$ (simultaneously for all directed
varieties). Let $Z\subset X_k$ be the base locus of all polynomials
$Q\in H^0(X,E_{k,m'}^\GG V^*\otimes \cO(-A))$ with $m'<m$. A priori, this defines
merely an algebraic set in the Green-Griffiths bundle $X_k^\GG=
(J_kV\ssm\{0\})/\bC^*$, but since the global polynomials $\varphi^*Q$
also enter the game, we know that the base locus is $\bG_k$-invariant,
and thus descends to~$X_k$. Let $f:(\bC,T_\bC)\to (X,V)$. By the 
induction hypothesis, we know that 
$f_{[k]}(\bC)\subset Z$. Therefore
$f$ can also be viewed as a entire curve drawn in the directed variety
$(Z,W)$ induced by $(X_k,V_k)$. By (7.36), we have a decomposition
$$
(\varphi^*P)(g_{[k]})=\sum_{\alpha\in\bN^k,\,|\alpha|_w=m}\varphi^{(\alpha)}(0)\,
P_\alpha(g_{[k]}),\quad
\hbox{with $\deg(P_\alpha)<\deg(P)$ for $\alpha\ne(m,0,\ldots,0)$},
$$
and since $P_\alpha(g_{[k]})=0$ for all germs of curves $g$ of $(Z,W)$
when $\alpha\ne(m,0,\ldots,0)$, we conclude that $P$ defines an
invariant jet differential when it is restricted to $(Z,W)$, in other words it
still defines a section of
$$
H^0\big(Z,(\cO_{X_k}(m)\otimes\pi^*_{k,0}\cO_X(-A))_{\restriction Z}\big).
$$
We can then apply the Ahlfors-Schwarz lemma in the way we did it
in \S8.C to conclude that $P(f_{[k]})=0$. \qed

\section{Morse inequalities and
the Green-Griffiths-Lang conjecture}

The goal of this section is to study the existence and properties of
entire curves $f:\bC\to X$ drawn in a complex irreducible
$n$-dimensional variety~$X$, and more specifically to show that they 
must satisfy certain global algebraic or differential equations 
as soon as $X$ is projective of general type. By means of
holomorphic Morse inequalities and a probabilistic analysis of
the cohomology of jet spaces, it is possible to prove a significant step 
of the generalized Green-Griffiths-Lang conjecture.
The use of holomorphic Morse inequalities 
was first suggested in [Dem07a], and then carried out in an algebraic context
by S.~Diverio in his PhD work ([Div08, Div09]). The general more analytic
and more powerful results presented here first appeared in [Dem11, Dem12].

\plainsubsection 9.A. Introduction|

Let $(X,V)$ be a directed variety. By defi\-nition, proving the algebraic degeneracy of an entire curve $f;(\bC,T_\bC)\to(X,V)$ means finding
a non-zero polynomial $P$ on $X$ such that $P(f)=0$. As already 
explained in \S$\,$8, all known methods of proof are based on establishing 
first the existence of certain algebraic differential equations 
$P(f\,;\,f',f'',\ldots,f^{(k)})=0$ of some order $k$, and then trying to find
enough such equations so that they cut out a proper algebraic
locus $Y\subsetneq X$. We use for this global sections of
$H^0(X,E^\GG_{k,m}V^*\otimes \cO(-A))$ where $A$ is ample, and
apply the fundamental vanishing theorem~8.15.
It is expected that the global sections of
$H^0(X,E^\GG_{k,m}V^*\otimes \cO(-A))$ are precisely those which
ultimately define the algebraic locus $Y\subsetneq X$ where the curve
$f$ should lie.  The problem is then reduced to (i) showing that there are 
many non-zero sections of $H^0(X,E^\GG_{k,m}V^*\otimes \cO(-A))$ and
(ii) understanding what is their joint base locus. 
The first part of this program is the main result of this section.

\claim 9.1. Theorem| Let $(X,V)$ be a directed projective variety
such that $K_V$ is big and let $A$ be an ample divisor. Then 
for $k\gg 1$ and $\delta\in\bQ_+$ small enough, $\delta\le c(\log k)/k$,
the number of sections
$h^0(X,E^\GG_{k,m}V^*\otimes\cO(-m\delta A))$ has maximal growth, 
i.e.\ is larger that $c_km^{n+kr-1}$ for some~$m\ge m_k$, 
where $c,\,c_k>0$, $n=\dim X$ and $r=\rank V$. In particular, entire curves
$f:(\bC,T_\bC)\to(X,V)$ satisfy $($many$)$ algebraic differential 
equations.
\endclaim

The statement is very elementary to check when $r=\rank V=1$, and 
therefore when $n=\dim X=1$. In higher dimensions  $n\ge 2$, only 
very partial results were known before Theorem 9.1 was obtained in [Dem11],
[and they dealt merely with the absolute case $V=T_X$]. In dimension~$2$,
Theorem~9.1 is a consequence of the Riemann-Roch calculation of 
Green-Griffiths [GrGr80], 
combined with a vanishing theorem due to Bogomolov [Bog79] -- the latter 
actually only applies to the top cohomology group $H^n$, and things 
become much more delicate when extimates of intermediate cohomology
groups are needed. In higher dimensions, Diverio [Div08, Div09] proved the 
existence of sections of $H^0(X,E^\GG_{k,m}V^*\otimes\cO(-1))$
whenever $X$ is a hypersurface of $\bP^{n+1}_\bC$ of high degree $d\ge
d_n$, assuming $k\ge n$ and $m\ge m_n$. More recently, Merker [Mer15]
was able to treat the case of arbitrary hypersurfaces of general type,
i.e.\ $d\ge n+3$, assuming this time $k$ to be very large.  The latter
result is obtained through explicit algebraic calculations of the
spaces of sections, and the proof is computationally very
intensive. B\'erczi [Ber15, Ber18] also obtained related results with a
different approach based on residue formulas, assuming e.g.\ $d\ge n^{9n}$.

All these approaches are algebraic in nature. Here,
however, our techniques are based on more elaborate curvature estimates in 
the spirit of Cowen-Griffiths [CoGr76]. They require holomorphic Morse 
inequalities (see 9.10 below) -- and we
do not know how to translate our method in an algebraic setting.
Notice that holomorphic Morse inequalities are essentially
insensitive to singularities, as we can pass to non-singular models
and blow-up $X$ as much as we want: if $\mu:\wt X\to X$ is a
modification then $\mu_*\cO_{\wt X}=\cO_X$ and for $q\geq 1$,
$R^q\mu_*\cO_{\wt X}$
is supported on a codimension $1$ analytic subset (even a codimension $2$
subset if $X$ is smooth). It~follows from the
Leray spectral sequence that the cohomology estimates for $L$ on $X$ or for
$\wt L=\mu^*L$ on $\wt X$ differ by negligible terms, i.e.
$$
h^q(\wt X,\wt L^{\otimes m})-h^q(X,L^{\otimes m})=O(m^{n-1}).\leqno(9.2)
$$
Finally, singular holomorphic Morse inequalities (in the form obtained 
by L.~Bonavero 
[Bon93]) allow us to work with singular Hermitian metrics $h$; this is the
reason why we will only require to have big line bundles rather than
ample line bundles. In the case of linear subspaces $V\subset T_X$, we
introduce singular Hermitian metrics as follows.

\claim 9.3. Definition| A singular Hermitian metric on a linear
subspace $V\subset T_X$ is a metric $h$ on the fibers of $V$ such that
the function $\log h:\xi\mapsto\log|\xi|_h^2$ is locally integrable 
on the total space of $V$.
\endclaim

Such a metric can also be viewed as a singular Hermitian metric on
the tauto\-logical line bundle $\cO_{P(V)}(-1)$ on the projectivized bundle
$P(V)=V\ssm\{0\}/\bC^*$, and therefore its dual metric $h^*$ defines a 
curvature current $\Theta_{\cO_{P(V)}(1),h^*}$ of type $(1,1)$ on
$P(V)\subset P(T_X)$, such that
$$
p^*\Theta_{\cO_{P(V)}(1),h^*} ={\ii\over 2\pi}\ddbar\log h,\qquad
\hbox{where $p:V\ssm\{0\}\to P(V)$}.\leqno(9.4)
$$
If $\log h$ is quasi-plurisubharmonic $($or quasi-psh, which means psh 
modulo addition of a smooth function$)$ on $V$, then $\log h$ is indeed
locally integrable, and we have moreover
$$\Theta_{\cO_{P(V)}(1),h^*}\ge -C\omega\leqno(9.5)$$
for some smooth positive $(1,1)$-form on $P(V)$ and some constant $C>0\;$;
conversely, if (9.5) holds, then $\log h$ is quasi-psh.

\claim 9.6. Definition| We will say that a singular Hermitian metric $h$ on
$V$ is {\rm admissible} if $h$ can be written as $h=e^\varphi h_{0|V}$ where
$h_0$ is a smooth positive definite Hermitian on $T_X$ and $\varphi$
is a quasi-psh weight with analytic singularities on $X$, as in 
Definition~$9.3$. 
Then $h$ can be seen as a singular Hermitian metric on $\cO_{P(V)}(1)$,
with the property that it induces a smooth positive definite metric
on a Zariski open set $X'\subset X\ssm \Sing(V)\,;$ we will denote by
$\Sing(h)$ the complement of the largest such Zariski
open set $X'$ $($so that $\Sing(h)\supset\Sing(V))$.
\endclaim

If $h$ is an admissible metric on~$V$, we define 
$\cO_h(V^*)$ to be the sheaf of germs of holomorphic sections 
of $V^*_{\restriction X\ssm\Sing(h)}$ which are $h^*$-bounded near $\Sing(h)$; 
by the assumption on the analytic singularities, 
this is a coherent sheaf (as the direct image of some coherent sheaf on
$P(V)$), and actually, since $h^*=e^{-\varphi}h_0^*$, it is 
a subsheaf of the
sheaf $\cO(V^*):=\cO_{h_0}(V^*)$ associated with a smooth positive
definite metric $h_0$  on $T_X$. If $r$ is the generic rank of $V$ and
$m$ a positive integer, we define similarly 
$$
\leqalignno{
\bddK_{V,h}^{[m]}={}&\hbox{sheaf of germs of holomorphic sections of 
$(\det V^*_{\restriction X'})^{\otimes m}=(\Lambda^rV^*_{\restriction X'})^{\otimes m}$}&(9.7)\cr
\noalign{\vskip-2pt}
&\hbox{which are $\det h^*$-bounded},\cr}
$$
so that $\bddK_V^{[m]}:=\bddK_{V,h_0}^{[m]}$ according to Def.~2.7.
For a given admissible Hermitian structure $(V,h)$, we define similarly
the sheaf $E^\GG_{k, m}V^*_h$ to be the sheaf of 
polynomials defined over $X\ssm\Sing(h)$ which are \hbox{``$h$-bounded''}.
This means that when they are viewed as polynomials $P(z\,;\,\xi_1,\ldots,
\xi_k)$ in terms of $\xi_j=(\nabla_{h_0}^{1,0})^jf(0)$ where
$\nabla_{h_0}^{1,0}$ is the $(1,0)$-component of the induced Chern 
connection on $(V,h_0)$, there is a uniform bound 
$$
\big|P(z\,;\,\xi_1,\ldots,\xi_k)\big|\le 
C\Big(\sum\Vert \xi_j\Vert_h^{1/j}\Big)^m\leqno(9.8)
$$
near points of $X\ssm X'$ (see section~2 for more details on
this). Again, by a direct image argument, one sees that
$E^\GG_{k, m}V^*_h$ is always a coherent sheaf. The sheaf
$E^\GG_{k, m}V^*$ is defined to be $E^\GG_{k, m}V^*_h$ 
when $h=h_0$ (it is actually independent of the choice of $h_0$, as
follows from arguments similar to those given in section~2). Notice
that this is exactly what is needed to extend the proof of the vanishing
theorem~8.15 to the case of a singular linear space $V\,$; the value 
distribution theory argument can only work when the functions 
$P(f\,;\;f',\ldots,f^{(k)})(t)$ do not exhibit poles, and this is guaranteed
here by the boundedness assumption.

Our strategy can be described as follows. We consider the
Green-Griffiths bundle of $k$-jets $X_k^\GG=J^kV\ssm\{0\}/\bC^*$,
which by (7.7) consists of a fibration in {\em weighted projective
spaces}, and its associated tautological sheaf 
$$L=\cO_{X_k^\GG}(1),$$
viewed rather as a virtual $\bQ$-line bundle $\cO_{X_k^\GG}(m_0)^{1/m_0}$
with $m_0=\lcm(1,2,\,...\,,k)$.
Then, if $\pi_k:X_k^\GG\to X$ is the natural projection, we have
$$
E^\GG_{k,m}=(\pi_k)_*\cO_{X_k^\GG}(m)\quad\hbox{and}\quad
R^q(\pi_k)_*\cO_{X_k^\GG}(m)=0~\hbox{for $q\ge 1$}.
$$
Hence, by the Leray spectral sequence we get for every invertible sheaf
$F$ on $X$ the isomorphism
$$
H^q(X,E^\GG_{k,m}V^*\otimes F)\simeq 
H^q(X_k^\GG,\cO_{X_k^\GG}(m)\otimes\pi_k^*F).\leqno(9.9)
$$
The latter group can be evaluated thanks to holomorphic Morse inequalities.
Let us recall the main statement.

\claim 9.10. Holomorphic Morse inequalities {\rm ([Dem85])}|Let $X$
be a compact complex  manifolds, $E\to X$ a holomorphic vector bundle of
rank $r$, and $(L,h)$ a Hermitian line bundle. The dimensions
$h^q(X,E\otimes L^m)$ of cohomology groups of the tensor powers 
$E\otimes L^m$ satisfy the following asymptotic estimates as 
$m\to +\infty\,:$

\noindent {\rm(WM)} Weak Morse inequalities$\,:$
$$h^q(X,E\otimes L^m)\le r {m^n\over n!}\int_{X(L,h,q)} (-1)^q \Theta_{L,h}^n + o(m^n),$$
where $X(L,h,q)$ denotes the open set of points $x\in X$ at which
the curvature form $\Theta_{L,h}(x)$ has signature $(q,n-q)\,;$\vskip2pt
\noindent {\rm(SM)} Strong Morse inequalities$\,:$
$$\sum_{0\le j\le q} (-1)^{q-j}h^j(X,E\otimes L^m) \le r {m^n\over n!}
\int_{X(L,h,\leq q)}(-1)^q\Theta_{L,h}^n+o(m^n),$$
where $X(L,h,\leq q)=\bigcup_{j\leq q}X(L,h,j)\,;$\vskip2pt
\noindent {\rm(RR)} Asymptotic Riemann-Roch formula$\,:$
$$\chi(X,E\otimes L^m) := \sum_{0\le j\le n} (-1)^jh^j(X,E\otimes L^m)
= r{m^n\over n!}\int_X \Theta_{L,h}^n + o(m^n)~.$$
\vskip-10pt
\endclaim

\noindent Moreover $($cf.\ Bonavero's PhD thesis {\rm [Bon93]}$)$, 
if $h=e^{-\varphi}$ is a singular Hermitian metric with analytic 
singularities of pole set $P=\varphi^{-1}(-\infty)$, the estimates 
still hold provided all cohomology groups
are replaced by cohomology groups 
$H^q(X,E\otimes L^m\otimes\cI(h^m))$ twisted with the 
corresponding $L^2$ multiplier ideal sheaves 
$$
\cI(h^m)=\cI(k\varphi)=\big\{f\in\cO_{X,x},\;\;\exists V\ni x,~
\int_V|f(z)|^2e^{-m\varphi(z)}d\lambda(z)<+\infty\big\},
$$
and provided the Morse integrals are computed on the regular locus of $h$, 
namely restricted to $X(L,h,q)\ssm\Sigma\,$:
$$
\int_{X(L,h,q)\ssm\Sigma}(-1)^q\Theta_{L,h}^n.
$$
The special case of 9.10~(SM) when $q=1$ yields a very useful
criterion for the existence of sections of large multiples of $L$.

\claim 9.11. Corollary|Let $L\to X$ be a holomorphic line bundle 
equipped with a singular Hermitian metric $h=e^{-\varphi}$ 
with analytic singularities of pole set $\Sigma=\varphi^{-1}(-\infty)$.
Then we have the following lower bounds
\vskip3pt
\plainitem{\rm(a)} at the $h^0$ level~$:$
$$
\eqalign{
h^0(X,E\otimes L^m)&\geq h^0(X,E\otimes L^m\otimes\cI(h^m))\cr
&\geq h^0(X,E\otimes L^m\otimes\cI(h^m))-h^1(X,E\otimes L^m\otimes\cI(h^m))\cr
&\geq r{k^n\over n!}\int_{X(L,h,\le 1)\ssm\Sigma}\Theta_{L,h}^n -o(k^n)~.\cr}
$$
Especially $L$ is big as soon as $\int_{X(L,h,\le 1)\ssm\Sigma}\Theta_{L,h}^n>0$
for some singular Hermitian metric $h$ on~$L$.
\vskip3pt
\plainitem{\rm(b)} at the $h^q$ level~$:$
$$
h^q(X,E\otimes L^m\otimes\cI(h^m))\geq r{k^n\over n!}
\sum_{j=q-1,q,q+1}(-1)^q\int_{X(L,h,j)\ssm\Sigma}\Theta_{L,h}^n -o(k^n)~.
$$
\endclaim

Now, given a directed manifold $(X,V)$, we can associate with any 
admissible metric $h$ on $V$ a metric
(or rather a natural family) of metrics on $L=\cO_{X_k^\GG}(1)$. 
The space $X_k^\GG$ always possesses quotient singularities if $k\ge 2$ 
(and even some more if $V$ is singular), but we do not really care since 
Morse inequalities still work in this setting thanks to Bonavero's
generalization. As we will see, it is then
possible to get nice asymptotic formulas as $m\to +\infty$. They appear to
be of a {\em probabilistic nature} if we take the components of the $k$-jet
(i.e.\ the successive derivatives $\xi_j=f^{(j)}(0)$, $1\le j\le k$) as 
random variables. This probabilistic behaviour was somehow  already
visible in the Riemann-Roch calculation of~[GrGr80].
In this way, assuming $K_V$ big, we produce a lot of sections
$\sigma_j=H^0(X_k^\GG,\cO_{X_k^\GG}(m)\otimes\pi_k^*F)$, corresponding
to certain divisors $Z_j\subset X_k^\GG$. The hard problem which
is left in order to complete a proof of the generalized
Green-Griffiths-Lang conjecture is to compute the base locus
$Z=\bigcap Z_j$ and to show that $Y=\pi_k(Z)\subset X$ must be a
proper algebraic variety. 

\plainsubsection 9.B. Hermitian geometry of weighted projective spaces|

The goal of this section is to introduce natural K\"ahler metrics on
weighted projective spaces, and to evaluate the corresponding volume forms.
Here we put $d^c={\ii\over 4\pi}(\dbar-\partial)$ so that
$dd^c={\ii\over 2\pi}\ddbar$. The normalization of the $d^c$ operator is
chosen such that we have precisely $(dd^c\log|z|^2)^n=\delta_0$ (the Dirac
mass at~$0$) for the Monge-Amp\`ere operator in $\bC^n$. Given a $k$-tuple
of ``weights'' $a=(a_1,\ldots,a_k)$, i.e.\ of integers $a_s>0$ with
\hbox{${\rm gcd}(a_1,\ldots,a_k)=1$}, we introduce the weighted projective
space $P(a_1,\ldots,a_k)$ to be the quotient of \hbox{$\bC^k\ssm\{0\}$} by
the corresponding weighted $\bC^*$ action:
$$
P(a_1,\ldots,a_k)=\bC^k\ssm\{0\}/\bC^*,\qquad
\lambda\cdot z =(\lambda^{a_1}z_1,\ldots,\lambda^{a_k}z_k),\quad
\lambda\in\bC^*.
\leqno(9.12)
$$
As is well known, this defines a toric $(k-1)$-dimensional algebraic variety 
with quotient singularities. On this variety, we introduce the possibly
singular (but almost everywhere smooth and non-degenerate) K\"ahler form 
$\omega_{a,p}$ defined by
$$
\pi_a^*\omega_{a,p}=dd^c\varphi_{a,p},\qquad
\varphi_{a,p}(z)={1\over p}\log\sum_{1\le s\le k}|z_s|^{2p/a_s},
\leqno(9.13)
$$
where $\pi_a:\bC^k\ssm\{0\}\to P(a_1,\ldots,a_k)$ is the canonical projection
and $p>0$ is a positive constant. It is clear that $\varphi_{p,a}$ is real
analytic on $\bC^k\ssm\{0\}$ if $p$ is an integer and a common multiple of 
all weights $a_s$, and we will implicitly pick such a $p$ later on 
to avoid any difficulty. Elementary calculations  give the following
well-known formula for the volume
$$
\int_{P(a_1,\ldots,a_k)}\omega_{a,p}^{k-1}={1\over a_1\ldots a_k}\leqno(9.14)
$$
(notice that this is independent of $p$, as it is obvious 
by Stokes theorem, since the cohomology class of $\omega_{a,p}$ does 
not depend on $p$).

Our later calculations will require a slightly more general setting.
Instead of looking at $\bC^k$, we consider the weighted $\bC^*$ action
defined by
$$
\bC^{|r|}=\bC^{r_1}\times\ldots\times\bC^{r_k},\qquad
\lambda\cdot z =(\lambda^{a_1}z_1,\ldots,\lambda^{a_k}z_k),\quad
\lambda\in\bC^*.
\leqno(9.15)
$$
Here $z_s\in\bC^{r_s}$ for some $k$-tuple $r=(r_1,\ldots,r_k)$ and
$|r|=r_1+\ldots+r_k$. This gives rise to a weighted projective space 
$$
\leqalignno{
&P(a_1^{[r_1]},\ldots,a_k^{[r_k]})=P(a_1,\ldots,a_1,\ldots,a_k,\ldots,a_k),\cr
&\pi_{a,r}:\bC^{r_1}\times\ldots\times\bC^{r_k}
\ssm\{0\}\longrightarrow P(a_1^{[r_1]},\ldots,a_k^{[r_k]}),&(9.16)\cr}
$$
obtained by repeating $r_s$ times each weight $a_s$. On this space, we
introduce the degenerate K\"ahler metric $\omega_{a,r,p}$ such that
$$
\pi_{a,r}^*\omega_{a,r,p}=dd^c\varphi_{a,r,p},\qquad
\varphi_{a,r,p}(z)={1\over p}\log\sum_{1\le s\le k}|z_s|^{2p/a_s},
\leqno(9.17)
$$
where $|z_s|$ stands now for the standard Hermitian norm 
$(\sum_{1\le j\le r_s}|z_{s,j}|^2)^{1/2}$ on~$\bC^{r_s}$. This metric
is cohomologous to the corresponding ``polydisc-like'' 
metric $\omega_{a,p}$ already defined, and therefore Stokes theorem implies
$$
\int_{P(a_1^{[r_1]},\ldots,a_k^{[r_k]})}
\omega_{a,r,p}^{|r|-1}={1\over a_1^{r_1}\ldots a_k^{r_k}}.
\leqno(9.18)
$$
Using standard results of integration theory (Fubini, change of variable 
formula...), one obtains:

\claim 9.19. Proposition| Let $f(z)$ be a bounded function on
$P(a_1^{[r_1]},\ldots,a_k^{[r_k]})$ which is continuous outside of
the hyperplane sections $z_s=0$. We also view $f$ as a $\bC^*$-invariant
continuous function on $\prod(\bC^{r_s}\ssm\{0\})$. Then
$$
\eqalign{
&\int_{P(a_1^{[r_1]},\ldots,a_k^{[r_k]})}f(z)\,\omega_{a,r,p}^{|r|-1}\cr
&={(|r|-1)!\over \prod_s a_s^{r_s}}\kern-2pt
\int_{(x,u)\in\Delta_{k-1}\times\prod S^{2r_s-1}}\kern-3pt
f(x_1^{a_1/2p}u_1,\ldots,x_k^{a_k/2p}u_k)\kern-4.5pt
\prod_{1\le s\le k}\kern-2.5pt
{x_s^{r_s-1}\over(r_s-1)!}\,dx\,d\mu(u)\kern-3pt
\cr}
$$
where $\Delta_{k-1}$ is the $(k-1)$-simplex $\{x_s\ge 0$, $\sum x_s=1\}$,
$dx=dx_1\wedge\ldots\wedge dx_{k-1}$ its standard measure, and where
$d\mu(u)=d\mu_1(u_1)\ldots d\mu_k(u_k)$ is the rotation invariant 
probability measure on the product $\prod_sS^{2r_s-1}$ of unit 
spheres in $\bC^{r_1}\times\ldots\times\bC^{r_k}$. As a consequence
$$
\lim_{p\to+\infty}
\int_{P(a_1^{[r_1]},\ldots,a_k^{[r_k]})}f(z)\,\omega_{a,r,p}^{|r|-1}=
{1\over \prod_s a_s^{r_s}}\int_{\prod S^{2r_s-1}}f(u)\,d\mu(u).
$$
\endclaim

Also, by elementary integrations by parts and induction on 
$k,\,r_1,\ldots,r_k$, it can be checked that
$$
\int_{x\in\Delta_{k-1}}
\prod_{1\le s\le k}x_s^{r_s-1}dx_1\ldots dx_{k-1}
={1\over (|r|-1)!}\prod_{1\le s\le k}(r_s-1)!~.
\leqno(9.20)
$$
This implies that $(|r|-1)!\prod_{1\le s\le k}
{x_s^{r_s-1}\over(r_s-1)!}\,dx$ is a probability measure on $\Delta_{k-1}$.

\plainsubsection 9.C. Probabilistic estimate of the curvature of $k$-jet bundles|

Let $(X,V)$ be a compact complex directed non-singular variety. To
avoid any technical difficulty at this point, we first assume that $V$
is a holomorphic vector subbundle of $T_X$, equipped with a smooth
Hermitian metric $h$.

According to the notation already specified in \S$\,$7, we
denote by $J^kV$ the bundle of $k$-jets of holomorphic curves
$f:(\bC,0)\to X$ tangent to $V$ at each point. Let us set $n=\dim_\bC
X$ and $r=\rank_\bC V$. Then $J^kV\to X$ is an algebraic fiber bundle
with typical fiber $\bC^{rk}$, and we get a projectivized $k$-jet bundle 
$$
X^\GG_k:=(J^kV\ssm\{0\})/\bC^*,\qquad \pi_k:X^\GG_k\to X,
\leqno(9.21)
$$
which is a $P(1^{[r]},2^{[r]},\ldots,k^{[r]})$ weighted projective 
bundle over $X$, and we have the direct image formula
$(\pi_k)_*\cO_{X^\GG_k}(m)=\cO(E_{k,m}^\GG V^*)$
(cf.\ Proposition~7.9). In the sequel, we do not make a direct 
use of coordinates, 
because they need not be related in any way to the Hermitian 
metric $h$ of $V$.  Instead, we choose a local holomorphic coordinate frame
$(e_\alpha(z))_{1\le\alpha\le r}$ of $V$ on a neighborhood $U$ of~$x_0$, 
such that
$$
\langle e_\alpha(z),e_\beta(z)\rangle =\delta_{\alpha\beta}+
\sum_{1\le i,j\le n,\,1\le\alpha,\beta\le r}c_{ij\alpha\beta}z_i\overline z_j+
O(|z|^3)\leqno(9.22)
$$
for suitable complex coefficients $(c_{ij\alpha\beta})$. It is a standard fact
that such a normalized coordinate system always exists, and that the 
Chern curvature tensor ${\ii\over 2\pi}D^2_{V,h}$ of $(V,h)$ at $x_0$ 
is then given by
$$
\Theta_{V,h}(x_0)=-{\ii\over 2\pi}
\sum_{i,j,\alpha,\beta}
c_{ij\alpha\beta}\,dz_i\wedge d\overline z_j\otimes e_\alpha^*\otimes e_\beta.
\leqno(9.23)
$$
Consider a local holomorphic connection $\nabla$ on $V_{\restriction U}$ (e.g.\ the 
one which turns $(e_\alpha)$ into a parallel frame), and take 
$\xi_k=\nabla^kf(0)\in V_x$ defined inductively 
by $\nabla^1 f=f'$ and $\nabla^sf=\nabla_{f'}(\nabla^{s-1}f)$. This
gives a local identification
$$
J_kV_{\restriction U}\to V_{\restriction U}^{\oplus k},\qquad
f\mapsto(\xi_1,\ldots,\xi_k)=(\nabla f(0),\ldots,\nabla f^k(0)),
$$
and the weighted $\bC^*$ action on $J_kV$ is expressed in this setting by
$$
\lambda\cdot(\xi_1,\xi_2,\ldots,\xi_k)=(\lambda\xi_1,
\lambda^2\xi_2,\ldots,\lambda^k\xi_k).
$$
Now, we fix a finite open covering 
$(U_\alpha)_{\alpha\in I}$ of~$X$ by open coordinate charts such that
$V_{\restriction U_\alpha}$ is trivial, along with holomorphic connections 
$\nabla_\alpha$ on $V_{\restriction U_\alpha}$. Let $\theta_\alpha$ be a partition of
unity of $X$ subordinate to the covering $(U_\alpha)$. Let us fix 
$p>0$ and small parameters $1=\varepsilon_1\gg\varepsilon_2\gg\ldots\gg
\varepsilon_k>0$. Then we define a global 
weighted Finsler metric on $J^kV$ by putting for any $k$-jet $f\in J^k_xV$
$$
\Psi_{h,p,\varepsilon}(f):=\Big(\sum_{\alpha\in I}
\theta_\alpha(x)\sum_{1\le s\le k}\varepsilon_s^{2p}\Vert\nabla^s_\alpha f(0)
\Vert_{h(x)}^{2p/s}\Big)^{1/p},
\leqno(9.24)
$$
where $\Vert~~\Vert_{h(x)}$ is the Hermitian metric $h$ of $V$ evaluated
on the fiber $V_x$, $x=f(0)$. The function $\Psi_{h,p,\varepsilon}$ satisfies
the fundamental homogeneity property 
$$
\Psi_{h,p,\varepsilon}(\lambda\cdot f)=\Psi_{h,p,\varepsilon}(f)\,|\lambda|^2
\leqno(9.25)
$$
with respect to the $\bC^*$ action on $J^kV$, in other words, it induces
a Hermitian metric on the dual $L^*$ of the tautological $\bQ$-line bundle
$L_k=\cO_{X_k^\GG}(1)$ over $X_k^\GG$. The curvature of $L_k$ is given by
$$
\pi_k^*\Theta_{L_k,\Psi^*_{h,p,\varepsilon}}=dd^c\log\Psi_{h,p,\varepsilon}.
\leqno(9.26)
$$
Our next goal is to compute precisely the curvature and to apply
holomorphic Morse inequalities to $L\to X_k^\GG$ with the above metric.
It might look a priori like an untractable problem, since the definition of
$\Psi_{h,p,\varepsilon}$ is a rather unnatural one. However, the ``miracle''
is that the asymptotic behavior of $\Psi_{h,p,\varepsilon}$ as
$\varepsilon_s/\varepsilon_{s-1}\to 0$ is in some sense uniquely defined 
and very natural.
It will lead to a computable asymptotic formula, which is moreover
simple enough to produce useful results.

\claim 9.27. Lemma| On each coordinate chart $U$ equipped with
a holomorphic connection $\nabla$ of $V_{\restriction U}$, let us define 
the components of a $k$-jet $f\in J^kV$ by $\xi_s=\nabla^sf(0)$,
and consider the rescaling transformation 
$$\rho_{\nabla,\varepsilon}(\xi_1,\xi_2,\ldots,\xi_k)=
(\varepsilon_1^1\xi_1,\varepsilon_2^2\xi_2,\ldots,
\varepsilon_k^k\xi_k)\quad
\hbox{on $J^k_xV$, $x\in U$},
$$
$($it commutes with the $\bC^*$-action but is otherwise unrelated and 
not canonically defined over $X$ as it depends on the choice of 
$\nabla)$. Then, if $p$ is a multiple of $\lcm(1,2,\ldots,k)$ and
$\varepsilon_s/\varepsilon_{s-1}\to 0$ for all $s=2,\ldots,k$, the
rescaled function $\Psi_{h,p,\varepsilon}\circ\rho_{\nabla,\varepsilon}^{-1}
(\xi_1,\ldots,\xi_k)$ converges towards
$$
\bigg(\sum_{1\le s\le k}\Vert \xi_s\Vert^{2p/s}_h\bigg)^{1/p}
$$
on every compact subset of $J^kV_{\restriction U}\ssm\{0\}$,
uniformly in $C^\infty$ topology.
\endclaim

\plainproof. Let $U\subset X$ be an open set on which $V_{\restriction U}$ is trivial
and equipped with some holomorphic connection $\nabla$. Let us pick
another holomorphic connection $\wt\nabla=
\nabla+\Gamma$ where $\Gamma\in H^0(U,\Omega^1_X\otimes
\Hom(V,V))$. Then $\wt\nabla^2f=\nabla^2f+\Gamma(f)(f')\cdot f'$, and
inductively we get
$$
\wt\nabla^sf=\nabla^sf+P_s(f\,;\,\nabla^1f,\ldots,\nabla^{s-1}f)
$$
where $P(x\,;\,\xi_1,\ldots,\xi_{s-1})$ is a polynomial of weighted degree
$s$ in $(\xi_1,\ldots,\xi_{s-1})$ with holomorphic 
coefficients in $x\in U$. In other words, the corresponding change
in the parametrization of $J^kV_{\restriction U}$ is given by a $\bC^*$-homogeneous
transformation
$$
\wt\xi_s=\xi_s+P_s(x\,;\,\xi_1,\ldots,\xi_{s-1}).
$$
Let us introduce the corresponding rescaled components
$$
(\xi_{1,\varepsilon},\ldots,\xi_{k,\varepsilon})=
(\varepsilon_1^1\xi_1,\ldots,\varepsilon_k^k\xi_k),\qquad
(\wt\xi_{1,\varepsilon},\ldots,\wt\xi_{k,\varepsilon})=
(\varepsilon_1^1\wt\xi_1,\ldots,\varepsilon_k^k\wt\xi_k).
$$
Then
$$
\eqalign{
\wt\xi_{s,\varepsilon}
&=\xi_{s,\varepsilon}+
\varepsilon_s^s\,P_s(x\,;\,\varepsilon_1^{-1}\xi_{1,\varepsilon},\ldots,
\varepsilon_{s-1}^{-(s-1)}\xi_{s-1,\varepsilon})\cr
&=\xi_{s,\varepsilon}+O(\varepsilon_s/\varepsilon_{s-1})^s\,
O(\Vert\xi_{1,\varepsilon}\Vert+\ldots+\Vert\xi_{s-1,\varepsilon}
\Vert^{1/(s-1)})^s\cr}
$$
and the error terms are thus polynomials of fixed degree with arbitrarily
small coefficients as $\varepsilon_s/\varepsilon_{s-1}\to 0$. Now, the 
definition of $\Psi_{h,p,\varepsilon}$ consists of glueing the sums
$$
\sum_{1\le s\le k}\varepsilon_s^{2p}\Vert\xi_k\Vert_h^{2p/s}=
\sum_{1\le s\le k}\Vert\xi_{k,\varepsilon}\Vert_h^{2p/s}
$$
corresponding to $\xi_k=\nabla_\alpha^sf(0)$ by means of the partition
of unity $\sum\theta_\alpha(x)=1$. We see that by using the rescaled
variables $\xi_{s,\varepsilon}$ the changes occurring when replacing a
connection $\nabla_\alpha$ by an alternative one $\nabla_\beta$ are
arbitrary small in $C^\infty$ topology, with error terms uniformly
controlled in terms of the ratios $\varepsilon_s/\varepsilon_{s-1}$ on
all compact subsets of $V^k\ssm\{0\}$. This shows that in $C^\infty$
topology,
$\Psi_{h,p,\varepsilon}\circ\rho_{\nabla,\varepsilon}^{-1}(\xi_1,\ldots,\xi_k)$
converges uniformly towards $(\sum_{1\le s\le k}\Vert\xi_k\Vert_h^{2p/s})^{1/p}$,
whatever the trivializing
open set $U$ and the holomorphic connection $\nabla$ used to evaluate
the components and to perform the rescaling are.\qed

Now, we fix a point $x_0\in X$ and a local holomorphic frame 
$(e_\alpha(z))_{1\le\alpha\le r}$ satisfying (9.22) on a neighborhood $U$ 
of~$x_0$. We introduce the rescaled components 
$\xi_s=\varepsilon_s^s\nabla^sf(0)$ on $J^kV_{\restriction U}$ and compute
the curvature of
$$
\Psi_{h,p,\varepsilon}\circ\rho_{\nabla,\varepsilon}^{-1}(z\,;\,\xi_1,\ldots,\xi_k)
\simeq\bigg(\sum_{1\le s\le k}\Vert\xi_s\Vert^{2p/s}_h\bigg)^{1/p},
$$
(by Lemma 9.27, the errors can be taken arbitrary small in 
$C^\infty$ topology). We write $\xi_s=\sum_{1\le\alpha\le r}\xi_{s\alpha}
e_\alpha$. By (9.22) we have
$$
\Vert \xi_s\Vert_h^2=
\sum_\alpha|\xi_{s\alpha}|^2+
\sum_{i,j,\alpha,\beta}c_{ij\alpha\beta}z_i\overline z_j\xi_{s\alpha}
\overline \xi_{s\beta}+O(|z|^3|\xi|^2).
$$
The question is to evaluate the curvature of the weighted metric defined by
$$
\eqalign{
\Psi(z\,;\,\xi_1,\ldots,\xi_k)
&=\bigg(\sum_{1\le s\le k}\Vert\xi_s\Vert^{2p/s}_h\bigg)^{1/p}\cr
&=\bigg(\sum_{1\le s\le k}\Big(\sum_\alpha|\xi_{s\alpha}|^2+\!\!
\sum_{i,j,\alpha,\beta}c_{ij\alpha\beta}z_i\overline z_j\xi_{s\alpha}
\overline\xi_{s\beta}
\Big)^{p/s}\bigg)^{1/p}\kern-4pt{}+O(|z|^3).\cr}
$$
We set $|\xi_s|^2=\sum_\alpha|\xi_{s\alpha}|^2$. A straightforward 
calculation yields
$$
\log\Psi(z\,;\,\xi_1,\ldots,\xi_k)
={1\over p}\log\sum_{1\le s\le k}|\xi_s|^{2p/s}+
\sum_{1\le s\le k}{1\over s}{|\xi_s|^{2p/s}\over \sum_t|\xi_t|^{2p/t}}
\sum_{i,j,\alpha,\beta}c_{ij\alpha\beta}z_i\overline z_j
{\xi_{s\alpha}\overline\xi_{s\beta}\over|\xi_s|^2}+O(|z|^3).
$$
By (9.26), the curvature form of $L_k=\cO_{X_k^\GG}(1)$ 
is given at the central point $x_0$ by the following formula.

\claim 9.28. Proposition| With the above choice of coordinates and with
respect to the rescaled components $\xi_s=\varepsilon_s^s\nabla^sf(0)$ at 
$x_0\in X$, we have the approximate expression
$$
\Theta_{L_k,\Psi^*_{h,p,\varepsilon}}(x_0,[\xi])\simeq
\omega_{a,r,p}(\xi)+{\ii\over 2\pi}
\sum_{1\le s\le k}{1\over s}{|\xi_s|^{2p/s}\over \sum_t|\xi_t|^{2p/t}}
\sum_{i,j,\alpha,\beta}c_{ij\alpha\beta}
{\xi_{s\alpha}\overline\xi_{s\beta}\over|\xi_s|^2}\,dz_i\wedge d\overline z_j,
$$
where the error terms  are
$O(\max_{2\le s\le k}(\varepsilon_s/\varepsilon_{s-1})^s)$ uniformly on
the compact variety $X_k^\GG$. Here $\omega_{a,r,p}$ is the $($degenerate$)$
K\"ahler metric associated with the weight $a=(1^{[r]},2^{[r]},\ldots,k^{[r]})$ 
of the canonical $\bC^*$ action on $J^kV$.
\endclaim

Thanks to the uniform approximation, we can (and will) neglect the error 
terms in the calculations below. Since $\omega_{a,r,p}$ is positive definite
on the fibers of $X_k^\GG\to X$ (at least outside of the axes $\xi_s=0$), 
the index of the $(1,1)$
curvature form $\Theta_{L_k,\Psi^*_{h,p,\varepsilon}}(z,[\xi])$ is equal
to the index of the $(1,1)$-form
$$
\gamma_k(z,\xi):={\ii\over 2\pi}
\sum_{1\le s\le k}{1\over s}{|\xi_s|^{2p/s}\over \sum_t|\xi_t|^{2p/t}}
\sum_{i,j,\alpha,\beta}c_{ij\alpha\beta}(z)
{\xi_{s\alpha}\overline\xi_{s\beta}\over|\xi_s|^2}\,dz_i\wedge d\overline z_j
\leqno(9.29)
$$
depending only on the differentials $(dz_j)_{1\le j\le n}$ on~$X$. The 
$q$-index integral of $(L_k,\Psi^*_{h,p,\varepsilon})$ on $X^\GG_k$ is 
therefore equal to
$$
\eqalign{
&\int_{X^\GG_k(L_k,q)}\Theta_{L_k,\Psi^*_{h,p,\varepsilon}}^{n+kr-1}\cr
&\qquad{}={(n+kr-1)!\over n!(kr-1)!}
\int_{z\in X}\int_{\xi\in P(1^{[r]},\ldots,k^{[r]})}
\omega_{a,r,p}^{kr-1}(\xi)\bOne_{\gamma_k,q}(z,\xi)\gamma_k(z,\xi)^n,\cr}
$$
where $\bOne_{\gamma_k,q}(z,\xi)$ is the characteristic function of the open
set of points where $\gamma_k(z,\xi)$ has signature $(n-q,q)$ in terms of
the $dz_j$'s. Notice that since $\gamma_k(z,\xi)^n$ is~a determinant, the
product $\bOne_{\gamma_k,q}(z,\xi)\gamma_k(z,\xi)^n$ gives rise to a continuous
function on~$X^\GG_k$. Formula (9.20) with \hbox{$r_1=\ldots=r_k=r$} and
$a_s=s$ yields the slightly more explicit
integral
$$
\eqalign{
&\int_{X^\GG_k(L_k,q)}\Theta_{L_k,\Psi^*_{h,p,\varepsilon}}^{n+kr-1}=
{(n+kr-1)!\over n!(k!)^r}~~\times\cr
&\qquad\int_{z\in X}\int_{(x,u)\in\Delta_{k-1}\times(S^{2r-1})^k}
\bOne_{g_k,q}(z,x,u)g_k(z,x,u)^n\,
{(x_1\ldots x_k)^{r-1}\over (r-1)!^k}\,dx\,d\mu(u),\cr}
$$
where $g_k(z,x,u)=\gamma_k(z,x_1^{1/2p}u_1,\ldots,x_k^{k/2p}u_k)$ is given by
$$
g_k(z,x,u)={\ii\over 2\pi}\sum_{1\le s\le k}{1\over s}x_s
\sum_{i,j,\alpha,\beta}c_{ij\alpha\beta}(z)\,
u_{s\alpha}\overline u_{s\beta}\,dz_i\wedge d\overline z_j
\leqno(9.30)
$$
and $\bOne_{g_k,q}(z,x,u)$ is the characteristic function of its $q$-index 
set. Here 
$$
d\nu_{k,r}(x)=(kr-1)!\,{(x_1\ldots x_k)^{r-1}\over (r-1)!^k}\,dx
\leqno(9.31)
$$
is a probability measure on $\Delta_{k-1}$, and we can rewrite
$$
\leqalignno{
&\int_{X^\GG_k(L_k,q)}\Theta_{L_k,\Psi^*_{h,p,\varepsilon}}^{n+kr-1}
={(n+kr-1)!\over n!(k!)^r(kr-1)!}~~\times\cr
&\qquad\int_{z\in X}
\int_{(x,u)\in\Delta_{k-1}\times(S^{2r-1})^k}
\bOne_{g_k,q}(z,x,u)g_k(z,x,u)^n\,d\nu_{k,r}(x)\,d\mu(u).&(9.32)\cr}
$$
Now, formula (9.30) shows that $g_k(z,x,u)$ is a ``Monte Carlo''
evaluation of the curvature tensor, obtained by averaging the curvature
at random points $u_s\in S^{2r-1}$ with certain positive weights $x_s/s\,$; 
we should then think of the \hbox{$k$-jet}
$f$ as some sort of random variable such that the derivatives 
$\nabla^kf(0)$ are uniformly distributed in all directions. Let us compute
the expected value of
$(x,u)\mapsto g_k(z,x,u)$ with respect to the probability measure
$d\nu_{k,r}(x)\,d\mu(u)$. Since 
$\int_{S^{2r-1}}u_{s\alpha}\overline u_{s\beta}d\mu(u_s)={1\over r}
\delta_{\alpha\beta}$ and $\int_{\Delta_{k-1}}x_s\,d\nu_{k,r}(x)={1\over k}$,
we find
$$
{\bf E}(g_k(z,\bu,\bu))={1\over kr}
\sum_{1\le s\le k}{1\over s}\cdot{\ii\over 2\pi}\sum_{i,j,\alpha}
c_{ij\alpha\alpha}(z)\,dz_i\wedge d\overline z_j.
$$
In other words, we get the normalized trace of the curvature, i.e.
$$
{\bf E}(g_k(z,\bu,\bu))={1\over kr}
\Big(1+{1\over 2}+\ldots+{1\over k}\Big)\Theta_{\det(V^*),\det h^*},
\leqno(9.33)
$$
where $\Theta_{\det(V^*),\det h^*}$ is the $(1,1)$-curvature form of
$\det(V^*)$ with the metric induced by~$h$. It is natural to guess that 
$g_k(z,x,u)$ behaves asymptotically as its expected value
${\bf E}(g_k(z,\bu,\bu))$ when $k$ tends to infinity. If we replace brutally 
$g_k$ by its expected value in (9.32), we get the integral
$$
{(n+kr-1)!\over n!(k!)^r(kr-1)!}{1\over (kr)^n}
\Big(1+{1\over 2}+\ldots+{1\over k}\Big)^n\int_X\bOne_{\eta,q}\eta^n,
$$
where $\eta:=\Theta_{\det(V^*),\det h^*}$ and $\bOne_{\eta,q}$ is the
characteristic function of its $q$-index set in~$X$. The leading constant is
equivalent to $(\log k)^n/n!(k!)^r$ modulo a multiplicative factor
$(1+O(1/\log k))$. By working out a more precise analysis
of the deviation, the following result has been proved in [Dem11] and 
[Dem12].

\claim 9.34. Probabilistic estimate|
Fix smooth Hermitian metrics $h$ on $V$ and
$\omega={\ii\over 2\pi} \sum\omega_{ij}dz_i\wedge d\overline z_j$ on $X$. 
Denote by $\Theta_{V,h}=-{\ii\over 2\pi}\sum
c_{ij\alpha\beta}dz_i\wedge d\overline z_j\otimes e_\alpha^*\otimes
e_\beta$ the curvature tensor of $V$ with respect to an $h$-orthonormal frame
$(e_\alpha)$, and put
$$
\eta(z)=\Theta_{\det(V^*),\det h^*}={\ii\over 2\pi}\sum_{1\le i,j\le n}\eta_{ij}
dz_i\wedge d\overline z_j,\qquad
\eta_{ij}=\sum_{1\le\alpha\le r}c_{ij\alpha\alpha}.
$$
Finally consider the $k$-jet line bundle $L_k=\cO_{X_k^\GG}(1)\to
X_k^\GG$ equipped with the induced metric $\Psi^*_{h,p,\varepsilon}$
$($as defined above, with $1=\varepsilon_1\gg\varepsilon_2\gg\ldots\gg
\varepsilon_k>0)$. When $k$ tends 
to infinity, the integral of the top power of the curvature of $L_k$ on its
$q$-index set $X^\GG_k(L_k,q)$ is given by
$$
\int_{X^\GG_k(L_k,q)}\Theta_{L_k,\Psi^*_{h,p,\varepsilon}}^{n+kr-1}=
{(\log k)^n\over n!\,(k!)^r}\bigg(
\int_X\bOne_{\eta,q}\eta^n+O((\log k)^{-1})\bigg)
$$
for all $q=0,1,\ldots,n$, and the error term $O((\log k)^{-1})$ can be 
bounded explicitly in terms of $\Theta_V$, $\eta$ and $\omega$. Moreover, the 
left hand side is identically zero for $q>n$.
\endclaim

The final statement follows from the observation that the curvature of
$L_k$ is positive along the fibers of $X_k^\GG\to X$, by the 
plurisubharmonicity of the weight (this is true even 
when the partition of unity terms are taken into account, since they
depend only on the base); therefore the $q$-index sets are empty for
$q>n$. It will be useful to extend the above estimates to the 
case of sections of
$$
L_k=\cO_{X_k^\GG}(1)\otimes
\pi_k^*\cO\Big({1\over kr}\Big(1+{1\over 2}+\ldots+{1\over k}\Big)F\Big),
\leqno(9.35)
$$
where $F\in\Pic_\bQ(X)$ is an arbitrary $\bQ$-line bundle on~$X$ and 
$\pi_k:X_k^\GG\to X$ is the natural projection. We assume here
that $F$ is also equipped with a smooth Hermitian metric $h_F$. In formulas
(9.32), (9.33) and estimate 9.34, the renormalized curvature
$\eta_k(z,x,u)$ of $L_k$ takes the form
$$
\eta_k(z,x,u)={1\over {1\over kr}(1+{1\over 2}+\ldots+{1\over k})}g_k(z,x,u)+
\Theta_{F,h_F}(z),
\leqno(9.36)
$$
and by the same calculations its  expected value is
$$
\eta(z):={\bf E}(\eta_k(z,\bu,\bu))=\Theta_{\det V^*,\det h^*}(z)+
\Theta_{F,h_F}(z).
\leqno(9.37)
$$
Then the variance estimate for $\eta_k-\eta$ is unchanged, and the
$L^p$ bounds for $\eta_k$ are still valid, since our forms are just shifted
by adding the constant smooth term $\Theta_{F,h_F}(z)$. The probabilistic
estimate 9.34 is therefore still true exactly in the same form, provided
we use (9.35 -- 9.37) instead of the previously defined $L_k$, $\eta_k$
and~$\eta$. An application of holomorphic Morse inequalities gives the 
desired cohomology estimates for 
$$
\eqalign{
h^q\Big(X,E_{k,m}^\GG V^*&{}\otimes
\cO\Big({m\over kr}\Big(1+{1\over 2}+\ldots+{1\over k}\Big)F\Big)\Big)\cr
&{}=h^q(X_k^\GG,\cO_{X_k^\GG}(m)\otimes
\pi_k^*\cO\Big({m\over kr}\Big(1+{1\over 2}+\ldots+{1\over k}\Big)F\Big)\Big),
\cr}
$$
provided $m$ is sufficiently divisible to give a multiple of $F$ which
is a $\bZ$-line bundle.

\claim 9.38. Theorem| Let $(X,V)$ be a directed manifold, $F\to X$ a
$\bQ$-line bundle, $(V,h)$ and $(F,h_F)$ smooth Hermitian structures on $V$ 
and $F$ respectively. We define
$$
\eqalign{
L_k&=\cO_{X_k^\GG}(1)\otimes
\pi_k^*\cO\Big({1\over kr}\Big(1+{1\over 2}+\ldots+{1\over k}\Big)F\Big),\cr
\eta&=\Theta_{\det V^*,\det h^*}+\Theta_{F,h_F},\cr}
$$
and let $X(\eta,q)$ be the open set of points $x\in X$ where
$\eta(x)$ has signature $(q,n-q)$. We also set $X(\eta,\leq q)=
\bigcup_{j\leq q}X(\eta,j)$.
Then for all $q\ge 0$ and all $m\gg k\gg 1$ such that 
m is sufficiently divisible, we have
$$\leqalignno{\kern20pt
h^q(X_k^\GG,\cO(L_k^{\otimes m}))&\le {m^{n+kr-1}\over (n+kr-1)!}
{(\log k)^n\over n!\,(k!)^r}\bigg(
\int_{X(\eta,q)}(-1)^q\eta^n+O((\log k)^{-1})\bigg),&\hbox{\rm(a)}\cr
h^0(X_k^\GG,\cO(L_k^{\otimes m}))&\ge {m^{n+kr-1}\over (n+kr-1)!}
{(\log k)^n\over n!\,(k!)^r}\bigg(
\int_{X(\eta,\le 1)}\eta^n-O((\log k)^{-1})\bigg),&\hbox{\rm(b)}\cr
\cr
\chi(X_k^\GG,\cO(L_k^{\otimes m}))&={m^{n+kr-1}\over (n+kr-1)!}
{(\log k)^n\over n!\,(k!)^r}\big(
c_1(V^*\otimes F)^n+O((\log k)^{-1})\big).&\hbox{\rm(c)}\cr
\cr}
$$
\vskip-4pt
\endclaim

Green and Griffiths [GrGr80] already checked the Riemann-Roch
calculation (9.38$\,$c) in the special case
$V=T_X^*$ and $F=\cO_X$. Their proof is much simpler since it relies only
on Chern class calculations, but it cannot provide any information on
the individual cohomology groups, except in very special cases where
vanishing theorems can be applied; in fact in dimension 2, the
Euler characteristic satisfies $\chi=h^0-h^1+h^2\le h^0+h^2$, and hence
it is enough to get the vanishing of the top cohomology group $H^2$
to infer $h^0\ge\chi\,$; this works for surfaces by means of a well-known
vanishing theorem of Bogomolov which implies in general
$$H^n\bigg(X,E_{k,m}^\GG T_X^*\otimes
\cO\Big({m\over kr}\Big(1+{1\over 2}+\ldots+{1\over k}\Big)F\Big)\Big)\bigg)=0
$$
as soon as $K_X\otimes F$ is big and $m\gg 1$.

In fact, thanks to Bonavero's singular holomorphic Morse inequalities 
[Bon93], everything works almost unchanged 
in the case where $V\subset T_X$
has singularities and $h$ is an admissible metric on $V$ (see Definition~9.6).
We only have to find a blow-up $\mu:\wt X_k\to X_k$ so that
the resulting pull-backs $\mu^*L_k$ and $\mu^*V$ are locally free,
and $\mu^*\det h^*$, $\mu^*\Psi_{h,p,\varepsilon}$ only have divisorial
singularities. Then $\eta$ is a $(1,1)$-current with logarithmic poles,
and we have to deal with smooth metrics on $
\mu^*L_k^{\otimes m}\otimes\cO(-mE_k)$ where $E_k$ is a certain effective 
divisor on $X_k$ (which, by our assumption in Definition~9.6,
does not project onto $X$). The cohomology groups involved are then
the twisted cohomology groups
$$
H^q(X_k^\GG,\cO(L_k^{\otimes m})\otimes\cJ_{k,m}),
$$
where $\cJ_{k,m}=\mu_*(\cO(-mE_k))$ is the corresponding multiplier ideal sheaf,
and the Morse integrals need only be evaluated in the complement of the 
poles, i.e., on $X(\eta,q)\ssm S$ with $S=\Sing(V)\cup\Sing(h)$. Since
$$
(\pi_k)_*\big(\cO(L_k^{\otimes m})\otimes\cJ_{k,m}\big)\subset
E_{k,m}^\GG V^*\otimes
\cO\Big({m\over kr}\Big(1+{1\over 2}+\ldots+{1\over k}\Big)F\Big)\Big)
$$
we still get a lower bound for the $h^0$ of the latter sheaf (or for the $h^0$
of the un-twisted line bundle $\cO(L_k^{\otimes m})$ on $X_k^\GG$).
If we assume that $K_V\otimes F$ is big, these considerations
also allow us to obtain a strong estimate in terms of the volume, by
using an approximate Zariski decomposition on a suitable blow-up of~$(X,V)$.
The following corollary implies Theorem~9.1 as a consequence.

\claim 9.39. Corollary|
If $F$ is an arbitrary $\bQ$-line bundle over~$X$, one has
$$
\eqalign{
h^0\bigg(&X_k^\GG,\cO_{X_k^\GG}(m)\otimes\pi_k^*\cO
\Big({m\over kr}\Big(1+{1\over 2}+\ldots+{1\over k}\Big)F\Big)\bigg)\cr
&\ge {m^{n+kr-1}\over (n+kr-1)!}
{(\log k)^n\over n!\,(k!)^r}\Big(
\Vol(K_V\otimes F)-O((\log k)^{-1})\Big)-o(m^{n+kr-1}),\cr}
$$
when $m\gg k\gg 1$, in particular there are many sections of the
$k$-jet differentials of degree $m$ twisted by the appropriate
power of $F$ if $K_V\otimes F$ is big.
\endclaim

\plainproof. The volume is computed here as usual, i.e.\ after performing a
suitable log-resolution $\mu:\wt X\to X$ which converts $K_V$ into 
an invertible sheaf. There is of course nothing to prove if
$K_V\otimes F$ is not big, so we can assume $\Vol(K_V\otimes F)>0$.
Let us fix smooth Hermitian metrics $h_0$ on $T_X$
and $h_F$ on~$F$. They induce a metric $\mu^*(\det(h_0^{-1})\otimes
h_F)$ on $\mu^*(K_V\otimes F)$ which, by our definition of $K_V$, is
a smooth metric (the divisor produced by the log-resolution gets simplified
with the degeneration divisor of the pull-back of the quotient metric on
$\det(V^*)$ induced by $\cO(\Lambda^rT_X^*)\to \cO(\Lambda^r V^*)$).
By the result of Fujita [Fuji94] on
approximate Zariski decomposition, for every $\delta>0$, one can find
a modification $\mu_\delta:\wt X_\delta\to X$ dominating
$\mu$ such that
$$
\mu_\delta^*(K_V\otimes F) =\cO_{\wt X_\delta}(A+E)
$$
where $A$ and $E$ are $\bQ$-divisors, $A$ ample and $E$ effective,
with 
$$\Vol(A)=A^n\ge \Vol(K_V\otimes F)-\delta.$$
If we take a smooth metric $h_A$ with positive definite curvature form
$\Theta_{A,h_A}$, then we get a singular Hermitian metric $h_Ah_E$ on
$\mu_\delta^*(K_V\otimes F)$ with poles along $E$, i.e.\ the quotient
$h_Ah_E/\mu_\delta^*(\det(h_0^{-1})\otimes h_F)$ is of the
form $e^{-\varphi}$ where
$\varphi$ is quasi-psh with log poles $\log|\sigma_E|^2$ 
(mod $C^\infty(\wt X_\delta))$ precisely given
by the divisor~$E$. We then only need to take the singular metric $h$
on $T_X$ defined by
$$
h=h_0e^{{1\over r}(\mu_\delta)_*\varphi}
$$
(the choice of the factor ${1\over r}$ is there to correct adequately 
the metric on $\det V$). By construction $h$ induces an 
admissible metric on $V$ and the resulting 
curvature current $\eta=\Theta_{K_V,\det h^*}+\Theta_{F,h_F}$ is such that
$$
\mu_\delta^*\eta = \Theta_{A,h_A} +[E],\qquad
\hbox{$[E]={}$current of integration on $E$.}
$$
Then the $0$-index Morse integral in the complement of the poles 
is given by
$$
\int_{X(\eta,0)\ssm S}\eta^n=\int_{\wt X_\delta}\Theta_{A,h_A}^n=A^n\ge
\Vol(K_V\otimes F)-\delta
$$
and Corollary~9.39 follows from the fact that $\delta$ can be taken arbitrary 
small.\qed
\medskip

The following corollary implies Theorem~0.12.

\claim 9.40. Corollary|Let $(X,V)$ be a projective directed manifold such that
$K_V^\bullet$ is big, and $A$ an ample $\bQ$-divisor on $X$ such that 
$K_V^\bullet\otimes\cO(-A)$ is still big. Then, if we put $r=\rank V$ and
\hbox{$\delta_k={1\over kr}(1+{1\over 2}+\ldots+{1\over k})$},
the space of global invariant jet differentials
$$
H^0(X,E_{k,m}V^*\otimes\cO(-m\delta_kA))
$$
has $($many$)$ non-zero sections for $m\gg k\gg 1$ and $m$ sufficiently 
divisible.
\endclaim

\plainproof. Corollary 9.39 produces a non-zero section 
$P\in H^0(E_{k,m}^\GG V^*\otimes\cO_X(-m\delta_kA))$ for $m\gg k\gg 1$, and
the arguments given in subsection 7.E (cf.\ (7.36)) yield
a non-zero section
$$Q\in H^0(E_{k,m'}V^*\otimes\cO_X(-m\delta_kA)),\quad m'\le m.$$
By raising $Q$ to some power $p$ and using
a section $\sigma\in H^0(X,\cO_X(dA))$, we obtain a section
$$
Q^p\sigma^{mq}\in H^0(X,E_{k,pm'}V^*\otimes\cO(-m(p\delta_k-qd)A)).
$$
One can adjust $p$ and $q$ so that $m(p\delta_k-qd)=pm'\delta_k$ and
$pm'\delta_kA$ is an integral divisor.\qed

\claim 9.41. Example| {\rm In some simple cases, the above estimates can 
lead to very explicit results. Take for instance $X$ to be a smooth
complete intersection of multidegree $(d_1,d_2,\ldots,d_s)$ in $\bP^{n+s}_\bC$
and consider the absolute case $V=T_X$. Then 
$K_X=\cO_X(d_1+\ldots+d_s-n-s-1)$ and one can check via
explicit bounds of the error terms (cf.\ [Dem11], [Dem12])
that a sufficient condition for the existence of
sections is
$$
k\ge\exp\Big(7.38\,n^{n+\frac{1}{2}}\Big({\sum d_j+1\over\sum d_j-n-s-a-1}\Big)^n\Big).
$$
This is good in view of the fact that we can cover arbitrary smooth 
complete intersections of general type. On the other hand, even when the
degrees $d_j$ tend to $+\infty$, we still get a large lower bound
$k\sim \exp(7.38\,n^{n+\smash{\frac{1}{2}}})$ on the order of jets, and this
is far from being optimal$\,$: Diverio [Div08, Div09] has shown e.g.\ that
one can take
$k=n$ for smooth hypersurfaces of high degree, using the 
algebraic Morse inequalities of Trapani [Tra95].
The next paragraph uses essentially the same idea, in our more analytic
setting.}
\endclaim

\plainsubsection 9.D. Non probabilistic estimate of the Morse integrals|

We assume here that the curvature tensor $(c_{ij\alpha\beta})$ satisfies a lower
bound
$$
\sum_{i,j,\alpha,\beta}c_{ij\alpha\beta}\xi_i\ol\xi_ju_\alpha\ol u_\beta\ge -\sum\gamma_{ij}\xi_i
\ol\xi_j\;|u|^2,
\qquad\forall\xi\in T_X,~u\in V,\leqno(9.42)
$$
for some semi-positive $(1,1)$-form $\gamma={i\over 2\pi}\sum\gamma_{ij}(z)\,
dz_i\wedge d\ol z_j$ on~$X$. This is the
same as assuming that the curvature tensor of $(V^*,h^*)$ satisfies
the semi-positivity condition
$$
\Theta_{V^*,h^*}+\gamma\otimes\Id_{V^*}\ge 0\leqno(9.42')
$$
in the sense of Griffiths, or equivalently $\Theta_{V,h}-\gamma\otimes\Id_V\le 0$.
Thanks to the compactness of~$X$, such a form $\gamma$ always
exists if $h$ is an admissible metric on~$V$. Now, instead of replacing $\Theta_V$
with its trace free part $\wt\Theta_V$ and exploiting a Monte Carlo convergence
process, we replace $\Theta_V$ with $\Theta_V^\gamma=\Theta_V-\gamma\otimes\Id_V\le 0$, 
i.e.\ $c_{ij\alpha\beta}$ by 
$c_{ij\alpha\beta}^\gamma=c_{ij\alpha\beta}+\gamma_{ij}\delta_{\alpha\beta}$. Also, we take
a line bundle $F=A^{-1}$ with $\Theta_{A,h_A}\ge 0$, i.e.\ $F$ semi-negative.
Then our earlier formulas in Prop.~9.28, and (9.35), (9.36) become instead
$$
\leqalignno{
&g_k^\gamma(z,x,u)={\ii\over 2\pi}\sum_{1\le s\le k}{1\over s}x_s
\sum_{i,j,\alpha,\beta}c_{ij\alpha\beta}^\gamma(z)\,
u_{s\alpha}\overline u_{s\beta}\,dz_i\wedge d\overline z_j\ge 0,
&(9.43)\cr
&L_k=\cO_{X_k^\GG}(1)\otimes
\pi_k^*\cO\Big(-{1\over kr}\Big(1+{1\over 2}+\ldots+{1\over k}\Big)A\Big),
&(9.44)\cr
&\Theta_{L_k}=\eta_k(z,x,u)={1\over {1\over kr}(1+{1\over 2}+\ldots+{1\over k})}
g_k^\gamma(z,x,u)-(\Theta_{A,h_A}(z)+r\gamma(z)).&(9.45)\cr}
$$
In fact, replacing $\Theta_V$ by $\Theta_V-\gamma\otimes\Id_V$ has the effect of
replacing $\Theta_{\det V^*}=\Tr\Theta_{V^*}$ by $\Theta_{\det V^*}+r\gamma$. The major
gain that we have is that $\eta_k=\Theta_{L_k}$ is now expressed as a difference
of semi-positive $(1,1)$-forms, and we can exploit the following simple lemma, which is
the key to derive algebraic Morse inequalities from their analytic form
(cf.\ [Dem94], Theorem~12.3).

\claim 9.46.~Lemma|Let $\eta=\alpha-\beta$ be a difference of semi-positive 
$(1,1)$-forms on an $n$-dimensional complex manifold~$X$, 
and let $\bOne_{\eta,\le q}$ be the characteristic function of the
open set where $\eta$ is non-degenerate with a number of negative eigenvalues 
at most equal to~$q$.
Then
$$
(-1)^q\bOne_{\eta,\le q}~\eta^n\le \sum_{0\le j\le q}(-1)^{q-j}\alpha^{n-j}\beta^j,
$$
in particular
$$
\bOne_{\eta,\le 1}~\eta^n\ge \alpha^n-n\alpha^{n-1}\wedge \beta\qquad\hbox{for $q=1$.}
$$
\endclaim

\plainproof. Without loss of generality, we can assume $\alpha>0$ positive definite, so that
$\alpha$ can be taken as the base Hermitian metric on~$X$. Let us denote by
$$
\lambda_1\ge\lambda_2\ge\ldots\ge\lambda_n\ge 0
$$
the eigenvalues of $\beta$ with respect to $\alpha$. The eigenvalues of $\eta=\alpha-\beta$
are then given by 
$$
1-\lambda_1\le\ldots\le 1-\lambda_q\le 1-\lambda_{q+1}\le\ldots\le 1-\lambda_n\,;
$$
hence the open set $\{\lambda_{q+1}<1\}$ coincides with the support of 
$\bOne_{\eta,\le q}$, except that it may also contain a part of 
the degeneration set $\eta^n=0$. On the other hand we have
$${n\choose j}\alpha^{n-j}\wedge\beta^j=\sigma_n^j(\lambda)\,\alpha^n,$$
where $\sigma_n^j(\lambda)$ is the $j$-th elementary symmetric function in the $\lambda_j$'s.
Thus, to prove the lemma, we only have to check that
$$\sum_{0\le j\le q}(-1)^{q-j}\sigma_n^j(\lambda)-
\bOne_{\{\lambda_{q+1}<1\}}(-1)^q\prod_{1\le j\le n}(1-\lambda_j)\ge 0.$$
This is easily done by induction on~$n$ (just split apart the parameter
$\lambda_n$ and write $\sigma_n^j(\lambda)=
\sigma_{n-1}^j(\lambda)+\sigma_{n-1}^{j-1}(\lambda)\,\lambda_n$).\qed
\medskip

We apply here Lemma 9.46 with
$$
\alpha=g_k^\gamma(z,x,u),\qquad\beta=\beta_k=
{1\over kr}\Big(1+{1\over 2}+\ldots+{1\over k}\Big)(\Theta_{A,h_A}+r\gamma),
$$
which are both semi-positive by our assumption. The analogue of (9.32) leads to
$$
\eqalign{
&\int_{X^\GG_k(L_k,\le 1)}\Theta_{L_k,\Psi^*_{h,p,\varepsilon}}^{n+kr-1}\cr
&\quad{}={(n+kr-1)!\over n!(k!)^r(kr-1)!}\int_{z\in X}
\int_{(x,u)\in\Delta_{k-1}\times(S^{2r-1})^k}
\bOne_{g_k^\gamma-\beta_k,\le 1}\;(g_k^\gamma-\beta_k)^n\,d\nu_{k,r}(x)\,d\mu(u)\cr
&\quad{}\ge{(n+kr-1)!\over n!(k!)^r(kr-1)!}\int_{z\in X}
\int_{(x,u)\in\Delta_{k-1}\times(S^{2r-1})^k}
((g_k^\gamma)^n-n(g_k^\gamma)^{n-1}\wedge\beta_k)\,d\nu_{k,r}(x)\,d\mu(u).\cr}
$$
The resulting integral now produces a ``closed formula'' which can be
expressed solely in terms of Chern classes (at least if we assume that
$\gamma$ is the Chern form of some semi-positive line bundle). It is
just a matter of routine to find a sufficient condition for the
positivity of the integral. One can first observe that $g_k^\gamma$ is
bounded from above by taking the trace of $(c_{ij\alpha\beta})$, in
this way we get
$$
0\le g_k^\gamma\le\bigglp2pt(\sum_{1\le s\le k}{x_s\over s}\biggrp2pt)\big(\Theta_{\det V^*}+r\gamma\big)
$$
where the right hand side no longer depends on $u\in (S^{2r-1})^k$. 
Also, $g_k^\gamma$ can be written as a sum of semi-positive $(1,1)$-forms
$$
g_k^\gamma =\sum_{1\le s\le k}{x_s\over s}\theta^\gamma(u_s),\qquad
\theta^\gamma(u)=\sum_{i,j,\alpha,\beta}c_{ij\alpha\beta}^\gamma
u_\alpha\ol u_\beta\,dz_i\wedge d\ol z_j,
$$
and hence for $k\ge n$ we have
$$
(g_k^\gamma)^n\ge n!\sum_{1\le s_1<\ldots<s_n\le k}
{x_{s_1}\ldots x_{s_n}\over s_1\ldots s_n}\,
\theta^\gamma(u_{s_1})\wedge \theta^\gamma(u_{s_2})\wedge\ldots\wedge \theta^\gamma(u_{s_n}).
$$
Since $\int_{S^{2r-1}}\theta^\gamma(u)\,d\mu(u)={1\over r}\Tr(\Theta_{V^*}+\gamma)=
{1\over r}\Theta_{\det V^*}+\gamma$, we infer from this
$$
\eqalign{
&\int_{(x,u)\in\Delta_{k-1}\times(S^{2r-1})^k}
(g_k^\gamma)^n\,d\nu_{k,r}(x)\,d\mu(u)\cr
&\qquad\ge 
n!\sum_{1\le s_1<\ldots<s_n\le k}
{1\over s_1\ldots s_n}\Big(\int_{\Delta_{k-1}}x_1\ldots x_n\,d\nu_{k,r}(x)\Big)
\Big({1\over r}\Theta_{\det V^*}+\gamma\Big)^n.\cr}
$$
By putting everything together, we conclude:

\claim 9.47.~Theorem|Assume that $\Theta_{V^*}\ge-\gamma\otimes\Id_{V^*}$ with a
semi-positive $(1,1)$-form $\gamma$ on $X$. Then the Morse integral of the line bundle
$$
L_k=\cO_{X_k^\GG}(1)\otimes
\pi_k^*\cO\Big(-{1\over kr}\Big(1+{1\over 2}+\ldots+{1\over k}\Big)A\Big),\qquad
A\ge 0
$$
satisfies for $k\ge n$ the inequality
$$
\leqalignno{
&{1\over (n+kr-1)!}\int_{X^\GG_k(L_k,\le 1)}\Theta_{L_k,\Psi^*_{h,p,\varepsilon}}^{n+kr-1}\cr
&\qquad{}\ge{1\over n!(k!)^r(kr-1)!}
\int_Xc_{n,r,k}\big(\Theta_{\det V^*}+r\gamma\big)^n-c'_{n,r,k}
\big(\Theta_{\det V^*}+r\gamma\big)^{n-1}\wedge
\big(\Theta_{A,h_A}+r\gamma\big)&(*)\cr}
$$
where
$$
\eqalign{
c_{n,r,k}&={n!\over r^n}\bigglp2pt(\sum_{1\le s_1<\ldots<s_n\le k}
{1\over s_1\ldots s_n}\biggrp2pt)\int_{\Delta_{k-1}}x_1\ldots x_n\,d\nu_{k,r}(x),\cr
c'_{n,r,k}&={n\over kr}\Big(1+{1\over 2}+\ldots+{1\over k}\Big)
\int_{\Delta_{k-1}}\bigglp2pt(\sum_{1\le s\le k}{x_s\over s}\biggrp2pt)^{n-1}\,d\nu_{k,r}(x).
\cr}
$$
Especially we have a lot of sections in $H^0(X_k^\GG,mL_k)$, $m\gg 1$, as soon
as the difference occurring in $(*)$ is positive.
\endclaim

The statement is also true for $k<n$, but then $c_{n,r,k}=0$ and the lower bound $(*)$ cannot
be positive. By Corollary 9.11, it still provides a non-trivial lower bound for
$h^0(X_k^\GG,mL_k)-h^1(X_k^\GG,mL_k)$, though. For $k\ge n$ we have $c_{n,r,k}>0$ and
$(*)$ will be positive if $\Theta_{\det V^*}$ is large enough. By Formula~9.20 we have
$$
c_{n,r,k}={n!\,(kr-1)!\over (n+kr-1)!}
\smash{\sum_{1\le s_1<\ldots<s_n\le k}{1\over s_1\ldots s_n}\ge
{(kr-1)!\over (n+kr-1)!}},\leqno(9.48)
$$
(with equality for $k=n$). On the other hand, for any multi-index
$(\beta_1,\ldots,\beta_k)\in\bN^k$ with $\sum\beta_s=p$,
the H\"older inequality implies
$$
\int_{\Delta_{k-1}}x_1^{\beta_1}\ldots x_k^{\beta_k}\,d\nu_{k,r}(x)
\leq \prod_{s=1}
\bigg(\int_{\Delta_{k-1}}x_s^p\,d\nu_{k,r}(x)\bigg)^{\beta_s/p}
=\int_{\Delta_{k-1}}x_1^p\,d\nu_{k,r}(x).
$$
An expansion of $\big(\sum_{1\le s\le k}{x_s\over s}\big)^{n-1}$ by means of the
multinomial formula then yields
$$
\int_{\Delta_{k-1}}\bigglp2pt(\sum_{1\le s\le k}{x_s\over s}\biggrp2pt)^{n-1}\,d\nu_{k,r}(x)
\leq
\int_{\Delta_{k-1}}\bigglp2pt(\sum_{1\le s\le k}{1\over s}\biggrp2pt)^{n-1}
x_1^{n-1}\,d\nu_{k,r}(x).
$$
On the other hand, it is obvious that 
$\int_{\Delta_{k-1}}\big(\sum_{1\le s\le k}{x_s\over s}\big)^{n-1}\,d\nu_{k,r}(x)\geq
\int_{\Delta_{k-1}}x_1^{n-1}\,d\nu_{k,r}(x)$, thus the error in the above 
upper bound is at most by a factor 
$(1+{1\over 2}+\ldots+{1\over k})^n\leq(1+\log k)^n$. 
From this, we infer again by Formula (9.20) that
$$
\leqalignno{
c'_{n,r,k}&\leq
{n\over kr}\Big(1+{1\over 2}+\ldots+{1\over k}\Big)^n
\int_{\Delta_{k-1}}x_1^{n-1}\,d\nu_{k,r}(x),\cr
&={n\over kr}\Big(1+{1\over 2}+\ldots+{1\over k}\Big)^n~{(n+r-2)!\over (r-1)!}{(kr-1)!\over (n+kr-2)!}.&(9.49)\cr}
$$
Since ${n+kr-1\over k}=r+{n-1\over k}\leq n+r-1$, our bounds
(9.48) and (9.49) imply
$$
\leqalignno{
{c'_{n,r,k}\over c_{n,r,k}}
&\leq {n\over k}\Big(1+{1\over 2}+\ldots+{1\over k}\Big)^n
~{(n+r-2)!\over r!}\,(n+kr-1),&(9.50)\cr
{c'_{n,r,k}\over c_{n,r,k}}&\leq n\,\Big(1+{1\over 2}+\ldots+{1\over k}\Big)^n
~{(n+r-1)!\over r!}.
&(9.51)}
$$
The right hand side of (9.51) increases with $r$.
For $r\leq n$, the Stirling formula yields
$$
{c'_{n,r,k}\over c_{n,r,k}}
<(1+\log k)^n\,{(2n)!\over 2\,n!}
<(1+\log k)^n\,{\sqrt{2n}\,({2n\over e})^{2n}\over
2\sqrt{n}\,({n\over e})^n}={1\over\sqrt{2}}
\big(4e^{-1}\,n(1+\log k)\big)^n.\leqno(9.52_n)
$$
Up to the constant $4e^{-1}$, this is essentially the same bound as
the one obtained in [Dem12], which, however, included a numerical
mistake, making unclear whether the constant $4e^{-1}>1$ could be dropped 
there, as would follow from the claimed estimate. 
We will later need the following particular values 
(cf.\ Formula~(9.20) and [Dem11, Lemma 2.20]):
$$
\leqalignno{
\noalign{\vskip-3pt}
&c_{2,2,2}={1\over 20},\kern4.8pt\qquad c'_{2,2,2}={9\over 16},\kern9.8pt\qquad
{c'_{2,2,2}\over c_{2,2,2}}={45\over 4},&(9.52_2)\cr
&c_{3,3,3}={1\over 990},\qquad c'_{3,3,3}={451\over 4860},\qquad
{c'_{3,3,3}\over c_{3,3,3}}={4961\over 54}.&(9.52_3)\cr
\noalign{\vskip10pt}}
$$

\section{Hyperbolicity properties of  hypersurfaces of high degree}

\plainsubsection 10.A. Global generation of the twisted tangent space of the universal family|

In [Siu02, Siu04], Y.T.~ Siu developed a new strategy to produce jet
differentials, involving meromorphic vector fields on the total space of
jet bundles -- these vector fields are used to differentiate the
sections of $E_{k,m}^\GG$ so as to produce new ones with less zeroes. The approach
works especially well on universal families of hypersurfaces in
projective space, thanks to the good positivity properties of the
relative tangent bundle, as shown by L.~Ein [Ein88, Ein91] and C.~Voisin [Voi96].
This allows at least to prove the hyperbolicity of generic surfaces
and generic 3-dimensional hypersurfaces of sufficiently high
degree. We reproduce here the improved approach given by [Pau08] for the
twisted global generation of the tangent space of the space of
vertical two jets. The situation of $k$-jets in arbitrary dimension $n$ is 
substantially more involved, details can be found in [Mer09].

Consider the universal hypersurface $\cX\subset\bP^{n+1}\times\bP^{N_d}$ of degree $d$ 
given by the equation
$$
\sum_{|\alpha|=d}A_{\alpha}\,Z^\alpha=0,
$$
where $[Z]\in\bP^{n+1}$, $[A]\in\bP^{N_d}$, $\alpha=(\alpha_0,\dots,\alpha_{n+1})\in\bN^{n+2}$
and
$$
N_d={n+d+1\choose d}-1.
$$
Finally, we denote by $\cV\subset\cX$ the vertical tangent space, i.e.\ the kernel of the
projection
$$
\pi:\cX\to U\subset\bP^{N_d}
$$
where $U$ is the Zariski open set parametrizing smooth hypersurfaces, and by $J_k\cV$ the 
bundle of $k$-jets of curves tangent to $\cV$, i.e.\ curves contained in the fibers
$X_s=\pi^{-1}(s)$. The goal is to describe certain meromorphic vector fields on the total space of $J_k\cV$. By an explicit calculation of vector fields in coordinates, according to Siu's strategy, P\u{a}un [Pau08] was able to prove:

\claim 10.1.~Theorem|The twisted tangent space 
$T_{J_2\cV}\otimes\cO_{\bP^3}(7)\otimes\cO_{\bP^{N_d}}(1)$
is generated over by its global sections over the complement $J_2\cV\ssm\cW$ of the Wronskian locus~$\cW$. Moreover, one can choose generating global sections that are invariant with respect to the action of $\bG_2$ on~$J_2\cV$.
\endclaim

By similar, but more computationally intensive arguments [Mer09], one
can investigate the higher dimensional case. The following result
strengthens the initial announcement of~[Siu04].

\claim 10.2.~Theorem|Let $J_k^{\rm vert}(\cX)$ be the space of
vertical $k$-jets of the universal hypersurface
$$
\cX\subset\bP^{n+1}\times\bP^{N_d}
$$
parametrizing all projective hypersurfaces $X\subset\bP^{n+1}$ of degree $d$.
Then for $k=n$, there exist constants $c_n$ and $c'_n$ such that the twisted tangent bundle
$$
T_{J^{\rm vert}_k(\cX)}\otimes\cO_{\bP^{n+1}}(c_n)\otimes\cO_{\bP^{N_d}}(c'_n)
$$
is generated by its global $\bG_k$-invariant sections outside a certain exceptional 
algebraic subset
$\Sigma\subset J^{\rm vert}_k(\cX)$. One can take either $c_n={1\over 2}(n^2+5n)$, $c'_n=1$ and
$\Sigma$ defined by the vanishing of certain Wronskians, or $c_n=n^2+2n$ and a
smaller set $\wt\Sigma\subset\Sigma$ defined by the vanishing of the $1$-jet part.
\endclaim

\plainsubsection 10.B. General strategy of proof|

Let again $\cX\subset\bP^{n+1}\times\bP^{N_d}$ be the universal hypersurface of degree $d$
in $\bP^{n+1}$.

\noindent
(10.3) {\em Assume that we can prove the existence of a non-zero polynomial differential operator
$$
P\in H^0(\cX,E_{k,m}^\GG T^*_\cX\otimes\cO(-A)),
$$
where $A$ is an ample divisor on $\cX$, at least over some Zariski open set $U$
in the base of the projection $\pi:\cX\to U\subset\bP^{N_d}$.}

Observe that we now have a lot of techniques to do this; the existence of $P$ over 
the family follows from lower semi-continuity in the Zariski topology, once we know 
that such a section $P$ exists on a generic fiber $X_s=\pi^{-1}(s)$. 
Let $\cY\subset\cX$ be the set of points $x\in\cX$ where $P(x)=0$, as an element in
the fiber of the vector bundle $E_{k,m}^\GG T^*_\cX\otimes\cO(-A))$ at $x$. Then
$\cY$ is a proper algebraic subset of $\cX$, and after shrinking $U$ we may assume
that $Y_s=\cY\cap X_s$ is a proper algebraic subset of $X_s$ for every~$s\in U$.

\noindent
(10.4) {\em Assume also, according to Theorems $10.1$ and $10.2$, that we have enough global 
holomorphic $\bG_k$-invariant vector fields $\theta_i$ on 
$J_k\cV$ with values in the pull-back of some ample divisor $B$ on $\cX$,
in such a way that they generate  $T_{J_k\cV}\otimes p^*_kB$ over 
the dense open set $(J_k\cV)^\reg$ of regular $k$-jets, i.e.\ 
$k$-jets with non-zero first derivative $($here $p_k:J_k\cV\to \cX$ 
is the natural projection\/$)$.}

Considering jet differentials $P$ as functions on $J_k\cV$, the idea is 
to produce new ones by taking differentiations
$$
Q_j:=\theta_{j_1}\ldots\theta_{j_\ell}P,\qquad 0\le \ell\le m,~
j=(j_1,\ldots,j_\ell).
$$
Since the $\theta_j$'s are $\bG_k$-invariant, they are in particular
$\bC^*$-invariant; thus
$$
Q_j\in H^0(\cX,E_{k,m}^\GG T^*_\cX\otimes\cO(-A+\ell B))
$$
(and $Q$ is in fact $\bG_k'$ invariant as soon as $P$ is). In order to
be able to apply the vanishing theorems of \S$\,$8, we need $(A-mB)$ to
be ample, so $A$ has to be large compared to $B$. If $f:\bC\to X_s$ is
an entire curve contained in some fiber $X_s\subset\cX$, its lifting
$j_k(f):\bC\to J_k\cV$ has to lie in the zero divisors of all sections
$Q_j$.  However, every non-zero polynomial of degree $m$ has at any
point some non-zero derivative of order $\ell\le m$. Therefore, at any
point where the $\theta_i$ generate the tangent space to $J_k\cV$, we
can find some non-vanishing section~$Q_j$.  By the assumptions on the
$\theta_i$, the base locus of the $Q_j$'s is contained in the union of
$p_k^{-1}(\cY)\cup (J_k\cV)^\sing$; there is of course no way of
getting a non-zero polynomial at points of $\cY$ where $P$
vanishes. Finally, we observe that $j_k(f)(\bC)\not\subset
(J_k\cV)^\sing$ (otherwise $f$ is constant).  Therefore
$j_k(f)(\bC)\subset p_k^{-1}(\cY)$ and thus $f(\bC)\subset\cY$, i.e.
$f(\bC)\subset Y_s=\cY\cap X_s$.

\claim 10.5.~Corollary|Let $\cX\subset\bP^{n+1}\times\bP^{N_d}$ be the
universal hypersurface of degree $d$ in $\bP^{n+1}$. If $d\ge d_n$ is
taken so large that conditions $(10.3)$ and $(10.4)$ are met with
$(A-mB)$ ample, then the generic fiber $X_s$ of the universal family
$\cX\to U$ satisfies the Green-Griffiths conjecture, namely all entire
curves $f:\bC\to X_s$ are contained in a proper algebraic subvariety
$Y_s\subset X_s$, and the $Y_s$ can be taken to form an algebraic
subset $\cY\subset\cX$.
\endclaim

This is unfortunately not enough to get the hyperbolicity of $X_s$,
because we would have to know that $Y_s$ itself is
hyperbolic. However, one can use the following simple observation due
to Diverio and Trapani [DT10]. The starting point is the following
general, straightforward remark. Let $\cE\to \cX$ be a holomorphic
vector bundle let $\sigma\in H^0(\cX,\cE)\ne 0$; then, up to
factorizing by an effective divisor $D$ contained in the common zeroes
of the components of $\sigma$, one can view $\sigma$ as a section
$$
\sigma\in H^0(\cX,\cE\otimes\cO_\cX(-D)),
$$
and this section now has a zero locus without divisorial
components. Here, when $n\ge 2$, a very generic fiber $X_s$ has
Picard number one by the Noether-Lefschetz theorem, and so, after
shrinking $U$ if necessary, we can assume that $\cO_\cX(-D)$ is the
restriction of $\cO_{\bP^{n+1}}(-p)$, $p\ge 0$ by the effectivity of
$D$. Hence $D$ can be assumed to be nef. After performing this
simplification, $(A-mB)$ is replaced by $(A-mB+D)$, which is still ample
if $(A-mB)$ is ample.  As a consequence, we may assume $\codim\cY\ge 2$,
and after shrinking $U$ again, that all $Y_s$ have $\codim Y_s\ge 2$.

\claim 10.6.~Additional statement|In corollary $10.5$, under the same
hypotheses $(10.3)$ and $(10.4)$, one can take all fibers $Y_s$ to
have $\codim Y_s\ge 2$.
\endclaim

This is enough to conclude that $X_s$ is hyperbolic if $n=\dim X_s\le
3$. In fact, this is clear if $n=2$ since the $Y_s$ are then reduced
to points. If $n=3$, the $Y_s$ are at most curves, but we know by Ein
and Voisin that a very generic hypersurface $X_s\subset\bP^4$ of degree
$d\ge 7$ does not possess any rational or elliptic curve.  Hence $Y_s$
is hyperbolic and so is $X_s$, for $s$ generic.\qed

\claim 10.7. Corollary|Assume that $n=2$ or $n=3$, and that
$\cX\subset\bP^{n+1}\times\bP^{N_d}$ is the universal hypersurface of
degree $d\ge d_n\ge 2n+1$ so large that conditions $(10.3)$ and
$(10.4)$ are met with $(A-mB)$ ample. Then the very generic hypersurface
$X_s\subset\bP^{n+1}$ of degree $d$ is hyperbolic.
\endclaim

\plainsubsection 10.C. Proof of the Green-griffiths conjecture for
generic hypersurfaces in $\bP^{n+1}$|

One of the first significant steps towards the Green-Griffiths conjecture
is the result of Diverio, Merker and Rousseau [DMR10],
confirming the statement when $X\subset\bP^{n+1}_\bC$ is a generic 
hypersurface of large degree~$d$. Their proof yields a non-optimal
lower bound $d\ge 2^{n^5}$ for the degree; it is based on an essential
way on Siu's strategy as detailed in~\S$\,$10.B, combined with the
earlier techniques of [Dem95]. Using our improved bounds from~\S$\,$9.D,
we obtain here a better estimate (actually, an estimate
$O(\exp(n^{1+\varepsilon}))$ of exponential order $1$ rather than $5$).
For the algebraic degeneracy of entire curves in open complements
$X=\bP^n\ssm H$, a better bound
$d\geq 5n^2n^n$ has been obtained by Darondeau [Dar14, Dar16b].

\vbox{\claim 10.8.~Theorem|A generic hypersurface $X\subset\bP^{n+1}$ of 
degree $d\ge d_n$ with
$$
d_2=286,\quad d_3=7316,\quad
d_n=\left\lfloor {n^4\over\sqrt{2}}\big(4e^{-1}\,n(1+\log n)\big)^n
\right\rfloor\quad\hbox{for $n\ge 4$},
$$
satisfies the Green-Griffiths conjecture.
\endclaim}

\plainproof. Let us apply Theorem~9.47 with $V=T_X$, $r=n$ and $k=n$. 
The main starting point is the well
known fact that $T^*_{\bP^{n+1}}\otimes\cO_{\bP^{n+1}}(2)$ is 
semi-positive (in fact, generated by its sections). Hence the exact
sequence
$$
0\to \cO_{\bP^{n+1}}(-d)\to T^*_{\bP^{n+1}|X}\to T^*_X\to 0
$$
implies that $T^*_X\otimes\cO_X(2)\ge 0$. We can therefore take 
$\gamma=\Theta_{\cO(2)}=2\omega$ where $\omega$ is the Fubini-Study metric.
Moreover $\det(V^*)=K_X=\cO_X(d-n-2)$ has curvature $(d-n-2)\omega$, and
thus $\Theta_{\det(V^*)}+r\gamma=(d+n-2)\omega$.
The Morse integral to be computed when $A=\cO_X(p)$ is
$$
\int_X \Big(c_{n,n,n}(d+n-2)^n-c'_{n,n,n}(d+n-2)^{n-1}(p+2n)\Big)\omega^n,
$$
so the critical condition we need is
$$d+n-2>{c'_{n,n,n}\over c_{n,n,n}}(p+2n).$$
On the other hand, Siu's differentiation technique requires 
${m\over n^2}(1+{1\over 2}+\ldots+{1\over n})A-mB$ to be ample, where
$B=\cO_X(n^2+2n)$ by Merker's result (Theorem~10.2). This ampleness
condition yields
$$
{1\over n^2}\Big(1+{1\over 2}+\ldots+{1\over n}\Big)p-(n^2+2n)>0,
$$
so one easily sees that it is enough to take $p=n^4-2n$ for $n\ge
3$. Our estimates $(9.52_n)$ give the expected 
bound~$d_n$.\qed
\medskip

Thanks to 10.6, one also obtains the generic hyperbolicity of $2$ and
$3$-dimensional hypersurfaces of large degree.

\claim 10.9.~Theorem|For $n=2$ or $n=3$, a generic hypersurface
$X\subset\bP^{n+1}$ of degree $d\ge d_n$ is Kobayashi hyperbolic.
\endclaim

By using more explicit calculations of Chern classes (and invariant jets
rather than Green-Griffiths jets) Diverio-Trapani [DT10] obtained the
better lower bound $d\ge d_3=593$ in dimension~$3$. In the case of
surfaces, P\u{a}un [Pau08] obtained $d\ge d_2=18$, using deep results of
McQuillan [McQ98].

One may wonder whether it is possible to use jets of order $k<n$ in
the proof of Theorems~10.8 and 10.9. Diverio [Div08] showed that the answer
is negative (his proof is based on elementary facts of representation
theory and a vanishing theorem of Br\"uckmann-Rackwitz [BR90]):

\claim 10.10.~Proposition {\rm([Div08])}|Let 
$X\subset\bP^{n+1}$ be a smooth hypersurface. Then 
$$
H^0(X,E_{k,m}^\GG T^*_X)=0
$$
for $m\ge 1$ and $1\le k<n$. More generally, if $X\subset\bP^{n+s}$ is
a smooth complete intersection of codimension~$s$, there are no global jet differentials
for $m\ge 1$ and $k<n/s$.
\endclaim

\section{Strong general type condition and the GGL conjecture}

\plainsubsection 11.A. A partial result towards the Green-Griffiths-Lang 
conjecture|

The main result of this section is a proof of the partial solution to 
the Green-Griffiths-Lang conjecture asserted in Theorem 0.15.
The following important ``induction step'' can be derived by 
Corollary~9.39. Here $D_k$ denotes  again the sequence of
``vertical divsors''  defined in (6.9).

\claim 11.1. Proposition|Let $(X,V)$ be a directed pair where
$X$ is projective algebraic. Take an irreducible algebraic subset
$Z\not\subset D_k$ of the associated $k$-jet Semple bundle
$X_k$ that projects onto~$X_{k-1}$, $k\ge 1$, and assume that the 
induced directed space $(Z,W)\subset(X_k,V_k)$ is of general type 
modulo $X_\bu\to X$, $\rank W\ge 1$. Then there exists a divisor 
$\Sigma\subset Z_\ell$
in a sufficiently high stage of the Semple tower $(Z_\ell,W_\ell)$
associated with $(Z,W)$, such that every non-constant holomorphic map
$f:\bC\to X$ whose $k$-jet defines a morphism $f_{[k]}:(\bC,T_\bC)\to (Z,W)$
also satisfies $f_{[k+\ell]}(\bC)\subset\Sigma$.
\endclaim

\plainproof. Our hypothesis is that we can find an
embedded resolution of singularities
$$
\mu_{\ell_0}:(\widehat Z_{\ell_0}\subset\widehat X_{k+\ell_0})
\to (Z_{\ell_0}\subset X_{k+\ell_0}),\quad \ell_0\geq 0
$$
and $p\in\bQ_{\geq 0}$ such that 
$$
K_{\widehat W_{\ell_0}}^\bullet\otimes
\cO_{\widehat Z_{\ell_0}}(p)_{\restriction \widehat Z_{\ell_0}}\quad\hbox{is big 
over~$\widehat Z_{\ell_0}$}.
$$
Since Corollary~9.39 and the related lower bound of $h^0$ are universal in the 
category of directed varieties, we can apply them by replacing $(X,V)$
with $(\widehat Z_{\ell_0},\widehat W_{\ell_0})$, $r$ with
$r_0=\rank W$, and $F$ by 
$$
F_{\ell_0}=\cO_{\widehat Z_{\ell_0}}(p)\otimes
\mu_{\ell_0}^*\pi_{k+\ell_0,0}^*\cO_X(-\varepsilon A),
$$
where $A$ is an ample divisor on $X$ and $\varepsilon\in\bQ_{>0}$.
The assumptions show that $K_{\widehat W_{\ell_0}}\otimes F_{\ell_0}$
is still big on $\widehat Z_{\ell_0}$ for $\varepsilon$ small enough,
therefore, by applying our theorem and taking $m\gg\ell\gg \ell_0$, 
we get a large number of (metric bounded) sections of
$$\leqalignno{
\cO_{\widehat Z_\ell}(m)&\otimes{\widehat \pi}_{k+\ell,k+\ell_0}^*
\cO\Big({m\over\ell r_0}\Big(
1+{1\over 2}+\ldots+{1\over\ell}\Big)F_{\ell_0}\Big)\cr
&=\cO_{{\widehat Z}_\ell}(m\abu)\otimes\mu_\ell^*\pi_{k+\ell,0}^*\cO\Big(
-{m\varepsilon\over \ell r_0}\Big(
1+{1\over 2}+\ldots+{1\over \ell}\Big)A\Big)_{\restriction {\widehat Z}_\ell}\cr
&\subset\cO_{{\widehat Z}_\ell}((1+\lambda)m)
\otimes\mu_\ell^*\pi_{k+\ell,0}^*\cO\Big(
-{m\varepsilon\over \ell r_0}\Big(
1+{1\over 2}+\ldots+{1\over \ell}\Big)A\Big)_{\restriction {\widehat Z}_\ell},&(11.2)\cr}
$$
where
$\mu_\ell:(\widehat Z_\ell\subset \widehat X_{k+\ell})\to
(Z_\ell\subset X_{k+\ell})$ is an embedded resolution dominating
$\widehat X_{k+\ell_0}$, and $\abu\in\bQ_+^{\ell'}$ a positive weight 
of the form $(0,\ldots,\lambda,\ldots,0,1)$ with some
non-zero component $\lambda\in\bQ_+$ at index~$\ell_0$.
Let $\widehat\Sigma\subset \widehat Z_\ell$ be the divisor of such
a section. We apply the fundamental vanishing theorem 8.9 to
lifted curves $\widehat f_{[k+\ell]}:\bC\to \widehat Z_\ell$ and
sections of (11.2), and conclude that
$\widehat f_{[k+\ell}(\bC)\subset\widehat\Sigma$.
Therefore $f_{[k+\ell]}(\bC)\subset\Sigma:=\mu_\ell(\widehat\Sigma)$
and Proposition 11.1 is proved.\qed

We now introduce the ad hoc condition that will enable us to check
the GGL conjecture.

\claim 11.3. Definition|Let $(X,V)$ be a directed pair where
$X$ is projective algebraic. We say that $(X,V)$ is ``strongly of
general type'' if it is of general type and for every irreducible 
algebraic set $Z\subsetneq X_k$, $Z\not\subset D_k$, that projects 
onto~$X$, the induced directed structure \hbox{$(Z,W)\subset(X_k,V_k)$}
is of general type modulo $X_\bu\to X$.
\endclaim

\claim 11.4. Example|\rm The situation of a product 
$(X,V)=(X',V')\times(X'',V'')$ described in (0.14) shows that $(X,V)$ can
be of general type without being strongly of general type. In fact,
if $(X',V')$ and $(X'',V'')$ are of general type, then 
$K_V=\pr^{\prime\,*}K_{V'}\otimes\pr^{\prime\prime\,*}K_{V''}$ is big, so $(X,V)$ is 
again of general type. However 
$$
Z=P(\pr^{\prime\,*}V')=X'_1\times X''\subset X_1
$$
has a directed structure $W=\pr^{\prime\,*}V'_1$ which does not possess
a big canonical bundle over $Z$, since the restriction of $K_W$ to
any fiber $\{x'\}\times X''$ is trivial. The higher stages $(Z_k,W_k)$
of the Semple tower of $(Z,W)$ are given by $Z_k=X'_{k+1}\times X''$ and
$W_k=\pr^{\prime\,*}V'_{k+1}$, so it is easy to see that $\GG_k(X,V)$ contains
$Z_{k-1}$. Since $Z_k$ projects onto $X$, we have here $\GG(X,V)=X$
(see [DR15] for more sophisticated indecomposable examples).
\endclaim

\claim 11.5. Hypersurface case|\rm Assume that $Z\ne D_k$ is an irreducible 
hypersurface of $X_k$ that projects onto $X_{k-1}$. To simplify things further,
also assume that $V$ is non-singular. Since the Semple jet-bundles $X_k$ 
form a tower of $\bP^{r-1}$-bundles, their Picard groups satisfy
$\Pic(X_k)\simeq\Pic(X)\oplus\bZ^k$ and we have 
$\cO_{X_k}(Z)\simeq\cO_{X_k}(\abu)\otimes\pi_{k,0}^*B$
for some $\abu\in\bZ^k$ and $B\in\Pic(X)$, where $a_k=d>0$ is the
relative degree of the hypersurface over $X_{k-1}$. Let $\sigma\in H^0(X_k,
\cO_{X_k}(Z))$ be the section defining $Z$ in $X_k$.
The induced directed variety $(Z,W)$ has $\rank W=r-1=\rank(V)-1$ and
Formula (7.25) yields $K_{V_k}=\cO_{X_k}(-(r-1)\onebu)\otimes\pi_{k,0}^*(K_V)$.
We claim that
$$
K_W\supset\big(K_{V_k}\otimes \cO_{X_k}(Z)\big)_{\restriction Z}
\otimes\cJ_S=
\big(\cO_{X_k}(\abu-(r-1)\onebu)\otimes\pi_{k,0}^*(B\otimes K_V)\big)_{\restriction Z}
\otimes\cJ_S\leqno(11.5.1)
$$
where $S\subsetneq Z$ is the set (containing $Z_{\rm sing}$) where $\sigma$
and $d\sigma_{\restriction V_k}$ both vanish, and $\cJ_S$ is the ideal locally 
generated by the coefficients of $d\sigma_{\restriction V_k}$ along $Z=\sigma^{-1}(0)$. In
fact, the intersection $W=T_Z\cap V_k$ is transverse on $Z\smallsetminus S\,$;
then (11.5.1) can be seen by looking at the morphism
$$
V_{k|Z}\build\llra{4ex}^{d\sigma_{\restriction V_k}}_{}\cO_{X_k}(Z)_{\restriction Z},
$$
and observing that the contraction by $K_{V_k}=\Lambda^rV_k^*$ provides a
metric bounded section of the canonical sheaf $K_W$. In order to investigate
the positivity properties of  $K_W$, one has to show that $B$ cannot be too
negative, and in addition to control the singularity set $S$. The second
point is a priori very challenging, but we get useful information for
the first point by observing that $\sigma$ provides a morphism
$\pi_{k,0}^*\cO_X(-B)\to\cO_{X_k}(\abu)$, whence a non-trivial morphism
$$
\cO_X(-B)\to E_{\abu}:=(\pi_{k,0})_*\cO_{X_k}(\abu)
$$
By [Dem95, Section~12], there exists a filtration on $E_{\abu}$ such that 
the graded pieces are irreducible representations of $\GL(V)$ contained
in $(V^*)^{\otimes \ell}$, $\ell\le|\abu|$. Therefore we get a non-trivial 
morphism
$$
\cO_X(-B)\to (V^*)^{\otimes \ell},\qquad \ell\le|\abu|.\leqno(11.5.2)
$$
If we know about certain (semi-)stability properties of $V$, this can be used
to control the negativity of $B$.\qed
\endclaim

\noindent
We further need the following useful concept that slightly
generalizes entire curve loci.

\claim 11.6. Definition|If $Z$ is an algebraic set contained
in some stage $X_k$ of the Semple tower of~$(X,V)$, we define its
``induced entire curve locus'' $\IEL_{X,V}(Z)\subset Z$ to be the Zariski closure 
of the union $\bigcup f_{[k]}(\bC)$ of all jets of entire curves
$f:(\bC,T_\bC)\to(X,V)$ such that $f_{[k]}(\bC)\subset Z$.
\endclaim

We~have of course $\IEL_{X,V}(\IEL_{X,V}(Z))=\IEL_{X,V}(Z)$ by definition.
It is not hard to check that modulo certain ``vertical divisors'' of $X_k$, 
the $\IEL_{X,V}(Z)$ locus is essentially the same as 
the entire curve locus $\ECL(Z,W)$ of the induced directed variety,
but we will not use this fact here. Notice that if $Z=\bigcup Z_\alpha$ is
a decomposition of $Z$ into irreducible components, then
$$\IEL_{X,V}(Z)=\bigcup_\alpha \IEL_{X,V}(Z_\alpha).$$
Since $\IEL_{X,V}(X_k)=\ECL_k(X,V)$,
proving the Green-Griffiths-Lang property amounts to showing that
$\IEL_{X,V}(X)\subsetneq X$ in the stage $k=0$ of the tower. The basic
step of our approach is expressed in the following statement.

\claim 11.7. Proposition|Let $(X,V)$ be a directed variety and 
$p_0\le n=\dim X$, $p_0\ge 1$.
Assume that there is an integer $k_0\ge 0$ such that for every $k\ge k_0$
and every irreducible algebraic set $Z\subsetneq X_k$, $Z\not\subset D_k$, 
such that $\dim \pi_{k,k_0}(Z)\ge p_0$, the induced directed structure
\hbox{$(Z,W)\subset(X_k,V_k)$} is of general type modulo $X_\bu\to X$.
Then $\dim\ECL_{k_0}(X,V)<p_0$.
\endclaim

\plainproof. We argue here by contradiction, assuming that 
$\dim\ECL_{k_0}(X,V)\ge p_0$. If 
$$p'_0:=\dim\ECL_{k_0}(X,V)>p_0$$ 
and if we can prove the result for $p'_0$, we will already get a contradiction.
Hence we can assume without loss of generality that $\dim\ECL_{k_0}(X,V)=p_0$.
The main argument
consists of producing inductively an increasing sequence of integers
$$k_0<k_1<\ldots<k_j<\ldots$$
and directed varieties $(Z^j,W^j)\subset(X_{k_j},V_{k_j})$ satisfying
the following properties~:
\vskip3pt
{\plainitemindent=14mm
\plainitem{(11.7.1)} $Z^0$ is one of the irreducible 
components of $\ECL_{k_0}(X,V)$ and $\dim Z^0=p_0\;$;
\vskip3pt
\plainitem{(11.7.2)} $Z^j$ is one of the irreducible 
components of $\ECL_{k_j}(X,V)$ and $\pi_{k_j,k_0}(Z^j)=Z^0\;$;
\vskip3pt
\plainitem{(11.7.3)} for all $j\ge 0$, $\IEL_{X,V}(Z^j)=Z^j$ and
$\rank W_j\ge 1\;$;
\vskip3pt
\plainitem{(11.7.4)} for all $j\ge 0$, the directed variety 
$(Z^{j+1},W^{j+1})$ is contained 
in some stage (of order $\ell_j=k_{j+1}-k_j$) of the Semple
tower of $(Z^j,W^j)$, namely 
$$
(Z^{j+1},W^{j+1})\subsetneq (Z^j_{\ell_j},W^j_{\ell_j})\subset
(X_{k_{j+1}},V_{k_{j+1}})
$$
and
$$
W^{j+1}=\overline{T_{Z^{j+1\,\prime}}\cap W^j_{\ell_j}}=
\overline{T_{Z^{j+1\,\prime}}\cap V_{k_j}}
$$
is the induced directed structure; moreover $\pi_{k_{j+1},k_j}(Z^{j+1})=Z^j$.
\vskip3pt
\plainitem{(11.7.5)} for all $j\ge 0$, we have $Z^{j+1}\subsetneq 
Z^j_{\ell_j}$ but $\pi_{k_{j+1},k_{j+1}-1}(Z^{j+1})=Z^j_{\ell_j-1}$.\vskip3pt}

\noindent
For $j=0$, we simply take $Z^0$ to be one of the irreducible components 
$S_\alpha$ of $\ECL_{k_0}(X,V)$ such that $\dim S_\alpha=p_0$, which exists by 
our hypothesis that $\dim\ECL_{k_0}(X,V)=p_0$. Clearly, $\ECL_{k_0}(X,V)$
is the union of the $\IEL_{X,V}(S_\alpha)$ and we have 
$\IEL_{X,V}(S_\alpha)=S_\alpha$ for all those components. Thus 
$\IEL_{X,V}(Z^0)=Z^0$ and $\dim Z^0=p_0$. Assume that $(Z^j,W^j)$
has been constructed. The subvariety $Z^j$ cannot
be contained in the vertical divisor $D_{k_j}$. In fact
no irreducible algebraic set $Z$ such that $\IEL_{X,V}(Z)=Z$ can be 
contained in a vertical divisor $D_k$, because $\pi_{k,k-2}(D_k)$ 
corresponds to stationary jets in $X_{k-2}\,$; as every non-constant 
curve $f$ has non-stationary points, its $k$-jet $f_{[k]}$ cannot
be entirely contained in $D_k\,$; also the induced directed structure
$(Z,W)$ must satisfy $\rank W\ge 1$, otherwise $\IEL_{X,V}(Z)\subsetneq Z$.
Condition (11.7.2) implies that
$\dim\pi_{k_j,k_0}(Z^j)\ge p_0$. Therefore $(Z^j,W^j)$ is of general type modulo
$X_\bu\to X$ by the assumptions of the proposition. Thanks to 
Proposition~2.5, we get an algebraic subset $\Sigma\subsetneq Z^j_\ell$
in some stage of the Semple tower $(Z^j_\ell)$ of $Z^j$ such that
every entire curve $f:(\bC,T_\bC)\to(X,V)$ satisfying
$f_{[k_j]}(\bC)\subset Z^j$ also satisfies 
$f_{[k_j+\ell]}(\bC)\subset\Sigma$. By definition, this implies
the first inclusion in the sequence
$$
Z^j=\IEL_{X,V}(Z^j)\subset\pi_{k_j+\ell,k_j}(\IEL_{X,V}(\Sigma))\subset
\pi_{k_j+\ell,k_j}(\Sigma)\subset Z^j
$$
(the other ones being obvious), so we have in fact an equality throughout.
Let $(S'_\alpha)$ be the irreducible
components of $\IEL_{X,V}(\Sigma)$. We have $\IEL_{X,V}(S'_\alpha)=S'_\alpha$ and
one of the components $S'_\alpha$ must satisfy
$$\pi_{k_j+\ell,k_j}(S'_\alpha)=Z^j=Z^j_0.$$
We take $\ell_j\in[1,\ell]$ to be 
the smallest order such that $Z^{j+1}:=\pi_{k_j+\ell,k_j+\ell_j}(S'_\alpha)
\subsetneq Z^j_{\ell_j}$, and set $k_{j+1}=k_j+\ell_j>k_j$. 
By definition of $\ell_j$, we have $\pi_{k_{j+1},k_{j+1}-1}(Z^{j+1})=
Z^j_{\ell_j-1}$, otherwise $\ell_j$ would not be minimal.
We then get $\pi_{k_{j+1},k_j}(Z^{j+1})=Z^j$ and thus
$\pi_{k_{j+1},k_0}(Z^{j+1})=Z^0$ by induction, and all properties 
$(11.7.1-11.7.5)$ follow easily. Now, by Observation~7.29, we have
$$
\rank W^j<\rank W^{j-1}<\ldots<\rank W^1<\rank W^0=\rank V.
$$
This is a contradiction because we cannot have such an infinite sequence.
Proposition~11.7 is proved.\qed

\noindent
The special case $k_0=0$, $p_0=n$ of Proposition~11.7 yields the following
consequence.

\claim 11.8. Partial solution to the generalized GGL conjecture|Let 
$(X,V)$ be a directed pair that is strongly of general type.
Then the Green-Griffiths-Lang conjecture holds true for $(X,V)$, namely
$\ECL(X,V)\subsetneq X\,$; in other words
there exists a proper algebraic variety $Y\subsetneq X$ such that
every nonconstant holomorphic curve $f:\bC\to X$ tangent to $V$
satisfies $f(\bC)\subset Y$.
\endclaim

\claim 11.9. Remark|\rm The proof is not very constructive, but it is 
however theoretically effective. By this we mean that if $(X,V)$ is
strongly of general type and is taken in a bounded family of 
directed varieties, i.e.\ $X$ is embedded
in some projective space $\bP^N$ with some bound $\delta$ on the degree, 
and  $P(V)$ also has bounded degree${}\le \delta'$ when viewed as 
a subvariety of $P(T_{\bP^N})$, then one could theoretically derive bounds 
$d_Y(n,\delta,\delta')$ for the degree of the locus~$Y$. Also, there 
would exist bounds $k_0(n,\delta,\delta')$ for the orders $k$ and 
bounds $d_k(n,\delta,\delta')$ for the degrees of 
subvarieties $Z\subset X_k$ that have to be checked in the 
definition of a pair of strong general type. In fact, [Dem11] produces
more or less explicit bounds for the order $k$ such that Corollary~9.39
holds true. The degree of the divisor $\Sigma$ is given by a section of
a certain twisted line bundle $\cO_{X_k}(m)\otimes\pi^*_{k,0}\cO_X(-A)$
that we know to be big by an application of holomorphic Morse inequalities 
-- and the bounds for the degrees of $(X_k,V_k)$ then provide bounds
for~$m$.
\endclaim

\claim 11.10. Remark|\rm The condition that $(X,V)$ is strongly of 
general type seems to be related to some sort of stability condition.
We are unsure what is the most appropriate definition, but here is one
that makes sense. Fix an ample divisor $A$ on~$X$. For every 
irreducible subvariety $Z\subset X_k$ that projects onto $X_{k-1}$ for
$k\ge 1$, and $Z=X=X_0$ for $k=0$, we define the slope $\mu_A(Z,W)$
of the corresponding directed variety $(Z,W)$ to be
$$
\mu_A(Z,W)={\inf\lambda\over\rank W},
$$
where $\lambda$ runs over all rational numbers such that there exists
$\ell\geq 0$, a modification $\widehat Z_\ell\to Z_\ell$ and $p\in\bQ_+$
for which 
$$
K_{\widehat W_\ell}\otimes\big(\cO_{\widehat Z_\ell}(p)
\otimes\pi_{k+\ell,0}^*\cO(\lambda A)
\big)_{\restriction \widehat Z_\ell}\quad\hbox{is big on $\widehat Z_\ell$}
$$
(again, we assume here that $Z\not\subset D_k$ for $k\ge 2$). Notice that
by definition $(Z,W)$ is of general type modulo $X_\bu\to X$ if and
only if $\mu_A(Z,W)<0$, and that $\mu_A(Z,W)=-\infty$ 
if $\cO_{\widehat Z_\ell}(1)$ is big for some~$\ell$. Also, the proof
of Lemma~7.24 shows that for any $(Z,W)$ we have
$\mu_A(Z_\ell,W_\ell)=\mu_A(Z,W)$ for all~$\ell\ge 0$.
We say that $(X,V)$ is {\it $A$-jet-stable} (resp.\ 
{\it $A$-jet-semi-stable})
if $\mu_A(Z,W)<\mu_A(X,V)$ (resp.\ $\mu_A(Z,W)\le\mu_A(X,V)$) for
all $Z\subsetneq X_k$ as above. It is then clear that if
$(X,V)$ is of general type and $A$-jet-semi-stable, then it is strongly
of general type in the sense of Definition~11.3. It would be useful 
to have a better understanding of this condition of stability 
(or any other one that would have better properties).
\endclaim

\plainsubsection 11.B. Algebraic jet-hyperbolicity implies Kobayashi hyperbolicity|

Let $(X,V)$ be a directed variety, where $X$ is an irreducible projective 
variety; the concept still makes sense when $X$ is singular, by embedding
$(X,V)$ in a projective space $(\bP^N,T_{\bP^N})$ and taking the linear space
$V$ to be an irreducible algebraic subset of $T_{\bP^n}$ that is contained 
in $T_X$ at regular points of~$X$.

\claim 11.11. Definition|Let $(X,V)$ be a directed variety. We say that 
$(X,V)$ is algebraically jet-hyperbolic if
for every $k\ge 0$ and every irreducible algebraic subvariety 
$Z\subset X_k$ that is not contained in the 
union $\Delta_k$ of vertical divisors, the induced directed 
structure $(Z,W)$ either satisfies 
\hbox{$W=0$}, or is of general type modulo $X_\bu\to X$,
i.e.\ there exists $\ell\geq 0$ and $p\in\bQ_{\geq 0}$ such that 
$K_{\widehat W_\ell}^\bullet\otimes\cO_{\widehat Z_\ell}(p)$
is big over~$\widehat Z_\ell$, for some modification
$(\widehat Z_\ell,\widehat W_\ell)$ of the $\ell$-stage of
the Semple tower of $(Z,W)$.
\endclaim

\noindent
Proposition 7.33 can be restated:

\claim 11.12. Proposition|If a projective directed variety $(X,V)$ is
such that $\cO_{X_\ell}(\abu)$ is ample for some $\ell\geq 1$ and
some weight $\abu\in\bQ_{>0}^\ell$, then $(X,V)$ is algebraically 
jet-hyperbolic.
\endclaim

\noindent
In a similar vein, one would prove that if $\cO_{X_\ell}(\abu)$ is big
and the ``augmented base locus''\break
\hbox{$B=\Bs(\cO_{X_\ell}(\abu)\otimes\pi_{l,0}^*A^{-1})$} projects onto 
a proper subvariety $B'=\pi_{\ell,0}(B)\subsetneq X$, then 
$(X,V)$ is strongly of general type. In general, Proposition~11.7 gives
the following:

\claim 11.13. Theorem|Let $(X,V)$ be an irreducible projective directed 
variety that is algebraically jet-hyperbolic in the sense of the above 
definition. Then $(X,V)$ is Brody $($or Kobayashi$\,)$ hyperbolic, i.e.\ 
$\ECL(X,V)=\emptyset$.
\endclaim

\plainproof. Here we apply Proposition 11.7 with $k_0=0$ and $p_0=1$. It is
enough to deal with subvarieties $Z\subset X_k$ such that 
$\dim\pi_{k,0}(Z)\ge 1\,$; otherwise $W=0$ and can reduce $Z$ to a smaller
subvariety by (2.2). Then we conclude that
$\dim\ECL(X,V)<1$. All entire curves tangent to $V$ have to be constant, 
and we conclude in fact that $\ECL(X,V)=\emptyset$.\qed

\section{Proof of the Kobayashi conjecture on generic
hyperbolicity}

We give here a simple proof of the Kobayashi conjecture, combining ideas
of Green-Griffiths [GrGr80], Nadel [Nad89], Masuda-Noguchi [MaNo96],
Demailly [Dem95], Siu-Yeung [SiYe96a], Shiffman-Zaidenberg [ShZa02],
Brotbek [Brot17], Ya Deng [Deng16],
in chronological order. Related ideas had been used earlier in [Xie15],
and then in [BrDa17], to establish Debarre's conjecture on the ampleness of
the cotangent bundle of generic complete intersections of codimension
at least equal to dimension.

\plainsubsection 12.A. General Wronskian operators|

This section follows closely the work of D.~Brotbek [Brot17].
Let $U$ be an open set of a complex manifold $X$, $\dim X=n$, and
$s_0,\ldots,s_k\in\cO_X(U)$ be holomorphic functions. To these
functions, we can associate a Wronskian operator of order $k$ defined
by
$$
W_k(s_0,\ldots,s_k)(f)=\left|
\plainmatrix{
  s_0(f) & s_1(f) & \ldots &s_k(f)\cr
  D(s_0(f)) & D(s_1(f)) & \ldots &D(s_k(f))\cr
  \noalign{\vskip4pt}  
  \vdots & \vdots &  &\vdots\cr
  \noalign{\vskip4pt}  
  D^k(s_0(f)) & D^k(s_1(f)) & \ldots &D^k(s_k(f))\cr}\right|\leqno(12.1)
$$
where $f:(\bC,0)\ni t\mapsto f(t)\in U\subset X$ is a germ of holomorphic curve
(or a $k$-jet of curve), and $D=\frac{d}{dt}$. For a biholomorphic change of
variable $\varphi:(\bC,0)\to(\bC,0)$, we find by induction on $\ell$ a
polynomial differential operator $p_{\ell,i}$ of order${}\leq\ell$ acting on
$\varphi$ satisfying
$$
D^\ell(s_j(f\circ \varphi))=\varphi^{\prime\ell}D^\ell(s_j(f))\circ\varphi
+\sum_{i<\ell}p_{\ell,i}(\varphi)D^i(s_j(f))\circ\varphi.
$$
It follows easily from this that
$$
W_k(s_0,\ldots,s_k)(f\circ\varphi)=(\varphi')^{1+2+\cdots+k}
W_k(s_0,\ldots,s_k)(f)\circ\varphi,
$$ 
and hence $W_k(s_0,\ldots,s_k)(f)$ is an invariant
differential operator of degree $k'=\frac{1}{2}k(k+1)$. Especially, we get
in this way a section that we denote somewhat sloppily
$$
W_k(s_0,\ldots,s_k)=
\left|
\plainmatrix{
  s_0 & s_1 & \ldots &s_k\cr
  D(s_0) & D(s_1) & \ldots &D(s_k)\cr
  \noalign{\vskip4pt}
  \vdots & \vdots &  &\vdots\cr
  \noalign{\vskip4pt}  
  D^k(s_0) & D^k(s_1) & \ldots &D^k(s_k)\cr}\right|\in
H^0(U,E_{k,k'}T^*_X).\leqno(12.2)
$$

\claim 12.3. Proposition|These Wronskian operators satisfy the following
properties.
\vskip2pt
\plainitem{\rm(a)} $W_k(s_0,\ldots,s_k)$ is $\bC$-multilinear and
alternate in $(s_0,\ldots,s_k)$.
\vskip2pt\plainitem{\rm(b)} For any $g\in\cO_X(U)$, we have
$$W_k(gs_0,\ldots,gs_k)=g^{k+1}W_k(s_0,\ldots,s_k).$$
\endclaim

\noindent
Property 12.3~(b) is an easy consequence of the Leibniz formula
$$
D^\ell(g(f)s_j(f))=\sum_{k=0}^{\ell}{\ell\choose k}D^k(g(f))
D^{\ell-k}(s_j(f)),
$$
by performing linear combinations of rows in the determinants. This property
implies in its turn that for any $(k+1)$-tuple of sections
$s_0,\ldots,s_k\in H^0(U,L)$ of a holomorphic line bundle $L\to X$,
one can define more generally an operator
$$
W_k(s_0,\ldots,s_k)\in H^0(U,E_{k,k'}T^*_X\otimes L^{k+1}).
\leqno(12.4)
$$
In fact, when we compute the Wronskian
in a local trivialization of $L_{\restriction U}$, Property 12.3~(b)
shows that the determinant is independent of the trivialization.
Moreover, if $g\in H^0(U,G)$ for some line bundle $G\to X$, we have
$$
W_k(gs_0,\ldots,gs_k)=g^{k+1}W_k(s_0,\ldots,s_k)\in H^0(U,E_{k,k'}T^*_X
\otimes L^{k+1}\otimes G^{k+1}).\leqno(12.5)
$$
We consider here a line bundle $L\to X$ possessing a linear system
$\Sigma\subset H^0(X,L)$ of global sections such that
$W_k(s_0,\ldots,s_k)\not\equiv 0$
for generic elements $s_0,\ldots,s_k\in\Sigma$. We can then view
$W_k(s_0,\ldots,s_k)$ as
a section of $H^0(X_k,\cO_{X_k}(k')\otimes\pi_{k,0}^*L^{k+1})$
on the $k$-stage $X_k$ of the Semple tower. Very roughly, the idea for
the proof of the Kobayashi conjecture is to produce many such Wronskians,
and to apply the fundamental vanishing theorem 8.15 to exclude the
existence of entire curves. However, the vanishing theorem only holds
for jet differentials in $H^0(X_k,\cO_{X_k}(k')\otimes\pi_{k,0}^*A^{-1})$
with $A>0$, while the existence of sufficiently many sections
$s_j\in H^0(X,L)$ can be achieved only when $L$ is ample, so the
strategy seems a priori unapplicable.
It turns out that one can sometimes arrange the Wronkian operator coefficients
to be divisible by a section $\sigma_\Delta\in H^0(X,\cO_X(\Delta))$
possessing a large zero divisor~$\Delta$, so that
$$
\sigma_\Delta^{-1}W_k(s_0,\ldots,s_k)\in
H^0\big(X_k,\cO_{X_k}(k')\otimes\pi_{k,0}^*(L^{k+1}\otimes\cO_X(-\Delta))\big),
$$
and we can then hope that $L^{k+1}\otimes\cO_X(-\Delta))<0$.
Our goal is thus to find a variety $X$ and linear systems 
$\Sigma\subset H^0(X,L)$ for which the associated Wronskians
$W_k(s_0,\ldots,s_k)$ have a very high divisibility. The study of the
base locus of line bundles
$\cO_{X_k}(k')\otimes\pi_{k,0}^*(L^{k+1}\otimes\cO_X(-\Delta))$
and their related positivity properties will be taken care of by using
suitable blow-ups.

\plainsubsection 12.B. Using a blow-up of the Wronskian ideal sheaf|

We consider again a linear system $\Sigma\subset H^0(X,L)$
producing some non-zero Wronskian sections $W_k(s_0,\ldots,s_k)$, so that
$\dim\Sigma\geq k+1$. As the Wronskian is alternate and multilinear in
the arguments $s_j$, we get a meromorphic map
$X_k\merto P(\Lambda^{k+1}\Sigma^*)$ by sending a $k$-jet
$\gamma=f_{[k]}(0)\in X_k$ to the point of projective coordinates
$[W_k(u_{i_0},\ldots,u_{i_k})(f)(0)]_{i_0,\ldots,i_k}$
where $(u_j)_{j\in J}$ is a basis of $\Sigma$ and
$i_0,\ldots,i_k\in J$ are in increasing order. This
assignment factorizes through the Pl\"ucker embedding into a
meromorphic map
$$
\Phi:X_k\merto \Gr_{k+1}(\Sigma)
$$
into the Grassmannian
of dimension $(k+1)$ subspaces of $\Sigma^*$ (or codimension $(k+1)$
subspaces of $\Sigma$, alternatively). In fact, if $L_{\restriction U}
\simeq U\times\bC$ is a trivialization of $L$ in a neighborhood of
a point $x_0=f(0)\in X$, we can consider the map
$\Psi_U:X_k\to \Hom(\Sigma,\bC^{k+1})$
given by
$$\pi_{k,0}^{-1}(U)\ni
f_{[k]}\mapsto\big(s\mapsto(D^{\ell}(s(f))_{0\leq\ell\leq k})\big),$$
and associate either the kernel $\Xi\subset\Sigma$ of
$\Psi_U(f_{[k]})$, seen as a point $\Xi\in\Gr_{k+1}(\Sigma)$, or 
$\Lambda^{k+1}\Xi^\perp\subset\Lambda^{k+1}\Sigma^*$, seen as a point
of $P(\Lambda^{k+1}\Sigma^*)$ 
(assuming that we are at a point where the rank
is equal to~$(k+1)$). Let $\cO_{\Gr}(1)$ be the tautological very ample line
bundle on $\Gr_{k+1}(\Sigma)$ (equal to the restriction of
$\cO_{P(\Lambda^{k+1}\Sigma^*)}(1))$. By construction, $\Phi$ is
induced by the linear system of sections 
$$
W_k(u_{i_0},\ldots,u_{i_k})\in
H^0(X_k,\cO_{X_k}(k')\otimes\pi_{k,0}^*L^{k+1}),
$$
and we thus get a natural isomorphism
$$
\cO_{X_k}(k')\otimes\pi_{k,0}^*L^{k+1}\simeq\Phi^*\cO_{\Gr}(1)\quad
\hbox{on $X_k\ssm B_k$},\leqno(12.6)
$$
where $B_k\subset X_k$ is the base locus of our linear system of Wronskians.
The presence of the indeterminacy set $B_k$ may create trouble in analyzing the positivity of our line bundles, so we are going to use an appropriate
blow-up to resolve the indeterminacies. For this purpose, we introduce the ideal
sheaf $\cJ_{k,\Sigma}\subset\cO_{X_k}$ generated by the linear system~$\Sigma$,
and take 
a modification $\mu_{k,\Sigma}:\widehat X_{k,\Sigma}\to X_k$ in such a way
that $\mu_{k,\Sigma}^*\cJ_{k,\Sigma}=\cO_{\widehat X_{k,\Sigma}}(-F_{k,\Sigma})$
for some divisor $F_{k,\Sigma}$ in
$\widehat X_{k,\Sigma}$. Then $\Phi$ is resolved into a morphism
\hbox{$\Phi\circ\mu_{k,\Sigma}:\widehat X_{k,\Sigma}\to\Gr_{k+1}(\Sigma)$},
and on~$\widehat X_{k,\Sigma}$, (12.6) becomes an everywhere defined
isomorphism
$$
\mu_{k,\Sigma}^*\big(\cO_{X_k}(k')\otimes\pi_{k,0}^*L^{k+1})\otimes
\cO_{\widehat X_{k,\Sigma}}(-F_{k,\Sigma})
\simeq(\Phi\circ\mu_{k,\Sigma})^*\cO_{\Gr}(1).\leqno(12.7)
$$
In fact, we can simply take $\widehat X_k$ to be the
normalized blow-up of $\cJ_{k,\Sigma}$, i.e.\ the normalization of the
closure $\Gamma\subset X_k\times\Gr_{k+1}(\Sigma)$ of the graph of $\Phi$
and $\mu_{k,\Sigma}:\widehat X_k\to X_k$ to be the
composition of the normalization map $\widehat X_k\to\Gamma$
with the first projection~$\Gamma\to X_k$. $\big[$The Hironaka
desingularization theorem would possibly allow us to replace $\widehat X_k$
by a nonsingular modification, and $F_{k,\Sigma}$ by a simple normal
crossing divisor
on the desingularization; we will avoid doing so here, as we would
otherwise need to show the existence of universal desingularizations when
$(X_t,\Sigma_t)$ is a family of linear systems of $k$-jets of sections
associated with a family of algebraic varieties$\big]$.  The
following basic lemma was observed by Ya Deng [Deng16].

\claim 12.8. Lemma|Locally over coordinate open sets $U\subset X$ on which
$L_{\restriction U}$ is trivial, there is a
maximal ``Wronskian ideal sheaf'' $\cJ^X_k\supset\cJ_{k,\Sigma}$ in $\cO_{X_k}$
achieved by linear systems
$\Sigma\subset H^0(U,L)$. It is attained globally on $X$ whenever the
linear system $\Sigma\subset H^0(X,L)$ generates $k$-jets of sections
of $L$ at every point. Finally, it is ``universal'' in the sense that is does
not depend on $L$ and behaves functorially under immersions: if $\psi:X\to Y$
is an immersion and $\cJ^X_k$, $\cJ^Y_k$ are the corresponding
Wronskian ideal sheaves in $\cO_{X_k}$, $\cO_{Y_k}$, then
$\smash{\psi_k^*\cJ^Y_k=\cJ^X_k}$ with respect to the induced
immersion $\psi_k:X_k\to Y_k$.
\endclaim

\proof The (local) existence of such a maximal ideal sheaf is merely a
consequence of the strong Noetherian property of coherent ideals.
As observed at the end of section 6.A, the bundle
$X_k\to X$ is a locally trivial tower of $\bP^{n-1}$-bundles, with a
fiber $\cR_{n,k}$ that is a rational \hbox{$k(n-1)$}-dimensional variety; over
any coordinate open set $U\subset X$ equipped with local coordinates
\hbox{$(z_1,\ldots,z_n)\,{\in}\,B(0,r)\,{\subset}\,\bC^n$}, it is
isomorphic to the
product $U\times\cR_{n,k}$, the fiber over a point $x_0\in U$ being identified
with the central fiber through a translation $(t\mapsto f(t))\mapsto
(t\mapsto x_0+f(t))$ of germs of curves. In this setting, $\cJ^X_k$ is
generated by the functions in $\cO_{X_k}$ associated with Wronskians
$$
X_{k\,\restriction U}\ni \xi=f_{[k]}\mapsto W_k(s_0,\ldots,s_k)(f)\in\cO_{X_k}(k')_{\restriction\cR_{n,k}},\quad s_j\in H^0(U,\cO_X),
$$
by taking local trivializations $\cO_{X_k}(k')_{\xi_0}\simeq\cO_{X_k,\xi_0}$
at points~$\xi_0\in X_k$. In fact, it is enough
to take Wronskians associated with {\it polynomials}
$s_j\in\bC[z_1,\ldots,z_n]$. To see this, one can e.g.\ invoke Krull's lemma
for local rings, which implies $\cJ^X_{k,\xi_0}=\bigcap_{\ell\geq 0}
(\cJ^X_{k,\xi_0}+\gm_{\xi_0}^{\ell+1})$, and to observe that $\ell$-jets of
Wronskians $W_k(s_0,\ldots,s_k)$ (mod $\gm_{\xi_0}^{\ell+1}$) depend only on
the $(k+\ell)$-jets of the sections $s_j$ in $\cO_{X,x_0}/\gm_{x_0}^{k+\ell+1}$,
where $x_0=\pi_{k,0}(\xi_0)$. Therefore, polynomial sections $s_j$ or
arbitrary holomorphic functions $s_j$ define the same $\ell$-jets of
Wronskians for any~$\ell$. Now, in the case of polynomials, it is clear that
translations $(t\mapsto f(t))\mapsto(t\mapsto x_0+f(t))$ leave
$\cJ_k^X$ invariant, hence $\cJ_k^X$ is the pull-back by the
second projection $X_{k\,\restriction U}\simeq U\times\cR_{n,k}\to \cR_{n,k}$
of its restriction to any of the fibers $\pi_{k,0}^{-1}(x_0)\simeq\cR_{n,k}$.
As the $k$-jets of the $s_j$'s at $x_0$ are sufficient to determine the
restriction
of our Wronskians to $\pi_{k,0}^{-1}(x_0)$, the first two claims of Lemma 5.3
follow. The universality property comes from the fact
that $L_{\restriction U}$ is trivial (cf.~(12.3~b)) and that germs of sections
of $\cO_X$ extend to germs of sections of $\cO_Y$ via the immersion~$\psi$.
(Notice that in this discussion, one may have to pick Taylor expansions of
order${}>k$ for $f$ to reach all points of the fiber $\pi_{k,0}^{-1}(x_0)$,
the order $2k-1$ being sufficient by [Dem95, Prop.~5.11], but this fact
does not play any role here). A consequence of universality is that $\cJ^X_k$
does not depend on coordinates nor on the geometry of~$X$.\qed

\noindent
The above discussion combined with Lemma 12.8 leads to the following statement.

\claim 12.9. Proposition|Assume that $L$ generates all $k$-jets of sections
$($e.g.\, take $L=A^p$ with $A$ very ample and $p\geq k)$, and
let $\Sigma\subset H^0(X,L)$ be a linear system that also generates
$k$-jets of sections at any point of $X$. Then we have a universal
isomorphism
$$
\mu_k^*\big(\cO_{X_k}(k')\otimes\pi_{k,0}^*L^{k+1})\otimes
\cO_{\widehat X_{k,\Sigma}}(-F_k)
\simeq(\Phi\circ\mu_k)^*\cO_{\Gr_{k+1}(\Sigma)}(1),
$$
where $\mu_k:\widehat X_k\to X_k$ is the normalized blow-up of
the $($maximal$\,)$ ideal sheaf $\cJ^X_k\subset\cO_{X_k}$ associated with 
order~$k$ Wronskians, and $F_k$ the universal divisor of $\widehat X_k$
resolving $\cJ^X_k$.
\endclaim

\plainsubsection 12.C. Specialization to suitable hypersurfaces|

Let $Z$ be a non-singular $(n+1)$-dimensional projective variety, and
let $A$ be a very ample divisor on~$Z\,$; the fundamental example is of
course $Z=\bP^{n+1}$ and $A=\cO_{\bP^{n+1}}(1)$. Our goal is to
show that a sufficiently general ($n$-dimensional) hypersurface
$X=\{x\in Z\,;\;\sigma(x)=0\}$ defined by a section $\sigma\in H^0(Z,A^d)$,
$d\gg 1$, is Kobayashi hyperbolic. A basic idea,
inspired by some of the main past contributions, such as
Brody-Green [BrGr77], Nadel [Nad89], Masuda-Noguchi [MaNo96],
Shiffman-Zaidenberg [ShZa02] and [Xie15],
is to consider hypersurfaces defined by special equations,
e.g.\ deformations of unions
of hyperplane sections $\tau_1\cdots\tau_d=0$ or of Fermat-Waring
hypersurfaces $\sum_{0\leq j\leq N}\tau_j^d=0$, for suitable sections
$\tau_j\in H^0(Z,A)$. Brotbek's main idea developed in [Brot17] is that
a carefully selected hypersurface may have enough Wronskian sections 
to imply the ampleness of some tautological jet line bundle --
a Zariski open property. Here, we take $\sigma \in H^0(X,A^d)$ equal to
a sum of terms
$$
\sigma=\sum_{0\leq j\leq N}a_jm_j^\delta,\quad
a_j\in H^0(Z,A^\rho),~~m_j\in H^0(Z,A^b),~~n<N\leq k,~~
d=\delta b+\rho,
\leqno(12.10)
$$
where $\delta \gg 1 $ and the $m_j$ are ``monomials'' of the same degree $b$,
i.e.\ products of $b$ ``linear'' sections $\tau_I\in H^0(Z,A)$, and the
factors $a_j$ are general enough. The integer $\rho$ is taken in the
range $[k,k+b-1]$, first to ensure that $H^0(Z,A^\rho)$ generates
$k$-jets of sections, and second, to allow $d$ to be an arbitrary
large integer (once $\delta\geq\delta_0$ has been chosen large enough).

The monomials $m_j$ will be chosen in such a way
that for suitable  \hbox{$c\in\bN$}, $1\leq c \leq N$,
any subfamily of $c$ terms $m_j$ shares one common
factor~$\tau_I\in H^0(X,A)$. To this end, we consider all subsets
$I\subset\{0,1,\ldots,N\}$ with
$\card I=c\,$; there are $B:={N+1\choose c}$
subsets of this type. For all such~$I$, we select sections
$\tau_I\in H^0(Z,A)$ such that $\prod_I \tau_I=0$ is
a simple normal crossing divisor in~$Z$ (with all of its components of
multiplicity~$1$). For $j=0,1,\ldots,N$ given, the number of
subsets $I$ containing $j$ is $b:={N\choose c-1}$. We put
$$
m_j=\prod_{I\ni j}\tau_I\in H^0(Z,A^b).\leqno(12.11)
$$
By construction, every family $m_{i_1},\ldots,m_{i_c}$ of sections shares
the common factor $\tau_I\in H^0(X,A)$ where $I=\{i_1,\ldots,i_c\}$.
The first step consists in checking that we can achieve $X$ to be smooth with
these constraints.

\claim 12.12. Lemma|Assume $N\geq c(n+1)$. Then,
for a generic choice of the sections $a_j\in H^0(Z,A^\rho)$ and
$\tau_I\in H^0(Z,A)$, the hypersurface $X=\sigma^{-1}(0)\subset Z$
defined by $(12.10\hbox{-}12.11)$ is non-singular. Moreover, under the same
condition for $N$, the intersection of $\prod \tau_I=0$ with~$X$
can be taken to be a simple normal crossing divisor in $X$.
\endclaim

\plainproof. As the properties considered in the Lemma are Zariski open
properties in terms of the $(N+B+1)$-tuple $(a_j,\tau_I)$, 
it is sufficient to prove the result for a specific choice
of the $a_j$'s: we fix here $a_j=\tilde\tau_j\tau_{I(j)}^{\rho-1}$ where
$\tilde\tau_j\in H^0(X,A)$, $0\leq j\leq N$ are new sections such that
$\prod\tilde\tau_j\prod\tau_I=0$ is a simple normal crossing divisor,
and $I(j)$ is any subset of cardinal $c$ containing~$j$. Let $H$ be the
hypersurface of degree $d$ of $\bP^{N+B}$
defined in homogeneous coordinates $(z_j,z_I)\in\bC^{N+B+1}$ by $h(z)=0$ where
$$
h(z)=\sum_{0\leq j\leq N}z_jz_{I(j)}^{\rho-1}\prod_{I\ni j}z_I^\delta,
$$
and consider the morphism $\Phi:Z\to \bP^{N+B}$ such that
$\Phi(x)=(\tilde\tau_j(x),\tau_I(x))$. With our choice of the $a_j$'s,
we have $\sigma=h\circ\Phi$. Now,
when the $\tilde\tau_j$ and $\tau_I$ are general enough, the map $\Phi$
defines an embedding of $Z$ into $\bP^{N+B}$ (for this, one needs
$N+B\geq 2(\dim Z)+1=2n+3$, which is the case by our assumptions).
Then, by definition, $X$ is isomorphic to the intersection
of $H$ with $\Phi(Z)$. Changing generically the $\tilde\tau_j$ and
$\tau_I$'s can be achieved by composing $\Phi$ with a generic automorphism
$g\in\Aut(\bP^{N+B})=\PGL_{N+B+1}(\bC)$ (as $\GL_{N+B+1}(\bC)$ acts
transitively on $(N+B+1)$-tuples of linearly independent linear forms).
As $\dim g\circ \Phi(Z)=\dim Z=n+1$,
Lemma~12.12 will follow from a standard Bertini argument if we
can check that $\Sing(H)$ has codimension at least $(n+2)$ in~$\bP^{N+B}$.
In fact, this condition implies $\Sing(H)\cap (g\circ\Phi(Z))=\emptyset$
for $g$ generic, while $g\circ\Phi(Z)$ can be chosen transverse to $\Reg(H)$.
Now, a sufficient condition for smoothness is that one of the
differentials $dz_j$, $0\leq j\leq N$, appears with a non-zero factor in
$dh(z)$ (just neglect the other differentials $*dz_I$ in this argument).
We infer from this and the fact that $\delta\geq 2$ that $\Sing(H)$ consists of
the locus defined by $\prod _{I\ni j}z_I=0$ for all $j=0,1,\ldots,N$.
It~is the union of the linear subspaces $z_{I_0}=\ldots=z_{I_N}=0$ for
all possible  choices of subsets $I_j$ such that $I_j\ni j$.
Since $\card I_j=c$, the equality $\bigcup I_j=\{0,1,\ldots,N\}$ implies
that there are at least $\lceil (N+1)/c\rceil$ distinct subsets $I_j$
involved in each of these linear subspaces, and the equality can be reached.
Therefore $\codim \Sing(H)=\lceil (N+1)/c\rceil\geq n+2$ as soon as
$N \geq c(n+1)$. By the same argument, we can assume that the intersection
of $Z$ with at least $(n+2)$ distinct hyperplanes $z_I=0$ is empty.
In order that $\prod\tau_I=0$ defines a normal crossing
divisor at a point $x\in X$,
it is sufficient to ensure that for any family $\cG$ of coordinate
hyperplanes $z_I=0$, $I\in \cG$, with $\card \cG\leq n+1$, we have
a ``free'' index $j\notin
\bigcup_{I\in \cG}I$ such that $x_I\neq 0$ for all $I\ni j$, so that
$dh$ involves a non-zero term $*\,dz_j$ independent of the $dz_I$, $I\in \cG$.
If this fails, there must be at least $(n+2)$ hyperplanes $z_I=0$
containing $x$, associated either with $I\in \cG$, or with other $I$'s covering
$\complement\big(\bigcup_{I\in \cG}I\big)$. The corresponding bad locus is
of codimension at least $(n+2)$ in~$\bP^{N+B}$ and can be avoided by
$g(\Phi(Z))$ for a generic choice of $g\in\Aut(\bP^{N+B})$. Then
$X\cap\bigcap_{I\in \cG}\tau_I^{-1}(0)$
is smooth of codimension equal to~$\card \cG$.\qed

\plainsubsection 12.D. Construction of highly divisible Wronskians|

To any families $s,\,\hat\tau$ of sections
$s_1,\ldots,s_r\in H^0(Z,A^k)$, $\hat\tau_1,\ldots,\hat\tau_r\in H^0(Z,A)$,
and to each subset $J\subset\{0,1,\ldots,N\}$ with $\card J=c$, we associate a
Wronskian operator of order $k$ (i.e.\ a $(k+1)\times(k+1)$-determinant)
$$
W_{k,s,\hat\tau,a,J}=W_k\big(s_1\hat\tau_1^{d-k},\ldots,s_r\hat\tau_r^{d-k},
(a_jm_j^\delta)_{j\in\complement J}\big),\quad r=k+c-N,~~|\complement J|=N-c.
\leqno(12.13)
$$
We assume here again that the $\hat\tau_j$ are chosen so that
$\prod\hat\tau_j\prod\tau_I=0$ defines a simple normal crossing divisor
in $Z$ and $X$. Since $s_j\hat\tau_j^{d-k},\,a_jm_j^\delta\in H^0(Z,A^d)$,
formula (12.4) applied with $L=A^d$ implies that
$$
W_{k,s,\hat\tau,a,J}\in H^0(Z,E_{k,k'}T^*_Z\otimes A^{(k+1)d}).
\leqno(12.14)
$$
However, we are going to see that $W_{k,s,\hat\tau,a,J}$ and its
restriction $W_{k,s,\hat\tau,a,J\restriction X}$ are divisible by
monomials $\hat\tau^\alpha \tau^\beta$ of very large degree,
where $\hat\tau$, resp.\ $\tau$, denotes the collection of sections
$\hat\tau_j$, resp.\ $\tau_I$ in $H^0(Z,A)$. In this way, we will see that
we can even obtain a negative exponent of $A$ after simplifying
$\hat\tau^\alpha\tau^\beta$ in $W_{k,s,\hat\tau,a,J\restriction X}$.
This simplification process is a generalization of techniques already
considered by [Siu87] and [Nad89] (and later [DeEG97]) in relation
with the use of meromorphic connections of low pole order.

\claim 12.15. Lemma|Assume that $\delta\geq k$.
Then the Wronskian operator $W_{k,s,\hat\tau,a,J}$, resp.\
$W_{k,s,\hat\tau,a,J\restriction X}$, is divisible by a monomial
$\hat\tau^\alpha\tau^\beta$,
resp.\ $ \hat\tau^\alpha\tau^\beta\tau_J^{\delta-k}$
$($with a multiindex notation
$\hat\tau^\alpha\tau^\beta=\prod\hat\tau_j^{\alpha_j}
\prod\tau_I^{\beta_I})$, and
$$
\alpha,\beta\geq 0,\quad
|\alpha|=r(d-2k),\quad
|\beta|=(N+1-c)(\delta-k)b.
$$
\endclaim

\plainproof. $W_{k,s,\hat\tau,a,J}$ is obtained as a determinant whose
$r$ first columns are the derivatives $D^\ell(s_j\hat\tau_j^{d-k})$ and the
last $(N+1-c)$ columns are the $D^\ell(a_jm_j^\delta)$, divisible respectively by
$\hat\tau_j^{d-2k}$ and $m_j^{\delta-k}$. As $m_j$ is of the form $\tau^\gamma$,
$|\gamma|=b$, this implies the divisibility of $W_{k,s,\hat\tau,a,J}$ by a
monomial of the form $\hat\tau^\alpha\tau^\beta$, as asserted.
Now, we explain why one can gain the additional factor
$\tau_J^{\delta-k}$ dividing the restriction
$W_{k,s,\hat\tau,a,J\restriction X}$. First notice that $\tau_J$ {\it does not}
appear as a factor in $\hat\tau^\alpha\tau^\beta$, precisely because
the Wronskian involves only terms $a_jm_j^\delta$ with $j\notin J$, and thus
these $m_j$'s do not contain $\tau_J$. Let us pick $j_0=\min(\complement J)
\in\{0,1,\ldots,N\}$. Since $X$ is defined by $\sum_{0\leq j\leq N}
a_jm_j^\delta=0$, we have identically
$$
a_{j_0}m_{j_0}^\delta=-\sum_{i\in J}a_im_i^\delta
-\sum_{i\in \complement J\ssm\{j_0\}}a_im_i^\delta
$$
in restriction to $X$, whence (by the alternate property of
$W_k({\scriptstyle\bullet})$)
$$
W_{k,s,\hat\tau,a,J\restriction X}=-\sum_{i\in J}
W_k\big(s_1\hat\tau_1^{d-k},\ldots,s_r\hat\tau_r^{d-k},
a_im_i^\delta,(a_jm_j^\delta)_{j\in\complement J\ssm\{j_0\}}
\big)_{\restriction X}.
$$
However, all terms $m_i$, $i\in J$, contain by definition the factor $\tau_J$,
and the derivatives $D^\ell({\scriptstyle\bullet})$ leave us a factor
$m_i^{\delta-k}$ at least. Therefore, the above restricted Wronskian is also
divisible by $\tau_J^{\delta-k}$, thanks to the fact that
$\prod\hat\tau_j\prod\tau_I=0$ 
forms a simple normal crossing divisor in~$X$.\qed

\claim 12.16. Corollary|For $\delta\geq k$, there
exists a monomial $\hat\tau^{\alpha_J}\tau^{\beta_J}$
dividing $W_{k,s,\hat\tau,a,J\restriction X}$
such that
$$
|\alpha_J|+|\beta_J|=(k+c-N)(d-2k)+(N+1-c)(\delta-k)b+(\delta-k),
$$
and we have
$$
\widetilde W_{k,s,\hat\tau,a,J\restriction X}:=
(\hat\tau^{\alpha_J}\tau^{\beta_J})^{-1}
W_{k,s,\hat\tau,a,J\restriction X}
\in H^0(X,E_{k,k'}T^*_X\otimes A^{-p})
$$
where
$$
p=|\alpha_J|+|\beta_J|-(k+1)d=(\delta-k)-(k+c-N)2k-(N+1+c)(kb+\rho).
\leqno(12.17)
$$
In particular, we have $p>0$ for $\delta$ large enough $($all other parameters
being fixed or bounded$\,)$, and under this assumption, the
fundamental vanishing theorem $8.15$ implies that all entire curves
$f:\bC\to X$ are annihilated by these Wronskian operators. 
\endclaim

\proof\ In fact,
$$
(k+1)d=(k+c-N)d+(N+1-c)d=(k+c-N)d+(N+1-c)(\delta b+\rho),
$$
and we get (12.17) by subtraction.\qed

\plainsubsection 12.E. Control of the base locus for sufficiently
general coefficients $a_j$ in $\sigma$|

The next step is to control more precisely the base locus of these
Wronskians and to find conditions on $N$, $k$, $c$, $d=b\delta+\rho$
ensuring that the base locus is empty for a generic choice of the
sections $a_j$ in $\sigma=\sum a_jm_j$. Although we will not formally use it,
the next lemma is useful to realize that the base locus is related
to a natural rank condition.

\claim 12.18. Lemma|Set $u_j:=a_jm_j^\delta$. The base locus in
$X_k^\reg$ of the above Wronskians
$W_{k,s,\hat\tau,a,J\restriction X}$, when $s,\,\hat\tau$ vary,
consists of jets $f_{[k]}(0)\in X_k^\reg$ such that the matrix 
$(D^\ell(u_j\circ f)(0))_{0\leq\ell\leq k,\,j\in\complement J}$ is not of
maximal rank $($i.e., of rank${}<\card\complement J=N+1-c)\,;$
if $\delta>k$, the base locus includes
all jets $f_{[k]}(0)$ such that $f(0)\in \bigcup_{I \neq J}\tau_I^{-1}(0)$.
When $J$ also varies, the base locus of all $W_{k,s,\hat\tau,a,J\restriction X}$
in the Zariski open set $X'_k:=X_k^\reg\ssm\bigcup_{|I|=c}\tau_I^{-1}(0)$
consists of all $k$-jets such that
$\rank(D^\ell(u_j\circ f)(0))_{0\leq\ell\leq k,\,0\leq j\leq N}\leq N-c$.
\endclaim

\plainproof. If $\delta>k$ and $m_j\circ f(0)=0$ for some $j\in J$, we have
in fact $D^\ell(u_j\circ f)(0)=0$
for all derivatives $\ell\leq k$, because the exponents involved in all
factors of the differentiated monomial $a_jm_j^\delta$ are at least equal
to $\delta-k>0$. Hence the rank of the matrix cannot be maximal. Now, assume that
$m_j\circ f(0)\neq 0$ for all $j\in\complement J$, i.e.\
$$
x_0:=f(0)\in X\ssm\bigcup_{j\in\complement J}m_j^{-1}(0)=
X\ssm \bigcup_{I \neq J}\tau_I^{-1}(0).\leqno(12.19)
$$
We take sections $\hat\tau_j$ so that $\hat\tau_j(x_0)\neq 0$, and then
adjust the $k$-jet of the sections $s_1,\ldots,s_r$ in order to
generate any matrix
of derivatives $(D^\ell(s_j(f)\hat\tau_j(f)^{d-k})(0))_{0\leq\ell\leq k,\,
j\in\complement J}$
(the fact that $f'(0)\neq 0$ is used for this!). Therefore, by expanding the
determinant according to the last $(N+1-c)$ columns, we see that the
base locus is defined by the equations
$$\det(D^\ell(u_j(f))(0))_{\ell\in L,\,j\in\complement J}=0,\qquad
\forall L\subset\{0,1,\ldots,k\},~~|L|=N+1-c,\leqno(12.20)
$$
equivalent to the non-maximality of the rank. The last assertion follows
by a simple linear algebra argument.\qed

For a finer control of the base locus, we adjust the family of coefficients
$$
a=(a_j)_{0\leq j\leq N}\in S:=H^0(Z,A^\rho)^{\oplus(N+1)}\leqno(12.21)
$$
in our section $\sigma=\sum a_jm_j^\delta\in H^0(Z,A^d)$, and denote by
$X_a=\sigma^{-1}(0)\subset Z$ the corresponding hypersurface. By Lemma~12.12,
we know that there is a Zariski open set $U\subset S$ such that
$X_a$ is smooth and $\prod \tau_I=0$
is a simple normal crossing divisor in $X_a$ for all $a\in U$.
We consider the Semple tower $X_{a,k}:=(X_a)_k$ of $X_a$,
the ``universal blow-up'' $\mu_{a,k}:\widehat X_{a,k}\to X_{a,k}$
of the Wronskian ideal sheaf $\cJ_{a,k}$ such that
$\mu_{a,k}^*\cJ_{a,k}=\cO_{\widehat X_{a,k}}(-F_{a,k})$ for some
``Wronskian divisor'' $F_{a,k}$ in~$\widehat X_{a,k}$. By the universality
of this construction (cf.~Lemma~12.8),
we can also embed $X_{a,k}$ in the Semple tower
$Z_k$ of $Z$, blow up the Wronskian ideal sheaf $\cJ^Z_k$ of $Z_k$ to
get a Wronskian divisor $F^Z_k$ in $\widehat Z_k$ where $\mu_k:\widehat Z_k\to
Z_k$ is the blow-up map. Then $F_{a,k}$ is
the restriction of $F^Z_k$ to $\widehat X_{a,k}\subset\widehat Z_k$.
Our section $\widetilde W_{k,s,\hat\tau,a,J\restriction X_a}$ 
is the restriction of a \emph{meromorphic} section defined on $Z$, namely
$$
(\hat\tau^{\alpha_J}\tau^{\beta_J})^{-1}W_{k,s,\hat\tau,a,J}=
(\hat\tau^{\alpha_J}\tau^{\beta_J})^{-1}W_k\big(s_1\hat\tau_1^{d-k},...\,,s_r\hat\tau_r^{d-k},(a_jm_j^\delta)_{j\in\complement J}\big).
\leqno(12.22)
$$
It induces over the Zariski open set $Z'=Z\ssm\bigcup_I\tau_I^{-1}(0)$
a holomorphic section
$$
\sigma_{k,s,\hat\tau,a,J}
\in H^0\big(\widehat Z'_k,
\mu_k^*(\cO_{Z_k}(k')\otimes\pi_{k,0}^*A^{-p})\otimes
\cO_{\widehat Z_k}(-F^Z_k)\big)
\leqno(12.23)
$$
(notice that the relevant factors $\hat\tau_j$ remain divisible on the whole
variety~$Z$). By construction, thanks to the divisibility property explained
in Lemma~12.15, the restriction of this section to $\widehat X'_{a,k}=
\widehat X_{a,k}\cap\widehat Z'_k$ extends holomorphically
to~$\widehat X_{a,k}$, i.e.\
$$
\sigma_{k,s,\hat\tau,a,J\restriction \widehat X_{a,k}}
\in H^0\big(\widehat X_{a,k},\mu_{a,k}^*(\cO_{X_{a,k}}(k')\otimes\pi_{k,0}^*A^{-p})
\otimes\cO_{\widehat X_{a,k}}(-F_{a,k})\big).
\leqno(12.24)
$$
(Here the fact that we took $\widehat X_{k,a}$ to be normal implies that
the divided section is indeed holomorphic on $\widehat X_{k,a}$,
as $\widehat X_{k,a}\cap
\mu_k^{-1}\big(\pi_{k,0}^{-1}\bigcap_{I\in\cG}\tau_I^{-1}(0)\big)$ has
the expected codimension${}=\card\cG$ for any family $\cG$).

\claim 12.25. Lemma|Let $V$ be a finite dimensional vector space over $\bC$,
$\Psi:V^p\to\bC$ a non-zero alternating multilinear form, and let
$m,c\in\bN$, $c<m\leq p$, $r=p+c-m\geq 0$. Then the subset
$T\subset V^m$ of vectors $(v_1,\ldots,v_m)\in V^m$ such that
$$
\Psi(h_1,\ldots,h_r,(v_j)_{j\in\complement J})=0\quad
\hbox{for all $J\subset\{1,\ldots,m\}$, $|J|=c$, and all $h_1,\ldots,h_r\in V$},
\leqno(*)
$$
is a closed algebraic subset of codimension${}\geq(c+1)(r+1)$.
\endclaim

\plainproof. A typical example is $\Psi=\det$ on a $p$-dimensional vector
space~$V$, then $T$ consists of $m$-tuples of vectors of rank${}<p-r$,
and the assertion concerning the codimension is well known (we will reprove
it anyway). In general, the algebraicity of $T$ is obvious. We argue by
induction on~$p$, the result being trivial for $p=1$ (the kernel of a non-zero
linear form is indeed of codimension${}\geq 1$). If $K$ is the kernel
of $\Psi$, i.e.\ the subspace of vectors $v\in V$
such that $\Psi(h_1,\ldots,h_{p-1},v)=0$ for all $h_j\in V$, then $\Psi$ induces
an alternating multilinear form $\overline\Psi$ on $V/K$, whose kernel
is equal to $\{0\}$. The proof is thus reduced
to the case when $\Ker\Psi=\{0\}$. Notice that we must have
$\dim V\geq p$, otherwise $\Psi$ would vanish. If $\card\complement J=m-c=1$,
condition $(*)$ implies that $v_j\in\Ker\Psi=\{0\}$ for all $j$, and hence
$\codim T=\dim V^m\geq mp=(c+1)(r+1)$, as desired. Now, assume $m-c\geq 2$,
fix $v_m\in V\ssm\{0\}$ and consider the
non-zero alternating multilinear form on $V^{p-1}$ such that
$$
\Psi'_{v_m}(w_1,\ldots,w_{p-1}):=\Psi(w_1,\ldots,w_{p-1},v_m).
$$
If $(v_1,\ldots,v_m)\in T$, then $(v_1,\ldots,v_{m-1})$ belongs to
the set $T'_{v_m}$ associated with the new data $(\Psi'_{v_m},p-1,m-1,c,r)$.
The induction hypothesis implies that $\codim T'_{v_m}\geq (c+1)(r+1)$,
and since the projection $T\to V$ to the first factor
admits the $T'_{v_m}$ as its fibers, we conclude that
$$
\codim T\cap((V\ssm\{0\})\times V^{m-1})\geq (c+1)(r+1).
$$
By permuting the arguments $v_j$, we also conclude that
$$
\codim T\cap(V^{k-1}\times(V\ssm\{0\})\times V^{m-k})\geq (c+1)(r+1)
$$
for all $k=1,\ldots,m$. The union 
$\bigcup_k(V^{k-1}\times(V\ssm\{0\})\times V^{m-k})\subset V^m$  leaves 
out only $\{0\}\subset V^m$  whose codimension
is at least $mp\geq (c+1)(r+1)$, so Lemma~12.25 follows.\qed

\claim 12.26. Proposition|Consider in $U\times\widehat Z'_k$ the set 
$\Gamma$ of pairs $(a,\xi)$ such that $\sigma_{k,s,\hat\tau,a,J}(\xi)=0$ for
all choices of $s$, $\hat\tau$ and $J\subset\{0,1,\ldots,N\}$ with
$\card J=c$. Then $\Gamma$ is an algebraic set of dimension
$$
\dim \Gamma\leq \dim S-(c+1)(k+c-N+1)+n+1+kn.
$$
As a consequence, if $(c+1)(k+c-N+1)>n+1+kn$, there exists
$a\in U\subset S$ such that the base locus of the family of sections
$\sigma_{k,s,\hat\tau,a,J}$ in $\widehat X_{a,k}$ lies over
$\bigcup_IX_a\cap\tau_I^{-1}(0)$.
\endclaim

\plainproof. The idea is similar to [Brot17, Lemma 3.8], but somewhat simpler
in the present context. Let us consider a point 
$\xi\in\widehat Z'_k$ and the
$k$-jet $f_{[k]}=\mu_k(\xi)\in Z'_k$, so that
$x=f(0)\in Z'=Z\ssm\bigcup_I\tau_I^{-1}(0)$. Let us take the $\hat\tau_j$ such
that $\hat\tau_j(x)\neq 0$. Then, we do not have to pay attention to the
non-vanishing factors $\hat\tau^{\alpha_J}\tau^{\beta_J}$, and
the $k$-jets of sections $m_j$ and $\hat\tau_j^{d-k}$ are invertible
near~$x$. Let $e_A$ be a local generator of $A$ near $x$ and
$e_\cL$ a local generator of the invertible sheaf
$$
\cL=\mu_k^*\cO_{Z_k}(k')\otimes\cO_{\widehat Z_k}(-F^Z_k)
$$
near $\xi\in\widehat Z'_k$. Let $J^k\cO_{Z,x}=\cO_{Z,x}/\gm_{Z,x}^{k+1}$
be the vector space of $k$-jets of functions on $Z$ at $x$.
By definition of the Wronskian ideal and of the associated divisor~$F^Z_k$,
we have a \emph{non-zero} alternating multilinear form
$$
\Psi:(J^k\cO_{Z,x})^{k+1}\to\bC,\qquad
(g_0,\ldots,g_k)\mapsto \mu_k^*W_k(g_0,\ldots,g_k)(\xi)/e_\cL(\xi).
$$
The simultaneous vanishing of our sections at $\xi$ is equivalent
to the vanishing of
$$
\Psi\big(
s_1\hat\tau_1^{d-k}e_A^{-d},\ldots,s_r\hat\tau_r^{d-k}e_A^{-d},
(a_jm_j^\delta e_A^{-d})_{j\in\complement J}\big)
\leqno(12.27)
$$
for all $(s_1,\ldots,s_r)$. Since $A$ is very ample and $\rho\geq k$, the
power $A^\rho$ generates $k$-jets at every point $x\in Z$,
and thus the morphisms
$$
H^0(Z,A^\rho)\to J^k\cO_{Z,x},\quad
a\mapsto am_j^\delta e_A^{-d}\quad\hbox{and}\quad
H^0(Z,A^k)\to J^k\cO_{Z,x},\quad s\mapsto s\hat\tau_j^{d-k}e_A^{-d}
$$
are surjective. Lemma~12.25 applied with $r=k+c-N$ and $(p,m)$ replaced
by $(k+1,N+1)$ implies that the codimension of families
$a=(a_0,\ldots,a_N)\in S=H^0(Z,A^\rho)^{\oplus(N+1)}$ for which
$\sigma_{k,s,\hat\tau,a,J}(\xi)=0$ for
all choices of $s$, $\hat\tau$ and $J$ is at least $(c+1)(k+c-N+1)$,
i.e.\ the dimension is at most $\dim S-(c+1)(k+c-N+1)$. When we let
$\xi$ vary over $\widehat Z_k'$ which has dimension $(n+1)+kn$
and take into account the fibration $(a,\xi)\mapsto\xi$,
the dimension estimate of Proposition~12.26 follows. Under
the assumption
$$
(c+1)(k+c-N+1)>n+1+kn\leqno(12.28)
$$
we have $\dim \Gamma<\dim S$, and so the image of the projection $\Gamma\to S$,
$(a,\xi)\mapsto a$, is a constructible algebraic subset distinct from~$S$.
This concludes the proof.\qed

Our final goal is to completely eliminate the base locus. Proposition 12.26
indicates that we have to pay attention to the intersections
$X_a\cap\tau_I^{-1}(0)$. For $x\in Z$, we let $\cG$ be the family of
hyperplane sections $\tau_I=0$ that contain $x$. We introduce the set
$P=\{0,1,\ldots,N\}\ssm\bigcup_{I\in\cG}I$ and the smooth intersection
$$
Z_\cG=Z\cap\bigcap_{I\in\cG}\tau_I^{-1}(0),
$$
so that $N'+1:=\card P\geq N+1-c\card\cG$ and $\dim Z_\cG=n+1-\card\cG$.
If $a\in U$ is such that $x\in X_a$, we also look at the intersection
$$
X_{\cG,a}=X_a\cap\bigcap_{I\in\cG}\tau_I^{-1}(0),
$$
which is a smooth hypersurface of $Z_\cG$. In that situation,
we consider Wronskians
$W_{k,s,\hat\tau,a,J}$ as defined above, but we now take 
$J\subset P$, $\card J=c$, $\complement J=P\ssm J$, $r'=k+c-N'$.

\claim 12.29. Lemma|In the above setting, if we assume $\delta>k$, the
restriction
$W_{k,s,\hat\tau,a,J\restriction X_{\cG,a}}$
is still divisible by a monomial $\hat\tau^{\alpha_J}\tau^{\beta_J}$
such that
$$
|\alpha_J|+|\beta_J|=(k+c-N')(d-2k)+(N'+1-c)(\delta-k)b+(\delta-k).
$$
Therefore, if
$$
p'=|\alpha_J|+|\beta_J|-(k+1)d=(\delta-k)-(k+c-N')2k-(N'+1+c)(kb+\rho)
$$
as in $(12.17)$, we obtain again holomorphic sections
\begin{align*}
&\widetilde W_{k,s,\hat\tau,a,J\restriction X_{\cG,a}}:=
(\hat\tau^{\alpha_J}\tau^{\beta_J})^{-1}
W_{k,s,\hat\tau,a,J\restriction X_{\cG,a}}
\in H^0(X_{\cG,a},E_{k,k'}T^*_X\otimes A^{-p'}),\\
\noalign{\vskip4pt}  
&\sigma_{k,s,\hat\tau,a,J\restriction \pi_{k,0}^{-1}(X_{\cG,a})}
\in H^0\big(\pi_{k,0}^{-1}(X_{\cG,a}),
\mu_{a,k}^*(\cO_{X_{a,k}}(k')\otimes\pi_{k,0}^*A^{-p'})
\otimes\cO_{\widehat X_{a,k}}(-F_{a,k})\big).
\end{align*}
\endclaim

\plainproof. The arguments are similar to those employed in
the proof of Lemma~12.15. Let $f_{[k]}\in X_{a,k}$ be a $k$-jet such that
$f(0)\in X_{\cG,a}$ (the $k$-jet need not be entirely contained
in $X_{\cG,a}$). Putting $j_0=\min(\complement J)$, we observe that we
have on $X_{\cG,a}$ an identity
$$
a_{j_0}m_{j_0}^\delta=-\sum_{i\in P\ssm\{j_0\}}a_im_i^\delta
=-\sum_{i\in J}a_im_i^\delta-\sum_{P\ssm (J\cup\{j_0\})}a_im_i^\delta,
$$
because $m_i=\prod_{I\ni i}\tau_I=0$ on $X_{\cG,a}$
for $i\in\complement P=\bigcup_{I\in\cG}I$ (one of the factors $\tau_I$ is
such that $I\in\cG$, so $\tau_I=0$). If we compose with a germ
$t\mapsto f(t)$ such that $f(0)\in X_{\cG,a}$ (even though $f$ does not
necessarily lie entirely in $X_{\cG,a}$), we get
$$
a_{j_0}m_{j_0}^\delta(f(t))=
-\sum_{i\in J}a_im_i^\delta(f(t))-\sum_{P\ssm (J\cup\{j_0\})}a_im_i^\delta(f(t))
+O(t^{k+1})
$$
as soon as $\delta>k$. Hence we have an equality for
all derivatives $D^\ell({\scriptstyle\bullet})$, $\ell\leq k$ at $t=0$, and
$$
W_{k,s,\hat\tau,a,J\restriction X_{\cG,a}}(f_{[k]})=-\sum_{i\in J}
W_k\big(s_1\hat\tau_1^{d-k},
\ldots,s_{r'}\hat\tau_{r'}^{d-k},
a_im_i^\delta,(a_jm_j^\delta)_{j\in P\ssm(J\cup\{j_0\})}
\big)_{\restriction X_{\cG,a}}(f_{[k]}).
$$
Then, again, $\tau_J^{\delta-k}$ is a new additional common factor of
all terms in the sum, and we conclude as in Lemma 12.15 and
Corollary 12.16.\qed

Now, we analyze the base locus of these new sections on
$$
\bigcup_{a\in U}\mu_{a,k}^{-1}\pi_{k,0}^{-1}(X_{\cG,a})\subset
\mu_{k}^{-1}\pi_{k,0}^{-1}(Z_\cG)\subset\widehat Z_k.
$$
As $x$ runs in $Z_\cG$ and $N'<N$,
Lemma 12.25 shows that (12.28) can be replaced by the less
demanding condition
$$
(c+1)(k+c-N'+1)>n+1-\card\cG+kn=
\dim \mu_{k}^{-1}\pi_{k,0}^{-1}(Z_\cG).\leqno(12.28')
$$
A proof entirely similar to that of Proposition 12.26 shows
that for a generic choice of $a\in U$, the base locus of these sections
on $\widehat X_{\cG,a,k}$ projects onto
$\bigcup_{I\in\complement\cG}X_{\cG,a}\cap\tau_I^{-1}(0)$.
Arguing inductively on $\card\cG$, the base locus can be shrinked step
by step down to empty set (but it is in fact sufficient to stop when
$X_{\cG,a}\cap\tau_I^{-1}(0)$ reaches dimension $0$).

\plainsubsection 12.F. Nefness and ampleness of appropriate tautological
line bundles|

At this point, we have produced a smooth family $\cX_S\to U\subset S$
of particular hypersurfaces in $Z$, namely
$X_a=\hbox{$\{\sigma_a(z)=0\}$}$, $a\in U$,
for which a certain ``tautological'' line bundle has an empty base locus
for sufficiently general coefficients:

\claim 12.30. Corollary|Under condition $(12.28)$ and the hypothesis
$p>0$ in $(12.17)$, the following properties hold.
\plainitem{\rm(a)} The line bundle
$$
\cL_a:=\mu_{a,k}^*(\cO_{X_{a,k}}(k')\otimes\pi_{k,0}^*A^{-1})
\otimes\cO_{\widehat X_{a,k}}(-F_{a,k})
$$
is nef on $\widehat X_{a,k}$ for general $a\in U$, i.e.\ for
$a\in U'\subset U$, where $U'$ is a dense Zariski open subset.
\plainitem{\rm(b)} Let $\Delta_a=\sum_{2\leq\ell\leq k}\lambda_\ell D_{a,\ell}$
be a positive rational combination of vertical divisors of the
Semple tower and $q\in\bN$, $q\gg 1$, an integer such that
$$
\cL'_a:=\cO_{X_{a,k}}(1)\otimes\cO_{a,k}(-\Delta_a)\otimes\pi_{k,0}^*A^q
$$
is ample on $X_{a,k}$. Then the $\bQ$-line bundle
$$
\cL'_{a,\varepsilon,\eta}:=
\mu_{a,k}^*(\cO_{X_{a,k}}(k')\otimes\cO_{X_{a,k}}(-\varepsilon\Delta_a)
\otimes\pi_{k,0}^*A^{-1+q\varepsilon})
\otimes\cO_{\widehat X_{a,k}}(-(1+\varepsilon\eta)F_{a,k})
$$
is ample on $\widehat X_{a,k}$ for $a\in U'$, for some $q\in\bN$ and
$\varepsilon,\eta\in\bQ_{>0}$ arbitrarily small.\vskip0pt
\endclaim

\plainproof. (a) This would be obvious if we had global sections generating
$\cL_a$ on the whole of $\widehat X_{a,k}$, but our sections are only defined on
a stratification of~$\widehat X_{a,k}$. In any case, if
$C\subset\widehat X_{a,k}$
is an irreducible curve, we take a maximal family $\cG$ such that
$C\subset X_{\cG,a,k}$. Then, by what we have seen, for
$a\in U$ general enough, we can find global sections of $\cL_a$ on
$\widehat X_{\cG,a,k}$ such that $C$ is not contained in their base locus.
Hence $\cL_a\cdot C\geq 0$ and $\cL_a$ is nef for $a$ in a dense Zariski
open set $U'\subset U$.
\medskip

\noindent
(b) The existence of $\Delta_a$ and $q$ follows from Proposition~7.19 and
Corollary~7.21, which even provide universal values for $\lambda_\ell$ and~$q$. 
After taking the blow up $\mu_{a,k}:\widehat X_{a,k}\to X_{a,k}$ 
(cf.\ (12.7)), we infer that
$$
\cL'_{a,\eta}:=\mu_{a,k}^*\cL'\otimes\cO_{\widehat X_{a,k}}
(-\eta F_{a,k})=\mu_{a,k}^*\big(\cO_{X_{a,k}}(1)\otimes\cO_{X_{a,k}}(-\Delta_a)
\otimes\pi_{k,0}^*A^q\big)\otimes\cO_{\widehat X_{a,k}}(-\eta F_{a,k})
$$
is ample for $\eta>0$ small. The result now follows by taking a combination
$$
\cL_{a,\varepsilon,\eta}=\cL_a^{1-\varepsilon/k'}\otimes(\cL'_{a,\eta})^\varepsilon.\eqno\square
$$

\claim 12.31. Corollary|Let $\cX\to\Omega$ be the universal
family of hypersurfaces $X_\sigma=\{\sigma(z)=0\}$, $\sigma\in \Omega$,
where $\Omega\subset P(H^0(Z,A^d))$ is the dense Zariski
open set over which the family is smooth. On the ``Wronskian blow-up''
$\widehat X_{\sigma,k}$ of $X_{\sigma,k}$, let us consider the line bundle
$$
\cL_{\sigma,\varepsilon,\eta}:=
\mu_{\sigma,k}^*(\cO_{X_{\sigma,k}}(k')\otimes\cO_{X_{\sigma,k}}
(-\varepsilon\Delta_\sigma)\otimes\pi_{k,0}^*A^{-1+q\varepsilon})
\otimes\cO_{\widehat X_{\sigma,k}}(-(1+\varepsilon\eta)F_{\sigma,k})
$$
associated with the same choice of constants as in Cor.~$12.30$.
Then $\cL'_{\sigma,\varepsilon,\eta}$
is ample on $\widehat X_{\sigma,k}$ for $\sigma$ in a dense Zariski
open set $\Omega'\subset \Omega$.
\endclaim

\proof By Cor.~12.30~(b), we can find $\sigma_0\in H^0(Z,A^d)$ such that
$X_{\sigma_0}=\sigma_0^{-1}(0)$ is smooth and
$\cL_{\sigma_0,\varepsilon,\eta}^m$ is an ample line bundle on
$\widehat X_{\sigma_0,k}$ ($m\in \bN^*$). As ampleness is a
Zariski open condition,
we conclude that $\cL_{\sigma,\varepsilon,\eta}^m$ remains ample for
a general section~$\sigma\in H^0(Z,A^d)$, i.e.\ for
$[\sigma]$ in some Zariski open set $\Omega'\subset\Omega$.
Since $\mu_{\sigma,k}(F_{\sigma,k})$ is contained in
the vertical divisor of $X_{\sigma,k}$, we conclude by Theorem~8.8
that $X_\sigma$ is Kobayashi hyperbolic for $[\sigma]\in\Omega$.\qed

\plainsubsection 12.G. Final conclusion and computation of degree bounds|

At this point, we fix our integer parameters to meet all conditions
that have been found. We must have $N\geq c(n+1)$ by Lemma 12.12,
and for such a large value of $N$, condition (12.28) can hold only
when $c\geq n$, so we take $c=n$ and $N=n(n+1)$.
Inequality (12.28) then requires $k$ large enough, $k=n^3+n^2+1$
being the smallest possible value. We find
$$
b={N\choose c-1}={n^2+n\choose n-1}=n\frac{(n^2+n)\ldots(n^2+2)}{n!}.
$$
We have $n^2+k=n^2(1+k/n^2)<n^2\exp(k/n^2)$ and 
by Stirling's formula, $n!>\sqrt{2\pi n}\,(n/e)^n$. Therefore
$$
b<\frac{n^{2n-1}\exp((2+\cdots+n)/n^2)}{\sqrt{2\pi n}\,(n/e)^n}<
\frac{e^{n+\frac{1}{2}+\frac{1}{2n}}}{\sqrt{2\pi}}\,n^{n-\frac{3}{2}}.
$$
Finally, we divide $d-k$ by $b$, get in this way $d-k=b\delta+\lambda$,
$0\leq \lambda<b$, and put $\rho=\lambda+k\geq k$. Then
$\delta+1\geq(d-k+1)/b$ and formula (12.17) yields
$$
\eqalign{
p&=(\delta-k)-(n^3+1)2k-(n^2+2n+1)(kb+\rho)\cr
&\geq(d-k+1)/b-1-(2n^3+3)k-(n^2+2n+1)(kb+k+b-1).\cr}
$$
Therefore $p>0$ is achieved as soon as
$$
d\geq d_n=k+b\big(1+(2n^3+3)k+(n^2+2n+1)(kb+k+b-1)\big),
$$
where
$$
k=n^3+n^2+1,\quad b={n^2+n\choose n-1}.
$$
The dominant term in $d_n$ is $k(n^2+2n+1)b^2\sim e^{2n+1}n^{2n+2}/2\pi$.
By more precise numerical calculations and Stirling's asymptotic expansion,
one can show in fact that $d_n\leq \lfloor(n+4)\,(en)^{2n+1}/2\pi\rfloor$
for $n\geq 4$ (which is also an equivalent and a close
approximation as~$n\to+\infty$), while $d_1=61$, $d_2=6685$,
$d_3=2825761$. We can now state the main result of this section.

\claim 12.32. Theorem|Let $Z$ be a projective $(n+1)$-dimensional manifold
and $A$ a very ample line bundle on $Z$. Then, for a general section
$\sigma\in H^0(Z,A^d)$ and $d\geq d_n$, the hypersurface
$X_\sigma=\sigma^{-1}(0)$ is Kobayashi hyperbolic, and in fact,
algebraically jet hyperbolic in the sense of Definition~$11.11$. 
The bound $d_n$ for the degree can be taken to be
$$
d_n=\lfloor(n+4)\,(en)^{2n+1}/2\pi\rfloor
\quad\hbox{for $n\ge 4$},
$$
and for $n\leq 3$, one can take $d_1=4$, $d_2=6685$, $d_3=2825761$.
\endclaim

A simpler (and less refined) choice is $\tilde d_n=\lfloor
\frac{1}{3}(en)^{2n+2}\rfloor$, which is valid for all~$n$.
These bounds are only slightly weaker than the ones found by Ya Deng 
in [Deng16, Deng17], namely $\tilde d_n=O(n^{2n+6})$.\medskip

\plainproof. The
bound $d_1=4$ (instead the insane value~$d_1=61$) can be obtained in
an elementary way by adjunction: sections of $A$ can be
used to embed any polarized surface $(Z,A)$ in $\bP^N$ (one
can always take $N=5$), and we have
$K_{X_\sigma}=K_{Z\restriction X_\sigma}\otimes A^d$,
along with a surjective morphism $\Omega^2_{\bP^N}\to K_Z$.
As $\Omega^2_{\bP^N}\otimes\cO(3)=\Lambda^{N-2}(T_{\bP^N}\otimes\cO(-1))$ is
generated by sections, this implies that $K_Z\otimes A^3$ is also
generated by sections, and hence $K_{X_\sigma}$ is ample for $d\geq 4$.\qed

\plainsubsection 12.H. Further comments|

\noindent
{\bf 12.33.} Our bound $d_n$ is rather large, but just as in Ya Deng's
effective approach of Brotbek's theorem [Deng17], the bound 
holds for a property that looks substantially stronger than 
hyperbo\-licity, namely the ampleness of the pull-back of some (twisted) 
jet bundle 
\hbox{$\mu_k^*\cO_{\widehat X_k}(\abu)\otimes\cO_{\widehat X_k}(-F'_k)$}.
Section 11.B provides much weaker conditions for hyperbolicity, but
checking them is probably more involved.
\medskip

\noindent
{\bf 12.34.} After these notes were written, Riedl and Yang [RiYa18] proved
the important and somewhat surprising result that the lower bound
estimates $d_\GG(n)$ and
$d_\Kob(n)$, respectively for the Green-Griffiths-Lang and Kobayashi
conjectures for general hypersurfaces in $\bP^{n+1}$, can be related
by $d_\Kob(n):=d_\GG(2n-2)$. This should be understood in the sense
that a solution of the generic
$(2n-2)$-dimensional Green-Griffiths conjecture for $d\geq d_\GG(2n-2)$ implies
a solution of the $n$-dimensional Kobayashi conjecture for the same
lower bound. We refer to [RiYa18] for the precise statement,
which requires an ad hoc assumption on the algebraic dependence of
the Green-Griffiths locus with respect to a variation of
coefficients in the defining polynomials. In combination with [DMR10],
this gives a completely new proof of the Kobayashi conjecture, and the
order~$1$ bound $d_\GG(n)=O(\exp(n^{1+\varepsilon}))$ of [Dem12] implies
a similar
bound $d_\Kob(n)=O(\exp(n^{1+\varepsilon}))$ for the Kobayashi conjecture --
just a little bit weaker than what our direct proof gave (Theorem 12.32).
In [MeTa19], Merker and Ta were able to improve the Green-Griffiths bound
to $d_\GG(n)=o((\sqrt{n}\log n)^n)$, using a strengthening of
Darondeau's estimates [Dar16a, Dar16b], along with very delicate
calculations. The Riedl-Yang result then implies
$d_\Kob(n)=O((n\log n)^{n+1})$, which is the best bound known at
this time.
\medskip

\noindent
{\bf 12.35.} In [Ber18], G.~B\'erczi stated a positivity
conjecture for Thom polynomials of Morin singularities (see also [BeSz12]),
and showed that it would imply a polynomial bound $d_n=2\,n^9+1$ for the
generic hyperbolicity of hypersurfaces.
\medskip

\noindent
{\bf 12.36.} In the unpublished preprint [Dem15], we introduced an
alternative strategy for the proof of the Kobayashi conjecture which
appears to be still incomplete at this point. We nevertheless hope that
a refined version could one day lead to linear bounds such as
$d_\Kob(n)=2n+1$. The rough idea was to establish a $k$-jet analogue of
Claire Voisin's proof [Voi96] of the Clemens conjecture. Unfortunately,
Lemma 5.1.18 as stated in
[Dem15] is incorrect -- the assertion concerning the $\Delta$ divisor
introduced there simply does not hold. It is however conceivable that a
weaker statement holds, in the form of a control of the degree of the
divisor $\Delta$, and in a way that would still be sufficient to imply
similar consequences for the generic positivity of tautological
jet bundles, as demanded e.g.\ in section 11.B.
\bigskip

\section*{References}
\vskip5pt

\begingroup
\fontsize{10pt}{12pt}\selectfont

\Bibitem[Ahl41]&Ahlfors, L.V.:& The theory of meromorphic curves.&
Acta Soc.\ Sci.\ Finn.\ Nova Ser.~A., {\bf 3} (1941), 1--31&

\Bibitem[ASS97]&Arrondo, E., Sols, I., Speiser, R.:& Global moduli
for contacts.& Ark.\ Mat.\ {\bf 35} (1997), 1--57&

\Bibitem[AzSu80]&Azukawa, K. and Suzuki, M.& Some examples of algebraic
degeneracy and hyperbolic manifolds.& Rocky Mountain
J.\ Math.\ {\bf 10} (1980), 655--659&

\Bibitem[Ber15]&B\'erczi, G.:& Towards the Green-Griffiths-Lang conjecture
via equivariant localisation.& arXiv: math.AG/1509.03406&

\Bibitem[Ber18]&B\'erczi, G.:& Thom polynomials and the 
Green-Griffiths-Lang conjecture for hypersurfaces with polynomial degree.&
Intern.\ Math.\ Res.\ Not.\ https://doi.org/10.1093/imrn/rnx332 (2018), 1--56&

\Bibitem[BeKi12]&B\'erczi, G., Kirwan, F.:& A geometric construction
for invariant jet differentials.& Surveys in Diff.\ Geom., Vol XVII (2012),
79--125&

\Bibitem[BeSz12]&B\'erczi, G., Szenes, A.:&
Thom polynomials of Morin singularities.& Annals of Math.\
{\bf 175} (2012), 567--629&

\Bibitem[Blo26a]&Bloch, A.:& Sur les syst\`emes de fonctions uniformes
satisfaisant \`a l'\'equation d'une vari\'et\'e alg\'ebrique dont
l'irr\'egularit\'e d\'epasse la dimension.& J.\ Math.\ >Pures Appl., {\bf 5}
(1926), 19--66&

\Bibitem[Blo26b]&Bloch, A.:& Sur les syst\`emes de fonctions holomorphes
\`a vari\'et\'es lin\'eaires lacunaires.& Ann.\ Ecole Normale, {\bf 43}
(1926), 309--362&

\Bibitem[Bog77]&Bogomolov, F.A.:& Families of curves on a surface of general
type.& Soviet Math.\ Dokl.\ {\bf 236} (1977), 1294--1297&

\Bibitem[Bog79]&Bogomolov, F.A.:& Holomorphic tensors and vector bundles on
projective manifolds.& Math.\ USSR Izvestija {\bf 13}/3 (1979), 499--555&

\Bibitem[Bon93]&Bonavero, L.:& In\'egalit\'es de Morse holomorphes
singuli\`eres.& C.~R.\ Acad.\ Sci.\ Paris S\'er.~I Math.\ {\bf 317} (1993)
1163-–1166&

\Bibitem[Bro78]&Brody, R.:& Compact manifolds and hyperbolicity.& Trans.\
Amer.\ Math.\ Soc.\ {\bf 235} (1978), 213--219&

\Bibitem[BrGr77]&Brody, R., Green, M.:& A family of smooth hyperbolic
hypersurfaces in $\bP^3$.& Duke Math.\ J.\ {\bf 44} (1977), 873--874&

\Bibitem[Brot17]&Brotbek, D.:& On the hyperbolicity of general hypersurfaces.&
Publications math\'ematiques de l'IH\'ES, {\bf 126} (2017), 1--34&

\Bibitem[BrDa17]&Brotbek, D., Darondeau, L.& Complete intersection varieties
with ample cotangent bundles.& Invent.\ Math.\ {\bf 212} (2018), 913--940&

\Bibitem[BR90]&Br\"uckmann P., Rackwitz, H.-G.:& $T$-Symmetrical
Tensor Forms on Complete Intersections.& Math.\ Ann.\ {\bf 288} (1990), 627--635&

\Bibitem[Bru02]&Brunella, M.:& Courbes enti\`eres dans les surfaces 
alg\'ebriques complexes (d'apr\`es McQuillan, Demailly-El Goul, $\ldots$).& 
S\'eminaire Bourbaki, Vol.\ 2000/2001. Ast\'erisque {\bf 282} (2002), Exp.\ 
No.~881, 39–-61&

\Bibitem[Bru03]&Brunella, M.:& Plurisubharmonic variation of the leafwise 
Poincar\'e metric.& Int.\ J.~Math. {\bf 14} (2003) 139--151&

\Bibitem[Bru05]&Brunella, M.:& On the plurisubharmonicity of the leafwise
Poincar\'e metric on projective manifolds.& J.~Math.\ Kyoto Univ.\ {\bf 45}
(2005) 381--390&

\Bibitem[Bru06]&Brunella, M.:& A positivity property for foliations on 
compact K\"ahler manifolds.& Internat.\ J.~Math.\ {\bf 17} (2006) 35--43&

\Bibitem[Can00]&Cantat, S.:& Deux exemples concernant une conjecture de 
Serge Lang.& C.~R.\ Acad.\ Sci.\ Paris S\'er.~I Math.\ {\bf 330} (2000), 
581--586&

\Bibitem[CaGr72]&Carlson, J.A., Griffiths, P.A.:& A defect relation for
equidimensional holomorphic mappings between algebraic varieties.&
Ann.\ Math.\ {\bf 95} (1972), 557--584&

\Bibitem[Cart28]&Cartan, H.:& Sur les syst\`emes de fonctions holomorphes
\`a vari\'et\'es lin\'eaires lacunaires et leurs applications, Th\`ese,
Paris.& Ann.\ Ecole Normale, {\bf 45} (1928), 255--346&

\Bibitem[Carl72]&Carlson, J.A.:& Some degeneracy theorems for entire
functions with values in an algebraic variety.& Trans.\ Amer.\ Math.\
Soc.\ {\bf 168} (1972), 273–-301&

\Bibitem[CFZ17]&Ciliberto, C., Flamini, F., Zaidenberg, M.:&
A remark on the intersection of plane curves.& 
math.AG, arXiv:1704.00320, 16p&

\Bibitem[CKM88]&Clemens, H., Koll\'ar, J., Mori, S.:& Higher dimensional
complex geometry.& Ast\'erisque {\bf 166}, 1988&

\Bibitem[Cle86]&Clemens, H.:& Curves on generic hypersurfaces.& Ann.\
Sci.\ Ec.\ Norm.\ Sup.\ {\bf 19} (1986), 629--636&

\Bibitem[ClR04]&Clemens, H., Ran, Z.:& Twisted genus bounds for subvarieties 
of generic hypersurfaces.& Amer.\ J.\ Math.\ {\bf 126} (2004), 89--120&

\Bibitem[CoKe94]&Colley, S.J., Kennedy, G.:& The enumeration of
simultaneous higher order contacts between plane curves.& Compositio
Math.\ {\bf 93} (1994), 171--209&

\Bibitem[Coll88]&Collino, A.:& Evidence for a conjecture of Ellingsrud
and Str\o mme on the Chow ring of ${\bf Hilb}_d(\bP^2)$.& Illinois J.\
Math.\ {\bf 32} (1988), 171--210&

\Bibitem[CoGr76]&Cowen, M., Griffiths, P.A.:& Holomorphic curves and
metrics of negative curvature.& J.~Analyse Math.\ {\bf 29} (1976),
93--153&

\Bibitem[Dar14]&Darondeau, L.:& Effective algebraic degeneracy of entire
curves in complements of smooth projective hypersurfaces.&
preprint Univ.\ Paris Sud Orday, arXiv:1402.1396&

\Bibitem[Dar16a]&Darondeau, L.:& Fiber integration on the Demailly tower.&
Ann.\ Inst.\ Fourier {\bf 66} (2016), 29--54&

\Bibitem[Dar16b]&Darondeau, L.:& On the logarithmic Green-Griffiths
conjecture.& Int.\ Math.\ Res.\ Not.\ (2016), no.~6, 1871--1923&

\Bibitem[DPP06]&Debarre, O., Pacienza, G., P{\u{a}}un, M.:& Non-deformability 
of entire curves in projective hypersurfaces of high degree.& Ann.\ Inst.\
Fourier (Grenoble) {\bf 56} (2006) 247--253&

\Bibitem[Dem85]&Demailly, J.-P.:& Champs magn\'etiques et in\'egalit\'es de
Morse pour la $d''$-coho\-mo\-logie.& Ann.\ Inst.\ Fourier
(Grenoble) {\bf 35} (1985) 189--229&

\Bibitem[Dem90a]&Demailly, J.-P.:& Cohomology of $q$-convex spaces in top degrees.& 
Math.\ Zeitschrift {\bf 203} (1990) 283--295&

\Bibitem[Dem90b]&Demailly, J.-P.:& Singular hermitian metrics on positive
line bundles.& Proceedings of the Bayreuth conference ``Complex algebraic
varieties'', April~2--6, 1990, edited by K.~Hulek, T.~Peternell,
M.~Schneider, F.~Schreyer, Lecture Notes in Math.\ ${\rm n}^\circ\,$1507,
Springer-Verlag (1992), 87--104&

\Bibitem[Dem92]&Demailly, J.-P.:& Regularization of closed positive
currents and Intersection Theory.& J.\ Alg.\ Geom.\ {\bf 1} (1992)
361--409&

\Bibitem[Dem94]&Demailly, J.-P.:& $L^2$ vanishing theorems for positive line bundles and adjunction theory.& alg-geom/9410022$\,$; Lecture Notes of the 
CIME Session ``Transcendental methods in Algebraic Geometry'', Cetraro, 
Italy, July 1994, Ed.\ F.~Catanese, C.~Cili\-berto, Lecture Notes in Math.,
Vol.~1646, CIME Found.\ Subser., Springer-Verlag, 1--97&

\Bibitem[Dem95]&Demailly, J.-P.:& Algebraic criteria for Kobayashi
hyperbolic projective varieties and jet differentials.& AMS Summer
School on Algebraic Geometry, Santa Cruz 1995, Proc.\ Symposia in
Pure Math., vol.\ {\bf 062.2}, ed.\ by J.~Koll\'ar and R.~Lazarsfeld,
Amer.\ Math.\ Soc., Providence, RI (1997), 285–-360&

\Bibitem[Dem97]&Demailly, J.-P.:& Vari\'et\'es hyperboliques et
\'equations diff\'erentielles alg\'ebriques.& Gaz.\ Math.\ {\bf 73}
(juillet 1997) 3--23&

\Bibitem[Dem07a]&Demailly, J.-P.:& Structure of jet differential rings and 
holomorphic Morse inequalities.& Talk at the CRM Workshop 
``The geometry of holomorphic and algebraic curves in complex algebraic
varieties'', Montr\'eal, May 2007&

\Bibitem[Dem07b]&Demailly, J.-P.:& On the algebraic structure of the
ring of jet differential operators.& Talk at the conference
``Effective aspects of complex hyperbolic varieties'', Aber Wrac'h,
France, September 10-14, 2007&

\Bibitem[Dem11]&Demailly, J.-P.:& 
Holomorphic Morse Inequalities and the Green-Griffiths-Lang Conjecture.&
Pure and Applied Math.\ Quarterly {\bf 7} (2011), 1165--1208&

\Bibitem[Dem12]&Demailly, J.-P.:& 
Hyperbolic algebraic varieties and holomorphic differential equations.&
expanded version of the lectures given at the annual meeting of VIASM,
Acta Math.\ Vietnam.\ {\bf 37} (2012), 441-–512&

\Bibitem[Dem14]&Demailly, J.-P.:& 
Towards the Green-Griffiths-Lang conjecture.&
Conference ``Analysis and Geometry'', Tunis,
March 2014, in honor of Mohammed Salah Baouendi, 
ed.\ by A.~Baklouti, A.\ El Kacimi, S.~Kallel, N.~Mir, Springer
Proc.\ Math.\ Stat., {\bf 127}, Springer-Verlag, 2015,
141--159&

\Bibitem[Dem15]&Demailly, J.-P.:& 
Proof of the Kobayashi conjecture on the hyperbolicity
of very general hypersurfaces.& arXiv:1501.07625, math.CV, unpublished&

\Bibitem[DeEG97]&Demailly, J.-P., El Goul, J.:& Connexions m\'eromorphes 
projectives et vari\'et\'es alg\'ebriques hyperboliques.& C.\ R.\ Acad.\ 
Sci.\ Paris, S\'erie~I, {\bf 324} (1997), 1385--1390&

\Bibitem DeEG00&Demailly, J.-P., El Goul, J.:& Hyperbolicity of
generic surfaces of high degree in projective 3-space.& Amer.\ J.\ 
Math.\ {\bf 122} (2000), 515--546&

\Bibitem[DeLS94]&Demailly, J.-P., Lempert L., Shiffman, B.:& Algebraic
approximation of holomorphic maps from Stein domains to projective
manifolds.& Duke Math.~J.\ {\bf 76} (1994) 333--363&

\Bibitem DePS94&Demailly, J.-P., Peternell, Th., Schneider, M.:&
Compact complex manifolds with numerically effective tangent bundles.&
Algebraic Geometry {\bf 3} (1994), 295--345&

\Bibitem[Deng16]&Deng, Ya.:& Effectivity in the hyperbolicity related
problems.& Chap.~4 of the PhD memoir ``Generalized Okounkov Bodies,
Hyperbolicity-Related and Direct Image Problems'' defended on June~26, 2017
at Universit\'e Grenoble Alpes, Institut Fourier, arXiv:1606.03831, math.CV&

\Bibitem[Deng17]&Deng, Ya.:& On the Diverio-Trapani Conjecture&
arXiv:1703.07560, math.CV&

\Bibitem[DGr91]&Dethloff, G., Grauert, H.:& On the infinitesimal deformation
of simply connected domains in one complex variable.& International Symposium
in Memory of Hua Loo Keng, Vol.\ II (Beijing, 1988), Springer, Berlin, (1991),
57--88&

\Bibitem[DLu01]&Dethloff, G., Lu, S.:& Logarithmic jet bundles and 
applications.& Osaka J.\ Math.\ {\bf 38} (2001), 185--237&

\Bibitem[DTH16a]&Dinh, Tuan Huynh:& Examples of hyperbolic
hypersurfaces of low degree in projective spaces.& 
Int.\ Math.\ Res.\ Not.\ {\bf 18} (2016) 5518–-5558&

\Bibitem[DTH16b]&Dinh, Tuan Huynh:& Construction of hyperbolic
hypersurfaces of low degree in $\bP^n(\bC)$.& Int.~J.\ Math.\ {\bf 27}
(2016) 1650059 (9 pages)&

\Bibitem[Div08]&Diverio, S.:& Differential equations on complex projective
hypersurfaces of low dimension.& Compos.\ Math.\ {\bf 144} (2008) 920-–932&

\Bibitem[Div09]&Diverio, S.:& Existence of global invariant jet
differentials on projective hypersurfaces of high degree.& Math.\
Ann.\ {\bf 344} (2009) 293-–315&

\Bibitem[DMR10]&Diverio, S., Merker, J., Rousseau, E.:& Effective
algebraic degeneracy.& Invent.\ Math.\ {\bf 180} (2010) 161--223&

\Bibitem[DR11]&Diverio, S., Rousseau, E.:& A survey on hyperbolicity of projective
hypersurfaces.& Publ.\ Mat.\ IMPA, Inst.\ Nac.\ Mat.\ Pura. Apl.,
Rio de Janeiro, 2011&

\Bibitem[DR15]&Diverio, S., Rousseau, E.:& The exceptional set and the
Green–Griffiths locus do not always coincide.& Enseign.\ Math.\
{\bf 61} (2015) 417--452&

\Bibitem[DT10]&Diverio, S., Trapani, S.:& A remark on the codimension
of the Green-Griffiths locus of generic projective hypersurfaces of
high degree.& J.\ Reine Angew.\ Math.\ {\bf 649} (2010) 55–-61&

\Bibitem[Dol81]&Dolgachev, I.:& Weighted projective varieties.&
Proceedings Polish-North Amer.\ Sem. on Group Actions and Vector Fields,
Vancouver, 1981, J.B.~Carrels editor, Lecture Notes in Math.\ {\bf 956},
Springer-Verlag (1982), 34--71&

\Bibitem[Duv04]&Duval J.:& Une sextique hyperbolique dans $\bP^3(\bC)$.& 
Math.\ Ann.\ {\bf 330} (2004) 473–-476&

\Bibitem[Duv08]&Duval J.:& Sur le lemme de Brody.& Invent.\ Math.\ {\bf 173} (2008), 
305--314&

\Bibitem[Ein88]&Ein L.:& Subvarieties of generic complete intersections.&
Invent.\ Math.\ {\bf 94} (1988), 163--169&

\Bibitem[Ein91]&Ein L.:& Subvarieties of generic complete intersections, II.&
Math.\ Ann.\ {\bf 289} (1991), 465--471&

\Bibitem[EG96]&El Goul, J.:& Algebraic families of smooth hyperbolic
surfaces of low degree in $\bP^3_\bC$.& Manuscripta Math.\ {\bf 90}
(1996), 521--532&

\Bibitem[EG97]&El Goul, J.:& Propri\'et\'es de n\'egativit\'e de courbure
des vari\'et\'es alg\'ebriques hyperboliques.& Th\`ese de Doctorat, Univ.\ 
de Grenoble~I (1997)&

\Bibitem[Fuj72]&Fujimoto H.:& On holomorphic maps into a taut complex space&
Nagoya Math.~J.\ {\bf 46} (1972), 49--61&

\Bibitem[Fuj01]&Fujimoto, H.:& A family of hyperbolic hypersurfaces in
the complex projective space.& Complex Variables, Theory Appl.\ {\bf 43},
no.~3-4, (2001) 273-283&

\Bibitem[Fuji94]&Fujita, T.:& Approximating Zariski decomposition of big line
bundles.& Kodai Math.\ J.\ {\bf 17} (1994) 1--3&

\Bibitem[Ghe41]&Gherardelli, G.:& Sul modello minimo della varieta degli
elementi differenziali del $2^\circ$ ordine del piano projettivo.&
Atti Accad.\ Italia.\ Rend., Cl.\ Sci.\ Fis.\ Mat.\ Nat.\ (7)
{\bf 2} (1941), 821--828&

\Bibitem[Gra89]&Grauert, H.:& Jetmetriken und hyperbolische Geometrie.&
Math.\ Zeitschrift {\bf 200} (1989), 149--168&

\Bibitem[GRe65]&Grauert, H., Reckziegel, H.:& Hermitesche Metriken und
normale Familien holomorpher Abbildungen.& Math.\ Zeitschrift {\bf 89}
(1965), 108--125&

\Bibitem[GrGr80]&Green, M., Griffiths, P.A.:& Two applications of algebraic
geometry to entire holomorphic mappings.& The Chern Symposium 1979,
Proc.\ Internal.\ Sympos.\ Berkeley, CA, 1979, Springer-Verlag, New York
(1980), 41--74&

\Bibitem[Gri71]&Griffiths, P.:& Holomorphic mappings into canonical
algebraic varieties.& Ann.\ of Math.\ {\bf 98} (1971), 439--458&

\Bibitem[Har77]&Hartshorne, R.:& Algebraic geometry.& Gard.\ Texts in
Math., {\bf 52}, Springer-Verlag, 1977&

\Bibitem[Hir64]&Hironaka, H.:& Resolution of singularities of an algebraic
variety over a field of characteristic zero.& Ann.\ of Math.\ {\bf 79}
(1964) 109--326&

\Bibitem[HVX17]&Huynh, D.T., Vu, D.V., Xie, S.Y.:& Entire holomorphic
curves into projective plane intersecting few generic algebraic curves.&
math.AG, arXiv:1704.03358, 9p&

\Bibitem[Kaw80]&Kawamata, Y.:& On Bloch's conjecture.& Invent.\ Math.\
{\bf 57} (1980), 97--100&

\Bibitem[Kob70]&Kobayashi, S.:& Hyperbolic manifolds and holomorphic
mappings.& Pure Appl.\ Math., {\bf 2}, Marcel Dekker, New York, NY, 1970&

\Bibitem[Kob75]&Kobayashi, S.:& Negative vector bundles and complex
Finsler structures.& Nagoya Math.\ J.\ {\bf 57} (1975), 153--166&

\Bibitem[Kob76]&Kobayashi, S.:& Intrinsic distances, measures and geometric
function theory.& Bull.\ Amer.\ Math.\ Soc.\ {\bf 82} (1976), 357--416&

\Bibitem[Kob80]&Kobayashi, S.:& The first Chern class and holomorphic
tensor fields.& J.\ Math.\ Soc.\ Japan {\bf 32} (1980), 325--329&

\Bibitem[Kob81]&Kobayashi, S.:& Recent results in complex differential
geometry.& Jber.\ dt.\ Math.-Verein.\ {\bf 83} (1981), 147--158&

\Bibitem[Kob98]&Kobayashi, S.:& Hyperbolic complex spaces.& Grundlehren der
Mathematischen Wissenschaften, volume 318, Springer-Verlag, Berlin, 1998&

\Bibitem[KobO71]&Kobayashi, S., Ochiai, T.:& Mappings into compact complex
manifolds with negative first Chern class.& J.\ Math.\ Soc.\ Japan
{\bf 23} (1971), 137--148&

\Bibitem[KobO75]&Kobayashi, S., Ochiai, T.:& Meromorphic mappings onto
compact complex spaces of general type.& Invent.\ Math.\ {\bf 31}
(1975), 7--16&

\Bibitem[KobR91]&Kobayashi, R.:& Holomorphic curves into algebraic
subvarieties of an abelian variety.& Internat.\ J.\ Math.\ {\bf 2} (1991),
711--724&

\Bibitem[LaTh96]&Laksov, D., Thorup, A.:& These are the differentials of
order $n$.& Trans.\ Amer.\ Math.\ Soc.\ {\bf 351} (1999), 1293–-1353&

\Bibitem[Lang86]&Lang, S.:& Hyperbolic and Diophantine analysis.&
Bull.\ Amer.\ Math.\ Soc.\ {\bf 14} (1986), 159--205&

\Bibitem[Lang87]&Lang, S.:& Introduction to complex hyperbolic spaces.&
Springer-Verlag, New York (1987)&

\Bibitem[Lu96]&Lu, S.S.Y.:& On hyperbolicity and the Green-Griffiths 
conjecture for surfaces.& Geometric Complex Analysis, ed.\ by J.~Noguchi 
et al., World Scientific Publishing Co.\ (1996), 401--408&

\Bibitem[LuMi95]&Lu, S.S.Y., Miyaoka, Y.:& Bounding curves in algebraic
surfaces by genus and Chern numbers.& Math.\ Research Letters {\bf 2}
(1995), 663-676&

\Bibitem[LuMi96]&Lu, S.S.Y., Miyaoka, Y.:& Bounding codimension one
subvarieties and a general inequality between Chern numbers.&
Amer.\ J.\ of Math. {\bf 119} (1997), 487--502&

\Bibitem[LuWi12]&Lu, S.S.Y., Winkelmann, J.:& Quasiprojective varieties
admitting Zariski dense entire holomorphic curves.& Forum Math.\
{\bf 24} (2012), 399--418&

\Bibitem[LuYa90]&Lu, S.S.Y., Yau, S.T.:& Holomorphic curves in surfaces
of general type.& Proc.\ Nat.\ Acad.\ Sci.\ USA, {\bf 87} (January 1990),
80--82&

\Bibitem[MaNo96]&Masuda, K., Noguchi, J.:& A construction of hyperbolic
hypersurface of $\bP^n(\bC)$.& Math.\ Ann.\ {\bf 304} (1996), 339--362&

\Bibitem[McQ96]&McQuillan, M.:& A new proof of the Bloch conjecture.&
J.\ Alg.\ Geom.\ {\bf 5} (1996), 107--117&

\Bibitem[McQ98]&McQuillan, M.:& Diophantine approximation and foliations.&
Inst.\ Hautes \'Etudes Sci.\ Publ.\ Math.\ {\bf 87} (1998), 121--174&

\Bibitem[McQ99]&McQuillan, M.:& Holomorphic curves on hyperplane sections 
of $3$-folds.& Geom.\ Funct.\ Anal.\ {\bf 9} (1999), 370--392&

\Bibitem[Mer08]&Merker, J.:& Jets de Demailly-Semple d'ordres 4 et 5
en dimension 2.& Int.~J.\ Contemp.\ Math.\ Sci.\ {\bf 3-18} (2008), 861--933&

\Bibitem[Mer09]&Merker, J.:& Low pole order frames on vertical jets of the 
universal hypersurface& Ann.\ Inst.\ Fourier (Grenoble), {\bf 59} (2009),
1077--1104&

\Bibitem[Mer10]&Merker, J.:&Application of computational invariant theory to
Kobayashi hyperbolicity and to Green–Griffiths algebraic degeneracy&
J.~of Symbolic Computation, {\bf 45} (2010) 986–1074&

\Bibitem[Mer15]&Merker, J.:& Complex projective hypersurfaces of general type: towards a conjecture of Green and Griffiths& arXiv:1005.0405; Algebraic
differential equations for entire holomophic curves in projective
hypersurfaces of general type: optimal lower degree bound. In:
Geometry and Analysis on Manifolds, in memory of Professor Shoshichi Kobayashi,
(eds.\ T.~Ochiai, T.~Mabuchi, Y.~Maeda, J.~Noguchi and A.~Weinstein), 
Progress in Mathematics, Birkh\"auser {\bf 308} (2015), 41--142&

\Bibitem[MeTa19]&Merker, J., Ta, The-Anh:&
Degrees $d\geq(\sqrt{n}\log n)^n$ and $d\geq (n\log n)^n$
in the Conjectures of Green-Griffiths and of Kobayashi.&
arXiv:1901.04042, math.AG&

\Bibitem[Mey89]&Meyer, P.-A.:& Qu'est ce qu'une diff\'erentielle d'ordre
$n\,$?& Expositiones Math.\ {\bf 7} (1989), 249--264&

\Bibitem[Miy82]&Miyaoka, Y.:& Algebraic surfaces with positive indices.
Classification of algebraic and analytic manifolds.&
Katata Symp.\ Proc.\ 1982, Progress in Math., vol.~39, Birkh\"auser,
1983, 281--301&

\Bibitem[MoMu82]&Mori, S., Mukai, S.:& The uniruledness of the moduli space
of curves of genus~$11$.& In: Algebraic Geometry Conference
Tokyo-Kyoto 1982, Lecture Notes in Math.\ 1016, 334--353&

\Bibitem[Nad89]&Nadel, A.M.:& Hyperbolic surfaces in $\bP^3$.&
Duke Math.\ J.\ {\bf 58} (1989), 749--771&

\Bibitem[Nog77a]&Noguchi, J.:& Holomorphic curves in algebraic varieties.&
Hiroshima Math.\ J.\ {\bf 7} (1977), 833--853&

\Bibitem[Nog77b]&Noguchi, J.:& Meromorphic mappings into a compact complex
space.& Hiroshima Math.~J.\ {\bf 7} (1977), 411--425&

\Bibitem[Nog81a]&Noguchi, J.:& A higher-dimensional analogue of Mordell's
conjecture over function fields.& Math.\ Ann.\ {\bf 258} (1981),
207-212&

\Bibitem[Nog81b]&Noguchi, J.:& Lemma on logarithmic derivatives and
holomorphic curves in algebraic varieties.& Nagoya Math.\ J.\ {\bf 83}
(1981), 213--233&

\Bibitem[Nog86]&Noguchi, J.:& Logarithmic jet spaces and extensions of
de Franchis' theorem.& Contributions to Several Complex Variables, 
Aspects of Math., E9, Vieweg, Braunschweig (1986), 227--249&

\Bibitem[Nog91]&Noguchi, J.:& Hyperbolic manifolds and diophantine geometry.&
Sugaku Expositiones {\bf 4}, Amer.\ Math.\ Soc.\ Providence, RI (1991),
63--81&

\Bibitem[Nog96]&Noguchi, J.:& Chronicle of Bloch's conjecture.&
Private communication&

\Bibitem[Nog98]&Noguchi, J.:& On holomorphic curves in semi-abelian
varieties.& Math.\ Z.\ {\bf 228} (1998), 713--721&

\Bibitem[NoOc90]&Noguchi, J.; Ochiai, T.:& Geometric function theory in
several complex variables.& Japanese edition, Iwanami, Tokyo, 1984$\,$;
English translation, xi${}+282$~p., Transl.\ Math.\ Monographs
{\bf 80}, Amer.\ Math.\ Soc., Providence, Rhode Island, 1990&

\Bibitem[NoWi13]&Noguchi, J., Winkelmann, J.:& Nevanlinna Theory in
Several Complex Variables and Diophantine Approximation.&
GrundLehren der Math.\ Wiss.\ Vol.~350, Springer,
Tokyo-Heidelberg-New York-Dordrecht-London, 2014&

\Bibitem[NWY07]& Noguchi, J., Winkelmann, J., Yamanoi, K.:&
Degeneracy of holomorphic curves into algebraic varieties.&
J.\ Math.\ Pures Appl.\ (9) {\bf 88}, No. 3, (2007), 293--306&

\Bibitem[NWY13]& Noguchi, J., Winkelmann, J., Yamanoi, K.:& Degeneracy of
holomorphic curves into algebraic varieties. II.& Vietnam J.\ Math.\
{\bf 41}, No. 4, (2013), 519--525&

\Bibitem[Och77]&Ochiai, T.:& On holomorphic curves in algebraic varieties
with ample irregularity.& Invent.\ Math. {\bf 43} (1977), 83--96&

\Bibitem[Pac04]&Pacienza, G.:& Subvarieties of general type on a general projective 
hypersurface.& Trans.\ Amer.\ Math.\ Soc.\ {\bf 356} (2004), 2649--2661&

\Bibitem[PaRo07]&Pacienza, G., Rousseau, E.:& On the logarithmic Kobayashi
conjecture.& J.~reine angew.~Math.\ {\bf 611} (2007), 221--235&

\Bibitem[Pau08]&P\u{a}un, M.:& Vector fields on the total space of 
hypersurfaces in the projective space and hyperbolicity.&
Math.\ Ann.\ {\bf 340} (2008) 875--892&

\Bibitem[RiYa18]&Riedl, E., Yang, D.:& Applications of a grassmannian technique in hypersurfaces.& June 2018, arXiv:1806.02364, math.AG&

\Bibitem[Rou06]&Rousseau, E.:& \'Etude des jets de Demailly-Semple en dimension $3$.&
Ann.\ Inst.\ Fourier (Grenoble) {\bf 56} (2006), 397--421&

\Bibitem[Roy71]&Royden, H.:& Remarks on the Kobayashi metric.& In:
Several Complex Variables II, Proc.\
Maryland Conference on Several Complex Variables, Lecture Notes,
Vol.~185, Springer-Verlag, Berlin (1971), 125--137&

\Bibitem[Roy74]&Royden, H.:& The extension of regular holomorphic maps.&
Proc.\ Amer.\ Math.\ Soc.\ {\bf 43} (1974), 306--310&

\Bibitem [Sei68]&Seidenberg, A.:& Reduction of the singularities of the
differential equation $Ady=Bdx$.& Amer.\ J.~of Math.\ {\bf 90} (1968), 
248--269&

\Bibitem [Sem54]&Semple, J.G.:& Some investigations in the geometry of
curves and surface elements.& Proc.\ London Math.\ Soc.\ (3) {\bf 4}
(1954), 24--49&

\Bibitem [ShZa02]&Shiffman, B., Zaidenberg, M.:& Hyperbolic hypersurfaces in
in $\bP^n$ of Fermat-Waring type.& Proc.\ Amer.\ Math.\ Soc. {\bf 130} (2002),
2031--2035&

\Bibitem[Shi98]&Shirosaki, M.:& On some hypersurfaces and holomorphic
mappings.& Kodai Math.\ J.~{\bf 21}, no.~1, (1998), 29--34&

\Bibitem [Siu76]&Siu, Y.T.:& Every Stein subvariety admits a
Stein neighborhood.& 
Invent.\ Math.\ {\bf 38} (1976) 89--100&

\Bibitem [Siu87]&Siu, Y.T.:& Defect relations for holomorphic maps between
spaces of different dimensions.& Duke Math.\ J.\ {\bf 55} (1987),
213--251&

\Bibitem[Siu93]&Siu, Y.T.:& An effective Matsusaka big theorem.&
Ann.\ Inst.\ Fourier (Grenoble), {\bf 43} (1993), 1387--1405&

\Bibitem[Siu97]&Siu, Y.T.:& A proof of the general Schwarz lemma using
the logarithmic derivative lemma.& Personal communication, April
1997&

\Bibitem[Siu02]&Siu, Y.T.:& Some recent transcendental techniques in
algebraic and complex geometry.& Proceedings of the International
Congress of Mathematicians, Vol.~I (Beijing, 2002), Higher
Ed.\ Press, Beijing (2002) 439--448&

\Bibitem[Siu04]&Siu, Y.T.:& Hyperbolicity in complex geometry.& The 
legacy of Niels Henrik Abel, Springer, Berlin (2004) 543--566&

\Bibitem[Siu15]&Siu, Y.T:& Hyperbolicity of generic high-degree hypersurfaces
in complex projective spaces.& Inventiones Math.\ {\bf 202} (2015)
1069--1166&

\Bibitem[SiYe96a]&Siu, Y.T., Yeung, S.K.:& Hyperbolicity of the complement of
a generic smooth curve of high degree in the complex projective plane.&
Invent.\ Math.\ {\bf 124} (1996), 573--618&

\Bibitem[SiYe96b]&Siu, Y.T., Yeung, S.K.:& A generalized Bloch's theorem
and the hyperbolicity of the complement of an ample divisor in an abelian
variety.& Math.\ Ann.\ {\bf 306} (1996), 743--758&

\Bibitem[SiYe97]&Siu, Y.T., Yeung, S.K.:& Defects for ample divisors of
Abelian varieties, Schwarz lemma and hyperbolic hypersurfaces of low degree.& 
Amer.~J.\ Math.\ {\bf 119} (1997), 1139--1172&

\Bibitem[Tra95]&Trapani, S.:& Numerical criteria for the positivity of the
difference of ample divisors.& Math.\ Z.\ {\bf 219} (1995), 387--401&

\Bibitem[Tsu88]&Tsuji, H.:& Stability of tangent bundles of minimal 
algebraic varieties.& Topology {\bf 27} (1988), 429--442&

\Bibitem[Ven96]&Venturini S.:& The Kobayashi metric on complex spaces.&
Math.\ Ann.\ {\bf 305} (1996), 25--44&

\Bibitem[Voi96]&Voisin, C.:& On a conjecture of Clemens on rational
curves on hypersurfaces.& J.\ Diff.\ Geom.\ {\bf 44} (1996) 200--213,
Correction: J.\ Diff.\ Geom.\ {\bf 49} (1998), 601--611&

\Bibitem[Voj87]&Vojta, P.:& Diophantine approximations and value
distribution theory.& Lecture Notes in Math.\ {\bf 1239},
Springer-Verlag, Berlin, 1987&

\Bibitem[Win07]&Winkelmann, J.:& On Brody and entire curves.& Bull.\ 
Soc.\ Math.\ France {\bf 135} (2007), 25–-46&

\Bibitem[Xie15]&Xie, S.-Y.:& On the ampleness of the cotangent bundles
of complete intersections.& Inventiones Math.\ {\bf 212} (2018) 941--996&

\Bibitem[Xu94]&Xu, G.:& Subvarieties of general hypersurfaces in projective
space.& J.\ Differential Geometry {\bf 39} (1994), 139--172&

\Bibitem[Zai87]&Zaidenberg, M.G.:& The complement of a generic
hypersurface of degree $2n$ in $\bC\bP^n$ is not hyperbolic.&
Siberian Math.\ J.\ {\bf 28} (1987), 425--432& 

\Bibitem[Zai93]&Zaidenberg, M.G.:& Hyperbolicity in projective spaces.&
International Symposium on Holomorphic mappings, Diophantine Geometry
and Related topics, 1992, 
R.I.M.S.\ Lecture Notes ser.\ {\bf 819},
R.I.M.S.\ Kyoto University (1993), 136--156&

\endgroup

\vskip15pt
\parindent=0cm
Version of June 19, 2019, proof corrected on November 28, 2019,\\
printed on \today, \timeofday\\
The original publication is available at {\tt www.springerlink.com}
(Japanese Journal of Mathematics)
\vskip15pt

Jean-Pierre Demailly\hfil\break
Universit\'e Grenoble Alpes, Institut Fourier (Math\'ematiques)\hfil\break
UMR 5582 du C.N.R.S., 100 rue des Maths, 38610 Gi\`eres, France\hfil\break
{\em e-mail:}\/ jean-pierre.demailly@univ-grenoble-alpes.fr

\end{document}